\documentclass[10pt,leqno]{amsart}

\newtheorem{theorem}{Theorem}

\newtheorem{lemma}{Lemma}

\newtheorem{corollary}{Corollary}
\newtheorem{definition}{Definition}
\newtheorem{example}{Example}

\newtheorem{convention}{Convention}
\newtheorem{notation}{Notation}

\usepackage{graphics,epsfig,amssymb,amsfonts,amscd,xypic}
\def\A{{\mathbb A}}
\def\B{{\mathbb B}}
\def\C{{\mathbb C}}

\def\I{{\mathbb I}}
\def\K{{\mathbb K}}
\def\L{{\mathbb L}}
\def\N{{\mathbb N}}

\def\R{{\mathbb R}}

\def\dl{[\![}
\def\dr{]\!]}

\begin{document}


\title[Geometry of generic submanifolds and reflection principle
]{
On the local geometry of generic submanifolds of
$\C^n$\\ and the analytic reflection principle\\
(Part~I)
}

\author{Jo\"el Merker}

\address{
CNRS, Universit\'e de Provence, LATP, UMR 6632, CMI, 
39 rue Joliot-Curie, 13453 Marseille Cedex 13, France}

\email{merker@cmi.univ-mrs.fr}



\date{\number\year-\number\month-\number\day}

\maketitle

\begin{center}
\begin{minipage}[t]{13.25cm}
\baselineskip 0.32cm
\bigskip

{\scriptsize
\centerline{\bf Table of contents~:}

\smallskip

{\bf Chapter~1.~General introduction for Part~I\dotfill 1.}

{\bf Chapter~2.~Geometry of complexified generic submanifolds and Segre 
chains \dotfill 2.}

{\bf Chapter~3.~Nondegeneracy conditions for generic submanifolds \dotfill 3.}

{\bf Chapter~4.~Nondegeneracy conditions for power series CR mappings 
\dotfill 4.}

}
\end{minipage}
\end{center}

\bigskip
\bigskip

\noindent
{\Large \bf Chapter~1: General introduction for Part~I}

\section*{\S1.1.~Resumed historical background}

\subsection*{1.1.1.~Local Lie groups and the no Riemann mapping theorem at 
the boundary} Inspired by the general idea that, in analogy with
\'E.~Galois's group theory of algebraic equations, group analysis of
differential equations should provide precious information about their
solvability, S.~Lie began around 1873--80 the classification of all
continuous local groups of transformation acting on $\C^n$. He quickly
succeeded for $n=1$ and achieved the case $n=2$~({\it see}~[18]), but
the unavoidable complexity and richness for $n=3$ exhausted his
efforts; moreover, after more than one century, the task has never
been achieved. Nevertheless, especially for $n=2$, Lie's
classificationf had the enormous power of providing any possible
application to the study of transformations preserving arbitrary types
of geometric structures. Thanks to the influence of G.~Darboux, the
works of S.~Lie were rapidly known to French mathematicians.  Based on
the general approach of S.~Lie, H.~Poincar\'e ({\it see}~[24])
discovered in 1907 that the automorphism groups of the two-dimensional
unit ball $\B^2:=\{(z_1,z_2)\in \C^2: \, \vert z_1 \vert^2+\vert
z_2\vert^2< 1\}$ and of the bidisc $\Delta^2:=\{(z_1,z_2)\in\C^2: \,
\vert z_1 \vert < 1, \, \vert z_2 \vert < 1\}$ are represented by
rational, but not isomorphic transformations and he deduced
immediately that $\B^2$ and $\Delta^2$ are not biholomorphically
equivalent.  This discovery was the starting point of the {\it no
Riemann mapping theorem} in several complex variables.

In the beginning of the twentieth century, the birth of pluricomplex
geometry also coincided with two other fundamental memoirs of
F.~Hartogs~[13] (1906) and of E.E.~Levi~[16] (1910). However, whereas
this direction had important ramifications in the years 1930-50,
especially with the works of W.~Osgood, of H.~Kneser, of R.~Fueter, of
E.~Martinelli, of K.~Behnke, of F.~Sommer, of S.~Bochner, and
culminated in the complete solution of the so-called {\sl problem of
Levi} given by K.~Oka in 1951--52, the direction initiated by
H.~Poincar\'e in 1907 lay dormant for approximatively sixty years,
with the major exception of four consecutive and historically isolated
memoirs of B.~Segre~[25], [26] and of \'E.~Cartan~[3], [4] in the
years 1931-32. Based on works of S.~Lie, of A.~Tresse (a french
student of S.~Lie), and of the young mathematician B.~Segre,
\'E.~Cartan (who also had defended his thesis under the direction of
S.~Lie) provided an essentially complete classification of all Levi
nondegenerate real analytic local hypersurfaces in $\C^2$, which
ultimately relies on S.~Lie's far reaching works about the
classification of second order ordinary differential equations.
Approximatively thirty years later, in 1974, \'E.~Cartan's equivalence
algorithm has been conducted in $\C^n$ for $n\geq 2$ by S.S.~Chern and
J.K.~Moser in~[5] to provide an {\it a priori} complete list of
differential invariants for Levi nondegenerate real analytic
hypersurface in $\C^n$ for $n\geq 2$, {\it see} also further
developments by A.G.~Vitushkin~[33] and
N.G.~Kruzhillin~[15]. However, the classification problem (in the
sense of S.~Lie) for real analytic Levi nondegenerate hypersurfaces in
$\C^n$ for $n\geq 3$ is essentially left incomplete by the analysis
in~[5], because the list of differential invariants does not provide
immediately a list of all possible automorphisms groups.

\subsection*{1.1.2.~Reflection principle and regularity of CR 
mappings} The real birth of Cauchy-Riemann geometry occured in the
beginning of the years 1970, especially when in 1974, C.~Fefferman
({\it see}~[10]) established that every biholomorphism between two
smoothly bounded strongly pseudoconvex domains extend smoothly as a CR
diffeomorphism between their boundaries. It is still conjectured, but
up to now unproved, that the result remains true without any
pseudoconvexity assumption.  Thus, the classification of bounded
domains up to biholomorphisms reduces to the classification of
boundaries up to CR diffeomorphisms. At the same time, S.~Pinchuk
discovered in~[21], [22], [23] an important local extension theorem 
for CR mappings
between real analytic hypersurfaces, in which both the Schwarz
reflection principle phenomenon and the Hartogs extension phenomenon
contribute to the analytic continuation of CR mappings. In 1977,
generalizing H.~Poincar\'e's grounding result, S.M. Webster
established in~[34] a general result according to which local CR mappings are
complex algebraic as soon as the source and target hypersurfaces are
algebraic. Thus, in the aim of generalizing Carath\'eodory's
theorems about the boundary regularity of conformal maps in the
complex plane, the grounding works of C.~Fefferman, of S.~Pinchuk and
of S.M.~Webster initiated a completely new subject about the
regularity (or the analytic continuation) of biholomorphic (or proper)
mappings (or of local CR mappings). Since then, this subject has been
very active during almost thirty years and a substantial amount of
efforts by numerous mathematicians has led to some remarkable
refinements of the original statements\footnotemark[1].

\footnotetext[1]{
However, one should remind of the historical bifurcation between the
classification problem and the reflection principle. It is probable
that too much emphasis has been put in the last decade on the
reflection principle, which occulted in part the original motivation
of classifying domains. Whereas this memoir is exclusively devoted to
the so-called {\it analytic reflection principle}, we believe that it
is time to come back to the original program of research hidden in the
mathematical treasures of S.~Lie and of \'E.~Cartan, as suggested for
instance in the recent works~[29],~[31].}

\section*{\S1.2.~Conceptional description of the topics adressed in 
Part~I of this memoir}

\subsection*{1.2.1.~Division in two parts}
This memoir is devoted to a synthetic exposition of some recent
results in the direction of the so-called {\it analytic reflection
principle}. This terminology justifies by the fact that most arguments
and proofs are based on Taylor series considerations.  The main topics
adressed in Part~I of this memoir is to study {\it ab initio} the
local geometry of arbitrary real algebraic or analytic submanifolds of
$\C^n$ which are {\sl generic}, namely which satisfy
$T_pM+iT_pM=T_p\C^n$ for every $p\in M$. Our main goal is to explain
how to go beyond the classical notion of Levi nondegeneracy, taking
account of the complexity due to arbitrary dimension and codimension,
in order to formulate appropriate generalizations of the reflection
principle.  In a forthcoming volume of the same collection, Part~II of
this memoir, which will be accompagnied by its own conceptional
introduction, will be specifically devoted to the study of the
analytic reflection principle. Thus the present Part~I is a kind of
thorough preparation, of which we can now present a quick description.

\subsection*{1.2.2.~Canonical pair of foliations attached to the
extrinsic complexification of a local generic submanifold} As will be
established in Theorem~2.1.22, a given real analytic generic
submanifold $M$ in $\C^n$ of codimension $d$ and of CR dimension
$m=n-d$ may be locally represented, in a neighborhood of one of its
points $p$, thanks to some appropriate coordinates
$t=(z,w)\in\C^m\times \C^d$ vanishing at $p$, by means of $d$ complex
defining fundamental equations of the form
\def\theequation{1.2.3}\begin{equation}
\bar w_j=\Theta_j(\bar z, z, w), \ \ \ \ \ \ j=1,\dots,d.
\end{equation}
Here, we assume that the Taylor series of the complex analytic
functions $\Theta_j(\bar z, z, w)=\sum_{\beta\in\N^m, \,
\gamma\in\N^m, \, \delta\in\N^d}\, \Theta_{j,\beta,\gamma,\delta}\,
\bar z^\beta \, z^\gamma\, w^\delta$ converge normally in some
polydisc centered at the origin in $\C^{m+m+d}$. As will be more
evident in the sequel, one finds many advantages to deal with complex
defining equations instead of real defining equations. Of course, the
conjugate defining equations $w_j=\overline{\Theta}_j(z,\bar z, \bar
w)
=\sum_{\beta\in\N^m, \,
\gamma\in\N^m, \, \delta\in\N^d}\, 
\overline{\Theta}_{j,\beta,\gamma,\delta}\,
z^\beta \, \bar z^\gamma\, \bar 
w^\delta$ must define the same generic real submanifold $M$, and the
ambiguity due to complex defining equations disappears thanks to the
existence of the following fundamental functional equations, obtained
in Theorem~2.1.32:
\def\theequation{1.2.4}\begin{equation}
\bar w_j\equiv
\Theta_j(\bar z, z, \overline{\Theta(\bar z, z, w)}), \ \ \ \ \
j=1,\dots,d.
\end{equation}  
We say that $M$ is {\sl algebraic} (in the sense of J. Nash) if the series
$\Theta_j(\bar z,z,w)$ are algebraic (a power series
$\varphi(x)\in\C\{x\}$ is (Nash) {\sl algebraic} if there exists a nonzero
polynomial $P(x,y)\in\C[x,y]$ such that $P(x,\varphi(x))\equiv 0$).

Following S.M. Webster's general philosophy ({\it cf.}~[37]), let $\tau=
(\zeta, \xi) \in \C^m \times \C^d$ be new independent coordinates
corresponding to $\bar t= (\bar z,\bar w)$ and define the {\it
extrinsic complexification} of $M$ to be the $d$-codimensional complex
submanifold $\mathcal{M}$ of $\C^{2n}$ defined by the equations
\def\theequation{1.2.5}\begin{equation}
\xi_j=\Theta_j(\zeta,z,w), \ \ \ \ \ \ j=1,\dots,d, 
\end{equation}
or equivalently by $w_j=\overline{\Theta}_j(z,\zeta,\xi)$,
$j=1,\dots,d$. Notice that the expressions $\Theta_j(\zeta,z,w)
=\sum_{\beta\in\N^m, \, \gamma\in\N^m, \, \delta\in\N^d}\,
\Theta_{j,\beta,\gamma,\delta}\, \zeta^\beta \, z^\gamma\, w^\delta$
are meaningful only because the $\Theta_j$ are converging power
series.  This submanifold $\mathcal{M}$ comes immediately equipped
with two foliations $\mathcal{F}=:\mathcal{M}\cap \{\tau=ct.\}$ and
$\underline{\mathcal{F}}= \mathcal{M}\cap \{t=ct.\}$ by
$m$-dimensional complex submanifolds, which were essentially
discovered by B. Segre in~[25], [26] ({\it see} also [3], [34]).  We
call the leaves of $\mathcal{F}$ the {\sl complexified Segre
varieties} and the leaves of $\underline{\mathcal{F}}$ the {\sl
conjugate complexified Segre varieties}. As we shall argue throughout
this memoir, the main features of the geometry of $M$ are hidden
behind the interweaving of this pair of foliations $(\mathcal{F},
\underline{\mathcal{F}})$ lying on its complexification $\mathcal{M}$.

Since we are mainly interested in the study of mappings, let $M'$ be a
second generic submanifold of codimension $d'$ in $\C^{n'}$ defined
similarly by complex defining equations $\bar w_{j'}'=
\Theta_{j'}'(\bar z',z',w')$, $j'=1,\dots,d'$, where $m'=n'-d'$ and
$t'=(z',w')\in\C^{m'}\times \C^{d'}$, and let a local mapping
$t'=h(t)=(h_1(t),\dots,h_{n'}(t))$ from $\C^n$ to $\C^{n'}$ be a {\sl
local power series CR mapping from $M$ to $M'$}. By this, we mean
precisely that there exists a $d'\times d$ matrix of power series
$a(t,\bar t)$ such that if we split $h(t)=(f(t),g(t))\in\C^{m'}\times
\C^{d'}$, then the following vectorial formal power series holds in
$\C\dl t,\bar t\dr^{d'}$:
\def\theequation{1.2.6}\begin{equation}
\bar g(\bar t)-\Theta'(\bar f(\bar t),f(t),g(t))\equiv
a(t,\bar t) \, [\bar w-\Theta(\bar z,z,w)].
\end{equation}
In this memoir, we shall always assume that $M$ and $M'$ are
real algebraic or analytic and we shall consider three different
regularity classes for $h$, namely either $h(t)$ is a purely formal
power series, or it is complex analytic, or it is complex algebraic.
By complexifying~\thetag{1.2.6}, we trivially obtain the following
identity in $\C\dl t,\tau\dr^{d'}$:
\def\theequation{1.2.7}\begin{equation}
\bar g(\tau)-\Theta'(\bar f(\tau),f(t),g(t))\equiv
a(t,\tau)\, [\xi-\Theta(\zeta,z,w)],
\end{equation}
which means precisely that the power series mapping $(t',\tau')=
(h(t),\bar h(\tau))$ maps the complexifications $\mathcal{M}$ into the
complexification $\mathcal{M}'$.
We shall denote by $h^c(t,\tau):=(h(t),\bar h(\tau))$ this
complexified mapping. A straightforward but crucial
observation is that $h^c$ stabilizes the two pairs of foliations, namely it
satisfies $h^c( \mathcal{F}) \subset \mathcal{ F}'$ and
$h^c(\underline{\mathcal{F}}) \subset \underline{ \mathcal{F}}'$. The
following symbolic figure is an attempt to illustrate this
stabilization property\footnotemark[2].

\footnotetext[2]{However, we warn the reader that the
representation is slightly incorrect, because the ambient codimensions
$d$ and $d'$ in $\mathcal{M}$ and in $\mathcal{M}'$ of the unions of
foliations $\mathcal{F}\cup \underline{\mathcal{F}}$ and
$\mathcal{F}'\cup \underline{\mathcal{F}}'$ is invisible in this
picture.  One should imagine for instance that $\mathcal{M}$ and
$\mathcal{M}'$ are three-dimensional spaces equipped with pairs of
foliations by horizontal orthogonal real lines.}

\begin{center}
\input stabilization-foliations.pstex_t
\end{center}

{\it Some strong rigidity properties are due to the fact that
$h^c=(h,\bar h)$ must respect these two pairs of foliations}.  For
instance, a theorem due to S.M.~Webster in~[34] states that every
local biholomorphism $h: M\to M'$ between two {\it Levi nondegenerate}
real algebraic hypersurfaces must be a complex algebraic mapping. This
theorem may be interpreted intuitively by thinking that $h^c$ (which
is {\it a priori} only complex analytic) is forced to be as smooth as
the two pairs of foliations $(\mathcal{F}, \underline{\mathcal{F}})$
and $(\mathcal{F}',\underline{\mathcal{F}}')$ are, namely to be
complex algebraic.

\subsection*{1.2.9.~Beyond Levi nondegeneracy: Minimality and finite 
nondegeneracy} In S.M.~Webster's theorem, behind Levi nondegeneracy
are hidden two highly different and independent concepts, the notion
of {\sl minimality} (in the sense of J.-M.~Tr\'epreau and
A.E.~Tumanov, following the general approach of H.J.~Sussmann in~[32],
{\it see} also~[2]) and the notion of {\sl finite nondegeneracy}
(introduced for the first time by K. Diederich 
and S.M. Webster in~[9], and by C.K.~Han in~[12] and then studied by
S.M.~Baouendi, P.~Ebenfelt and L.P.~Rothschild in~[1]).

The first main concept of minimality is of geometric nature and may
easily be described in terms of the pair of foliations $(\mathcal{F},
\underline{\mathcal{F}})$. Let $z_1\in \C^m$.  We denote by
$\underline{\Gamma}_1(z_1)$ the point located in the (vertical)
$m$-dimensional complex leaf $\underline{\mathcal{F}}_0$ passing
through the origin which lies at distance $z_1$ from the origin, {\it
see} {\sc Figure~1.2.8}.  Of course, $\underline{\Gamma}_1(z_1)$
belongs to $\mathcal{M}$. In other words, we move vertically from the
origin up to distance $z_1\in\C^m$. Let $z_2\in \C^m$. From this point
$\underline{\Gamma}_1(z_1)$, we then move horizontally up to distance
$z_2$, namely following the $m$-dimensional complex leaf
$\mathcal{F}_{\underline{\Gamma}_1(z_1)}$. We denote by
$\underline{\Gamma}_2(z_{(2)})$ the resulting point, {\it see} again
{\sc Figure~1.2.8}, where we use the notation
$z_{(2)}:=(z_1,z_2)\in\C^{2m}$. Of course, the point
$\underline{\Gamma}_2(z_{(2)})$ also belongs to $\mathcal{M}$.  Let
$z_3\in\C^m$. We further move vertically up to distance $z_3$ and we
denote the resulting point by $\underline{\Gamma}_3(z_{(3)})$, where
$z_{(3)}=(z_1,z_2,z_3)\in\C^{3m}$. More generally, by following
alternately the two foliations $\underline{\mathcal{F}}$ and
$\mathcal{F}$, we may define for every $k\in\N$ a point
$\underline{\Gamma}_k(z_{(k)})$, where $z_{(k)}=(z_1,\dots,
z_k)\in\C^{km}$, which belongs to $\mathcal{M}$. It is easy
to see that the mapping $z_{(k)}\mapsto \underline{\Gamma}_k(
z_{(k)})$ satisfies $\underline{\Gamma}_k(0)=0$ and has the same
regularity as $\mathcal{M}$, namely it is complex algebraic or
analytic. We call this map the {\sl $k$-th conjugate Segre chain}.
The precise construction of $\underline{\Gamma}_k$ is presented in
Chapter~2 and there are combinatorial formulas which yield the
complete expression of $\underline{\Gamma}_k(z_{(k)})$ by means of the
fundamental power series $\Theta_j(\zeta,z,w)$.

The complexification $\mathcal{M}$ is then called {\sl minimal at the
origin} if there exists an integer $k$ such that for every
neighborhood $\mathcal{V}_k$ of the origin in $\C^{km}$, its image
$\underline{\Gamma}_k(\mathcal{V}_k)$ contains a neighborhood of the
origin in $\mathcal{M}$. Intuitively, the concept of minimality says
that one can reach every point in a neighborhood of the origin in
$\mathcal{M}$ by following alternately the two canonical foliations.
Since the two foliations $\mathcal{F}$ and $\underline{\mathcal{F}}$
are biholomorphically invariant, 
the notion of minimality at one point $p\in M$ so defined is 
independent of the choice of coordinates vanishing at $p$.  It is
elementary to see that a Levi nondegenerate real analytic hypersurface
in $\C^n$ ($n\geq 2$) is minimal at every point.

The second main concept of finite nondegeneracy is of analytic nature and
it may be easily described by means of a development in power series
of the defining equations of $\mathcal{M}$:
\def\theequation{1.2.10}\begin{equation}
\xi_j=\sum_{\beta\in\N^m}\, 
\zeta^\beta \, \Theta_{j,\beta}(t).
\end{equation}
Here, the $\Theta_{j,\beta}(t)=\Theta_{j,\beta}(z,w)$ are complex
algebraic or analytic power series converging normally in a uniform
polydisc centered at the origin. Let $k\in\N$. By the {\sl $k$-th
Segre mapping} we mean the local complex algebraic or analytic mapping
\def\theequation{1.2.11}\begin{equation}
\mathcal{Q}_k: 
\C^n\ni t \longmapsto \left(
\Theta_{j,\beta}(t)
\right)_{1\leq j\leq d, \, 
\vert \beta \vert \leq k}\in\C^{N_{d,n,k}},
\end{equation}
where the integer $N_{d,n,k}$ denotes the total number of $k$-th jets
of a $d$-vectorial mapping of $n$ independent variables, namely
$N_{d,n,k}=d \, {(n+k)!\over n! \ k!}$. The generic submanifold $M$ is
then called {\sl finitely nondegenerate at the origin} if there exists
an integer $k$ such that the $k$-th Segre mapping is of (maximal
possible) rank $n$ at the origin. Although the mapping $\mathcal{Q}_k$
is defined in terms of a system of coordinates and although it seems
to depend on the choice of complex defining equations for $M$, it may
be established that its properties are essentially biholomorphically
and invariantly attached to $M$, and in particular, the notion of
finite nondegeneracy at a point $p\in M$ so defined is independent of
the choice of coordinates vanishing at $p$.  One can show that Levi
nondegeneracy of $M$ at the origin (in the sense that the kernel of
the vector-valued Levi form of $M$ is zero) is equivalent to the fact
that the mapping $\mathcal{Q}_1$ is of rank $n$ at the origin, hence
the notion of finite nondegeneracy is a generalization of the notion
of Levi nondegeneracy. More generally, $M$ is called {\sl
holomorphically nondegenerate at the origin} (in the sense of
N.~Stanton, {\it cf.}~[28]) if there exists an integer $k$ such that
the generic rank of $\mathcal{Q}_k$ is equal to $n$. Further study of
nondegeneracy conditions on the mapping $\mathcal{Q}_k$ are presented
in Chapter~3. Since this has been suggested in~[8] and~[9], we also
endeavour a self-contained study of {\sl jets of Segre varieties}, a
fundamental topic for which we know no complete background
reference.

\subsection*{1.2.12.~Local geometry at a Zariski-generic point}
Why are minimality and finite nondegeneracy adequate concepts from the
point of view of local Cauchy-Riemann geometry~? Firstly, because it
may be established that attached to a given arbitrary connected real
agebraic or analytic generic submanifold $M$ in $\C^n$, there exists
an invariant integer $d_{2,M}$ and a proper real algebraic or analytic
subvariety $E\subset M$ such that for every point $p\in M\backslash
E$, there exists a neighborhood $V_p$ of $p$ in $\C^n$ and a system of
complex algebraic or analytic coordinates $(t_1,\dots,t_n)$ centered
at $p$ such $M\cap V_p$ is contained in the transverse intersection of
$d_{2,M}$ Levi flat hypersurfaces defined by $\{\bar t_1=t_1,\dots,
t_{ d_{2,M}}=t_{ d_{2,M}}\}$ and such that, moreover, for every
constant $(c_1, \dots,c_{ d_{2,M}}) \in \R^{ d_{2,M}}$, the
intersection $M_c:= M \cap \{ t_1= c_1, \dots, t_{ d_{2,M}}=c_{
d_{2,M}}\} \cap V_p$ is minimal at every point (Corollary~2.8.5).
Here, the $M_c$ are elementary ``bricks'' and there is no ``complex
link'' between them. Hence one may think that minimality is a
good ``general'' assumption.

Secondly, it may be furthermore established that there exists an
invariant integer $n_M$ with $d\leq n_M \leq n$ and another proper
subvariety $F\subset M$ such that for every point $p\in M\backslash
F$, there exists a neighborhood $V_p$ of $p$ in $\C^n$ and a system of
coordinates centered at $p$ in which $M\cap V_p$ is a product
$\underline{M}_p'\times \Delta^{n-n_M}$ of a $d$-codimensional generic
submanifold $\underline{M}_p'$ in $\C^{n_M}$ by a complex polydisc
$\Delta^{n-n_M}$, such that moreover $\underline{M}_p'$ is finitely
nondegenerate at its central point (Theorem~3.5.48). In particular,
$M$ is holomorphically nondegenerate if and only if $n=n_M$, in which
case $M$ is finitely nondegenerate at every point of $M\backslash F$.
Generally speaking, from the point of view of CR geometry where
``Complex'' and ``Real'' concepts should be truly associated, the
factor $\Delta^{n-n_M}$ is essentially superfluous, hence one should
think that finite nondegeneracy is a good ``general'' assumption.

Whereas minimality and finite nondegeneracy do not impose dimensional
restrictions, it is well known that the assumption of Levi
nondegeneracy requires that $d\leq m^2$. In addition, there exist some
classes of hypersurfaces in $\C^3$ whose Levi form is of rank one at
every point and which are finitely nondegenerate at every point ({\it
see} Examples~3.2.15 and~3.2.20). In sum, we believe that minimality
and finite nondegeneracy are adequate assumptions.

\subsection*{1.2.13.~Nondegeneracy conditions for power series CR mappings}
In Chapter~4 of part~I of this memoir, we shall introduce various
nondegeneracy conditions for power series CR mappings.  As in \S1.2.2,
let $h$ be a power series CR mapping from $M$ to $M'$, whose
complexification $h^c=(h,\bar h)$ satisfies the fundamental
identities~\thetag{1.2.7}, which yields after replacing $\xi$ by
$\Theta(\zeta,t)$ the following formal identities in $\C\dl
\zeta,t\dr$:
\def\theequation{1.2.14}\begin{equation}
\bar g_{j'}(\zeta,\Theta(\zeta,t))\equiv 
\Theta_{j'}(\bar f(\zeta,\Theta(\zeta,t)),h(t)), \ \ \ \ \
j'=1,\dots,d'.
\end{equation}
We consider the following pairwise
commuting $m$ vector fields tangent to $\mathcal{M}$
\def\theequation{1.2.15}\begin{equation}
\underline{\mathcal{L}}_k:={\partial \over \partial \zeta_k}+
\sum_{j=1}^d\, {\partial \Theta_j\over \partial \zeta_k}
(\zeta,t)\,{\partial \over \partial \xi_j},
\end{equation}
which span the leaves of $\underline{ \mathcal{F}}$ at every point.
For every $\beta=(\beta_1, \dots, \beta_m) \in \N^m$, we introduce the
derivation $\underline{\mathcal{L}}^\beta:= \underline{ \mathcal{
L}}_1^{\beta_1}\cdots \underline{\mathcal{L}}_m^{\beta_m}$ of order
$\vert \beta\vert$, which we apply to the
equations~\thetag{1.2.14}. After some computations, this yields an
expression of the form $R_{j', \beta}'(t, \tau, ( \partial_\tau^\alpha
\bar h( \tau))_{\vert \alpha \vert \leq \vert \beta \vert}: \, h(t))$,
where $R_{j',\beta}'$ is a certain analytic expression in its
variables.  Based on the properties of the infinite collection of
functions $R_{j',\beta}'$, we shall formulate five technical
nondegeneracy conditions about $h$. For further intuitive explanation,
we refer to the beginning of the conceptional introduction of the
forthcoming Part~II of this memoir.

\smallskip
\noindent
{\tt Note to the Russian translator(s):} Since my English has probably
some deficiencies, please, do not hesitate to arrange the translation
in classical Russian style. For any question about the meaning of a
phrase or of a paragraph which would be difficult to understand and
difficult to translate, please, do not hesitate to contact me by
e-mail, or by regular mail, asking me to rewrite phrases or
paragraphs in a better, more understandable style, which I would
diligently do.

\bigskip

\noindent
{\Large \bf Chapter~2: Geometry of complexified generic submanifolds and Segre chains}

\section*{2.1.~Elementary local geometry of Cauchy-Riemann submanifolds} 

\subsection*{2.1.1.~Formal, Analytic, Algebraic}
As we shall essentially deal in the two parts of this memoir with
local power series centered at the origin, we start with classical
definitions.  Let $\K$ be the field $\R$ of real numbers or the field
$\C$ of complex numbers.  Let $n\in\N$ be a positive integer.  Let
$x=(x_1,\dots,x_n)\in \K^n$. We denote $\vert x \vert:=\max\{\vert x_1
\vert, \dots, \vert x_n\vert\}$. Let $\K\dl x\dr$ denote the local
ring of formal power series in the $n$ variables $(x_1,\dots,x_n)$. By
definition, an element $\varphi(x)\in \K\dl x\dr$ writes in the form
$\varphi(x)=\sum_{\alpha\in \N^n}\, \varphi_\alpha\, x^\alpha$, where
$x^\alpha=x_1^{\alpha_1}\cdots x_n^{\alpha_n}$ and $\varphi_\alpha\in
\K$ for every multiindex $\alpha=(\alpha_1,\dots,\alpha_n)\in\N^n$.
We say that $\varphi$ is a {\sl $\K$-formal} power series.  Such a
power series $\varphi(x)= \sum_{\alpha\in\N^n}\, \varphi_\alpha\,
x^\alpha$ is {\sl identically zero} if all its coefficients
$\varphi_\alpha$ are zero. We write this property $\varphi(x)\equiv
0$. If $\alpha=(\alpha_1,\dots,\alpha_n)$ is a multiindex in $\N^n$,
we denote its {\sl length} by $\vert \alpha \vert :=
\alpha_1+\cdots+\alpha_n$ and the corresponding partial derivative of
a power series by $\partial_x^\alpha\varphi(x):=
\partial_{x_1}^{\alpha_1}\cdots
\partial_{x_n}^{\alpha_n}\varphi(x)$. Sometimes, we use also the
equivalent notation $\partial^{\vert \alpha \vert}\varphi(x)/ \partial
x_1^{\alpha_1}\cdots \partial x_n^{\alpha_n}$. We have
$\varphi_\alpha=[1/\alpha!]\, \partial_x^\alpha
\varphi(x)\vert_{x=0}$.  If the coefficients satisfy a Cauchy estimate
like $\vert \varphi_\alpha\vert \leq C \, \rho^{-\vert
\alpha\vert}$, where $C>0$ and $\rho>0$, the series {\sl converges
normally} in the polydisc $\Delta_n(\rho)=\{x\in \K^n: \, 
\vert x \vert < \rho\}$. We say that $\varphi$ is
{\sl $\K$-analytic} and we write $\varphi(x) \in\K \{x\}$. If there
exists moreover a nonzero polynomial $P(X_1, \dots,X_n, \Phi)\in
\K[X_1, \dots, X_n, \Phi]$ such that $P(x_1,\dots,x_n, \varphi (x_1,
\dots, x_n))\equiv 0$ for all $x\in\ \Delta_n (\rho)$, we say that
$\varphi$ is {\sl $\K$-algebraic} (in the sense of J.~Nash) and we write
$\varphi(x) \in \mathcal{A}_\K \{x\}$.  By classical elimination theory
it follows that if, more generally, a power series $\varphi(x) \in
\K\{x\}$ satisfies a polynomial equation $P(\varphi(x)) \equiv 0$, where
$P(\Phi)\in \mathcal{A}_\K\{x\}[T]$ is a polynomial in the
indeterminate $\Phi$ with coefficients in $\mathcal{A}_\K\{x\}$, then
$\varphi(x)$ is $\K$-algebraic. We have the following inclusions:
\def\theequation{2.1.2}\begin{equation}
\K\dl x\dr \supset \K\{x\} \supset
\mathcal{A}_\K\{x\},
\end{equation}
which are all strict. The three rings $\K\dl x\dr$, $\K\{x\}$ and
$\mathcal{A}_\K\{x\}$ are local, noetherian and they
enjoy the Weierstrass division property.

\subsection*{2.1.3.~Composition, differentiation and implicit 
function theorem} Furthermore, the rings $\K\dl x\dr$, $\K\{x\}$ and
$\mathcal{A}_\K\{x\}$ are stable under elementary
algebraic operations, under composition, under differentiation and the
implicit function theorem holds true.  Only $\mathcal{A}_\K\{x\}$ is
dramatically unstable under integration. The following known theorem,
that we shall admit, summarizes these properties.

\def\thetheorem{2.1.4}\begin{theorem}
The following three statements hold true:
\begin{itemize}
\item[{\bf (1)}]
Let $n$ and $d$ be positive integers, let $x\in\K^n$, let $y\in\K^d$, let
$\varphi(x)$ belong to $\K\dl x\dr$, to
$\K\{x\}$ or to $\mathcal{A}_\K\{x\}$, 
let $\psi_1(y),\dots,\psi_n(y)$ belong to $\K\dl y\dr$, to $\K\{y\}$, or to
$\mathcal{A}_\K\{y\}$ and vanish at the origin.  Then
$\varphi(\psi_1(y),\dots,\psi_n(y))$ belongs $\K\dl y\dr$, to $\K\{y\}$, or to
$\mathcal{A}_\K(y)$.
\item[{\bf (2)}]
Let $n$ be a positive integer and let $x\in\K^n$.  If a power series
$\varphi(x)$ belongs to $\K\dl x\dr$, to $\K\{x\}$, or to
$\mathcal{A}_\K\{x\}$, then for every multiindex $\alpha\in\N^n$, the
partial derivative $\partial_x^\alpha\varphi(x)$ also belongs to
$\K\dl x\dr$, to $\K\{x\}$, or to $\mathcal{A}_\K\{x\}$.
\item[{\bf (3)}]
Let $n$ and $d$ be positive integers, let $x\in\K^n$, $y\in\K^d$
and let $H_1(x,y),\dots, H_d(x,y)$ be a collection of formal, analytic or
algebraic power series vanishing at the origin, namely the $H_j(x,y)$
belongs to $\K\dl x,y\dr$, to $\K\{x,y\}$, or to $\mathcal{A}_\K\{x,y\}$
and they satisfy $H_j(0,0)=0$ for $j=1,\dots,d$. Assume that the functional
determinant $(\partial H_{j_1}/\partial y_{j_2}(0))_{1\leq j_1,j_2\leq
d}$ does not vanish. Then there exists a unique $\K^d$-valued power
series mapping $\varphi(x)= (\varphi_1(x),\dots, \varphi_d(x))$, where
the $\varphi_j(x)$ belong to $\K\dl x\dr$,
to $\K\{x\}$, or to $\mathcal{A}_\K\{x\}$ and vanish at the origin, such that
$H_j(x,\varphi(x))\equiv 0$ for $j=1,\dots,d$.
\end{itemize}
\end{theorem}

\subsection*{2.1.5.~Local submanifolds and their mappings}
By definition, a {\sl local submanifold} $M$ of $\K^n$ is identified with
the data of $d\leq n$ power series $r_1(x),\dots,r_d(x)$ vanishing at
the origin and which belong to $\K\dl x\dr$, to $\K\{x\}$, or to
$\mathcal{A}_\K\{x\}$ such that the linear forms $dr_1(0),\dots,
dr_d(0)$ are linearly independent.  Two data $r(x)=(r_1(x),
\dots,r_d(x))$ and $\widehat{r} (x)=(\widehat{ r}_1(x), \dots,
\widehat{r}_d(x))$ define the {\sl same submanifold} if there exists
an invertible $d\times d$ matrix $a(x)=(a_{j_1, j_2}(x))_{ 1\leq
j_1,j_2\leq d}$ of power series in $\K\dl x\dr$, in $\K\{x \}$, or in
$\mathcal{A}_\K\{x\}$, such that $\widehat{r}_j(x)\equiv
\sum_{l=1}^d\, a_{j,l}(x)\, r_l(x)$, or in matrix notation $\widehat{
r}(x)\equiv a(x)\, r(x)$. Clearly this defines an equivalence relation
between $d$-tuples of power series $r(x)=(r_1(x),\dots,r_d(x))$ whose
differentials are independent at the origin. A submanifold identifies
with an equivalence class.  We shall write $M:
r_1(x)=\cdots=r_d(x)=0$, keeping in mind that the identification of
$M$ with its ``zero set'' is meaningless in the formal category. We
call $d$ the {\sl codimension} of $M$.  Let $x'=\Phi(x)$ be a formal,
algebraic or analytic invertible change of coordinates centered at the
origin and let $x=\Phi'(x')$ denote its inverse.  The {\sl transformed
submanifold} $M':= \Phi(M)$ is defined by the collection $r'(x'):=(r_1
(\Phi'(x')), \dots, r_d(\Phi'(x')))$. It follows from the formal,
analytic or algebraic implicit function theorem, that there exists a
local invertible transformation $x'=\Phi(x)$ such that $r_1
(\Phi'(x')) =x_1', \dots, r_d(\Phi'(x')) =x_d'$, so the image
$M':=\Phi (M)$ writes $M': x_1'=\dots=x_d'=0$.

Let $n$ and $n'$ be two positive integers. A formal, analytic
or algebraic {\sl local mapping from $\K^n$ to $\K^{n'}$} consists of the
datum of a $n'$-tuple $h(x)=(h_1(x),\dots,h_{n'}(x))$ of power series
$h_{i'}(x)$ in $\K\dl x\dr$, in $\K\{x\}$, or in $\mathcal{A}_\K\{x\}$,
with $h_{i'}(0)=0$ for $i'=1,\dots,n'$.  We write $x'=h(x)$. If $n''$
is a third positive integer and if $x''=g(x')$ is a second formal,
analytic or algebraic mapping, the {\sl composition} $x''=g(h(x))$ is
the collection of power series $(g_1(h(x)),\dots,g_{n''}(h(x)))$, which is a local
mapping from $\K^n$ to $\K^{n''}$. If $\widetilde{x}=\Phi(x)$ and
$\widetilde{x}'=\Psi(x')$ are changes of coordinates in $\K^n$ and in
$\K^{n'}$, the {\sl transformed mapping} $\widetilde{h}$ is the
mapping $\widetilde{x}'=\widetilde{h}(\widetilde{x})$, where
$\widetilde{h}(\widetilde{x}):=\Psi(h(\Phi^{-1}(\widetilde{x})))$. The
{\sl rank at the origin} of $h$ is the rank of the Jacobian matrix
$(\partial h_{i'}(0)/\partial x_i)_{1\leq i'\leq n', \, 1\leq i\leq
n}$. We denote it by ${\rm rk}_0(h)$. 
The {\sl generic rank} of $h$ is
the largest integer $e\leq \min (n,n')$ such that there exists an
$e\times e$ minor of the Jacobian matrix ${\rm Jac}\, h(x)$ which does
not vanish identically, but all $(e+1)\times (e+1)$ minors do vanish
identically. We denote it by ${\rm genrk}_\K(h)$.

Let now $d$ and $d'$ be two positive integers and let $M:
r_1(x)=\dots=r_d(x)=0$ and $M': r_1'(x')=\dots=r_{d'}'(x')=0$ be two
formal, analytic or algebraic submanifolds.  We say that {\sl $h$ maps $M$
into $M'$} if there exists a $d'\times d$ matrix $b(x)=
(b_{j',j}(x))_{1\leq j'\leq d',\, 1\leq j\leq d}$ of formal, analytic
or algebraic power series such that $r_{j'}'(h(x))\equiv \sum_{j=1}^d\,
b_{j',j}\, r_j(x)$, or in matrix notation $r'(h(x))\equiv b(x)\,
r(x)$. This definition is meaningful, since if $\widehat{r}(x)=a(x)\,
r(x)$ and $\widehat{r}'(x')= a'(x')\, r'(x')$ denote two equivalent
defining formal, analytic or algebraic defining power series for $M$
and for $M'$, then $\widehat{r}'(h(x))\equiv a'(h(x)) \, r'(h(x))\equiv
a'(h(x))\, b(x)\, [a(x)]^{-1} \ \widehat{r}(x)$,
so we have $\widehat{r}' (h(x)) \equiv \widehat{ b}(x)\, \widehat{r}
(x)$ with $\widehat{b} (x):= a'(h(x))\, b(x)\, [a(x)]^{-1}$.

\subsection*{2.1.6.~Cauchy-Riemann submanifolds of $\C^n$}
We want to study some aspects of the geometry of real submanifolds of
$\C^n$. Most often in this memoir, we shall mainly be concentrated on
the local study of pieces of submanifolds centered at one point. {\it
However, we stress that we shall never use the language of germs,
because it might sometimes be confusing}. Hence we have to work within
precise neighborhoods of central points.

Without loss of generality, we can assume that the center point is the
origin in suitable coordinates $z=x+iy\in \C^n$.  Thus, we consider a
(local) {\it real}\, $d$-codimensional submanifold $M$ of
$\C^n\cong\R^{2n}$ passing through the origin which defined by
equations $r_1(x,y)=\dots= r_d(x,y)=0$ where the differentials
$dr_1,\dots,dr_d$ are linearly independent at the origin. If
$z=x+iy\in\C$ or equivalently $(x,y)\in\R^{2n}$, we use the cube norms
$\vert x \vert = \max_{1\leq k\leq n} \, \vert x_k \vert$, $\vert y
\vert = \max_{1\leq k\leq n} \, \vert y_k \vert$ and the polydisc norm
$\vert z\vert=\max_{1\leq k\leq n} \, \vert z_k\vert$, where $\vert
z_k\vert= (x_k^2+y_k^2)^{1/2}$. For given $\nu\in\N$ with $\nu\geq 1$
and $\rho\in\R$ with $\rho>0$ we denote by $\I_\nu(\rho)$ the real
cube $(-\rho,\rho)^\nu$ in $\R^\nu$. If $\rho>0$, we denote by
$\Delta_n(\rho)=\{z\in\C^n: \, \vert z\vert < \rho\}$ the open
polydisc of radius $\rho$ centered at the origin. {\it Throughout
this memoir, we shall always work with cubes and polydiscs}. 

For useful and complete background about Cauchy-Riemann (CR for
short) structures, we refer the reader to~[1], [7]. Here, we only
give quick definitions for the purpose of being self-contained.  Let
$J$ denote the complex structure of $T\C^n$, acting on real vectors as
if it were multiplication by $\sqrt{-1}$, hence satisfying $J^2=-{\rm
Id}$. Let $M$ be a connected local real algebraic or analytic
submanifold of $\C^n$ of codimension $d$. For $p\in M$, the smallest
$J$-invariant subspace of the tangent space $T_pM$ is given by
$T_p^cM:=T_pM\cap JT_pM$ and is called the {\sl complex tangent space to
$M$ at $p$}.

\def\thedefinition{2.1.7}\begin{definition}
{\rm
The submanifold $M$ is called
\begin{itemize}
\item[{\bf (1)}]
{\sl Holomorphic} if $T_p^cM=T_pM$ at every point $p\in M$;
\item[{\bf (2)}]
{\sl Totally real} if $T_p^cM=\{0\}$ at every point $p\in M$;
\item[{\bf (3)}]
{\sl Generic} if $T_pM+JT_pM=T_p\C^n$ at every point $p\in M$;
\item[{\bf (4)}]
{\sl Cauchy-Riemann} (CR for short) if
the dimension of $T_p^cM$ is equal to a fixed constant at
every point $p\in M$.
\end{itemize} 
}
\end{definition}

In particular, holomorphic and totally real submanifolds are obviously
CR. The generic submanifolds are also CR (and in fact of minimal
possible CR dimension), because by the dimension formula $\dim_\R
(E+F)=\dim_\R E+ \dim_\R F - \dim_\R (E \cap F)$ for real vector
subspaces, we deduce from $\dim_\R(T_pM +JT_pM)=2n$ that 
$\dim_\R (T_pM\cap J T_pM)=2n-2d$, which is constant. We shall
remember that for generic submanifolds, the CR dimension is given 
by $m=n-d$.

By means of the dimension formula, we also see that if $M$ is totally
real, then $\dim_\R\, M\leq n$; also, if $M$ is generic, then $\dim_\R
\, M\geq n$. If $M$ is totally real and generic, then $\dim_\R \,
M=n$. In this case, we call $M$ {\sl maximally real}.

The two $J$-invariant spaces $T_pM\cap JT_pM$ and $T_pM+ JT_pM$ are
clearly of even real dimension.  We denote by $m_p$ the integer
${1\over 2}\dim_\R (T_pM\cap JT_pM)$ and call it the {\sl CR dimension
of $M$ at $p$}.  We denote by $c_p$ the integer $n-{1\over 2} \dim_\R(
T_pM+JT_pM)$ and call it the {\sl holomorphic codimension of $M$ at
$p$}. Of course, we have $c_p=d-n+m_p$.  In terms of these two
integers $m_p$ and $c_p$, we may rephrase the above definition as follows.

\def\thedefinition{2.1.8}\begin{definition} 
{\rm
The $d$-codimensional real submanifold $M\subset \C^n$ is
\begin{itemize}
\item[{\bf (1')}]
Holomorphic 
if $2n-d=\dim_\R\, M=2m_p$ at every point $p\in M$;
\item[{\bf (2')}]
Totally real if $m_p=0$ at every point $p\in M$;
\item[{\bf (3')}]
Generic if $m_p=n-d$ at every point $p\in M$; in this case, $m_p$ is
as small as possible and we call the integer $m:=n-d$ the {\sl CR
dimension} of $M$;
\item[{\bf (4')}]
CR if $m_p$ is equal to a fixed constant $m$ at every point $p\in M$;
in this case, we call the integer $m$ the {\sl CR dimension} of
$M$; also, it follows that the holomorphic codimension
$c_p:=d-n+m_p=d-n+m$ is constant and we call it the {\sl holomorphic
codimension} of $M$.
\end{itemize}
}
\end{definition}

For the proof of the following local graph representation theorem, we
refer to~[1], [7].

\def\thetheorem{2.1.9}\begin{theorem}
Let $M$ be a real algebraic or analytic submanifold of
$\C^n$ of codimension $d$.
\begin{itemize}
\item[{\bf (1)}]
Assume that $M$ is holomorphic, let $m={1\over 2}\, \dim_\R M$ be the
CR dimension of $M$ and let $d_1:={1\over 2}\, d$. Then for every
point $p_0\in M$, there exist local complex algebraic or analytic
coordinates $(z,w)\in\C^m\times \C^{d_1}$ vanishing at $p_0$ and there
exists $\rho_1>0$ such that $M\cap \Delta_n(\rho_1)$ is given by the
$d_1$ complex equations $w_j=0$, $j=1,\dots,d_1$ or equivalently by
the $d$ real equations ${\rm Re}\, w_j={\rm Im}\, w_j=0$,
$j=1,\dots,d_1$.
\item[{\bf (2)}]
Assume that $M$ is totally real, let $c=d-n\geq 0$ be the holomorphic
codimension of $M$ and let $d_1:=d-2c$.  Then for every point $p_0\in
M$, there exist complex algebraic or analytic coordinates
$(w,v)\in\C^{d_1}\times \C^c$ centered at $p_0$ and there exists
$\rho_1>0$ such that $M\cap \Delta_n(\rho_1)$ is given by the $d$ real
equations
\def\theequation{2.1.10}\begin{equation}
\left\{
\aligned
{\rm Im}\, w_j=
& \ 
0, \ \ \ \ \ \ \ \ \ \ \ \ \ \ \ \ \ \ \ \
j=1,\dots,d_1,\\
{\rm Re}\, v_k= 
& \
{\rm Im}\, v_k=0, \ \ \ 
\ \ \ \ \ \ \, k=1,\dots,c. 
\endaligned\right.
\end{equation}
\item[{\bf (3)}]
Assume that $M$ is generic and let $m=d-n$ be the CR dimension of $M$.
Then for every point $p_0\in M$ and
for every choice of complex affine coordinates $(z,w)\in
\C^m\times \C^d$ centered at $p_0$ such that $T_{p_0}^cM\cap
\{w=0\}=\{0\}$, there exists $\rho_1>0$ and there exist uniquely defined
real algebraic or analytic functions $\varphi_j$, $j=1,\dots,d$,
converging normally in the cube $\I_{2m+d}(2\rho_1)$ and
vanishing at the origin such that
$M\cap \Delta_n(\rho_1)$ is given by the $d$ real equations
\def\theequation{2.1.11}\begin{equation}
{\rm Im}\, w_j=\varphi_j({\rm Re}\, z,{\rm Im}\, z,
{\rm Re}\, w), \ \ \ \ \ 
j=1,\dots,d.
\end{equation}
We can in addition choose the coordinates in order that $T_0M$ is
given by the equations ${\rm Im}\, w_j=0, \, j=1,\dots, d$, in which
case we have moreover $d\varphi_j(0)=0$, for $j=1,\dots,d$.
\item[{\bf (4)}]
Assume that $M$ is CR, let $m$ be the CR dimension of $M$, let
$c=d-n+m$ be the holomorphic codimension of $M$ and let $d_1:=d-2c\geq
0$. Then for every point $p_0\in M$, there exist local complex
algebraic or analytic coordinates $(z,w,v)
\in\C^m\times \C^{d_1}\times
\C^c$ centered at $p_0$ with $T_{p_0}^cM\cap \{w=v=0\}=\{0\}$
and there exist real algebraic or analytic functions $\varphi_j$
converging normally in the cube $\I_{2m+d_1}(2\rho_1)$ for some
$\rho_1>0$ and vanishing at the origin such that $M\cap
\Delta_n(\rho_1)$ is given by the $d$ real equations
\def\theequation{2.1.12}\begin{equation}
\left\{
\aligned
{\rm Im}\, w_j= 
& \
\varphi_j({\rm Re}\, z,{\rm Im}\, z,{\rm Re}\, w), \ \ \ \ \ j=1,\dots,d_1, \\
{\rm Re}\, v_k= 
& \
{\rm Im}\, v_k=0, \ \ \ \ \ \ \ \ \ \, k=1,\dots,c. 
\endaligned\right.
\end{equation}
In particular, $M$ is contained and generic in the complex linear
subspace $(\C^m\times\C^{d_1}\times \{0\})\cap \Delta_n(\rho_1)$, 
which we call the {\sl
intrinsic complexification} of $M$. We can in addition choose the
coordinates in order that $T_0M$ is given 
by the equations ${\rm Im}\, w_j=0$, $j=1,\dots, d_1$, 
${\rm Re}\, v_k={\rm Im}\, v_k=0$, $k=1,\dots,c$, in which case
we have moreover $d\varphi_j(0)=0$ for $j=1,\dots,d_1$.
\end{itemize}
\end{theorem}

\subsection*{2.1.13.~Complex defining equations}
We now consider a real algebraic or analytic {\it generic}\,
submanifold $M$ of $\C^n$ given as in Theorem~2.1.9 by real defining
equations $v_j= \varphi_j(x,y,u)$, $j=1,\dots,d$, where
$(z,w)=(x+iy,u+iv)\in \C^m\times \C^d$.  Unless the contrary is
explicitely mentioned, {\it our generic submanifolds will always be of
positive codimension $d\geq 1$ and of positive CR dimension $m\geq
1$}. Without loss of generality, we can assume that $d\varphi_j(0)=0$
for $j=1,\dots,d$. Replacing $x$ by $(z+\bar z)/2$, $y$ by $(z-\bar
z)/2i$, $u$ by $(w+\bar w)/2$ and $v$ by $(w-\bar w)/2i$ in the
defining equations of $M$, which yields
\def\theequation{2.1.14}\begin{equation}
{w_j-\bar w_j\over 2i}=
\varphi_j\left(
{z+\bar z\over 2}, {z-\bar z\over 2i}, 
{w+\bar w\over 2}
\right), \ \ \ \ \ 
\end{equation}
for $j=1,\dots, d$, then by means of the algebraic or analytic
implicit function theorem, we can solve the $\bar w_j$ in terms of
$(\bar z, z, w)$, which yields
\def\theequation{2.1.15}\begin{equation}
\bar
w_j=\Theta_j(\bar z,z,w), 
\ \ \ \ \ \
j=1,\dots,d, 
\end{equation} 
for some complex algebraic
or analytic functions $\Theta_j$, which vanish at the origin
 and which are defined in
a neighborhood of the origin in $\C^{2m+d}$. Shrinking $\rho_1>0$ if
necessary, we can assume that the $\Theta_j$ converge normally in
$\Delta_{2m+d}(2\rho_1)$. We call these new equations {\sl complex
defining equations} for $M$ and we want to compare them with the real
defining equations.

Generally speaking, given an arbitrary series
$\Phi(t)=\sum_{\gamma\in\N^n}\, \Phi_\gamma\, t^\gamma$ with complex
coefficients $\Phi_\gamma\in\C$, we are led to define the series
$\overline{\Phi}(t):= \sum_{\gamma\in\N^n}\, \overline{\Phi}_\gamma\,
t^\gamma$ by conjugating only its complex coefficients.  With this
definition, the conjugation operator (overline) can be applied
independently over functions and over variables, as shown by the
functional equation $\overline{\Phi(t)}\equiv \overline{\Phi}(\bar
t)$. We shall use this property very frequently.

Let $(z,w)\in M$. Conjugating the defining equations of $M$, 
we get $w_j=\overline{\Theta}_j(z,\bar z,\bar w)$ and replacing the
$\bar w_l$ by their value $\Theta_l$, we get the following
equation, valuable for all $(z,w)$ belonging to $M$:
\def\theequation{2.1.16}\begin{equation}
w_j= \overline{\Theta}_j(z,\bar z, 
\Theta(\bar z, z, w)), \ \ \ \ \
j=1,\dots, d.
\end{equation}
But as we may write $(z,w)=(z,u+i\varphi(x,y,u))\in M$, where
$u=(u_1,\dots,u_d)={\rm Re}\, w$, we can replace in~\thetag{2.1.16},
which yields a power series identity in terms of the variables $(x,y,
u)$ for all $(x,y,u)\in\I_{2m+d}(\rho_1)$. As the $(2m+d)$-dimensional
real algebraic or analytic submanifold $\{(x,y, u+i\varphi(x,y,u))\}$
of $\C^{2m+d}$ is maximally real, by an application of the generic
uniqueness principle, we get the power series
identity
\def\theequation{2.1.17}\begin{equation}
w_j\equiv \overline{\Theta}_j(z,\bar z, 
\Theta(\bar z, z, w)), \ \ \ \ \
j=1,\dots, d.
\end{equation}
in $\C\{z,\bar z, \bar w\}$ or for all 
$(\bar z, z, w)\in\Delta_{2m+d}(\rho_1)$.

Conversely, suppose that these power series
identities~\thetag{2.1.17} holds. By the implicit function
theorem, there exists unique complex algebraic or analytic solutions
$\varphi_j((z+\zeta)/2,(z-\zeta)/2i,w)$, 
$j=1,\dots,d$, $z\in \C^m$, $\zeta\in \C^m$,
$w\in \C^d$ of the functional equations
\def\theequation{2.1.18}\begin{equation}
\left\{
\aligned
w_j-i\varphi_j((z+\zeta)/2,
& \,
(z-\zeta)/2i,\zeta,w)\equiv \\
&
\equiv \Theta_j(\zeta,z,w+i\varphi((z+\zeta)/2,(z-\zeta)/2i,w)),
\endaligned\right.
\end{equation}
for $j=1,\dots,d$.  We then claim that $\varphi_j(x,y,u)$ is real for
$x=(x_1,\dots,x_m)\in\R^m$, $y=(y_1,\dots,y_m)\in\R^m$ and
$u=(u_1,\dots,u_d)\in \R^d$. Indeed, by replacing first $w$ by
$u+i\varphi((z+\bar z)/2,(z-\bar z)/2i,u)$ in
the functional equations~\thetag{2.1.17}, we get
\def\theequation{2.1.19}\begin{equation}
\left\{
\aligned
u_j+i\varphi_j((z+\bar z)/2,
& \,
(z-\bar z)/2i,u)\equiv\\
&
\overline{\Theta}_j(z,\bar z,
\Theta(\bar z,z,u+i\varphi((z+\bar z)/2,(z-\bar z)/2i,u))),
\endaligned\right.
\end{equation}
for $j=1,\dots,d$. Using the implicit equations~\thetag{2.1.18} which
define $\varphi$ with $\zeta$ replaced by $\bar z$ and $w$ replaced by
$u$, we may then simplify the terms behind $\Theta$ in~\thetag{2.1.19}, wich
yields
\def\theequation{2.1.20}\begin{equation}
u_j+i\varphi_j(x,y,u)\equiv
\overline{\Theta}_j(z,\bar z, 
u-i\varphi(x,y,u)), \ \ \ \ \ 
j=1,\dots,d.
\end{equation}
Conjugating these identities, we get
\def\theequation{2.1.21}\begin{equation}
u_j-i\overline{\varphi}_j(x,y,u)\equiv
\Theta_j(\bar z,z, 
u+i\overline{\varphi}(x,y,u)), \ \ \ \ \ 
j=1,\dots,d.
\end{equation}
Comparing with the implicit equations~\thetag{2.1.18} with $\zeta$
replaced by $\bar z$ and $w$ replaced by $u$, we see that
$\varphi(x,y,u)$ and $\overline{\varphi}(x,y,u)$ are solutions of the
same implicit equations. By uniqueness in the implicit function
theorem, we obtain $\overline{\varphi}(x,y,u)\equiv \varphi(x, y,u)$,
as claimed.  Finally, the identities $u_j-i\varphi_j(x, y,u)\equiv
\Theta_j(z,\bar z,u+i\varphi(x,y,u))$ show that the set of points
$(z,w)$ satisfying $\bar w_j=\Theta_j(\bar z,z,w)$, $j=1,\dots,d$,
coincides with the real algebraic or analytic generic submanifold of
equations $v_j=\varphi_j(x,y,u)$, for $j=1,\dots d$.  In conclusion,
we have established the following important theorem which we shall use
very often in the sequel.

\def\thetheorem{2.1.22}\begin{theorem} 
Let $M$ be a real algebraic or analytic generic submanifold of
codimension $d\geq 1$ and of CR dimension $m=n-d\geq 1$ in $\C^n$.  Then
for every point $p_0\in M$, and for every choice of complex affine
coordinates $t=(z,w)\in\C^m\times \C^d$ centered at $p_0$ such that
$T_{p_0}^cM\cap \{w=0\}=\{0\}$, there exists $\rho_1>0$ and there
exist uniquely defined complex algebraic or analytic functions
$\Theta_j$, $j=1,\dots,d$, vanishing at the origin, defined and
converging normally in $\Delta_{2m+d}(2\rho_1)$ such that $M\cap
\Delta_n(\rho_1)$ is given by the $d$ complex defining equations
\def\theequation{2.1.23}\begin{equation}
\bar w_j=\Theta_j(z,\bar z, w), \ \ \ \ \ 
j=1,\dots,d,
\end{equation}
or equivalently by the $d$ conjugate complex
defining equations 
\def\theequation{2.1.24}\begin{equation}
w_j=\overline{\Theta}_j(z,\bar z,\bar w), \ \ \ \ \ 
j=1,\dots,d.
\end{equation}
Here, the vector-valued mapping $\Theta:=(\Theta_1,\dots,\Theta_d)$
satisfies the two conjugate vectorial functional equations
\def\theequation{2.1.25}\begin{equation}
\left\{
\aligned
\bar w \equiv 
& \
\Theta(\bar z, z, \overline{\Theta}(z,\bar z, \bar w)), \\
w\equiv 
& \
\overline{\Theta}(z,\bar z, 
\Theta(\bar z, z, w)).
\endaligned\right.
\end{equation}
Conversely, given a collection $\Theta=(\Theta_1,\dots,\Theta_d)$ of
complex algebraic or analytic functions vanishing at the origin,
converging normally in $\Delta_{2m+d}(2\rho_1)$ for some $\rho_1>0$
and satisfying the functional equations~\thetag{2.1.25}, then the set $M:=
\{(z,w)\in \Delta_n(\rho_1): \, \bar w_j=\Theta_j(\bar z,z,w), \,
j=1,\dots,d\}$ is a real generic submanifold of codimension $d$.
Finally, with these equations, a basis of $(0,1)$ vector
fields tangent to $M$ is given for $k=1,\dots,m$ by
\def\theequation{2.1.26}\begin{equation}
\overline{L}_k:=
{\partial \over \partial \bar z_k}+\sum_{j=1}^d\, 
{\partial \Theta_j \over \partial \bar z_k}
(\bar z, z, w)\,  
{\partial \over \partial \bar w_j}.
\end{equation}
\end{theorem}

Let $\tau=(\zeta,\xi)\in\C^m\times \C^d$ be new independent complex
variables. As in \S2.2.10 below, we define the {\sl extrinsic complexification
$\mathcal{M}$ of $M$} to be the complex analytic or algebraic
$d$-codimensional submanifold of $\C^n\times \C^n$ defined by the
equations $\xi_j-\Theta_j(\zeta,t)=0$, $j=1,\dots,d$.  The following
lemma, which is equivalent to the functional 
equations~\thetag{2.1.25}, will also be very useful.

\def\thelemma{2.1.27}\begin{lemma}
There exists an invertible $d\times d$ matrix $a(t,\tau
)$ of algebraic or analytic power series such that
\def\theequation{2.1.28}\begin{equation}
\xi-\Theta(\zeta,t)\equiv 
a(t,\tau)\,[w-\overline{\Theta}(w,\tau)].
\end{equation}
\end{lemma}

\proof
We consider the involution $\sigma$ defined by $\sigma(t,\tau):=(
\tau,t)$.  Let us say that an ideal $\mathcal{J}$ of
$\mathcal{A}_\C\{t,\tau\}$ or of $\C\{t,\tau\}$ is {\sl invariant
under the involution} $\sigma$ if for every element
$\psi(t,\tau) \in\mathcal{J}$, we have $\overline{\psi}
(\sigma(t,\tau))=\overline{\psi} (\tau,t) \in\mathcal{J}$.  Then the ideal
$\mathcal{J}$ generated by the functions
$(w_j-\xi_j) /2i-\varphi_j ((z+\zeta)/2, (z-\zeta)/2i,(w+\xi)/2)$, for
$j=1, \dots,d$, is clearly invariant under the involution $\sigma$,
since the $\varphi_j$ are real functions. By the implicit function
theorem, we can solve as above with respect to $\xi$, and we have
\def\theequation{2.1.29}\begin{equation}
\mathcal{J}=
\left<\xi_j-\Theta(\zeta,t)\right>_{1\leq j\leq d}=
\sigma_*(\mathcal{J})=
\left<w_j-\overline{\Theta}(z,\tau)\right>_{
1\leq j\leq d},
\end{equation}
from which the existence of the matrix $a(t,\tau)$ follows.
\endproof

\subsection*{2.1.30.~Existence of normal coordinates}
For notational convenience, it is often more appropriate in the real
algebraic or analytic categories to consider power series not with
respect to $(x,y)\in \R^m\times \R^m$ but with respect to $(z,\bar
z)\in\C^m\times \C^m$, which is equivalent because $(x,y)=((z+\bar z)/2, 
(z-\bar z)/2i)$ and $(z,\bar z)=(x+iy, x-iy)$.

\def\theconvention{2.1.31}\begin{convention}
{\rm 
As the local generic submanifold $M$ is algebraic, analytic or formal,
we shall write its defining equations $v_j=\varphi_j(z,\bar z,u)$, 
for $j=1,\dots,d$, in
coordinates $(z,w)=(x+iy,u+iv)\in \C^m\times \C^d$, where the
$\varphi_j$ are power series with respect to $(z,\bar z,u)$ centered
at the origin. In such a representation, we mix the real and the
complex variables.
}
\end{convention}

We can now state the existence of {\sl normal coordinates}, which are
coordinates in which the conditions~\thetag{2.1.34} below are satisfied. 
Such coordinates are not unique, as may easily be verified.

\def\thetheorem{2.1.32}\begin{theorem}
Let $M$ be as in Theorem~2.1.22. Then there exists a complex algebraic or
analytic change of coordinates $t'=h(t)$ of the special form
\def\theequation{2.1.33}\begin{equation}
z'=z, \ \ \ \ \ 
w'=g(z,w),
\end{equation}
such that the image $M':=h(M)$ has real defining equations of the form
$v_j'= \varphi_j'(z',\bar z',u')$, $j=1,\dots,d$
and complex defining equations of the form $\bar w_j'=\Theta_j'(\bar
z',z',w')$, $j=1,\dots,d$, satisfying
\def\theequation{2.1.34}\begin{equation}
\left\{
\aligned
{}
&
\varphi_j'(0,\bar z', u')\equiv 
\varphi_j'(z',0, u')\equiv 0,\\
&
\Theta_j'(0,z',w')\equiv \Theta_j'(\bar z',0,w')\equiv w_j'.
\endaligned\right. 
\end{equation}
\end{theorem}

\proof
After a linear transformation of the form~\thetag{2.1.33},
we can assume that $T_0M=\{v=0\}$, 
hence $d\varphi_j(0)=0$, $j=1,\dots,d$. Next,
the local transformation defined by $z=z'$, $w=w'+i\varphi(0,0,w')$
straightens the maximally real $d$-dimensional submanifold
\def\theequation{2.1.35}\begin{equation}
\{(0,v+i\varphi(0,0,v)): \, 
j=1,\dots,d\}\subset \{0\}\times \Delta_d(\rho_1)
\end{equation} 
to the $d$-dimensional plane $\{(0,v')\}$. Thus, we can also assume
that $\varphi_j(0,0,v)\equiv 0$, $j=1,\dots,d$. It follows that
$\Theta_j(0,0,w)\equiv w_j$. We continue the proof with the 
complex defining equations of $M$.

But before proceeding further, we remark firstly that by reality of the
power series $\varphi_j'$, we have $\varphi_j'(z',\bar z',u')\equiv
\overline{\varphi}_j'(\bar z',z',u')$, whence the collection of
relations $\varphi_j'(0,\bar z',u')\equiv 0$, $j=1,\dots,d$, is
equivalent to the collection of relations $\varphi_j'(z',0,u')\equiv
0$, $j=1,\dots,d$.  Secondly, using the functional
relations~\thetag{2.1.25}, we see immediately that the collection of
relations $\Theta_j'(0,z',w')\equiv w_j'$, $j=1,\dots,d$, is also
equivalent to the collection of relations $\Theta_j'(\bar
z',0,w')\equiv w_j'$, $j=1,\dots,d$. Thirdly, by inspecting the way
how the real and the complex defining equations of $M'$ are related
({\it see} especially the proof of Theorem~2.1.22),
we observe easily that the collection of relations in the first line
of~\thetag{2.1.34} is equivalent to the collection of relation in the
second line of~\thetag{2.1.34}. Consequently, it suffices
to find a change of coordinates of the form~\thetag{2.1.33} such that
$\Theta_j'(\bar z',0,w')\equiv w_j'$ for $j=1,\dots,d$.

We then claim that the transformation $(z',w'):=(z, \Theta(0,z,w))$ is
appropriate. Indeed, working with the extrinsic complexifications
$\mathcal{ M}$ and $\mathcal{ M}'$, we have $(z,w, \zeta,\xi)\in
\mathcal{ M}$ if and only if $(z, \Theta (0,z,w), \zeta, \overline{
\Theta} (0,\zeta,\xi))\in \mathcal{M}'$, which yields
\def\theequation{2.1.36}\begin{equation}
\overline{\Theta}(0,\zeta, \xi)=\xi'=\Theta'(
\zeta, z, \Theta(0,z,w)),
\end{equation}
again for $(z,w,\zeta,\xi)\in \mathcal{M}$. 
Replacing $\xi$ by its value $\Theta(\zeta,z,w)$ on 
$\mathcal{M}$ and setting $z=0$, we obtain the power series identity
\def\theequation{2.1.37}\begin{equation}
\overline{\Theta}(0,\zeta,\Theta(\zeta,0,w))\equiv
\Theta'(\zeta,0,\Theta(0,0,w)).
\end{equation}
But we remember that the functional equations~\thetag{2.1.25} hold,
which enables us to simplify the left hand side and we remember that
we have already the relation $\Theta(0,0,w)\equiv w$, which enables us to
simplify the right hand side and we obtain the desired power series
identity
\def\theequation{2.1.38}\begin{equation}
w\equiv \Theta'(\zeta,0,w).
\end{equation}
This completes the proof of Theorem~2.1.32.
\endproof

\subsection*{2.1.39.~The formal case}
All the previous computations are meaningful in the purely formal
case. Especially, Theorems~2.1.22 and~2.1.32 hold true in the formal
category.

\subsection*{2.1.40.~Conclusion}
As we shall observe and confirm in the sequel, the representation of
$M$ by complex defining equations is substantially more convenient and
more tractable than the representation by real defining equations.
We remind that, unless the contrary is explicitely mentioned,
{\it our generic submanifolds will always be of positive
codimension $d\geq 1$ and of positive CR dimension $m\geq 1$}.

\def\thenotation{2.1.41}\begin{notation}
{\rm 
Throughout this memoir, we shall fix the following notations:
\begin{itemize}
\item[{\bf (1)}]
The generic submanifold $M$ of $\C^n$ will be of codimension $d\geq 1$
and of CR dimension $m=n-d\geq 1$. The coordinates on $\C^n$ will
be denoted by 
\def\theequation{2.1.42}\begin{equation}
t=(t_1,\dots,t_n)=(z_1,\dots,z_k,w_1,\dots,w_d)\in
\C^m\times \C^ d=\C^n
\end{equation} 
and the complex defining equations of $M$ by 
$\bar w_j=\Theta_j(\bar z,z,w)$, $j=1,\dots,d$. 
\item[{\bf (2)}]
The index $i\in \N$ will run from $1$ to $n$, namely $i=1,\dots,n$,
for instance in the denotation of a vector field $L=\sum_{i=1}^n\, a_i(t)\,
\partial_{t_i}$. The letter $i$ will also be used to denote
$\sqrt{-1}$.
\item[{\bf (3)}]
The index $j\in\N$ will run from $1$ to $d$, namely $j=1,\dots,d$.
\item[{\bf (4)}]
The index $k\in\N$ will run from $1$ to $m$, namely $k=1,\dots,m$.
The letter $k$ will also often be used to
denote another integer varying in $\N$.
\item[{\bf (5)}]
Also, 
we denote $z_k=x_k+iy_k$, $w_j=u_j+iv_j$, 
$x=(x_1,\dots,x_m)$, $y=(y_1,\dots,y_m)$,
$u=(u_1,\dots,u_d)$ and $v=(v_1,\dots,v_d)$.
\end{itemize}
}
\end{notation}

\subsection*{2.1.43.~Precise definition of local generic submanifolds}
Throughout this memoir, we shall often need to localize our geometric
constructions. It is therefore necessary to formulate once for all
times a firm and precise choice of local representation.

\def\thedefinition{2.1.44}\begin{definition}
{\rm 
A {\sl local generic submanifold $M$ of $\C^n$} of codimension
$d\geq 1$ and of CR dimension $m=n-d\geq 1$ is defined
in coordinates $t=(z,w)=(x+iy,u+iv)\in \C^m\times \C^d$
vanishing at a point $p_0\in M$ as a graph
\def\theequation{2.1.45}\begin{equation}
M=\{(z,w)\in \Delta_n(\rho_1): \, 
v_j=\varphi_j(x,y,u), \, j=1,\dots,d\},
\end{equation}
where the functions $\varphi_j$ are real algebraic or analytic
for $\vert (x,y,u)\vert < 2\rho_1$.  We also require that for all
$\rho$ with $0\leq \rho \leq \rho_1$, we have $\vert \varphi(x,y,u)
\vert < \rho$ if $\vert (x,y,u) \vert < \rho$, namely $M$ is a ``good
graph'', as shown in {\sc Figure~2.1.47} below. Of course, after
peharps shrinking $\rho_1>0$, this condition is automatically
satisfied if we adjust the coordinates in order that $T_0M=\{v=0\}$.
In fact, we shall often
prefer the representation by complex defining equations
\def\theequation{2.1.46}\begin{equation}
M=\{(z,w)\in\Delta_n(\rho_1): \, 
\bar w_j=\Theta_j(\bar z, z, w), \, 
j=1,\dots,d\},
\end{equation}
where $M$ is again a good 
graph and the $\Theta_j$ converge normally 
for $\vert (\bar z, z, w)\vert < 2\rho_1$.
}
\end{definition}

\begin{center}
\input admitted.pstex_t
\end{center}

The {\sl process of localization} consists in choosing subsequent
smaller polydiscs $\Delta_n(\rho_2)$, $\Delta_n(\rho_3)$,
$\Delta_n(\rho_4)$, $\dots$, with $0 < \cdots < \rho_4 < \rho_3 <
\rho_2 < \rho_1$, where the choice of smaller radii $\rho_2, \rho_3,
\rho_4$ depends on the construction of further geometric objects
related to $M$. We say that $p_0$ is the {\sl central point}. This
process may be illustrated symbolically as follows:

\begin{center}
\input local-gen-subm.pstex_t
\end{center}

If a globally defined connected generic submanifold $M$ is given, for
every point $p_0\in M$, we can obviously {\sl localize $M$ at $p_0$}
by choosing complex affine coordinates vanishing at $p_0$ such that
$M$ in a neighborhood of $p_0$ is represented as in Definition~2.1.44.

\section*{\S2.2.~Segre varieties and extrinsic complexification}

\subsection*{2.2.1.~Reality condition}
Although we shall mainly work with the complex defining equations, it
will be useful on occasion to work with arbitrary real defining
equations. So we consider an arbitrary set of $d$ real defining
equations for $M$ which we denote by $\rho_j(t,\bar t)=0$,
$j=1,\dots,d$, for instance $\rho_j(t,\bar t):= v_j
-\varphi_j(z,\bar z, u)$, with $\rho_j(0)=0$. Here, we
assume that the complex differentials $\partial \rho_1, \dots,
\partial \rho_d$ are linearly independent at the origin, so that $M$
is generic. By reality of the $\rho_j$, we have $\rho_j(t,\bar
t)\equiv \overline{\rho_j(t,\bar t)}$. Developping the $\rho_j$ in
power series, we may write $\rho_j(t,\bar t)\equiv \sum_{\mu,\,
\nu\in\N^n}\, \rho_{j,\mu,\nu}\, t^\mu\, \bar t^\nu$, with
$\rho_{j,\mu,\nu}\in \C$.  From this functional equation, we deduce
that $\rho_{j,\mu,\nu}=\overline{\rho_{j,\nu,\mu}}$ for all
$j,\mu,\nu$. Conversely, any such converging power series with complex
coefficients satisfying $\rho_{j,\mu,\nu}=\overline{\rho_{j,\nu,\mu}}$
takes only real values.
As an application, we may write $\overline{\rho_j(t,\bar t)}\equiv
\bar \rho_j(\bar t,t)$, so the reality condition on $\rho_j$ is simply
\def\theequation{2.2.2}\begin{equation}
\rho_j(t,\bar t)\equiv \bar \rho_j(t,\bar t), \ \ \ \ \ \
j=1,\dots,d.
\end{equation}
Now, let $\tau\in \C^n$ be a new independent variable corresponding to
the extrinsic complexification of the variable $\bar t$. We shall
write symbolically $\tau:=(\bar t)^c$, where the letter ``c'' stands
for the word ``complexified''. As~\thetag{2.2.2} is equivalent to
$\rho_{j,\mu,\nu}=\overline{\rho_{j,\nu,\mu}}$, we observe that the
complexified series $\rho_j(t,\tau)$ satisfy the important symmetry
functional equation
\def\theequation{2.2.3}\begin{equation}
\rho_j(t,\tau)\equiv \bar \rho_j(\tau,t),  \ \ \ \ \ \
j=1,\dots,d.
\end{equation}
which is simply obtained by replacing $\bar t$ by $\tau$ in~\thetag{2.2.2}.
We can summarize these observations.

\def\thelemma{2.2.4}\begin{lemma}
As the defining functions $\rho_j(t,\bar t)=
\sum_{\mu,\nu\in\N^n}\, \rho_{j,\mu,\nu}\, 
t^\mu\, \bar t^\nu$, $j=1,\dots,d$ are real power series, we have
$\rho_{j,\mu,\nu}=\overline{\rho_{j,\nu,\mu}}$ for all
$j,\mu,\nu$ and
\begin{itemize}
\item[{\bf (1)}]
$\rho_j(t,\tau)\equiv \bar \rho_j(\tau,t)$.
\item[{\bf (2)}]
$\rho_j(t,\tau)=0$ if and only if $\rho_j(\bar \tau,\bar t)=0$.
\end{itemize}
\end{lemma}

Property {\bf (2)} follows trivially from {\bf (1)} and will be useful
later.

\subsection*{2.2.5.~Classical Segre varieties and conjugate Segre varieties}
Let $\rho_1>0$ such that the $\rho_j$ converge normally in
$\Delta_{2n}(2\rho_1)$ and consider the zero-set
$M:=\{t\in\Delta_n(\rho_1): \, \rho_j(t,\bar t)=0\}$. Let
$\rho_j'(t,\bar t)=0$ be an another choice of defining equations for
the same generic submanifold. It follows that there exists an
invertible $d\times d$ matrix $a(t,\bar t)$ of real power series (of
the same regularity as $M$, namely algebraic
or analytic) such that $\rho'(t,\bar t)\equiv a(t,\bar
t)\, \rho(t,\bar t)$. With this relation, we observe easily that the
(classical) {\sl Segre variety} associated to a point
$p\in\Delta_n(\rho_1)$ with coordinates $t_p=(t_{1p},\dots,t_{np})\in\C^n$,
defined as in~[34] by
\def\theequation{2.2.6}\begin{equation}
S_{\bar t_p}:=\{t\in\Delta_n(\rho_1): \, 
\rho_j(t,\bar t_p)=0, \, 
j=1,\dots,d\},
\end{equation}
does not depend on the choice of defining equations for $M$, namely we
also have $S_{\bar t_p}=\{t\in\Delta_n(\rho_1): \, \rho'(t,\bar
t_p)=0\}$.  In the litterature, this Segre variety is usually denoted
by $Q_p$, {\it cf.}~[1], [8], [9], [12], [22], [23], [27], [30], [31], [34],
[35], [36].  Here, we choose instead the letter ``$S$'', because it is
the initial of the name Segre. More importantly, we stress the
notation $S_{\bar t_p}$ or $S_{\bar p}$, and not $S_p$, {\it with the
bar of complex conjugation over $t_p$}, as in the expresson
$\rho_j(t,\bar t_p)$.

In fact, for reasons of symmetry, we are also led to define the
{\sl conjugate Segre variety} by
\def\theequation{2.2.7}\begin{equation}
\overline{S}_{t_p}:=
\{\bar t\in \Delta_n(\rho_1): \, 
\rho_j(t_p,\bar t)=0, \, 
j=1,\dots,d\}.
\end{equation}
To the author's knowledge, conjugate Segre varieties are not
considered in the literature.  As a matter of fact, if for an
arbitrary subset $E\subset \C^n$ we define the set of conjugate points
of $E$ by $\overline{E}:=\{\bar t: \, t\in E\}$, it follows that
$\overline{S}_{t_p}$ is just the set of conjugate points of $S_{\bar
t_p}$, as the reader may verify thanks to Lemma~2.2.4.  It follows that
we can write
\def\theequation{2.2.8}\begin{equation}
\overline{S}_{t_p}
=\overline{S}_{\bar{\bar t}_p}=
\overline{S_{\bar t_p}} \ \ \ \ \ 
{\rm and} \ \ \ \ \
S_{\bar t_p}=\overline{\overline{S}}_{\bar t_p}=
\overline{\overline{S}_{t_p}}.
\end{equation}
with the complex conjugation operator acting separately as an
involution over the letter $S$ and over its argument $t_p$.  Finally,
we would like to observe that the third (tempting) definition $\{t\in
\Delta_n(\rho_1): \, \rho_j(t_p,\bar t)=0, \, 1,\dots,d\}$ (instead
of~\thetag{2.2.7}) does not provide us with the correct definition of
conjugate Segre variety, because Lemma~2.2.4 implies that this set
coincides in fact with $S_{\bar t_p}$.

We claim that both Segre and conjugate Segre varieties are biholomorphic
invariants of $M$. Indeed, let $t'=h(t)$ be a local biholomorphic
change of coordinates and denote by $t=h'(t')$ its inverse, by
$M':=h(M)$ and by $\rho_j'(t',\bar t'):= \rho_j(h'(t'),\bar h'(\bar
t'))$ the defining equations of $M'$, for
$j=1,\dots,d$. Since $h$ maps $M$ into $M'$, according to \S2.1.5,
there exists an invertible $d\times d$ matrix $a(t,\bar t)$ of power
series such that $\rho'(h(t),\bar h(\bar t))\equiv a(t,\bar t)\,
\rho(t,\bar t)$. From this relation, it follows easily that $h(S_{\bar
t_p})=S_{\bar h(\bar t_p)}'$ and $h(\overline{S}_{ t_p})=
\overline{S}_{ \bar h(\bar t_p)}'$, which proves the claim.
Finally, we collect some classical properties.

\def\thelemma{2.2.9}\begin{lemma}
The following four properties are satisfied:
\begin{itemize}
\item[{\bf (1)}]
$q\in S_{\bar t_p}$ if and only if $p\in S_{\bar t_q}$.
\item[{\bf (2)}]
$p\in S_{\bar t_p}$ if and only if $p\in M$.
\item[{\bf (3)}]
$\bar q\in \overline{S}_{t_p}$ if and only if $\bar p\in \overline{S}_{t_q}$.
\item[{\bf (4)}]
$\bar p\in \overline{S}_{t_p}$ if and only if $p\in M$.
\end{itemize}
\end{lemma}

\proof
Indeed, applying Lemma~2.2.4 {\bf (2)}, we have $\rho(t_q,\bar t_p)=0$
if and only if $\rho(t_p,\bar t_q)=0$, which yields {\bf (1)} and
{\bf (3)}. Also, we have $\rho(t_p,\bar t_p)=0$ if and
only if $t_p\in M$, which yields {\bf (2)} and
{\bf (4)}.
\endproof

\subsection*{2.2.10.~Extrinsic complexification}
Now, let $\zeta\in\C^m$ and $\xi\in\C^d$ denote some new independent
coordinates corresponding to the complexification of the variables
$\bar z$ and $\bar w$, which we denote symbolically by $\zeta:=(\bar
z)^c$ and $\xi:= (\bar w)^c$, where the letter ``c'' stands for the
word ``complexified''. We also write $\tau:=(\bar t)^c$, so
$\tau=(\zeta,\xi)\in\C^n$. The {\sl extrinsic complexification}\,
$\mathcal{M}:=(M)^c$ of $M$ is the complex $d$-codimensional
submanifold defined precisely by
\def\theequation{2.2.11}\begin{equation}
\mathcal{M}:=\{(z,w,\zeta,\xi)\in\Delta_n(\rho_1)\times
\Delta_n(\rho_1):\, 
\xi_j=\Theta_j(\zeta,z,w), \ 
j=1,\dots,d\}.
\end{equation} 
Let $\sigma$ denote the antiholomorphic involution defined by
$\sigma(t,\tau):=(\bar \tau,\bar t)$.  Since by Lemma~2.1.27, there
exists an invertible $d\times d$ matrix $a(t,\tau)$ of power series
such that $w-\overline{\Theta}(z,\zeta,\xi)\equiv
a(t,\tau)[\xi-\Theta(\zeta,z,w)]$, we see that $\sigma$ maps
$\mathcal{M}$ bi-antiholomorphically onto $\mathcal{M}$.  In this
chapter, we shall essentially deal with $\mathcal{M}$ instead of
dealing with $M$. In fact, $M$ clearly imbeds in $\mathcal{M}$ as the
intersection of $\mathcal{M}$ with the antiholomorphic diagonal
defined by $\underline{\Lambda}:=\{(t,\tau)\in\C^n\times \C^n: \,
\tau=\bar t\}$. Also, we shall very frequently use the fact that
$\mathcal{M}$ can be represented by the following two equivalent
families of $d$ complex defining equations:
\def\theequation{2.2.12}\begin{equation}
\mathcal{M}: \ \ 
w_j=\overline{\Theta}_j(z,\zeta,\xi),  \ \ \ \ \ j=1,\dots,d, \ \ \ \ \ 
{\rm or} \ \ \ \ \ 
\xi_j=\Theta_j(\zeta,z,w), \ \ \ \ \ j=1,\dots,d.
\end{equation}

\section*{\S2.3.~Complexified Segre varieties and 
complexified CR vector fields}

\subsection*{2.3.1.~Complexified Segre varieties}
Next, for $\tau_p\in\Delta_n(\rho_1)$
fixed, we define the associated {\sl complexified Segre variety} by 
\def\theequation{2.3.2}\begin{equation}
\mathcal{S}_{\tau_p}:=  
\{(t,\tau)\in \Delta_{2n}(\rho_1): \, 
\tau=\tau_p, \, 
w_j=\overline{\Theta}_j(z,\tau_p), \, 
j=1,\dots,d\}.
\end{equation}
We shall write symbolically $\mathcal{S}_{\tau_p}= (S_{\bar
t_p})^c$. Clearly, $\mathcal{S}_{\tau_p}$ is an $m$-dimensional
submanifold contained in $\mathcal{M}$ and it coincides in fact with
the intersection of $\mathcal{M}$ with the horizontal slice
$\{(t,\tau): \tau=\tau_p\}$.  Analogously, for
$t_p\in\Delta_n(\rho_1)$ fixed, we define the associated {\sl
conjugate complexified Segre variety} by
\def\theequation{2.3.3}\begin{equation}
\underline{\mathcal{S}}_{t_p}:=
\{(t,\tau)\in\Delta_{2n}(\rho_1): \, 
t=t_p, \, 
\xi_j=\Theta_j(\zeta,t_p), \, 
j=1,\dots,d\}.
\end{equation}
Clearly again, $\underline{\mathcal{S}}_{t_p}$ is an $m$-dimensional
submanifold contained in $\mathcal{M}$ and it coincides in fact with
the intersection of $\mathcal{M}$ with the vertical slice $\{(t,\tau):
t=t_p\}$.

\begin{center}
\input segre-blowing.pstex_t
\end{center}

It is very important to notice that the ambient Segre varieties
$S_{\bar t_p}$ and $\overline{S}_{t_p}$ are extrinsic to $M$: they lie
in general outside $M$, even in the Levi-flat case.  Moreover, the
union of $\cup_{p\in\Delta_n(\rho_1)} \, S_{\bar p}$ never makes a
foliation by $m$-dimensional submanifolds. These assertions may easily
be checked by inspecting the Levi-flat hyperplane $\{{\rm Im}\, w=0\}$
in $\C^n$ and the Heisenberg sphere ${\rm Im}\, w=\vert
z_1\vert^2+\cdots+ \vert z_{n-1}\vert^2$ in $\C^n$. Fortunately, by
the strange miracle of extrinsic complexification, {\it we blow-up the
two unions $\cup_p\, S_{\bar p}$ and $\cup_p\, \overline{S}_p$ in a
double foliation of $\mathcal{M}$ by complex $m$-dimensio\-nal Segre
varieties} (as explained in Theorem~2.3.9 below). 
This geometric observation is of utmost importance, is
illustrated symbolically in {\sc Figure}~2.3.4 and will be explained
more closely in the next subparagraphs.

\subsection*{2.3.5.~Complexified CR vector fields}
We consider the following ``natural'' basis of $(1,0)$ vector fields
tangent to $M$:
\def\theequation{2.3.6}\begin{equation}
L_k:=
{\partial\over \partial z_k}+
\sum_{j=1}^d\, 
{\partial \overline{\Theta}_j\over
\partial z_k}(z,\bar z,\bar w)\, 
{\partial \over \partial w_j}, 
\ \ \ \ \ k=1,\dots,m.
\end{equation}
One verifies immediately that $L_k(w_j-\overline{\Theta}_j
(z,\bar z,\bar w))\equiv 0$  for $k=1,\dots,m$ and $j=1,\dots,d$.
We also consider the conjugates of these vector fields, 
which form a basis of the $(0,1)$ vector fields tangent to $M$:
\def\theequation{2.3.7}\begin{equation}
\overline{L}_k:=
{\partial \over \partial \bar z_k}+
\sum_{j=1}^d\,  
{\partial \Theta_j\over
\partial \bar z_k}(\bar z,z,w)\, 
{\partial \over \partial \bar w_j}, 
\ \ \ \ \ k=1,\dots,m.
\end{equation}
Again obviously, we verify that $\overline{L}_k(\bar w_j-\Theta_j
(\bar z,z,w))\equiv 0$ for $k=1,\dots,m$ and $j=1,\dots,d$.
Of course, this second system of relations is the conjugate 
of the first.

By complexification, the vector fields behave as follows:
we write $[\chi(t,\bar{t})]^c=\chi(t,\tau)$, if $\chi(t,\bar t)$ is 
a real analytic function of $(t,\bar t)$ 
and $\left[\sum_{j=1}^{n} \,a_j(t,\bar{t}) \,\partial/\partial t_j+
\sum_{j=1}^n \, b_j(t,\bar{t}) \,\partial/\partial \bar{t}_j\right]^c:=
\sum_{j=1}^n \, a_j(t,\tau)\,
\partial/\partial t_j+ \sum_{j=1}^n \, b_j(t,\tau)\,
\partial/\partial\tau_j$. It follows that
$(L\chi)^c=L^c\chi^c$.

Consequently, we can complexify the pair of conjugate
generating families of CR vector fields tangent to $M$ given 
by~\thetag{2.3.6} and~\thetag{2.3.7}, namely the vector fields
$L_1,\dots,L_m$ and their conjugates
$\overline{L}_1,\dots,\overline{L}_m$ above. Their complexification
yields a pair of collections of $m$ vector fields defined 
explicitely over
$\Delta_n(\rho_1)\times \Delta_n(\rho_1)$ by
\def\theequation{2.3.8}\begin{equation}
\left\{
\aligned
\mathcal{L}_k:= & \ 
{\partial \over\partial z_k}+
\sum_{j=1}^d\,
{\partial\overline{\Theta}_j\over \partial z_k}
(z,\zeta,\xi)\,
\, {\partial\over\partial w_j}, \ \ \ \ 
k=1,\dots,m,\\
\underline{\mathcal{L}}_k:= & \
{\partial\over\partial \zeta_k}+\sum_{j=1}^d\,
{\partial \Theta_j\over \partial \zeta_k}(\zeta,z,w)\,
{\partial\over\partial \xi_j}, \ \ \ \
k=1,\dots,m.
\endaligned\right.
\end{equation}
We write $\mathcal{L}_k=(L_k)^c$ and $\underline{\mathcal{L}}_k=
(\overline{L}_k)^c$.  The reader may check directly that
$\mathcal{L}_k(w_j- \overline{\Theta}_j (z,\zeta,\xi)) \equiv 0$ (this
relation also holds by complexication), which shows that the vector
fields $\mathcal{L}_k$ are tangent to $\mathcal{M}$.  Similarly,
$\underline{\mathcal{L}}_k (\xi_j- \Theta_j(\zeta,z,w)) \equiv 0$, so
the vector fields $\underline{\mathcal{L}}_k$ are also tangent to
$\mathcal{M}$. Of course, all of this
is obvious, but we prefer to start 
slowly. Furthermore, we may check the commutation relations
$[L_{k_1},L_{k_2}]=0$, $[\overline{L}_{k_1},\overline{L}_{k_2}]=0$,
$[\mathcal{L}_{k_1},\mathcal{L}_{k_2}]=0$ and
$[\underline{\mathcal{L}}_{k_1},\underline{\mathcal{L}}_{k_2}]=0$ for
all $k_1, k_2=1,\dots,m$. By the theorem of Frobenius, it follows that
the two $m$-dimensional distributions spanned by the two collections
of $m$ vector fields $\{\mathcal{L}_k\}_{1\leq k\leq m}$ and
$\{\underline{\mathcal{L}}_k\}_{1\leq k\leq m}$ have the integral
manifold property. This is not astonishing, due to the fact that the
vector fields $\mathcal{L}_k$ are just the vector fields tangent to
the intersection of $\mathcal{M}$ with the sets $\{\tau=\tau_p=ct.\}$,
which are the $m$-dimensional complexified Segre varieties
$\mathcal{S}_{\tau_p}$ already defined above. Similarly, the
$\underline{\mathcal{L}}_k$ have the conjugate complexified Segre
varieties $\underline{\mathcal{S}}_{t_p}$ as integral manifolds.
Hence in fact, we do not have to appeal to the theorem  of
Frobenius.

All the geometric observations which we have done so far may be gathered in 
the following statement just below. We shall frequently use the abbreviations
$\mathcal{L}=\{\mathcal{L}_k\}_{1\leq k\leq m}$ and
$\underline{\mathcal{L}}=\{\underline{\mathcal{L}}_k\}_{
1\leq k\leq m}$. We denote by $\pi_t: (t,\tau)\mapsto t$ and
$\pi_\tau: (t,\tau)\mapsto \tau$ the two canonical projections.

\def\thetheorem{2.3.9}\begin{theorem}
Let $\mathcal{ M}=(M)^c$ be as above and let $\mathcal{ L}_k$, $k=1,
\dots,m$ be a basis of complexified $(1,0)$ vector fields tangent to
$M$ and let $\underline{ \mathcal{L}}_k$, $k=1,\dots,m$, be their
complexified conjugates. Recall that $\{\mathcal{ L}_k\}_{ 1\leq k\leq
m}$ and $\{ \underline{ \mathcal{ L}}_k\}_{ 1\leq k\leq m}$ are
Frobenius-integrable. Then the following four properties hold true:
\begin{itemize}
\item[{\bf (1)}]
$\mathcal{L}$ and $\underline{\mathcal{L}}$ 
induce naturally two local flow
foliations $\mathcal{F}_{\mathcal{L}}$ and 
${\mathcal{F}}_{\underline{\mathcal{L}}}$ of
$\mathcal{M}$.
\item[{\bf (2)}] 
If $\sigma(t,\tau):=(\bar \tau,\bar t)$, then
$\sigma(\mathcal{F}_{\mathcal{L}})=\mathcal{F}_{
\underline{\mathcal{L}}}$ and 
their two leaves passing through a point $p^c=(t_p,\bar t_p)\in 
\C^n\times \C^n$
satisfy $\mathcal{F}_{\mathcal{L}}(p^c)\cap 
\mathcal{F}_{\underline{\mathcal{
L}}}(p^c)=p^c$.
\item[{\bf (3)}]
The fibers of the projections $\pi_t$ and $\pi_{\tau}$ also coincide
with the leaves of the flow foliations 
$\mathcal{F}_{\underline{\mathcal{L}}}$ and $\mathcal{
F}_{\mathcal{L}}$, respectively.
\item[{\bf (4)}]
The leaves of the foliation 
$\mathcal{F}_{\mathcal{L}}$ are the Segre varieties
$\mathcal{S}_{\tau_p}$ and the leaves of the foliation $\mathcal{
F}_{\underline{\mathcal{L}}}$ are the conjugate Segre varieties
$\underline{\mathcal{S}}_{t_p}$\text{\rm :}
\def\theequation{2.3.10}\begin{equation}
\mathcal{F}_{\mathcal{L}}=
\bigcup_{\tau_p\in\Delta_n(\rho_1)}
{\mathcal{S}}_{\tau_p} \ \ \ \ \text{\rm and} \ \ \ \ \
\mathcal{F}_{\underline{\mathcal{L}}}=
\bigcup_{t_p\in\Delta_n(\rho_1)}
\underline{\mathcal{S}}_{t_p}.
\end{equation}
\end{itemize}

\noindent
In other words, the leaves of these two flow foliations are the two families
of complexified $($conjugate$)$ Segre varieties. In symbolic 
representation, for these two foliations, we have the 
correspondence\text{\rm \,:}
\def\theequation{2.3.11}\begin{equation}
\text{\it CR-flow foliations of} \ \ \mathcal{M}
\ \ \Longleftrightarrow \ \ \text{\it Foliations by complexified
Segre varieties}. 
\end{equation}
\end{theorem}

\subsection*{2.3.12.~Conclusion}
The following symbolic picture summarizes this geometrical
theorem. However, we warn the reader that the codimension $d\geq 1$ of
the union of the two foliations $\mathcal{F}_{\mathcal{L}}$ and
$\mathcal{F}_{\underline{\mathcal{L}}}$ in $\mathcal{M}$ is not
rendered visible in this two-dimensional figure. A three-dimensional
{\sc Figure~2.3.3} will be provided below.

\begin{center}
\input complexificationbig.pstex_t
\end{center}

\section*{2.4.~Multiple flows and Segre chains}

\subsection*{2.4.1.~Pair of complex flows}
Now, we introduce the ``multiple'' flows of the two collections of
conjugate vector fields $(\mathcal{L}_k)_{1\leq k\leq m}$ and
$(\underline{\mathcal{L}}_k)_{1\leq k\leq m}$. This multiple flow will
be used frequently throughout the next chapters of Part~II.
Precisely, for an arbitrary point $p=( w_p, z_p, \zeta_p, \xi_p) \in
\mathcal{M}$ and for an arbitrary complex ``multitime'' parameter
$z_1=( z_{1,1}, \dots, z_{1,m}) \in \C^m$, we define
\def\theequation{2.4.2}\begin{equation}
\left\{
\aligned
\mathcal{L}_{z_1}(z_p,w_p,\zeta_p,\xi_p) & \ :=
\exp(z_1\mathcal{L})(p):=
\exp(z_{1,1}\mathcal{L}_1(\cdots(\exp(z_{1,m}\mathcal{L}_m(
p)))\cdots)):=\\
& \
:=\left(z_p+z_1,\, \overline{\Theta}(z_p+z_1,\zeta_p,\xi_p),\,
\zeta_p,\,\xi_p\right).
\endaligned\right.
\end{equation}
With this formal definition, there exists a maximal connected open
subset $\Omega$ of $\mathcal{M}\times \C^m$ containing
$\mathcal{M}\times\{0\}$ such that
$\mathcal{L}_{z_1}(p)\in\mathcal{M}$ for all
$(z_1,p)\in\Omega$. Analogously, for $(\zeta_1,p)$ running in a
similar open subset $\underline{\Omega}$, we may define the map
\def\theequation{2.4.3}\begin{equation}
\underline{\mathcal{L}}_{\zeta_1}(z_p,w_p,\zeta_p,\xi_p):=
\left(z_p,\,w_p,\,\zeta_p+\zeta_1,\,\Theta(\zeta_p+\zeta_1,z_p,w_p)\right).
\end{equation} 
Of course, the
two mappings~\thetag{2.4.2} and~\thetag{2.4.3} are of the same
regularity as $M$, namely they are algebraic
or analytic.

\subsection*{2.4.4.~Segre chains}
Now, let us start from the point $p$ being the origin and let us move
alternately in the (horizontal) direction of $\mathcal{ F}_{
\mathcal{L}}$ (namely the direction of $\mathcal{ S}$) and in the
(vertical) direction of $\mathcal{ F}_{ \underline{ \mathcal{L}}}$
(namely the direction of $\underline{\mathcal{ S}}$). More precisely,
we consider the two maps $\Gamma_1 (z_1) := \mathcal{ L}_{ z_1}(0)$
and $\underline{ \Gamma}_1(z_1):= \underline{ \mathcal{ L}}_{
z_1}(0)$, where $z_1\in \C^m$.  Next, we start from these endpoints
and we move in the other direction. More precisely, we consider the
two maps
\def\theequation{2.4.5}\begin{equation}
\Gamma_2(z_1,z_2):=\underline{\mathcal{L}}_{z_2}(\mathcal{L}_{z_1}(0)), 
\ \ \ \ \
\underline{\Gamma}_2(z_1,z_2):=
\mathcal{L}_{z_2}(\underline{\mathcal{L}}_{z_1}(0)),
\end{equation}
where $z_1,\,z_2\in\C^m$. Also, we define $\Gamma_3(z_1,z_2,z_3):=
\mathcal{L}_{z_3}(\underline{\mathcal{L}}_{z_2}(\mathcal{L}_{z_1}(0)))$,
{\it etc.} For the sake of concreteness, let us exhibit the complete
expressions of $\Gamma_1$, $\Gamma_2$ and $\Gamma_3$, which follows by
a repeated application of formulas~\thetag{2.4.2} and~\thetag{2.4.3}:
\def\theequation{2.4.6}\begin{equation}
\left\{
\aligned
{}
&
\mathcal{L}_{z_1}(0)=\left(z_1, \, \overline{\Theta}
(z_1,0,0), \, 0, \, 0\right).\\
&
\underline{\mathcal{L}}_{z_2}(
\mathcal{L}_{z_1}(0))=\left(z_1, \, \overline{\Theta}(z_1,0,0), \, 
z_2, \, \Theta(z_2,z_1,\overline{\Theta}(z_1,0,0))\right).\\
&
\mathcal{L}_{z_3}(
\underline{\mathcal{L}}_{z_2}(
\mathcal{L}_{z_1}(0))))=\left(z_1+z_3, \,
\overline{\Theta}(z_1+z_3,z_2,\Theta(z_2,z_1,\overline{\Theta}
(z_1,0,0))),\,  z_2,\, \Theta(z_2,z_1,\overline{\Theta}(z_1,0,0) )\right).\\
\endaligned\right.
\end{equation}
By induction, for every integer $k\in\N$ with 
$k\geq 1$, we obtain two maps $\Gamma_k(z_1,\dots,z_k)$ and
$\underline{\Gamma}_k(z_1,\dots,z_k)$, where 
$z_1,\dots,z_k\in\C^m$. 
Clearly, there are precise combinatorial formulas
generalizing~\thetag{2.4.6}.
In the sequel, we shall often
use the notation $z_{(k)}:=(z_1,\dots,z_k)\in\C^{mk}$.
We shall call the map $\Gamma_k$ the {\sl $k$-th Segre chain} and
the map $\underline{\Gamma}_k$ the {\sl conjugate
$k$-th Segre chain}. Since
$\Gamma_k(0)=\underline{\Gamma}_k(0)=0$, for every $k\in\N_*$, there
exists a sufficiently small open polydisc $\Delta_{mk}(\delta_k)$
centered at the origin in $\C^{mk}$ with $\delta_k>0$ such that
$\Gamma_k(z_{(k)})$ and $\underline{\Gamma}_k(z_{(k)})$ belong to
$\mathcal{M}$ for all $z_{(k)} \in \Delta_{mk}(\delta_k)$.

We also exhibit a simple link between the maps $\Gamma_k$ and
$\underline{\Gamma}_k$. Let $\sigma$ be the antiholomorphic involution
defined by $\sigma(t,\tau):=(\bar\tau,\bar t)$.  Since
$w=\overline{\Theta}(z,\zeta,\xi)$ if and only if
$\xi=\Theta(\zeta,z,w)$, this involution maps $\mathcal{M}$ onto
$\mathcal{M}$ and it also fixes the antidiagonal $\underline{\Lambda}$
pointwise.  Using the definitions~\thetag{2.4.2} and~\thetag{2.4.3}, we
see readily that $\sigma(\mathcal{L}_{z_1}(0)) =
\underline{\mathcal{L}}_{\bar z_1}(0)$.  It follows generally that $\sigma(
\Gamma_k(z_{(k)}))= \underline{\Gamma}_k( \overline{z_{(k)}})$.
To give a concrete illustration, we may compute the explicit
expressions of $\underline{\Gamma}_1$, $\underline{\Gamma}_2$ 
and $\underline{\Gamma}_3$ and compare with~\thetag{2.4.6}:
\def\theequation{2.4.7}\begin{equation}
\left\{
\aligned
{}
&
\underline{\mathcal{L}}_{z_1}(0)=\left(0, \, 0, \, z_1, \,
\Theta(z_1,0,0)
\right).\\
& 
\mathcal{L}_{z_2}(\underline{\mathcal{L}}_{z_1}(0))=
\left(
z_2,\, \overline{\Theta}(z_2,z_1,\Theta(z_1,0,0)),\, 
z_1, \, 
\Theta(z_1,0,0)
\right).\\
& 
\underline{\mathcal{L}}_{z_3}(\mathcal{L}_{z_2}
(\underline{\mathcal{L}}_{z_1}(0)))=
\left(
z_2,\, \overline{\Theta}(z_2,z_1,\Theta(z_1,0,0)),\, 
z_1+z_3, \, \Theta(z_1+z_3,z_2,\overline{\Theta}(z_2,z_1,
\Theta(z_1,0,0)))\right)
\endaligned\right.
\end{equation}
Also, we observe that $\Gamma_{k+1}(z_{(k)},0)=\Gamma_k(z_{(k)})$,
since $\mathcal{L}_0$ and $\underline{\mathcal{L}}_0$ coincide with
the identity map by~\thetag{2.4.2} and~\thetag{2.4.3}. 
So, for $k$ increasing, the ranks at the origin of
the maps $\Gamma_k$ are increasing. We now introduce the
following important definition.

\def\thedefinition{2.4.8}\begin{definition}
{\rm The generic submanifold $M$ is said to 
be {\it minimal}\, at $p$ if the maps $\Gamma_k$ are 
of (maximal possible) rank equal to $2m+d={\rm dim}_\C\,\mathcal{M}$ at
the origin in $\Delta_{mk}(\delta_k)$ for all $k$ large enough.} 
\end{definition}

In other words, $M$ is minimal at $p$ if and only if sufficiently high
order Segre chains are submersive.  Equivalently, sufficiently high
order conjugate Segre chains are submersive. In the next Section~2.5,
we shall show that minimality is characterized by the fact that the
maps are of generic rank equal to $2m+d$ for all $k$ large enough.
First of all, in order to enlighten this definition, we shall prove a
general {\sl local orbit theorem} in the spirit of H.J.~Sussmann~[32].

\subsection*{2.4.9.~Local orbit theorem in the
$\K$-analytic category} Let $\K=\R$ or $\C$ and
let $\Delta:=\{x\in \K: \, \vert x\vert < 1\}$
and $r\Delta:=\{x\in \K: \, \vert x \vert <r\}$. Let $n\in \N$ with
$n\geq 1$.  In $\Delta^n$ equipped with coordinates
$x=(x_1,\dots,x_n)$, we consider the origin as a center point. Let
$\L=\{L^a\}_{1\leq a \leq A}$, $A\geq 1$, be a {\it finite}\, system
of nonzero vector fields defined all over $\Delta^n$. {\it We do not
require that this set is stable under taking linear combinations with
coefficients being analytic or algebraic functions over
$\Delta^n$}. Let $L\in \L$. As previously, we shall simply denote the
flow map of $L$ by $(t,x)\mapsto L_t(x)\equiv \exp(tL)(x)$.

We recall the defining properties of the flow map: $L_0(x)=x$ and
$\frac{d}{dt} (L_t(x))= L( L_t(x))$, where $L(x')$ denotes the value
of $L$ at $x'$. As is known, the $\K$-algebraic case is
exceptional because the flow of a vector field $L$ having
$\K$-algebraic coefficients is in general transcendent.  This is why
we shall in fact make two kinds of precise regularity assumptions:

\smallskip
\noindent
{\tt (Algebraic):}
The coefficients of all elements of $\L$ are
$\K$-algebraic {\it and moreover}\, their flow if also 
$\K$-algebraic. In fact, the $\K$-algebraicity of the flow
implies the $\K$-algebraicity of the coefficients.

\smallskip
\noindent
{\tt (Analytic):}
The coefficients of all elements of $\L$ are $\K$-analytic power
series centered at the origin converging in $\Delta^n$, whence their
flow is $\K$-analytic.

\smallskip

Choose now $r$ with $0<r\leq 1/2$. We first define finite
concatenations of flow mappings of vector fields in $\L$ as follows.
If $k\in \N_*$, $L=(L^1,\ldots,L^k)\in \L^k$, $t=(t_1,\ldots,t_k)\in
\K^k$ and $x\in (r\Delta)^n$, we use again the notation
$L_t(x)=L_{t_k}^k(\cdots(L_{t_1}^1(x))\cdots)$ whenever the
composition is defined. Anyway, after bounding $k\leq 3n$, it is clear
that there exists $\delta >0$ such that all maps $(t,x)\mapsto L_t(x)$
are well-defined for $t\in (2\delta \Delta)^k$,
$x\in(\frac{r}{2}\Delta)^n$ and satisfy $L_t(x)\in (r\Delta)^n$. 
By definition, the point $x'=L_{t_k}^k\circ \cdots \circ
L_{t_1}^1(x)$ is the endpoint of a piecewise smooth 
algebraic or analytic curve with origin
$x$: it consists in following $L^1$ during time $t_1$, 
following $L^2$ during time $t_2$ $\ldots$ and
following $L^k$ during time $t_k$.

Next, we shall say that an embedded small piece of $\K$-manifold
$N\subset \Delta^n$ passing through the origin (which is either
$\K$-algebraic or $\K$-analytic) is a {\sl weak $\L$-integral
manifold}\, if $T_xN \supset \L(x)$ for all $x\in N$. In the formal
case, however, this condition is meaningless. Equivalently, we mean
that for each $L\in \L$, $L\vert_N$ is tangent to $N$. Now, this new
condition makes sense in the formal case. We shall in fact consider
the formal case afterwards, as a generalization of
the algebraic and analytic cases. In particular,
it clearly ensues from the tangency of $L\vert_N$ to $N$ that any
integral curve of an element $L\in \L$ with origin a point $x\in N$
which belongs to $N$ stays in $N$.

Now, we introduce the following special definitions.  The {\sl
$\L$-orbit of $0$ in $\Delta^n$}, denoted by $\mathcal{ O}_{\L} (0)$,
is the set of all points $L_t(0)\in (r \Delta)^n$ for all $t\in
(\delta \Delta)^k$, $k\leq 3n$. The reason why we bound $k\leq 3n$
will be clear afterwards. We shall say that the open set $(r
\Delta)^n$ is {\sl $\L$-minimal at $0$}\, if $\mathcal{O}_{\L}(0)$
contains a polydisc $(\varepsilon \Delta)^n$, where $\varepsilon >0$.
If $L=(L^1, \dots,L^k) \in \L^k$ with $k \leq 3n$, we shall denote by
$\Gamma_L (t)$ the mapping $t \mapsto L_{t_k}^k( \cdots(L_{ t_1}^1
(0)) \cdots)$.

We can now state the {\sl local orbit theorem}. By induction, we
shall construct a special sequence of vector fields $L^{*k} :=
(L^{*1}, \ldots, L^{*k}) \in \L^k$, $k \in \N_*$. We state a
long but progressively explained theorem with the purpose of
exhibiting all the relevant informations.

\def\thetheorem{2.4.10}\begin{theorem} 
As above, let $\L=\{L^a\}_{1\leq a\leq A}$, $A\in \N_*$, be a finite nonempty
set of vector fields which are defined over $\Delta^n$ and satisfy
one of the two regularity assumptions {\tt (Algebraic)} or {\tt
(Analytic)}. Then there exists an integer $e\geq 1$ and a
multiplet of vector fields $L^*=(L^{*1},\ldots,L^{*e})\in \L^{e}$
such that the following seven
properties hold\text{\rm \,:}
\begin{itemize}
\item[{\bf (1)}] 
For every $k=1,\dots,e$, the map $(t_1,\dots,t_k)\mapsto
L_{t_k}^{*k}(\cdots (L_{t_1}^{*1}(0))\cdots)$ is of generic rank equal
to $k$.
\item[{\bf (2)}] 
For every arbitrary element $L'\in \L$, the 
map $(t_1,\dots,t_e,t')\mapsto 
L_{t'}'(L_{t_e}^{*e}(\cdots (L_{t_1}^{*1}(0))\cdots))$ is
of generic rank $e$, hence $e$ is the maximal possible generic rank.
\item[{\bf (3)}] 
There exists an element $t^*\in \Delta^e$ arbitrarily close to the
origin which is of the special form $(t_1^*,\ldots,t_{e-1}^*,0)$,
namely with $t_{e}^*=0$, and there exists an open connected neighborhood 
$\omega^*$ of
$t_*$ in $\Delta^e$ such that the map $\Gamma_{L^*}: t\mapsto
L_{t_e}^{*e}(\cdots(L_{t_1}^{*1}(0)))$ is of constant rank $e$ over
$\omega^*$.
\item[{\bf (4)}]
After setting $L^*:=(L^{*1}, \ldots, L^{*e})$,
$K^*:=(L^{*e-1},\ldots,L^{*1})$ and $s^* := (-t_{e-1}^*, \ldots,
-t_1^*)$, we then have $K^*_{ s^*} \circ L_{
t^*}^*(0)=0$. Furthermore, the map $\psi: \omega^* \to \Delta^n$
defined by $\psi: t\mapsto K_{s^*}^*\circ L_t^*(0)$ is of constant
rank equal to $e$ over the domain $\omega^*$.
\item[{\bf (5)}] 
Then the image $\psi
(\omega^*)$ is a piece of $\K$-manifold passing through
the origin which is either $\K$-algebraic or $\K$-analytic, because
the flows of the elements of $\L$ are $\K$-algebraic or $\K$-analytic.
\item[{\bf (6)}]
This piece of $\K$-manifold is a weak $\L$-integral manifold.
Furthermore, very weak $\L$-integral manifold passing
through $0$ must contain $\mathcal{O}_\L(0)$ in a neighborhood of
$0$.
\item[{\bf (7)}]
There exists $\varepsilon>0$ such that $\mathcal{O}_\L(0)\cap 
(\varepsilon\Delta)^n=\psi(\omega^*)$.
\end{itemize} 
In conclusion, the local orbit $\mathcal{O}_\L(0)$ is represented by
the small $\K$-manifold $\psi(\omega^*)$ and its dimension $e$ is
characterized by the generic rank properties {\bf (1)} and {\bf (2)}.
\end{theorem}

\proof
If all vector fields in $\L$ vanish at the origin, then
$\mathcal{O}_\L(0)=\{0\}$. We now exclude this possibility. Choose a
vector field $L^{*1}\in\L$ which does not vanish at $0$.  The map
$t_1\mapsto L_{t_1}^{*1}(0)$ is of generic rank one.  If there exists
$L'\in \L$ such that the map $(t_1,t')\mapsto
L_{t'}'(L_{t_1}^{*1}(0))$ is of generic rank two, we choose one such
$L'$ and we denote it by $L^{*2}$.  Continuing this process, we get
vector fields $L^{*1}, \dots,L^{*e}$ satisfying properties {\bf (1)}
and {\bf (2)}. Since the generic rank of the map $\Gamma_{L^*}:
(t_1,\dots,t_e)\mapsto L_{t_e}^{*e}(\cdots(L_{t_1}^{*1}(0))\cdots)$
equals $e$, and since this map is either $\K$-algebraic 
or $\K$-analytic, there exist elements $t^*\in\Delta^e$ arbitrarily
close to the origin such that its rank at $t_*$ equals $e$. We claim
that we can moreover choose $t^*$ to be of the special form
$(t_1^*,\dots,t_{e-1}^*,0)$, {\it i.e.}  with $t_e^*=0$.
This is a consequence of the following lemma.

\def\thelemma{2.4.11}\begin{lemma}
Let $n\in\N_*$, $e\in \N_*$, $t\in \K^e$ and
$t\mapsto \varphi(t)=(\varphi_1(t),\dots,\varphi_n(t))\in\K^n$
be a mapping of generic rank equal to $(e-1)$ which is either
$\K$-algebraic or $\K$-analytic. Let $L'\in \L$ and
assume that the mapping $\psi: (t,t')\mapsto L_{t'}'(\varphi(t))$ has
generic rank $e$. Then there exists a point $(t^*,0)$ at
which the rank of $\psi$ is equal to $e$.
\end{lemma}

\proof
Suppose on the contrary that at all points $(t^*,0)$, the map
$(t,t')\mapsto L_{t'}'(\varphi(t))$ has rank $\leq e-1$. Choose $t^*$
arbitrarily close to zero at which $\varphi$ has maximal, hence
locally constant rank $(e-1)$. By the rank theorem, there exists a
neighborhood $\omega^*$ of $t^*$ in $\K^e$ such that
$N:=\varphi(\omega^*)$ is a small piece of 
$\K$-manifold which is either $\K$-algebraic 
or $\K$-analytic. Let $L'\in \L$. Since by assumption
the rank of $(t,t')\mapsto L_{t'}'(\varphi(t))$ equals 
$(e-1)$ at every point $(t^*,0)\in \omega^*\times \K$, 
it follows necessarily that $L'$ is tangent 
to $N$. Consequently, the flow of $L'$ stabilizes $N$. 
Finally, for $\varphi(t)\in N$, we have $L_{t'}'(\varphi(t))\in N$ 
also, whence the rank of $(t,t')\mapsto L_{t'}'(\varphi(t))$ 
is less than or equal to $\dim_\K N=e-1$ at every point in a 
neighborhood of $(t^*,0)$ in $\K^n\times \K$. 
We have proved that the mapping $(t,t')\mapsto L_t'(\varphi(t))$ is
of generic rank $(e-1)$ in a neighborhood of $(t^*,0)$
in $\K^e\times \K$, hence everywhere by the principle of
analytic continuation, which contradicts the assumption that
it has generic rank equal to $e$. This completes the proof.
\endproof

\smallskip
\noindent
{\it End of the proof of Theorem~2.4.9.}  Now, we choose a point
$(t_1^*,\dots,t_{e-1}^*,0)\in \Delta^e$ arbitrarily close to $0$ at
which the rank of $t\mapsto
L_{t_e}^{*e} (\cdots (L_{t_1}^{*1} (0))\cdots)$ is maximal and equals
$e$, so we get {\bf (3)}.  In {\bf (4)}, the property $K_{s^*}^*\circ
L_{t^*}^*(0)$ is obvious, since the mapping:
\def\theequation{4.5.7}\begin{equation}
L_{-t_1^*}^{*1}\circ
\cdots \circ L_{-t_{e-1}^*}^*\circ
L_{0}^{*e}\circ L_{t_{e-1}^*}^*\circ
\cdots \circ
L_{t_1^*}^*(\cdot)={\rm Id}
\end{equation}
is the indentity. By {\bf (2)}, the mapping $(t,s)\mapsto K_s^*\circ
L_t^*(0)$ is also of generic rank $e$, whence its restriction to
$\omega^*\times \{s^*\}$ is of constant rank $e$ since the mapping
$K_{s^*}^*(\cdot)$ is a local diffeomorphism. We get {\bf (4)} and
then {\bf (5)} obviously.  So we have constructed a piece $N$ of
$\K$-manifold passing through the origin. Let $L'\in \L$. We claim
that $L'$ is tangent to $N$. Otherwise, if $L'$ would not be tangent,
the mapping $(t,s,t')\mapsto L_{t'}'(K_s(L_t(0)))$ would be of generic
rank $\geq e+1$, contrarily to the definition of $e$.
This completes the proof of Theorem~2.4.9.
\endproof

\section*{\S2.5.~Segre type and Segre multitype}

We now apply the general considerations of Theorem~2.4.9 to the
specific situation where $\K=\C$, where $\Delta^n$ is replaced by
$\mathcal{M}$ and where the collection $\L$ of vector fields is our
previous complexified vector fields $\{\mathcal{L}_k,
\underline{\mathcal{L}}_k\}_{1\leq k\leq m}$.  We also bring some
refinements.  

\subsection*{2.5.1.~Increasing generic ranks}
Let us denote by ${\rm genrk}_\C(\Phi)$ the {\sl generic
rank} of a $\C$-algebraic or $\C$-analytic map $\Phi: X\to Y$ of
connected complex manifolds. Here of course, we have ${\rm genrk}_\C
(\Gamma_1)= {\rm genrk}_\C (\underline{ \Gamma}_1)=m$ and ${\rm
genrk}_\C (\Gamma_2)= {\rm genrk}_\C (\underline{ \Gamma}_2)=2m$,
which is evident in equations~\thetag{2.4.6} and~\thetag{2.4.7}.  We
set $e_1:=m$ and $e_2:=m$. Next, we set $e_3:= {\rm
genrk}_\C(\Gamma_3)-2m$ and, by induction $e_{k+1}:={\rm
genrk}_\C(\Gamma_{k+1})-e_3-\cdots -e_k-2m$, whence ${\rm
genrk}_\C(\Gamma_k)=2m+e_3+\cdots+e_k$ if $k\geq 3$, and similarly, we
can define the sequence $\underline{e}_k$ for
$\underline{\Gamma}_k$. We notice at once that we have
$\underline{e}_k=e_k$, since $\sigma (\Gamma_k (z_{(k)}))=
\underline{\Gamma}_k (\overline{ z_{(k)}})$.

We claim that $e_l=0$ for all $l\geq k+1$ if $e_{k+1}=0$ and $e_k\neq
0$. In other words, the generic rank enjoys a stabilization
property. Indeed, we first choose a point $z_{(k)}^*$ arbitrarily
close to the origin in $\C^{mk}$ such that $\Gamma_k$ has (necessarily
locally constant) rank equal to $2m+e_3+\cdots+e_k$ at $z_{(k)}^*$ and
we set $q:=\Gamma_k(z_{(k)}^*)\in \mathcal{M}$. Then by the rank
theorem, the image $\mathcal{H}$ of a neighborhood $\mathcal{W}^*$ of
$z_{(k)}^*$ is a submanifold of $\mathcal{M}$ of dimension
$2m+e_3+\cdots+e_k$.  We claim that the vector fields $\mathcal{L}_k$
and $\underline{\mathcal{L}}_k$ are all tangent to $\mathcal{H}$.  For
instance, to fix ideas, we assume that $k$ is even (the odd case will
be similar). Thus we can write $\Gamma_k(z_{(k)})=
\underline{\mathcal{L}}_{z_k} (\cdots (\mathcal{L}_{z_1}(0))\cdots)$,
{\it i.e.} the chain $\Gamma_k$ ends-up with a
$\underline{\mathcal{L}}$. This shows that $\mathcal{H}$ is fibered by
the leaves of $\mathcal{F}_{\underline{\mathcal{L}}}$, so the
$\underline{\mathcal{L}}_k$ are already tangent to $\mathcal{H}$ at
every point. By differentiating $\Gamma_{k+1}=\mathcal{L}_{z_{k+1}}
(\Gamma_k(z_{(k)}))$ with respect to $z_{k+1}$ at $z_{k+1}=0$, we
obtain the $m$-dimensional space $\mathcal{L}(\Gamma_k(z_{(k)}))$,
namely the tangent space to the foliation $\mathcal{F}_{\mathcal{L}}$
at the point $\Gamma_k(z_{(k)})$.  Then the assumption $e_{k+1}=0$
entails that this space $\mathcal{L}(\Gamma_k(z_{(k)}))$ is
necessarily contained in the tangent space to $\mathcal{H}$ at
$\Gamma_k(z_{(k)})$, which proves the claim.  Finally, as the
$\mathcal{L}_k$ and the $\underline{\mathcal{L}}_k$ are all tangent to
$\mathcal{H}$, it follows that their local flow at $q$ is contained in
$\mathcal{H}$, whence the range of the subsequent $\Gamma_l$, $l\geq
k+1$, is locally contained in $\mathcal{H}$.  Because they are either
algebraic or analytic, this shows that their generic rank does not go
beyond $2m+e_3+\cdots+e_k$, which proves the claim.

In conclusion, there exists a well-defined integer $\mu_0\geq 2$
with $\mu_0\leq d+2$ such that $e_3>0,\dots,e_{\mu_0}>0$ and
$e_l=0$, for all $l\geq \mu_0+1$. We call the integer
$\mu_0$ the {\sl Segre type} of $\mathcal{M}$ at the origin and
we call the $\mu_0$-tuple $(m,m,e_3,\dots,e_{\mu_0})$ the {\sl
Segre multitype} of $\mathcal{M}$ at the origin. This {\sl Segre multitype}
simply recollects all the jumps of generic ranks of the
$\Gamma_k$. It is clear that Segre type and multitype are
biholomorphic invariants, because the Segre foliations defined by
$\mathcal{L}$ and $\underline{\mathcal{L}}$ are so. To summarize, we
have\,:
\begin{itemize}
\item[{\bf (a)}]
$\text{\rm genrk}_{\C} (\Gamma_k)=
2m+e_3+\cdots+e_k= \text{\rm genrk}_{\C} 
(\underline{\Gamma}_k)$, for $2\leq k \leq \mu_0$,
\item[{\bf (b)}]
$\text{\rm genrk}_{\C} (\Gamma_k)=2m+e_3+\cdots+e_{\mu_0}=
\text{\rm genrk}_{\C} (\underline{\Gamma}_k)$, for 
$k\geq \mu_0$.
\end{itemize}
The main advantage of dealing with $\mathcal{M}$, $\mathcal{L}$,
$\underline{\mathcal{L}}$, $\Gamma_k$ and $\underline{\Gamma}_k$ lies
in the fact that all these objects may be understood in a
coordinate-free way. Even the two projections $\pi_t$ and $\pi_{\tau}$
can be defined abstractly, because their fibers are the leaves of the
Segre foliations. As \S2.4.9, We may define the orbit
$\mathcal{O}_{\mathcal{L}, \underline{\mathcal{L}}}(\mathcal{M},0)$
and we have the following theorem.

\def\thetheorem{2.5.2}\begin{theorem}
Assume that $M$ is a real algebraic or analytic generic
submanifold of $\C^n$, let $p_0\in M$, let $\mathcal{M}:=(M)^c$ be its
complexification and let $(p_0)^c:=(p_0,\bar p_0)$, identified
with the origin in coordinates $(t,\tau)$.  Let
$\mathcal{L}_k$ and $\underline{\mathcal{L}}_k$, $k=1,\dots,n$, be
vector fields generating the pair of Segre foliations.  Then there exists
$z_{(\mu_0)}^*\in\C^{m\mu_0}$ arbitrarily close to the origin of the
form $z_{(\mu_0)}^*=(z_1^*,\dots,z_{\mu_0-1}^*,0)$ and a small
neighborhood $\mathcal{W}^*$ of $z_{(\mu_0)}^*$ in
$\Delta_{m\mu_0}(\delta_{\mu_0})$ such that, if we denote
$\omega_{(\mu_0-1)}^*:=(-z_{\mu_0-1}^*,\dots, -z_1^*)$, then we
have\text{\rm :}
\begin{itemize}
\item[{\bf (c)}]
The complex algebraic
or analytic map $\Gamma_{\mu_0}$ is of rank $2m+e_3+\cdots+e_{\mu_0}$ at
$z_{(\mu_0)}^*$.
\item[{\bf (d)}]
$\Gamma_{2\mu_0-1}(z_{(\mu_0)}^*,\omega_{(\mu_0-1)}^*)=0$.
\item[{\bf (e)}]
The restricted map $\Gamma_{2\mu_0-1}: 
\mathcal{W}^*\times \omega_{(\mu_0-1)}^* \to \mathcal{M}$
is of constant rank equal to $2m+e_3+\cdots+e_{\mu_0}$.
\item[{\bf (f)}]
The image $\Gamma_{2\mu_0-1} (\mathcal{W}^*\times
\omega_{(\mu_0-1)}^*)$ is a complex algebraic or analytic
manifold-piece, denoted by $\mathcal{O}_{\mathcal{L},
\underline{\mathcal{L}}}(\mathcal{M}, 0)$ and called the {\sl
$\{\mathcal{L}, \underline{\mathcal{L}}\}$-orbit of the origin} with
the property that the vector fields $\mathcal{L}_k$ and
$\underline{\mathcal{L}}_k$ are tangent to it.
\item[{\bf (g)}]
This integer $2m+e_3+\cdots+e_{\mu_0}$ is equal to $\dim_{\C} \,
\mathcal{O}_{\mathcal{L}, \underline{\mathcal{L}}}(\mathcal{M}, 0)$.
\end{itemize}
Of course, the same statement holds with
$\underline{\Gamma}_{2\mu_0-1}$ instead of $\Gamma_{2\mu_0-1}$.
\end{theorem}

\proof
The proof is similar to the proof of
Theorem~2.4.10, with minor modifications.
According to {\bf (b)}, $\Gamma_{\mu_0}$ is of generic rank
$2m+e_3+\cdots+e_{\mu_0}$. Consequently, for every point
$z_{(\mu_0)}^*\in\Delta_{m\mu_0}(\delta_{\mu_0})$ outside of some
proper complex subvariety, the map $\Gamma_{\mu_0}$ is of rank
$2m+e_3+\cdots+e_{\mu_0}$ {\it at} \, $z_{(\mu_0)}^*$. In fact, we
claim that we can even choose such a $z_{(\mu_0)}^*$ of the form
$(z_1^*,\dots,z_{\mu_0-1}^*,0)$, {\it i.e.} with
$z_{\mu_0}^*=0$. Indeed, as $\Gamma_{(k)}(z_{(k)}) =[\mathcal{L}
\,\text{\rm or} \, \underline{\mathcal{
L}}]_{z_k}(\Gamma_{(k-1)}(z_{(k-1)}))$, we can apply the following
lemma, which follows directly form Lemma~2.4.11.

\def\thelemma{2.5.3}\begin{lemma}
Let $z\in\C^m$, $z'\in\C^{m'}$, let $\Gamma(z')\in\mathcal{M}$ be a
formal, algebraic or analytic map of $z'$ with $\Gamma(0)=0$ and let
$\varphi(z,z'):= {\mathcal{ L}}_z(\Gamma(z'))$ or $\varphi(z,z'):=
\underline{\mathcal{L}}_z(\Gamma(z'))$. Then $\varphi$ attains its
maximal generic rank at some points of the form $(0,{z'}^*)\in \C^m\times
\C^{m'}$. 
\end{lemma}

\noindent
Now, we fix such a $z_{(\mu_0)}^*$ of the form
$(z_1^*,\dots,z_{\mu_0-1}^*,0)$, which satisfies {\bf (c)} and we
check that it satisfies the other claims.  Let
$\omega_{(\mu_0-1)}^*:=(-z_{\mu_0-1}^*,\dots,-z_1^*)$.  First, {\bf
(d)} is easy\,: suppose for instance $\mu_0$ is even, then we have
$\Gamma_{2\mu_0-1}(z_{ (\mu_0)}^*,\omega_{(\mu_0 -1)}^*)=
\mathcal{L}_{-z_1^*}\circ \cdots\circ \mathcal{
L}_{-z_{\mu_0-1}^*}\circ \underline{\mathcal{L}}_0\circ
\mathcal{L}_{z_{ \mu_0-1}^*}\circ \cdots \circ
\mathcal{L}_{z_1^*}(0)=0$, because $\underline{ \mathcal{L}}_0(q)= q$,
$\mathcal{L}_{-w}\circ \mathcal{L}_w= \hbox{Id}$ and $\underline{
\mathcal{L}}_{-\zeta}\circ \underline{\mathcal{ L}}_{
\zeta}=\hbox{Id}$. The odd case is similar.  Now, we prove {\bf
(e)}. The restricted map $z_{ (\mu_0)} \mapsto \mathcal{ L}_{-z_1^*}
\circ \cdots \circ \mathcal{ L}_{-z_{ \mu_0-1}^*}\circ \Gamma_{
\mu_0}(z_{ (\mu_0)})$ (again written in case $\mu_0$ is even), is
clearly of rank $2m+e_3+ \cdots+e_{\mu_0}$ at the point
$z_{(\mu_0)}^*$, because the map $q \mapsto \mathcal{ L}_{-z_1^*}
\circ\cdots\circ \mathcal{ L}_{-z_{ \mu_0 -1}^*}(q)$ is a local
biholomorphism, by definition of flows. Notice that
$\Gamma_{2\mu_0-1}$ is then of {\it constant rank} equal to $2m+e_3+
\cdots+e_{\mu_0}$ in a neighborhood of $(z_{ (\mu_0)}^*,
\omega_{(\mu_0-1)}^*)$ in $\mathcal{W}^*\times \omega_{ (\mu_0-
1)}^*$, since, by {\bf (b)}, $2m+e_3+ \cdots+e_{ \mu_0}$ is already
the maximum value of all the generic ranks of the $\Gamma_k$. This
proves {\bf (e)}. By definition, the orbit $\mathcal{O}_{\mathcal{L},
\underline{\mathcal{L}}}(\mathcal{M},0)$ is the union of the ranges of
the maps $\Gamma_k$ and $\underline{\Gamma}_k$. It is easy to check
that this double union coincides in fact with the union of only the
$\Gamma_k$ (or of only the $\underline{\Gamma}_k$), simply because,
setting $z_1=0$, we have $\Gamma_k(0,z_2,\dots,z_k)\equiv
\underline{\Gamma}_{k-1}(z_2,\dots,z_k)$.  Thanks to the constant rank
property {\bf (e)}, we already know that this orbit contains the
$(2m+e_3+\cdots+e_{\mu_0})$-dimensional manifold-piece passing through
$0$\,: $\mathcal{N}:=\Gamma_{2\mu_0-1} (\mathcal{
W}^*\times\omega_{(\mu_0-1)}^*)$. Because by {\bf (b)} the next
generic ranks for $k\geq 2\mu_0-1$ do not increase and because of the
principle of analytic continuation, we then deduce that all the ranges
of the subsequent $\Gamma_k$ are contained in this manifold
piece $\mathcal{N}$ and it follows that
$\mathcal{L}$ and $\underline{\mathcal{L}}$ are
tangent to this manifold-piece. 
In conclusion, we get {\bf (f)} and {\bf (g)}, which completes
the proof of Theorem~2.5.2.
\endproof

In the hypersurface case, we have the following
simple criterion of minimality, left to the reader.

\def\thecorollary{2.5.4}\begin{corollary}
If $M$ is a real algebraic
or analytic hypersurface, 
\text{\rm i.e.} if $d=1$, then
\begin{itemize}
\item[{\bf (h)}] 
$M$ is minimal at $0$
\ $\Longleftrightarrow$ \ $\mu_0=3$.
\item[{\bf (i)}]
$M$ is nonminimal at $0$ \ $\Longleftrightarrow$ \ $\mu_0=2$.
\end{itemize} 
\end{corollary}

\subsection*{2.5.7.~Example}
We now give a simple example in the hypersurface case which illustrates
statements {\bf (e)} and {\bf (f)} of Theorem~2.5.2 in a very concrete
way. We let $M$ be the hypersurface of $\C^2$ of equation 
$z=\bar z+iw^2\bar w^2$. We choose $p=0$ and here $2\mu_0-1=5$. 
We compute\,:
\def\theequation{2.5.8}\begin{equation}
\left\{
\aligned
{}
&
\Gamma_1(z_1)=
(z_1, \, 0, \, 0, \, 0)\\
&
\Gamma_2(z_1,z_2)=
(z_1, \, 0, \, z_2, \, -i z_1^2 z_2^2)\\\
&
\Gamma_3(z_1,z_2,z_3)=(z_1+z_3, \,iz_2^2[z_3^2+2z_1z_3], \,
z_2, \, -i z_1^2 z_2^2)\\
&
\Gamma_4(z_1,z_2,z_3,z_4)=(z_1+z_3, \, 
iz_2^2[z_3^2+2z_1z_3], \, z_2+z_4, \, \\
& 
\ \ \ \ \ \ \ \ \ \ \ \ \ \ \ \ \ \ \ \ \ \ \ \ \ \ \ \ \ \
iz_2^2[z_3^2+2z_1z_3]-i
[(z_2+z_4)(z_1+z_3)]^2)\\
&
\Gamma_5(z_1,z_2,z_3,z_4,z_5)=(z_1+z_3+z_5, \,
iz_2^2[z_3^2+2z_1z_3]-\\
&
\ \ \ \ \ \ \ \ \ \ \ \ \ \ \ \
-i[(z_2+z_4)(z_1+z_3)]^2)+i[(z_1+z_3+z_5)(z_2+z_4)]^2, \, \\
&
\ \ \ \ \ \ \ \ \ \ \ \ \ \ \ \
z_2+z_4, \, iz_2^2[z_3^2+2z_1z_3]-i
[(z_2+z_4)(z_1+z_3)]^2).
\endaligned\right.
\end{equation}
The maps $\Gamma_k$ have range in $\mathcal{M}$, on which either the
coordinates $(z,w,\zeta)$ or $(\zeta,\xi,z)$ can be chosen. We do the
first choice for $k$ even and the second choice for $k$ odd. Thus, we
view $\Gamma_5$ as a map $\C^5\to \C_{(z,\zeta,\xi)}^3$, {\it i.e.} we
forget the second $w$-coordinate in the above expression of $\Gamma_5$.
Now, computing the $3\times 5$ Jacobian matrix of $\Gamma_5$ at the point
$(z_{(3)}^*,\omega_{(2)}^*)$ as in Theorem~2.5.2 which is necessarily of
the form $(z_1^*,z_2^*,0,-z_2^*,-z_1^*)$, and for which we clearly
have $\Gamma_5(z_1^*,z_2^*,0,-z_2^*,-z_1^*)=0$, we see that the
determinant of the first $3\times 3$ submatrix is equal to
$2iz_1^*(z_2^*)^2$. Thus, it is nonzero for an arbitrary choice of
$z_1^*\neq 0$ and $z_2^*\neq 0$. By the way, the question arises
whether the integer $(2\mu_0-1)$ in Theorem~2.5.2 is optimal. Incidentally, 
this example shows that it is optimal. Indeed, if we ask whether there
exists $z_{(4)}^*=(z_1^*,z_2^*,z_3^*,z_4^*)$ such that
$\Gamma_4(z_{(4)}^*)=0$ and the rank at $z_{(4)}^*$ of the
differential of $\Gamma_4$ equals $3$ (the dimension of $\mathcal{M}$), then
looking at eqs~\thetag{2.5.8}, we get first $z_1^*+z_3^*=0$,
$(z_2^*)^2z_3^*[z_3^*+2z_1^*]=0$ and $z_2^*+z_4^*=0$, thus $z_{(4)}^*$
is necessarily of the form $(0,z_2^*,0,-z_2^*)$ or
$(z_1^*,0,-z_1^*,0)$. Viewing now $\Gamma_4$ as a map $\C^4\to
\C_{(z,w,\zeta)}^3$, and computing its $3\times 4$ Jacobian matrix at
such points, one sees that it is of rank $2$, which proves
the claim.

\subsection*{2.5.9.~Segre geometry in the formal category}
Replacing the complex defining equations~\thetag{2.1.15} by purely
formal power series, the reader may verify that all the previous
results are meaningful and may be expressed in terms of purely formal power
series.

\section*{\S2.6.~Extrinsic complexification of CR orbits}

\subsection*{2.6.1.~Intrinsic CR orbits and their smoothness}
On $M$ the complex tangent bundle $T^cM$ is generated by the $2m$
sections ${\rm Re}\, L_1, {\rm Im}\, L_1, \dots, {\rm Re}\, L_m, {\rm
Im}\, L_m$. In the formal case, their flow is formal.
In the analytic case, their flow is analytic. 
However, a very subtle point occurs in the algebraic
category. We have seen that the complex flows of the complexified
vector fields $\mathcal{L}_k$, given by~\thetag{2.4.2} and~\thetag{2.4.3}, are
algebraic. This is untrue about the real flows of the
real and imaginary parts of the vector fields $L_k$, as shows
the following example.

\def\theexample{2.6.2}\begin{example}
{\rm 
Let $M$ be the real analytic hypersurface passing
through the origin in $\C^2$ defined by 
${\rm Im}\, w=\sqrt{1+z\bar z}-1$. The vector field
$\overline{L}:=\partial_{\bar z}+
iz(1+z\bar z)^{-1/2}\, \partial_{\bar w}$ generates
$T^{0,1}M$. Its real and imaginary parts are given by 
\def\theequation{2.6.3}\begin{equation}
\left\{
\aligned
2\,{\rm Re}\, \overline{L}= & \
\partial_x-y(1+x^2+y^2)^{-1/2}\, \partial_u+
x(1+x^2+y^2)^{-1/2}\, \partial_v,\\
2\,{\rm Im}\, \overline{L}= & \
-\partial_y+x(1+x^2+y^2)^{-1/2}\, 
\partial_u+y(1+x^2+y^2)^{-1/2}\, \partial_v.
\endaligned\right.
\end{equation}
We claim that the flow of $2\, {\rm Re}\,\overline{L}$ is not
algebraic.  Indeed, let $s$ denote a real time parameter and let
$(x(s),y(s),u(s),v(s))$ be the unique integral curve of $2\, {\rm
Re}\, \overline{L}$ with $x(0)=x_0$, $y(0)=y_0$, $u(0)=u_0$ and
$v(0)=v_0$ with $(x_0+iy_0, u_0+iv_0)\in M$. We compute first
$x(s)=x_0+s$, $y(s)=y_0$, $v(s)=(1+y_0^2+(x_0+s)^2)^{1/2}-1$ and then
$u(s)$ satisifies the ordinary differential equation
\def\theequation{2.6.4}\begin{equation}
du(s)/ds=
-y_0(1+y_0^2(x_0+s)^2)^{-1/2},
\end{equation}
which may be integrated as
\def\theequation{2.6.5}\begin{equation}
u(s)=u_0-y_0\left[
{\rm Arcsh}\, \left(
{x_0+s\over \sqrt{1+y_0^2}}
\right)-
{\rm Arcsh}\,
\left(
{x_0\over \sqrt{1+y_0^2}}
\right)
\right]
\end{equation}
Consequently, the flow of $2\,{\rm Re}\, \overline{L}$ is not
algebraic. However, {\it we stress that the complex flow of the
complexified vector field $\underline{\mathcal{L}}=
\partial_\zeta+iz(1+z\zeta)^{-1/2}\,\partial_\xi$ is complex
algebraic}, as shown by the general expression~\thetag{2.4.3}, which
yields in this particular case $\underline{\mathcal{L}}_\zeta
(z_p,w_p,\zeta_p,\xi_p)= (z_p,w_p, \zeta_p+\zeta, w_p+2i(\sqrt{1+
z\zeta}-1)))$: indeed, this expression is clearly algebraic~!
}
\end{example}

We now consider the set $\L:=\{{\rm Re}\, L_1, {\rm Im}\, L_1, \dots,
{\rm Re}\, L_m, {\rm Im}\, L_m\}$ of $2m$ real vector fields
generating $T^cM$. Applying Theorem~2.4.10, 
we may construct the {\sl local CR orbits of
points $p$ in $M$}, which we denote by $\mathcal{O}_{CR}(M,p)$.
Since they are weak $T^cM$-integral 
manifolds, they are of the same CR dimension as $M$.
According to Example~2.6.2 just above, we only know that
$\mathcal{O}_{CR}(M,p)$ is a real analytic submanifold passing through
$p$. Algebraicity seems to be lost, but naturally, according to
Leibniz' principle of sufficient reason, one expects however that
$\mathcal{O}_{CR}(M,p)$ is algebraic if $M$ is algebraic

The above Example~2.6.2 shows that the algebraicity of the CR
orbits of a real algebraic generic manifold cannot be established
by means of sections of the complex tangent bundle $T^cM={\rm Re}\, 
T^{0,1}M$. Fortunately, 
by passing to the extrinsic complexification, we may avoid this
difficulty.

\def\thetheorem{2.6.6}\begin{theorem}
There is a one-to-one correspondence
between the CR orbits and their
extrinsic complexifications, namely for all 
$p$ in a neighborhood of the origin, we have\text{\rm \,:}
\begin{itemize}
\item[{\bf (1)}] 
$\mathcal{O}_{CR}(M,p)^c=\mathcal{O}_{\mathcal{L}, \underline{\mathcal{
L}}}(\mathcal{M}, p^c)$, and 
\item[{\bf (2)}] 
$\mathcal{O}_{CR}(M,p)=\pi_t(\underline{\Lambda}\cap
\mathcal{O}_{\mathcal{L}, \underline{\mathcal{L}}} (\mathcal{M}, p^c))$,
\end{itemize}
where we remind the notation $p^c=(p,\bar p)\in\C^n\times \C^n$.  In
particular, if $M$ is a real algebraic generic submanifold, the
complexified orbit $\mathcal{O}_{\mathcal{L},
\underline{\mathcal{L}}}(\mathcal{M},p^c)$ is complex algebraic,
whence by {\bf (2)} the CR orbit $\mathcal{O}_{CR}(M,p)$ is real
algebraic.
\end{theorem}

\proof
By Theorem~2.4.10, $\mathcal{O}_{CR}(M,p)$ is a real analytic
closed submanifold of $M$ passing through $p$. Even if $M$ is real analytic,
the flows of elements of $\L$ are not algebraic in general, as shows
Example~2.6.2, so we do not know more than the analyticity of
$\mathcal{O}_{CR}(M,p)$. Thus, let $\mathcal{O}$ be a small open
connected manifold-piece of $\mathcal{O}_{CR}(M,p)$ through $p$, and
let $\mathcal{O}^c$ be its extrinsic complexification. Because
$L_k|_{\mathcal{O}}$ and $\bar{L}_k|_{\mathcal{O}}$ are tangent to
$\mathcal{O}$, the {\it generic uniqueness principle} ({\it via}
$\mathcal{O}\subset \underline{\Lambda}$, where $\underline{\Lambda}$
is maximally real, so $\mathcal{O}$ is generic in $\mathcal{O}^c$)
entails that $\mathcal{L}_k|_{\mathcal{O}^c}$ and
$\underline{\mathcal{L}}_k|_{\mathcal{O}^c}$ are tangent to
$\mathcal{O}^c$. Therefore $\mathcal{O}^c$ is an integral manifold for
$\{\mathcal{L}, \underline{\mathcal{ L}}\}$ through $p^c$, whence
$\mathcal{O}^c \supset \mathcal{O}_{\mathcal{L},
\underline{\mathcal{L}}}(\mathcal{M}, p^c)$, since a characterizing
property of the orbit $\mathcal{O}_{\mathcal{L},
\underline{\mathcal{L}}}(\mathcal{M}, p^c)$ is to say that it is the
{\it smallest integral manifold-piece} for $\{\mathcal{L},
\underline{\mathcal{ L}}\}$ through $p^c$.

Conversely, Let $\mathcal{N}$ be a manifold-piece of
$\mathcal{O}_{\mathcal{L},
\underline{\mathcal{L}}}(\mathcal{M}, p^c)$ through $p^c$. We have
just shown that $\mathcal{N}\subset \mathcal{O}^c$, hence to finish the
proof, we want to show that $\mathcal{N}\supset \mathcal{O}^c$. We
claim that we have $\sigma(\mathcal{N})=\mathcal{N}$ in a neighborhood of
$p^c$. Indeed, By definition, the orbit is the following set of
endpoints of concatenations of flows of $\mathcal{L}$ and of flows of
$\underline{\mathcal{L}}$ (notice that because $\mathcal{L}_{z_2}
\circ \mathcal{L}_{z_1}= \mathcal{L}_{z_1+z_2}$ and
$\underline{\mathcal{L}}_{\zeta_2} \circ
\underline{\mathcal{L}}_{\zeta_1}=
\underline{\mathcal{L}}_{\zeta_1+\zeta_2}$ but $\mathcal{L}$ and
$\underline{\mathcal{L}}$ do not commute, there can be only two
different kinds of concatenated flow maps\,; we do not use the
abbreviated notation $\Gamma_k$ here)\,:
\def\theequation{2.6.7}\begin{equation}
\left\{
\aligned
{}
&
\mathcal{O}_{\mathcal{L},
\underline{\mathcal{L}}}(\mathcal{M}, p^c)= 
\{\mathcal{L}_{z_k} \circ \cdots \circ \mathcal{L}_{z_2}
\circ
\underline{\mathcal{L}}_{\zeta_1} \circ \mathcal{L}_{z_1} (p^c) : \\
&
\ \ \ \ \ \ \ \ \ \ \ \ \ \ \ \ \ \ \ \ \ \ \ \ \ \ \ \ \ \ 
: \left. z_1, \zeta_1, z_2,\dots,z_k\in \C \ \text{\rm small}, 
k\in \N\right\} \
\bigcup\\
&
\ \ \ \ \ \ \ \ \ \ \ \ \ \ \ \
\bigcup \ \{\underline{\mathcal{L}}_{\zeta_k}
\circ \cdots \circ\underline{\mathcal{L}}_{\zeta_2} \circ \mathcal{L}_{z_1}
\circ \underline{\mathcal{L}}_{\zeta_1} (p^c) :\\
&
\ \ \ \ \ \ \ \ \ \ \ \ \ \ \ \ \ \ \ \ \ \ \ \ \ \ \ \ \ \ \
: \zeta_1,z_1,\zeta_2,\dots,\zeta_k\in \C \ \text{\rm small}, k\in \N\}:=
E\cup F.
\endaligned\right.
\end{equation}
Now an examination of~\thetag{2.4.2} and~\thetag{2.4.3} 
shows that we have $\sigma(\mathcal{
L}_w(q))=\underline{\mathcal{L}}_{\bar{w}}(\sigma(q))$ and
$\sigma(\underline{\mathcal{L}}_{\zeta}(q))=
\mathcal{L}_{\bar{\zeta}}(\sigma(q))$, 
for each $q\in \mathcal{M}$. Consequently\,:
\def\theequation{2.6.8}\begin{equation}
\left\{
\aligned
{}
&
\sigma(\mathcal{L}_{z_k}\circ \underline{\mathcal{L}}_{\zeta_{k-1}} \circ
\cdots \circ \mathcal{L}_{z_1}(p^c))= \underline{
\mathcal{L}}_{\bar{w}_k}(\sigma( 
{\underline{\mathcal{L}}_{\zeta_{k-1}}}\circ
\cdots \circ \mathcal{L}_{z_1} (p^c)))=\\
&
=\underline{\mathcal{L}}_{\bar{w}_k} \circ
\mathcal{L}_{\bar{\zeta}_{k-1}}\circ 
\cdots \circ\sigma(\mathcal{L}_{z_1}(p^c))=
\underline{\mathcal{L}}_{\bar{w}_k} \circ 
\mathcal{L}_{\bar{\zeta}_{k-1}}\circ \cdots \circ \underline{
\mathcal{L}}_{\bar{w}_1}(p^c),
\endaligned\right.
\end{equation}
since $\sigma(p^c)=p^c$. This proves $F=\sigma(E)$, hence
$\sigma(\mathcal{O}_{\mathcal{L}, \underline{ \mathcal{L}}}
(\mathcal{M}, p^c))= \mathcal{O}_{ \mathcal{L}, \underline{
\mathcal{L}}} (\mathcal{M}, p^c)$, so we have
$\sigma(\mathcal{N})=\mathcal{N}$ as announced.
By Theorem~2.4.10, $\mathcal{N}$ is smooth at $p^c$ and
satisfies $\sigma( \mathcal{N})=\mathcal{N}$. To conclude, we need the
following lemma.

\def\thelemma{2.6.9}\begin{lemma}
There is a one-to-one correspondence between real analytic subsets
$\Sigma\subset M$ and complex analytic subvarieties $\Sigma_1$ of
$\mathcal{M}$ satisfying $\sigma(\Sigma_1)=\Sigma_1$ given by
$\Sigma\mapsto \Sigma^c=:\Sigma_1$, with inverse $\Sigma_1 \mapsto
\pi_t(\Sigma_1 \cap \underline{\Lambda})=:\Sigma$.
Furthermore, $\Sigma$ is a smooth submanifold if and only if
$\Sigma_1$ is smooth.
\end{lemma}

\proof
Let $\Sigma\subset M$ be given by real analytic equations
$\chi_l(t,\bar t)=0$, $l=1,\dots,c$.  If $\Sigma$ is smooth, we can
assume that $d\rho_1\wedge \cdots \wedge d\rho_d\wedge d\chi_1\wedge
\cdots\wedge d\chi_l(0)\neq 0$.  Let $\Sigma^c\subset \mathcal{M}$ be
defined by $\chi_l(t,\tau)=0$. Clearly, $\Sigma^c$ is again 
smooth and $\Sigma=\pi_t(\Sigma^c \cap \underline{\Lambda})$.

Conversely, let $\Sigma_1\subset \mathcal{M}$ be given by
$\chi_l(t,\tau)=0$, $l=1,\dots,c$. If $\Sigma_1$ is smooth, we can
assume that $d\rho_1\wedge \cdots \wedge d\rho_d\wedge d\chi_1\wedge
\cdots\wedge d\chi_l(0)\neq 0$. By assumption, $\Sigma_1$ is
$\sigma$-invariant, {\it i.e.} we have $\chi_l(t,\tau)=0$,
$l=1,\dots,c$ if and only if $\chi_l(\bar \tau,\bar t)=0$,
$l=1,\dots,c$. This implies that if we set $\Sigma:=\{t\in M: \,
\chi_l(t,\bar t)=0, \, l=1,\dots,c\}=\pi_t(\Sigma_1\cap
\underline{\Lambda})$, then $\Sigma$ is real, {\it i.e.}  satisfies
$t\in\Sigma$ if and only if $\bar t\in\Sigma$ and satisfies
$(\Sigma)^c=\Sigma_1$.  Finally, $\Sigma$ is smooth if $\Sigma_1$ is.
\endproof

Continuing the proof of Theorem~2.6.6, we now know that there exists
$N:= \pi_t(\mathcal{N} \cap \underline{\Lambda})$ a unique piece of a
real analytic submanifold $N\subset M$ passing through $p$ such that
$N^c=\mathcal{N}$.

Let us denote $\mathcal{N}=\{\rho(t,\tau)=0, \chi(t,\tau)=0\}$, so
that $N=\{\rho(t,\bar{t})=0, \chi(t,\bar{t})=0\}$. Then
$\mathcal{L}_k\rho=0$, $\underline{\mathcal{L}}_k\rho=0$, $\mathcal{L}_k
\chi=0$, $\underline{\mathcal{L}}_k \chi =0$ on $\{\rho=\chi=0\}$, since
$\mathcal{N}$ is an $\{\mathcal{L},
\underline{\mathcal{L}}\}$-integral manifold. Therefore, after
restriction to the antidiagonal
$\{\tau=\bar{t}\}=\underline{\Lambda}$, we have $L_k\rho=0$,
$\bar{L}_k\rho=0$, $L_k\chi=0$ and $\bar{L}_k\chi=0$ on
$\{\rho(t,\bar{t})=0, \chi(t,\bar{t})=0\}=N$, so that $N$ is an
$\{L,\bar{L}\}$-integral manifold. Thus by the minimality property of
CR-orbits, we have $N\supset \mathcal{O}$ as germs at $p$. By
complexifying, we get $\mathcal{N}\supset \mathcal{O}^c$, as desired.
\endproof

Thanks to Theorem~4, an equivalent definition of
minimality is as follows ({\it cf.}~Definition~2.4.8):

\def\thedefinition{2.6.10}\begin{definition}
{\rm 
The generic submanifold $M$ is called {\sl minimal at $p\in M$} if the
CR orbit $\mathcal{O}_{CR}(M,p)$ is of maximal dimension equal to
${\rm dim}_\R \, M$.
}
\end{definition}

\subsection*{2.6.11.~Segre type of $M$}
Now, let us define the maps
$\psi^1(z_{1}):=\pi_t( \Gamma_1 (z_{1}))$,
$\psi^2(z_1,z_2):=\pi_{\tau}(\Gamma_2 (z_1,z_2))$ and more
generally\,:
\def\theequation{2.6.12}\begin{equation}
\psi^{2j}(z_{(2j)}):=
\pi_{\tau}(\Gamma_{2j}(z_{(2j)}))
 \ \ \ \text{\rm and} \ \ \psi^{2j+1}(z_{(2j+1)}):=\pi_t(
\Gamma_{2j+1}(z_{(2j+1)})).
\end{equation}
Notice that by the 
definitions~\thetag{2.4.2} and~\thetag{2.4.3}), the action
of the flow of $\mathcal{L}$ leaves unchanged the
$(\zeta,\xi)$-coordinates, and vice versa, the action of the flow of
$\underline{\mathcal{L}}$ leaves unchanged the
$(z,w)$-coordinates.
Similarly also, we can define the maps $\underline{ \psi}^k$ by
$\underline{ \psi}^{2j}( z_{(2j)}):=\pi_t( \underline{ \Gamma}_{
2j}(z_{(2j)}))$ and $\underline{ \psi}^{ 2j+1}( z_{ (2j+1)}):= \pi_{
\tau}( \underline{ \Gamma}_{ 2j+1}( z_{(2j+1)}))$.  We need the
following lemma (of course, a similar statement also holds with
$\underline \Gamma_{ k+2}$ and $\underline{ \psi}^{ k+1}$ instead):

\def\thelemma{2.6.13}\begin{lemma}
For $0\leq k \leq \mu_0$, we have\text{\rm \,:}
\def\theequation{2.6.14}\begin{equation}
m+ \text{\rm genrk}_{\C} (\psi^{k+1})= \text{\rm
genrk}_{\C} (\Gamma_{k+2})=2m+e_3+\cdots+e_k,
\end{equation}
and $\text{\rm genrk}_{\C} (\psi^{k+1})=
m+e_3+\cdots+e_{\mu_0}$ for $k\geq \mu_0$.
\end{lemma}
 
\proof
For $k=0$, we have $\psi^1(z_1)=(z_1, i \bar \Theta(z_1,0,0))$,
whence $\text{\rm genrk}_{\C} (\psi^1)=m$
obviously. Recall that, by~\thetag{2.4.6}, we have
\def\theequation{2.6.15}\begin{equation}
\Gamma_2(z_1,z_2)=
\left(z_1, \, \overline{\Theta}(z_1,0,0), \, 
z_2, \,
\Theta(z_2,z_1,\overline{\Theta}(z_1,0,0))\right),
\end{equation}
so $m+\text{\rm genrk}_{\C} (\psi^1)= \text{\rm genrk}_{\C}
(\Gamma_2)=2m$. More generally, for $k=2j$, we have:
\def\theequation{2.6.16}\begin{equation}
\left\{
\aligned
{}
&
\mathcal{L}_{z_{2j+1}}(\Gamma_{2j}(z_{(2j)}))=
\mathcal{L}_{z_{2j+1}}(
z(z_{(2j)}), w(z_{(2j)}),\zeta(z_{(2j)}),\xi(z_{(2j)}))\\
&
\ \ \ \ \ \ \ \ \ \ \ \ \ \ \ \ \ \ \ \ \ \ \ \ 
=\left(z_{2j+1}+z(z_{(2j)}), \, \bar{\Theta}(z_{2j+1}+
z(z_{(2j)}), \zeta(z_{(2j)}),
\xi(z_{(2j)})), \, \zeta(z_{(2j)}), \, \xi(z_{(2j)})\right).
\endaligned\right.
\end{equation}
We choose the coordinates $(z,\zeta,\xi)$ on 
$\mathcal{M}$, whence we consider the map $\Gamma_{2j+1}(z_{(2j+1)})$ 
in~\thetag{2.6.16} to have range in
$\C_{(z,\zeta,\xi)}^{2m+d}$. It is then the map
$(z_{(2j)},z_{2j+1})\mapsto \left(z_{2j+1}+z(z_{(2j)}),
\zeta(z_{(2j)}), \, \xi(z_{(2j)})\right)$. It follows
immediately that
\def\theequation{2.6.17}\begin{equation}
\text{\rm genrk}_{\C} (\Gamma_{2j+1})=
m+\text{\rm genrk}_{\C} [z_{(2j)}\mapsto (\zeta(z_{(2j)}),
\xi(z_{(2j)}))]= m+\text{\rm genrk}_{\C} \psi^{2j}. 
\end{equation}
This completes the proof of Lemma~2.6.13.
\endproof

We now define the {\it Segre type of $M$ at $p\in M$} (not to be
confused with $\mu_0$) to be the smallest integer $\nu_0$ satisfying
$\text{\rm genrk}_{\C} (\psi^{\nu_0})=\text{\rm genrk}_{\C}
(\psi^{\nu_0+1})$. By~\thetag{2.6.17}, we readily observe that in
fact, we have $\nu_0=\mu_0-1$. The Segre type {\it of $M$} can be
related to its CR orbits as will be explained in the next
paragraph.

\subsection*{2.6.18.~Intrinsic complexification of CR-orbits}
By the {\sl intrinsic complexification $N^{i_c}$ of a real CR manifold
$N$}, we understand the smallest complex algebraic or analytic
manifold containing $N$ in $\C^n$. It exists and satisfies $\dim_\C
N^{i_c}={\rm CRdim}\, N+{\rm Hcodim}\, N$ ({\sl holomorphic codimension},
{\it cf.}~\S2.1.6). Let $\mathcal{O}$ denote a manifold-piece of
$\mathcal{O}_{CR}(M,p)$ through $p$ and let $\mathcal{O}^{i_c}$ be its
{\sl intrinsic complexification}, namely 
the smallest complex manifold of the ambient space $\C^n$ 
containing $\mathcal{O}$. By construction, the ranges of the
$\psi^{2j}$ are contained in $\C_\tau^n$, but we will prefer to work
in $\C_t^n$ (although it is equivalent in principle to work in
$\C_\tau^n$), hence we shall consider the
$\underline{\psi}^{2j}$ instead. We can now establish that $\text{\rm
genrk}_{\C} (\underline{\psi}^{\nu_0})=\text{\rm dim}_{\C}
\mathcal{O}^{i_c}$ and that the range of $\underline{\psi}^{2\nu_0}$
contains a manifold-piece of $\mathcal{ O}^{i_c}$ through $p$.

\def\thetheorem{2.6.19}\begin{theorem}
There exist some points 
$\underline{z}_{(2\nu_0)}^*\in\C^{m2\nu_0}$ arbitrarily close to
the origin and small neighborhoods $\mathcal{V}^*$ of 
$\underline{z}_{(2\nu_0)}^*$ in
$(\delta\Delta^m)^{2\nu_0}$ such that we have\text{\rm \,:}
\begin{itemize}
\item[{\bf (j)}]
$\underline{\psi}^{2\nu_0}(\underline{z}_{(2\nu_0)}^*)=p$.
\item[{\bf (k)}]
The map $\underline{\psi}^{2\nu_0}$ is of constant rank 
$m+e_3+\cdots+e_{\mu_0}$ in $\mathcal{V}^*$.
\item[{\bf (l)}]
$\underline{\psi}^{2\nu_0} (\mathcal{V}^*)$ is a 
manifold-piece $\mathcal{O}^{i_c}$ of the
intrinsically complexified CR orbit of $M$ through $p$.
\item[{\bf (m)}]
$m+e_3+\cdots+e_{\mu_0}=\text{\rm
dim}_{\C} \, \mathcal{O}^{i_c}=
\text{\rm CRdim} \, \mathcal{O}+{\rm Hcodim} \, \mathcal{O}$.
\end{itemize}
\end{theorem}

\proof
Recall that in view of Theorem~2.5.2, there exists
$\underline{z}_{(2\mu_0-1)}^*\in (\delta\Delta^m)^{2\mu_0-1}$ with
$\underline{\Gamma}_{2\mu_0-1}(\underline{z}_{(2\mu_0-1)}^*)=p^c$,
such that $\underline{\Gamma}_{2\mu_0-1}$ is of rank $2m+e_3+\cdots+e_{\mu_0}$ at
$\underline{z}_{(2\mu_0-1)}^*$. Looking again at~\thetag{2.6.17}
(for $k=2j+1$ odd, which we have not written, but the corresponding
equation is similar), and using the chain rule, we deduce that
$\underline{\psi}^{2\mu_0-2}$ is of rank $m+e_3+\cdots+e_{\mu_0}$ at
the point $\underline{z}_{2\nu_0}^*:=(\underline{z}_1^*,\dots,
\underline{z}_{2\mu_0-2}^*)$ and that
$\underline{\psi}^{2\nu_0}(\underline{z}_{(2\nu_0)}^*)=p$ (recall
$\nu_0=\mu_0-1$). This yields {\bf (j)} and {\bf (k)}. For reasons of
dimension, we already know that $\dim_\C \, \mathcal{O}^{i_c}$ must be
equal to $m+e_3+\cdots+e_{\mu_0}$, since ${\rm CRdim}\, \,
\mathcal{O}=m$ and $\dim_\C \, \mathcal{O}^c=m+\dim_\C \,
\mathcal{O}^{i_c}$. This yields {\bf (m)}. Finally, to
deduce {\bf (l)}, we claim that it can be observed that the range of
$\underline{\psi}^{2\nu_0}$ is {\it a priori} contained in
$\mathcal{O}^{i_c}$, and afterwards for dimensional reasons, the image
$\underline{\psi}^{2\nu_0}(\mathcal{V}^*)$ will necessarily be a
manifold-piece of $\mathcal{O}^{i_c}$ through $p$. To complete this
observation, we introduce holomorphic coordinates $(z,w_1,w_2)\in
\C^m\times \C^{e_3+\cdots+e_{\mu_0}}\times
\C^{n-m-e_3-\cdots-e_{\mu_0}}$ vanishing at $p$ in which the equation
of $\mathcal{O}^{i_c}$ is $\{w_2=0\}$, which is possible. Using the
assumption that $M\cap \{w_2=0\}$ is smooth and of CR dimension $m$,
one shows that the equations of $\mathcal{M}$ can then be written in
the form $w_1=\overline{\Theta}_1(z,\zeta,\xi_1,\xi_2)$ and
$w_2=\overline{\Theta}_2(z,\zeta,\xi_1,\xi_2)]$ with
$\overline{\Theta}_2(z,\zeta,\xi_1,0)\equiv 0$. Then an inspection of
the inductive construction of the maps $\underline{\Gamma}_k$ shows
that they have range contained in $\{w_2=0, \xi_2=0\}$, whence the
maps $\underline{\psi}_{2j}$ have range in $\{w_2=0\}$, as announced.
The proof of Theorem~2.6.19 is complete.
\endproof

\def\theexample{2.6.20}\begin{example}  
{\rm
Looking at the map $\Gamma_4$ in~\thetag{2.5.8},
we see that the integer $2\nu_0=2\mu_0-2$ satisfying the assertions
{\bf (j)} and {\bf (k)} of Theorem~2.6.19 is in general optimal.
}
\end{example}

\section*{\S2.7.~Segre chains and Segre sets in ambient space}

\subsection*{2.7.1.~Segre chains as $k$-th orbit chains of vector fields
} In this section, we shall define certain subsets of $\mathcal{M}$
which are the images of the Segre chains $\Gamma_k$.  These last
results will close up our presentation of the general theory of Segre
chains. Although we shall not use them in the sequel, their definition
seems to be interesting geometrically speaking. Here is an illustration:

\begin{center}
\input S01S02S03.pstex_t
\end{center}

At first, we come back
to the concatenated flow maps in~\thetag{2.4.2} and~\thetag{2.4.3}.
For each $k\in \N$, we choose in advance $\delta_k>0$ so small that
$[\mathcal{L}\, {\rm or}\, \underline{\mathcal{L}}]_{z_k}(\cdots(
\mathcal{L}_{z_1}(p^c))\cdots)$ belongs to $\Delta_{2n}(\rho_1)$ for
all $z_{(k)}\in \Delta_{mk}(\delta_k)$ and all $(t_p,\tau_p)\in
\mathcal{M}$ with $(t_p,\tau_p)\in \Delta_{2n}(\rho_1/2)$.  Looking 
at~\thetag{2.5.2} and~\thetag{2.5.3}, we see that (up to a
shrinking) the complexified Segre varieties of a point
$(t_p,\tau_p)\in\mathcal{M}$ can be defined by $\mathcal{S}_{\tau_p}:=
\{\mathcal{L}_{z_1}(t_p,\tau_p)\in \Delta_{2n}(\rho_1) : z_1\in
\Delta_m(\delta_1)\}$ and $\underline{\mathcal{S}}_{t_p} :=
\{\underline{\mathcal{L}}_{z_1}(t_p,\tau_p)\in\Delta_{2n}(\rho_1) :
z_1\in \Delta_m(\delta_1)\}$. At order $k=2$, we can define\,:
\def\theequation{2.7.3}\begin{equation}
\left\{
\aligned
{}
&
\mathcal{S}_{\tau_p}^2=\{\underline{\mathcal{L}}_{z_2}(
\mathcal{L}_{z_1} (t_p,\tau_p))\in \Delta_{2n}(\rho_1) :
(z_1, z_2)\in \Delta_{2m}(\delta_2)\},\\
&
\underline{\mathcal{S}}_{t_p}^2=\{\mathcal{L}_{z_2}
(\underline{\mathcal{L}}_{z_1}(t_p,\tau_p))\in \Delta_{2n}(\rho_1)
: (z_1,z_2)
\in \Delta_{2m}(\delta_2)\}.
\endaligned\right.
\end{equation}
More generally, we want to define the sets $\mathcal{S}_{\tau_p}^k$ and
$\underline{\mathcal{S}}_{t_p}^k$. By a slight abuse of language, we
shall also call these sets the $k$-th {\sl Segre chain} and the $k$-th
{\sl conjugate Segre chain}. Since we shall only use the mappings
$\Gamma_k$ and $\underline{\Gamma}_k$ and not their images, there will
be no risk of confusion.

\noindent 
At first, we remind
that, because only two ``starting actions''
$\mathcal{L}_{z_1}(t_p,\tau_p)$ and $\underline{\mathcal{
L}}_{\zeta_1}(t_p,\tau_p)$ can make a difference in a concatenation of
flows of $\mathcal{L}$ and of $\underline{\mathcal{L}}$, there
can exist only two different families of $k$-th Segre chains.
The precise definition of $\mathcal{S}_{\tau_p}^k$ and of
$\mathcal{S}_{t_p}^k$ is as follows:
\def\theequation{2.7.4}\begin{equation}
\left\{
\aligned
{}
&
\mathcal{S}_{\tau_p}^{2j}:=\{\underline{\mathcal{L}}_{z_{2j}}
\circ\cdots
\circ\mathcal{L}_{z_1}(t_p,\tau_p): 
\, (z_1,\dots,z_{2j})\in 
\Delta_{2mj}(\delta_{2j})\},\\
&
\mathcal{S}_{\tau_p}^{2j+1}:=\{\mathcal{L}_{z_{2j+1}}
\circ\cdots\circ 
\mathcal{L}_{z_1}(t_p,\tau_p): \, (z_1,\dots,z_{2j+1})\in
\Delta_{2mj+m}(\delta_{2j+1})\},\\
&
\underline{\mathcal{S}}_{t_p}^{2j}:=\{\mathcal{L}_{z_{2j}}\circ
\cdots\circ \underline{\mathcal{L}}_{z_1}(t_p,\tau_p) : \, 
\, (z_1,\dots,z_{2j})\in 
\Delta_{2mj}(\delta_{2j})\},\\
&
\underline{\mathcal{S}}_{t_p}^{2j+1}:=
\{\underline{\mathcal{L}}_{z_{2j+1}}\cdots\circ
\underline{\mathcal{L}}_{z_1}(t_p,\tau_p) : \, (z_1,\dots,z_{2j+1})\in
\Delta_{2mj+m}(\delta_{2j+1})\},
\endaligned\right.
\end{equation}
for $k=2j$ or $k=2j+1$, where $j\in \N$. Clearly, we
have $\mathcal{S}_{\tau_p}^k\subset \mathcal{M}$ and
$\underline{\mathcal{ S}}_{t_p}^k\subset \mathcal{M}$. As
$\sigma(\mathcal{L}_w(q))= \underline{\mathcal{L}}_{\bar
w}(\sigma(q))$, we have $\sigma({\mathcal{S}}_{\tau_p}^k)=
\underline{\mathcal{S}}_{\bar t_p}^k$. 

\subsection*{2.7.5.~Segre sets in ambient space}
We can now define Segre sets in ambient space as certain projections
of Segre chains.
The following picture gives an idea of the definition of Segre sets
as unions of Segre varieties in the case of a minimal hypersurface
in $\C^2$.

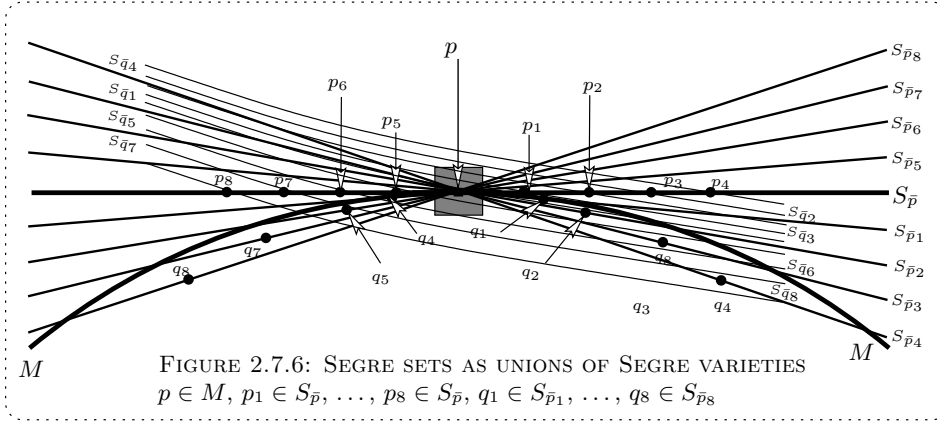
\begin{center}
\input segreset.pstex_t
\end{center} 

The sets
$S_{\bar{t}_p}^{2j+1}:=\pi_t(\mathcal{S}_{\bar{t}_p}^{2j+1}) \subset
\Delta_n(\rho_1)$, $\overline{S}_{t_p}^{2j+1}:=\pi_{\tau}
(\underline{\mathcal{ S}}_{t_p}^{2j+1})\subset \Delta_n(\rho_1)$,
$S_{\bar{t}_p}^{2j}:=\pi_{\tau}( \mathcal{S}_{\bar{t}_p}^{2j})\subset
\Delta_n(\rho_1)$ and $\overline{S}_{t_p}^{2j}:=
\pi_t(\underline{\mathcal{S}}_{t_p}^{2j})\subset \Delta_n(\rho_1)$ will be
called the {\sl $k$-th Segre sets and the $k$-th conjugate Segre
sets}, with $k=2j$ or $k=2j+1$. Notice that by the 
definitions~\thetag{2.4.2} and~\thetag{2.4.3}), the action
of the flow of $\mathcal{L}$ leaves unchanged the
$(\zeta,\xi)$-coordinates, and vice versa, the action of the flow
$\underline{\mathcal{L}}$ leaves unchanged the
$(z,w)$-coordinates. This is why in the definition of Segre sets,
we alternately project in the $\C_t^n$-space and in the
$\C_\tau^n$-space. 

An equivalent, purely set-theoretic definition of Segre sets is as
follows.  We define\,: $S_{\bar{t}_p}^0:=\{\bar{t}_p\}$, and
$S_{\bar{t}_p}^1:=S_{\bar{t}_p}= \bigcup_{\bar{t}\in S_{\bar{t}_p}^0}
S_{\bar{t}}$, $S_{\bar{t}_p}^2= \bigcup_{t\in S_{\bar{t}_p}^1}
\overline{S}_t$, and then inductively, for $j\in \N_*$,
$S_{\bar{t}_p}^{2j} =\bigcup_{t\in S_{\bar{t}_p}^{2j-1}}
\overline{S}_t$ and $S_{\bar{t}_p}^{2j+1}=\bigcup_{\bar{t}\in
S_{\bar{t}_p}^{2j}} S_{\bar{t}}$. On the other hand, we also define
$\overline{S}_{t_p}^0:=\{t_p\}$, and
$\overline{S}_{t_p}^1:=\overline{S}_t=\bigcup_{t\in
\overline{S}_{t_p}^0} \overline{S}_t$,
$\overline{S}_{t_p}^2:=\bigcup_{\bar{t}\in \overline{S}_{t_p}^1}
S_{\bar{t}}$, and inductively, for $j\in \N_*$,
$\overline{S}_{t_p}^{2j}:=\bigcup_{\bar{t}\in
\overline{S}_{t_p}^{2j-1}} S_{\bar{t}}$, and
$\overline{S}_{t_p}^{2j+1} = \bigcup_{t\in \overline{S}_{t_p}^{2j}}
\overline{S}_t$.
Finally, we mention the following two elementary properties\,:

\begin{itemize}
\item[{\bf (1)}]
$\overline{S_{\bar{t}_p}^k}=\overline{S}_{t_p}^k$ and
$S_{\bar{t}_p}^k=\overline{\overline{S}_{t_p}^k}$, $k\in \N$.
\item[{\bf (2)}]
$\overline{h}(S_{\bar{t}_p}^{2j})=
S_{\overline{h}(\bar{t}_p)}^{'2j}$,
$h(S_{\bar{t}_p}^{2j+1})=S_{\overline{h}(\bar{t}_p)}^{'2j+1}$,
$h(\overline{S}_{t_p}^{2j})=\overline{S}_{h(t_p)}^{'2j}$,
$\overline{h}(\overline{S}_{t_p}^{2j+1})=
\overline{S}_{h(t_p)}^{'2j+1}$.
\end{itemize}

\section*{\S2.8.~Generic Segre multitype}

\subsection*{2.8.1.~Segre chains with varying base point}
As before, let $M$ be a connected generic real algebraic or
analytic submanifold of $\C^n$. Let $p_0\in M$.  In a neighborhood of
$p_0$, we can consider the pair of Segre foliations of the local
complexification $\mathcal{M}$ of $M$. Let $p=(t_p,\tau_p)\in
\mathcal{M}$. Identifying the point $p_0$ with the origin in some system
of coordinates, we have denoted by $\Gamma_k(z_{(k)})$ the mapping
$[\mathcal{L}\, {\rm or} \, \underline{\mathcal{L}}]_{z_k}\circ \cdots
\circ \mathcal{L}_{z_1}(p_0,\bar p_0)$, namely the Segre chain with
base point $(p_0)^c=(0,0)\in\C^n\times
\C^n$. We want to let the base point vary, so we need a
new notation.  For $p=(t_p,\tau_p)$ in a neighborhood of $(p_0)^c$ in
$\mathcal{M}$, we define 
\def\theequation{2.8.2}\begin{equation}
\Gamma_k(z_{(k)},t_p,\tau_p):=
[\mathcal{L}\, {\rm or} \, 
\underline{\mathcal{L}}]_{z_k}\circ \cdots
\circ \mathcal{L}_{z_1}(t_p,\tau_p).
\end{equation}
Similarly, we define $\underline{\Gamma}_k (z_{(k)},t_p,\tau_p)$.
In order to indicate the dependence of the Segre type with respect to 
$p$, we shall denote it by $\mu_p$. Also, we shall denote the
Segre multitype at $p$ by $(m,m,e_{3,p},\dots,
e_{\mu_p,p})$.

By Theorem~2.5.2., the generic rank of $\Gamma_k$ stabilizes when $k\geq
\mu_p$. We know that
$\mu_p\leq d+2$ for all $p$ in a neighborhood of $p_0$ and that the
mapping $z_{(k)}\mapsto \Gamma_k(z_{(k)},t_p,\tau_p)$ provides an open
piece of the $\{\mathcal{L}, \underline{\mathcal{L}}\}$-orbit through
$p\in \mathcal{M}$ for all $k\geq 2\mu_p-1$. As $2\mu_p-1\leq 2d+3$,
and as $2d+3\leq 3(m+d)=3n$, because $m\geq 1$ and $d\geq 1$ by assumption, we
observe that the mapping $\Gamma_{3n}(z_{(3n)}, t_p,\tau_p)$ with the
uniform integer $k=3n$ suffices to construct the $\{\mathcal{L},
\underline{\mathcal{L}}\}$-orbits of all points $p$ in a neighborhood
of $(p_0)^c$ in $\mathcal{M}$.

Thanks to this observation and thanks to the algebraicity or the
analyticity of the mapping $\Gamma_{3n}(z_{(3n)}, t_p,\tau_p)$, it is
easy to see that there is a proper complex algebraic or complex
analytic subvariety $\mathcal{E}$ of $\mathcal{M}$ with the property
that the Segre type and multitype of $\mathcal{M}$ are constant at
every point $p\in \mathcal{M}\backslash \mathcal{E}$. We shall denote
these constants by $(m,m,e_{3,M},\dots,e_{\mu_M,M})$, where $\mu_M$ is
the constant Segre type of $\mathcal{M}$ outside $\mathcal{E}$.  In
particular, the orbit dimension $2m+e_{1,M}+\cdots+ e_{\mu_M,M}$ is
constant in a neighborhood of $p$. Moreover, it also follows from the
algebraicity or analyticity of the mappings $\Gamma_{3n}(z_{(3n)},
t_p,\tau_p)$ that the functions $p\mapsto \mu_p$ and $p\mapsto
e_{3,p}, \dots, p\mapsto e_{\mu_p,p}$ are lower semi-continuous.
Finally, using the $\sigma$-invariance of the CR orbits of
complexified points $p^c=(t_p,\bar t_p)\in \mathcal{M}\cap
\underline{\Lambda}$, with $p\in M$, we get the following theorem,
which states that the Segre geometry possesses constant invariants
over a Zariski open subset of $M$.

\def\thetheorem{2.8.3}\begin{theorem}
Let $M$ be a connected real algebraic or analytic generic
submanifold of $\C^n$ of codimension $d\geq 1$ and of CR dimension
$m=n-d\geq 1$. Then there is a proper real algebraic or analytic
subvariety $E$ of $M$ such that for every point 
$p\in M\backslash E$, the Segre type of $M$ at 
$p$ is constant equal to an integer $\nu_M=\mu_M-1\leq d+1$  and
the Segre multitype of $M$ at $p$ is also constant equal to the
multiplet $(m,e_{3,M},\dots,e_{M,\mu_M})$. In particular, 
the CR-orbit dimension 
$\dim_\R \, \mathcal{O}_{CR}(M,p)$ 
is constant equal to $2m+e_{3,M}+\cdots+
e_{\mu_M,M}$ for all $p\in M\backslash E$.
\end{theorem}

We call $\mu_M$ the {\sl generic Segre type of $\mathcal{M}$} and the
multiplet $(m,m,e_{3,M},\dots,e_{\mu_M,M})$ the
{\sl generic Segre multitype of $\mathcal{M}$}.
Let 
\def\theequation{2.8.4}\begin{equation}
d_M:= e_{3,M}+\cdots+e_{\mu_M,M}.
\end{equation}
We call the integer $2m+d_M$ the {\sl Zariski-generic orbit dimension}
of $M$.  We call the integer $d-d_M$ the {\sl Zariski-generic orbit
codimension} of $M$. Then using again the mapping
$\Gamma_{3n}(z_{(3n)}, t_p,\tau_p)$, we can derive the following
algebraic or analytic CR foliation theorem which shows that $d-d_M$
coincides with the Zariski-generic holomorphic codimension of the
intrinsic complexification of CR orbits.

\def\thecorollary{2.8.5}\begin{corollary}
Let $p\in M\backslash E$ and set $d_{2,M}:=d-e_{3,M}-\cdots-e_{\mu_M,M}$.
Then a neighborhood of $p$ in $M$ is real algebraically or
analytically foliated by CR orbits, namely there exist $d_{2,M}$ complex
algebraic or analytic functions $h_1,\dots,h_{d_{2,M}}$ with
$\partial h_1\wedge \cdots \wedge \partial h_{d_{2,M}}(p)\neq 0$ such that
\begin{itemize}
\item[{\bf (1)}]
$M$ is contained in $\{h_1=\bar h_1,\, \dots, \, h_{d_{2,M}}=\bar
h_{d_{2,M}}\}$. In other words, $M$
is contained in a transverse intersection of
$d_{2,M}$ Levi flat hypersurfaces in 
general position.
\item[{\bf (2)}]
For every $c=(c_1,\dots,c_{d_{2,M}})\in\R^{d_{2,M}}$,
the manifold $M_c:=M\cap \{h_1=c_1, \, \dots, h_{d_{2,M}}=c_{d_{2,M}}\}$ is a
CR orbit of $M$.
\end{itemize}
\end{corollary}

\proof
Thanks to the mapping $\Gamma_{3n}(z_{(3n)}, t_p,\tau_p)$, we find
real algebraic or analytic functions $h_1,\dots,h_{d_{2,M}}$ with
linearly independent real differentials such that the level sets
$\{h_1=c_1, \, \dots, h_{d_{2,M}}=c_{d_{2,M}}\}$
are the CR orbits of $M$ in a
neighborhood of $p$. Since the functions $h_1,\dots,h_{d_{2,M}}$ are
constant in each CR orbit, they are in particular trivially CR.  By
the Severi-Tomassini extension theorem, they extend complex
algebraically or analytically to a neighborhood of $p$ in $\C^n$. This
proves the corollary.
\endproof

Taking the functions $h_1,\dots,h_{d_{2,M}}$ as part of a system of
complex coordinates and applying Theorem~2.1.32, we deduce:

\def\thecorollary{2.8.6}\begin{corollary} 
For every point $p\in M\backslash E$, there exist complex algebraic or
analytic local normal coordinates $(z,w_1,w_2)\in\C^m\times
\C^{d-d_{2,M}}\times \C^{d_{2,M}}$ vanishing at $p$ such that the
complex defining equations of $M$ are of the form
\def\theequation{2.8.7}\begin{equation}
\left\{
\aligned
0
& \
=\bar w_2-w_2,\\
0
& \
=\bar w_1-\Theta_1(\bar z,z,w_1,w_2),
\endaligned\right.
\end{equation}
where $\Theta_1(0,z,w_1,w_2)\equiv w_1$, where for $u_{2,q}\in
\R^{d_{2,M}}$ sufficiently small, the sets $M\cap \{w_2=u_{2,q}=ct.\}$
coincide with the local CR orbit of the points $q=(0, 0,u_{2,q})\in M$
and where the generic submanifold of $\C^{m+d-d_{2,M}}$ defined by the
equations
\def\theequation{2.8.8}\begin{equation}
0 = \bar w_1-\Theta_1(\bar z, z, w_1,u_{2,q})
\end{equation}
is minimal at $(z,w_1)=(0,0)$ for every $u_{2,q}$.
\end{corollary}

\section*{\S2.9.~Local representation of nonminimal generic submanifolds}

As a conclusion, we can now produce a general summary of important
results which we will use constantly in the sequel.  In advance, we
formulate them in the most appropriate way for later use.  As in
Definition~2.1.44, let $M\subset \C^n$ be a local generic submanifold.

\subsection*{2.9.1.~Minimal generic submanifolds}
The following theorem is a corollary of Theorem~2.5.2 and
Theorem~2.6.19 in the minimal case
where $2m+e_2+\dots+e_{\mu_0}=2m+d$.  For later applications to the
study of CR mappings, it is more convenient to state it with the
conjugate Segre chain $\underline{ \Gamma}_{ 2\nu_0}$.

\def\thetheorem{2.9.2}\begin{theorem}
If $M$ is minimal at $p_0$, there exists a positive integer $\nu_0\leq
d+1$, the {\sl Segre type of $M$ at $p_0$}, there exists an element
$\underline{z}_{ (2\nu_0)}^*\in \C^{ 2m\nu_0}$ arbitrarily close to
the origin, there exists a $n$-dimensional complex affine subspace
$H^*$ passing through $\underline{z}_{ (2\nu_0)}^*$ in $\C^{
2m\nu_0}$ and there exists a complex affine
parametrization $s\mapsto \underline{z}_{(2\nu_0)}(s)$ of $H^*$ with
$\underline{z}_{(2\nu_0)}(0)= \underline{z}_{(2\nu_0)}^*$ such that
the mapping defined by composing the projection onto the
first factor with the $(2\nu_0)$-th Segre chain, namely
the mapping
\def\theequation{2.9.3}\begin{equation}
\C^n\ni s\longmapsto \pi_t
(\underline{\Gamma}_{2\nu_0}( \underline{ z}_{
(2\nu_0)}(s)))\in
\C^n
\end{equation}
is of rank $n$ and vanishes at $s=0$.
\end{theorem}

We shall use this formulation very frequently in Part~II of this
memoir.

\subsection*{2.9.4.~General generic submanifolds}
In the case where $M$ is not necessarily minimal, the holomorphic
codimension in $\C^n$ of the local CR orbit $\mathcal{O}_{CR}(M,p_0)$
is an arbitrary integer $d_2$ with $0\leq d_2 \leq d$ and $d_2=0$ if
and only if $M$ is minimal at $p_0$. We set $d_1:=d-d_2$, so
$\mathcal{O}_{CR}(M,p_0)$ is of dimension $2m+d_1$.  By
Theorem~2.6.20, the intrinsic complexification
$[\mathcal{O}_{CR}(M,p_0)]^{i_c}$ is a complex algebraic or analytic
CR submanifold of $\C^n$ passing through $p_0$ which is of complex
codimension $d_2$.  After straightening it, we can assume that in
coordinates $(z,w_1,w_2)\in\C^m\times \C^{d_1} \times \C^{d_2}$, it
coincides with $\{z_2=0\}$, so there are local defining equations for
$M$ of the form
\def\theequation{2.9.5}\begin{equation}
\left\{
\aligned
{}
&
\bar w_{1,j_1}
& \
=
& \
\Theta_{1,j_1}(\bar z, z, w_1,w_2),  \ \ \ \ \ 
j_1=1,\dots,d_1,\\
&
\bar w_{2,j_2}
& \
=
& \
\Theta_{2,j_2}(\bar z, z, w_1, w_2), \ \ \ \ \
j_2=1,\dots,d_2, 
\endaligned\right.
\end{equation}
where $\Theta_{2,j_2}(\bar z,z,w_1,0)\equiv 0$ and where the 
generic submanifold $M_1$ of $\C^{m+d_1}$ defined by
\def\theequation{2.9.6}\begin{equation}
M_1:=M\cap \{w_2=0\} 
\end{equation}
is minimal at the origin, with Segre type equal to $\nu_0$.
Complexifying $(M)^c:=\mathcal{M}$, we obtain the equations
\def\theequation{2.9.7}\begin{equation}
\left\{
\aligned
{}
&
\xi_{1,j_1}
& \
=
& \
\Theta_{1,j_1}(\zeta, z, w_1,w_2),  \ \ \ \ \ 
j_1=1,\dots,d_1,\\
&
\xi_{2,j_2}
& \
=
& \
\Theta_{2,j_2}(\zeta, z, w_1, w_2), \ \ \ \ \
j_2=2,\dots,d_2, 
\endaligned\right.
\end{equation}
for $k=1,\dots,m$ 
and the complexified $(1,0)$ vector fields
\def\theequation{2.9.8}\begin{equation}
\left\{
\aligned
\mathcal{L}_k:=
{\partial \over \partial z_k}+
& \
\sum_{j_1=1}^{d_1}\, 
{\partial \Theta_{1,j_1}\over \partial z_k}
(\zeta,z,w_1,w_2)\, 
{\partial \over \partial w_{1,j_1}}+\\
& \
+\sum_{j_2=1}^{d_2}\, 
{\partial \Theta_{2,j_2}\over \partial z_k}
(\zeta,z,w_1,w_2)\, 
{\partial \over \partial w_{2,j_2}},
\endaligned\right.
\end{equation} 
for $k=1,\dots,m$ 
and the complexified $(1,0)$ vector fields
\def\theequation{2.9.9}\begin{equation}
\left\{
\aligned
\underline{\mathcal{L}}_k:=
{\partial \over \partial \zeta_k}+
& \
\sum_{j_1=1}^{d_1}\, 
{\partial \overline{\Theta}_{1,j_1}\over \partial \zeta_k}
(z,\zeta,\xi_1,\xi_2)\, 
{\partial \over \partial \xi_{1,j_1}}+\\
& \
+\sum_{j_2=1}^{d_2}\, 
{\partial \overline{\Theta}_{2,j_2}\over \partial \zeta_k}
(z,\zeta,\xi_1,\xi_2)\, 
{\partial \over \partial \xi_{2,j_2}},
\endaligned\right.
\end{equation} 
In the ambient space $\C^{2n}$ of the complexification 
$\mathcal{M}$, we shall denote
the six coordinates in $\C^m\times \C^{d_1}\times \C^{d_2}\times
\C^m\times \C^{d_1} \times \C^{d_2}$ by  
\def\theequation{2.9.10}\begin{equation}
(z,w_1,w_2,\zeta,\xi_1,\xi_2).
\end{equation} 
In $\C^{2n}$, the set
\def\theequation{2.9.11}\begin{equation}
\mathcal{T}:=\{(0,0,w_2,0,\Theta_1(0,0,0,w_2),
\Theta_2(0,0,0,w_2),w_2)\} 
\end{equation}
is a transversal in $\mathcal{M}$ to 
the complexification $\mathcal{M}_1:=(M_1)^c$ given by 
\def\theequation{2.9.12}\begin{equation}
\mathcal{M}_1: \ \ 
w_2=\xi_2=0, \ \ \ \ \ 
\xi_1=\Theta_1(\zeta,
w_1,0),
\end{equation}
namely we have $T_0\mathcal{M}_1\oplus T_0\mathcal{T}=
T_0\mathcal{M}$. Of course, this transversal depends on the choice of
coordinates. To simplify a bit the expression of a choice of
$\mathcal{T}$, we can (without loss of generality) assume that the
coordinates $(z,w_1,w_2)$ are normal, as described in Theorem~2.1.32,
hence $\Theta_1(0,z,w_1,w_2)\equiv w_1$ and
$\Theta_2(0,z,w_1,w_2)\equiv w_2$. Then
\def\theequation{2.9.13}\begin{equation}
\mathcal{T}=\{(0,0,w_2,0,0,w_2)\}.
\end{equation}

With this choice, we may now generalize the definition of Segre chains
by including the transversal parameter $w_2$ as follows. 
Firstly, for $z_{(1)}\in \C^m$, we set
\def\theequation{2.9.14}\begin{equation}
\left\{
\aligned
\underline{\Gamma}_1(z_{(1)}: w_2):=
& \
\underline{\mathcal{L}}_{z_1}(0,0,w_2,0,0,w_2)\\
= 
& \
(0,0,w_2,z_1,\Theta_1(z_1,0,0,w_2),
\Theta_2(z_1,0,0,w_2)).
\endaligned\right.
\end{equation}
Secondly, for $z_{(2)}=(z_1,z_2)\in \C^{2m}$, we set 
\def\theequation{2.9.15}\begin{equation}
\left\{
\aligned
{}
&
\underline{\Gamma}_2(z_{(2)}: w_2)
& \
:=
& \
\mathcal{L}_{z_2}(\Gamma_1(z_1: w_2)), \\
&
\underline{\Gamma}_3(z_{(3)}: w_2)
& \
:=
& \
\underline{\mathcal{L}}_{z_3}(
\underline{\Gamma}_2(z_{(2)}: w_2)),
\endaligned\right.
\end{equation}
and so on by induction. As a slight generalization of Theorem~2.5.2,
we have the following theorem which describes the local Segre chain
geometry in a neighborhood of an arbitrary point $p_0$ of $M$, without
any nondegeneracy condition on $M$, in the most general setting.

\def\thetheorem{2.9.16}\begin{theorem}
If $d_2$ denotes the holomorphic codimension of the CR orbit of
$p_0$ in $\C^n$ and if the coordinates $(z,w_1,w_2)\in\C^m \times
\C^{d_1}\times \C^{d_2}$ are chosen such that $M_1:=M\cap \{w_2=0\}$
is the CR orbit of $p_0$, if the integer $\nu_0$ with $\nu_0\leq
d_2+1$ denotes the Segre type of $M$ at $p_0$, then there exists an
element $\underline{z}_{ (2\nu_0)}^*\in \C^{ 2m\nu_0}$ arbitrarily
close to the origin, there exists a $(n-d_2)$-dimensional complex
affine subspace $H^*$ passing through $\underline{z}_{ (2\nu_0)}^*$ in
$\C^{ 2m\nu_0}$, there exists a complex affine parametrization
$s\mapsto \underline{z}_{(2\nu_0)}(s)$ of $H^*$ with
$\underline{z}_{(2\nu_0)}(0)= \underline{z}_{(2\nu_0)}^*$ such that
the projection of the conjugate Segre chains with origin the
transversal $\mathcal{T}$ to the complexification $\mathcal{M}_1$,
namely
\def\theequation{2.9.17}\begin{equation}
\C^{n-d_2}\times \C^{d_2}\ni (s,w_2)\longmapsto \pi_t
(\underline{\Gamma}_{2\nu_0}( \underline{ z}_{
(2\nu_0)}(s): w_2
))=:t\in
\C^n
\end{equation}
is of rank $n$ and vanishes at $(s,w_2)=(0,0)$.
\end{theorem}

In particular, if $M$ is real algebraic or real analytic, if the local
CR orbits of points $p$ varying in a neighborhood of $p_0$ are all of
holomorphic codimension equal to $d_2$, it follows from
Corollary~2.8.6 that we can represent $M$ in a neighborhood of $p_0$ by
the equations
\def\theequation{2.9.18}\begin{equation}
\left\{
\aligned
{}
&
\bar w_{1,j_1}
& \
=
& \
\Theta_{1,j_1}(\bar z,z,w_1,w_2), \ \ \ \ \
& \ 
j_1=1,\dots,d_1, \\
&
\bar w_{2,j_2}
& \ 
=
& \
w_{2,j_2}, \ \ \ \ \ \ \ \ \ \ \ \ \ \ \ \ \ \ 
& \
j_2=1,\dots,d_2,
\endaligned\right.
\end{equation}
where the last two equations represent the transversal intersection of
$d_2$ Levi-flat hyperplanes in general position. 

\def\thecorollary{2.9.19}\begin{corollary}
If $M$ is real algebraic or analytic and if the orbit codimension
is constant in a neighborhood of $p_0$, then for every $u_2\in\R^{d_2}$ fixed,
the image of the mapping
\def\theequation{2.9.20}\begin{equation}
\C^{n-d_2}\ni s\longmapsto \pi_t
(\underline{\Gamma}_{2\nu_0}( \underline{ z}_{
(2\nu_0)}(s): u_2
))=:t\in
\C^n
\end{equation}
covers a local piece of the intrinsic complexification
of the CR orbit of the point in $M$ with coordinates 
$(0,0,u_2)$.
\end{corollary}

\bigskip

\noindent
{\Large \bf Chapter~3: Nondegeneracy conditions for generic submanifolds}

\section*{\S3.1.~Segre mapping}

\subsection*{3.1.1.~Definition}
Let $M$ be a connected generic submanifold of $\C^n$ of codimension
$d\geq 1$ and of CR dimension $m=n-d\geq 1$ and let $p_0\in M$.  As
provided by Theorem~2.1.9, we choose coordinates
$t=(t_1,\dots,t_n)=(z_1,\dots,z_m,w_1,\dots,w_d)\in\C^m\times \C^d$
vanishing at $p_0$ in which $M$ is represented by the $d$ complex
defining equations 
\def\theequation{3.1.2}\begin{equation} \bar w_j=\Theta_j(\bar
z,t)=\Theta_j(\bar z,z,w), \ \ \ \ \ j=1,\dots,d.
\end{equation}
We remind that for every choice of coordinates $(z,w)$ vanishing at
$p_0$ such that $T_{p_0}^cM\cap (\{0\}\times\C_w^d)=\{0\}$, there exists a
unique collection power series $\Theta_j(\bar z,t)$ such that $M$ is
represented by~\thetag{3.1.2}.  Here, we shall assume that the powers
series $\Theta_j(\bar z,z,w)$ are complex algebraic or analytic,
namely they belong to $\mathcal{A}_\C\{\bar z,z,w\}$ or to $\C\{\bar
z,z,w\}$. In this section, we shall only work at the central point $p_0$,
which is the origin in these coordinates.

By developing the series $\Theta_j(\bar z,t)$ in powers of $\bar z$,
we may write $\bar w_j=\sum_{\beta\in\N^m}\, \bar z^\beta\,
\Theta_{j,\beta}(t)$. In terms of such a development, the {\sl
infinite Segre mapping of $M$} is defined to be the mapping
\def\theequation{3.1.3}\begin{equation}
\mathcal{Q}_\infty: \ \ \C^n\ni t \longmapsto
(\Theta_{j,\beta}(t))_{1\leq j\leq d,\, 
\beta\in\N^m}\in \C^\infty.
\end{equation}
Let $k\in\N$. For finiteness reasons, 
it is convenient to truncate this infinite
collection and to define the 
{\sl $k$-th Segre mapping of $M$} by
\def\theequation{3.1.4}\begin{equation}
\mathcal{Q}_k: \ \ \C^n\ni t \longmapsto
(\Theta_{j,\beta}(t))_{1\leq j\leq d,\, 
\vert\beta\vert\leq k}\in \C^{N_{d,n,k}},
\end{equation}
where the integer $N_{d,n,k}$ denotes the number of $k$-th jets of a
$d$-vectorial mapping of $n$ independent variables $(t_1,\dots,t_n)$,
namely $N_{d,n,k}=d{(n+k)!\over n! \ k!}$. Let $k_2\geq k_1$ and let
$\pi_{k_2,k_1}$ denote the canonical projection $\C^{N_{d,n,k_2}}\to
\C^{N_{d,n,k_1}}$. Then we obviously have
$\pi_{k_2,k_1}(\mathcal{Q}_{k_2}(t))= \mathcal{Q}_{k_1}(t)$.

We shall see that these Segre mappings $\mathcal{Q}_k$ and
$\mathcal{Q}_\infty$ are of utmost importance among the biholomorphically 
invariant objects attached to a real algebraic or analytic generic
submanifold $M$.

\subsection*{3.1.5.~Transformation of the Segre mapping under a change of 
coordinates} Apparently, the definition of the mappings
$\mathcal{Q}_k$ strongly depends on the choice of coordinates and so
the $\mathcal{Q}_k$ do not seem to represent an invariant
analytico-geometric concept. However, we shall establish some
canonical transformation rules which will show that all the
definitions provided in this chapter are biholomorphically
invariant.

The necessary ingredients for a biholomorphic transformation are as
follows. Let $t'=h(t)= (h_1(t),\dots,h_n(t))$ be an invertible
transformation, where the series $h_i(t)$ belong to
$\mathcal{A}_\C\{t\}$ or to $\C\{t\}$, satisfy $h_i(0)=0$ and ${\rm
det} \, ([\partial h_{i_1}/\partial t_{i_2}](0)_{1\leq i_1,i_2\leq
n})\neq 0$. Let $t=h'(t')$ denote the inverse mapping and split the
mapping $h'=(f',g')= (f_1',\dots,f_m',g_1',\dots,g_d')$ according to
the splitting $t=(z,w)$ of the coordinates $t$. Furthermore,
substitute $t$ by $h'(t')$ and $\bar t$ by $\bar h'(\bar t')$
in~\thetag{3.1.2}, which yields
\def\theequation{3.1.6}\begin{equation}
\bar g_j'(\bar t')=\Theta_j(\bar f'(\bar t'),h'(t')), \ \ \ \ \
j=1,\dots,d.
\end{equation}
If necessary, we renumber the coordinates in order that after
splitting $t'=(z_1',\dots,z_m',w_1',\dots,w_d')$, after applying the
algebraic or the analytic implicit function theorem, we can solve
$\bar w'$ in terms of $(\bar z',z',w')$ in the following form which is
analogous to~\thetag{3.1.2}:
\def\theequation{3.1.7}\begin{equation}
\bar w_j'=\Theta_j'(\bar z',t')=
\Theta_j'(\bar z',z',w'), \ \ \ \ \
j=1,\dots,d.
\end{equation}
This yields an algebraic or an analytic
generic submanifold $M'$ of $\C^n$.  Finally, we develope these series in
powers of $\bar z'$, which yields $\bar w_j'=\sum_{\beta\in\N^m}\,
(\bar z')^\beta\, \Theta_{j,\beta}'(t')$, and we define the {\sl transformed
$k$-th Segre mapping of $M'$} by
\def\theequation{3.1.8}\begin{equation}
\mathcal{Q}_k': \ \ \C^n\ni t' \longmapsto
(\Theta_{j,\beta}'(t'))_{1\leq j\leq d,\, 
\vert\beta\vert\leq k}\in \C^{N_{d,n,k}},
\end{equation}
and we also define the {\sl transformed infinite Segre mapping
$\mathcal{Q}_\infty'$ of $M'$}.  With such notations at hand, we can
now state the following important transformation rules which we will
prove in Section~3.6 below.

\def\thetheorem{3.1.9}\begin{theorem}
For every $j=1,\dots,d$ and every $\beta\in\N^m$, there exists a
mapping $R_{j,\beta}$ which is complex algebraic or analytic in its
variables such that 
\def\theequation{3.1.10}\begin{equation}
\left\{
\aligned
\Theta_{j,\beta}'(h(t))+\sum_{\gamma\in\N^m\backslash \{0\}}\,
{(\beta+\gamma)!\over \beta! \ \gamma!}\,
(\bar f(0,\Theta(0,t)))^\gamma\, 
\Theta_{j,\beta+\gamma}'(h(t))\equiv \\
\equiv
Q_{j,\beta}(\{
\Theta_{j_1,\beta_1}(t)
\}_{1\leq j_1\leq d,\, \vert \beta_1\vert \leq \vert
\beta \vert}).
\endaligned\right.
\end{equation}
Here, the left hand side is the power series developement of
${1\over \beta!}\,[\partial_{\zeta'}^\beta\Theta_{j,\beta}'](
\bar f(0,\Theta(0,t)),h(t))$. A collection of
relations equivalent to the collection~\thetag{3.1.10}
is as follows
\def\theequation{7.1.11}\begin{equation}
\left\{
\aligned
\Theta_{j,\beta}'(h(t))\equiv 
& \
\sum_{\gamma\in\N^m}\, 
(-1)^\gamma\,
{(\beta+\gamma)! \over \beta ! \ \gamma!}\, 
(\bar f(0,\Theta(0,t)))^\gamma\, 
Q_{j,\beta+\gamma}(\{
\Theta_{j_1,\beta_1}(t)
\}_{1\leq j_1\leq d, \, 
\vert \beta_1 \vert \leq \vert \beta \vert+
\vert \gamma \vert}
)\\
\equiv 
& \
R_{j,\beta}(\{\Theta_{j_1,\beta_1}(t)\}_{1\leq 
j_1\leq d, \, \beta_1\in\N^m}).
\endaligned\right.
\end{equation}
for every $j=1,\dots,d$ and every $\beta\in\N^m$.
\end{theorem}

Admitting this theorem, we shall verify in Section~3.3 below that every
nondegeneracy condition on $M$ which is defined in terms of the Segre
mapping $\mathcal{Q}_\infty$ is invariant under changes of coordinates,
namely such a condition is satisfied for $(M, \mathcal{Q}_\infty)$ if and
only if it is satisfied for $(M', \mathcal{Q}_\infty')$. First of all, we
present five such nondegeneracy conditions.

\section*{3.2.~Five pointwise nondegeneracy conditions}

\subsection*{3.2.1.~Levi nondegeneracy} 
To begin with, we shall say that $M$ is {\sl Levi nondegenerate at}
$p_0$ if the first Segre mapping $\mathcal{Q}_1$ is of rank $n$ at the
origin.

We verify that this definition coincides with the usual one, firstly
in the case $d=1$. After diagonalizing the Hermitian matrix of its
Levi form, a given hypersurface $M$ passing through the origin in
$\C^n$ may be represented by the equation $\bar w= w+i\sum_{k=1}^r\,
\varepsilon_k\, \vert z_k\vert+ {\rm O}(3)$, where $0\leq r\leq n-1$
is the rank of the Levi form of $M$ at the origin and where
$\varepsilon_k=\pm 1$.  Then $\mathcal{Q}_1(t)=(w, i \varepsilon_1
z_1,\dots, i\varepsilon_r z_r,0,\dots,0)+{\rm O}(2)$. Clearly, the
rank of $\mathcal{Q}_1$ at $0$ equals $n$ if and only if $r=n-1$, as
announced.

Secondly, in the case $d\geq 2$, a generic submanifold of codimension
$d$ may be represented by equations $\bar w_j=w_j +i\, \sum_{ k_1,
k_2=1}^m\, a_{j,k_1,k_2} \, z_{k_1}\, \bar z_{k_2} +{\rm O}(3)$, where
$a_{j,k_1,k_2}=\overline{ a_{j,k_2, k_1}}\in\C$.  Introducing the
Hermitian forms $\left< A_j(z), \bar z \right>:=\sum_{ k_1,k_2 =1}^m\,
a_{j, k_1, k_2} \, z_{k_1}\, \bar z_{k_2}$, where the $A_j(z)= (\sum_{
k_1=1}^m\, a_{j,k_1,1} \, z_{k_1}, \dots, \sum_{k_1=1}^m\,
a_{j,k_1,m}\, z_{k_1})$ are complex linear endomorphisms of $\C^m$, we
know by definition that $M$ is Levi nondegenerate at the origin if and
only if $\bigcap_{j=1}^d {\rm Ker} \, A_j=\{0\}$. On the other hand,
we may write
\def\theequation{3.2.2}\begin{equation}
\mathcal{Q}_1(t)=
\left(w_j,iA_j(z)\right)_{1\leq j\leq d}+{\rm O}(2),
\end{equation}
hence the rank of $\mathcal{Q}_1$ at the origin equals $n$ if and only
if $\bigcap_{j=1}^d {\rm Ker} \, A_j=\{0\}$ again.  In conclusion, the
two definitions of Levi nondegeneracy are equivalent.

Equivalently, we remind that Levi nondegeneracy of $M$ can be expressed
directly by means of a collection of $d$ defining equations
$r_1(t,\bar t)= \cdots=r_d(t,\bar t)=0$ for $M$ in a neighborhood of
$p_0$, simply by the condition that the intersection of the kernels of
the Levi forms of the defining functions $r_1,\dots,r_d$ at the origin
reduces to $\{0\}$. In Lemma~3.2.6 below, this criterion is
generalized.

\subsection*{3.2.3.~Finite nondegeneracy} 
More generally, we shall say that $M$ is {\sl finitely nondegenerate
at} $p_0$ if there exists an integer $k$ such that the $k$-th Segre
mapping $\mathcal{Q}_k$ is of rank $n$ at the origin. Of course, this
implies that $\mathcal{Q}_l$ is also of rank $n$ for all $l\geq k$. If
$\ell_0$ denotes the smallest integer $k$ such that $\mathcal{Q}_k$ is of
rank $n$ at the origin, we shall say that $M$ is {\sl
$\ell_0$-finitely nondegenerate at} $p_0$. Evidently, $M$ is
$1$-finitely nondegenerate at $p_0$ if and only if it is Levi
nondegenerate. For the moment, we shall admit that $\ell_0$-nondegeneracy
is independent of the choice of coordinates and of defining equations
(this will be proved in Section~3.3 below and it follows in addition
from the proof of Lemma~3.2.6 below). By the definition of
$\mathcal{Q}_k$, the following characterization is immediate.

\def\thelemma{3.2.4}\begin{lemma} 
The generic submanifold $M$ is $\ell_0$-finitely nondegenerate at the origin if
and only if there exists multiindices
$\beta_*^1,\dots,\beta_*^n\in\N^m$ with $\vert \beta_*^i\vert \leq
\ell_0$, $i=1,\dots,n$, and integers $j_*^1,\dots,j_*^n$ with $1\leq
j_*^i\leq d$, $i=1,\dots,n$, such that the rank at the origin of the
mapping
\def\theequation{3.2.5}\begin{equation}
t\mapsto 
(\Theta_{j_*^1,\beta_*^1}(t),\dots,
\Theta_{j_*^n,\beta_*^n}(t))
\end{equation}
is equal to $n$, but such a property is impossible with
multiindices satisfying $\vert \beta_*^i \vert \leq \ell_0-1$.
\end{lemma}

Finite nondegeneracy of $M$ at the origin can be expressed directly by
means of a collection of $d$ defining equations $r_1(t,\bar
t)=\cdots= r_d(t,\bar t)=0$ for $M$ in a neighborhood of
$p_0$. Indeed, let $\overline{L}_1,\dots,\overline{L}_m$ be a basis of
$(0,1)$ vector fields tangent to $M$ in a neighborhood of $p_0$. If
$\beta=(\beta_1,\dots,\beta_m)\in\N^m$, we use the abbreviated
notation $\overline{L}^\beta$ for the derivation of order $\vert \beta
\vert$ defined by $\overline{L}_1^{\beta_1}\cdots
\overline{L}_m^{\beta_m}$. Let $\nabla_t(r_j)(t,\bar t)$ denote
the complex gradient of $r_j(t,\bar t)$, namely 
$([\partial r_j/\partial t_1](t,\bar t),\dots,
[\partial r_j/\partial t_n](t,\bar t))$.

\def\thelemma{3.2.6}\begin{lemma}
The generic submanifold $M$ is $\ell_0$-finitely nondegenerate at
$p_0$ if and only if the complex linear span of all derivatives
$\overline{L}^\beta[ \nabla_t(r_j)(t,\bar t)]\vert_{t=p_0}$, for
$\vert \beta\vert \leq \ell_0$ and for $1\leq j\leq d$, generates
$\C^n$, a property which we may write symbolically as follows
\def\theequation{3.2.7}\begin{equation}
{\rm Span}_\C \{
\overline{L}^\beta [\nabla_t(r_j)(t,\bar t)]\vert_{t=0}:\,
\beta\in\N^m, \, 
\vert \beta\vert \leq \ell_0, \, 
j=1,\dots,d\}=\C^n,
\end{equation}
but such a property cannot be satisfied 
with multiindices satisfying 
$\vert \beta \vert \leq \ell_0-1$.
\end{lemma}

\proof 
Firstly, we can check directly that this new condition does not
depend on the choice of a basis of $(0,1)$ vector fields tangent to
$M$. Indeed, if $\overline{L}_1',\dots, \overline{L}_m'$ is another
basis, there exists an invertible $m\times m$ matrix of power series
$b(t,\bar t)$ such that $\overline{L}_k'=\sum_{l=1}^m\, b_{k,l}(t,\bar
t)\, \overline{L}_l$. Then computing $(\overline{L}')^\beta\,
[\nabla_t(r_j)(t,\bar t)]\vert_{t=0}$, we obtain a linear combination
with constant coefficients of the vectors $(\overline{L})^{\beta_1}\,
[\nabla_t(r_j)(t,\bar t)]\vert_{t=0}$, where $\vert \beta_1\vert \leq
\vert \beta \vert$.  It follows that we have the inclusion relation
\def\theequation{3.2.8}\begin{equation}
\left\{
\aligned
{\rm Span}_\C \{
\overline{L'}^\beta [\nabla_t(r_j)(t,\bar t)]\vert_{t=0}:\, 
\vert \beta\vert \leq \ell_0, \, 
j=1,\dots,d\}\subset \\
\subset
{\rm Span}_\C \{
\overline{L}^\beta [\nabla_t(r_j)(t,\bar t)]\vert_{t=0}:\, 
\vert \beta\vert \leq \ell_0, \, 
j=1,\dots,d\}.
\endaligned\right.
\end{equation}
As $b$ is invertible, we can reverse the r\^oles of the vector fields
$\overline{L}_k$ and of the vector fields $\overline{L}_k'$, which yields the
opposite inclusion relation, hence an equality in~\thetag{3.2.8}.

Secondly, we verify that the new condition~\thetag{3.2.7} neither depends
on the choice of defining equations for $M$, nor on the choice of
coordinates. This is a little bit tedious, but the principle of proof
is also quite simple. Indeed, let $t'=h(t)$ be a change of coordinates
with $h(0)=0$ and assume that the image $M':=h(M)$ is represented by
equations $r_1'(t',\bar t')=\dots=r_d'(t',\bar
t')=0$. Equivalently, there exists an invertible $d\times d$ matrix
$a(t,\bar t)$ of complex algebraic or analytic power series such that
$r'(h(t),\bar h(\bar t))\equiv a(t,\bar t)\, r(t,\bar t)$.  As
we have already checked that the condition~\thetag{3.2.7} does not depend
on the choice of a basis of $(0,1)$ vector fields tangent to $M'$, we
may choose the basis $\overline{L}_k':= (\bar h)_*(\overline{L}_k)$,
namely $\overline{L}_k':= \sum_{i=1}^n\, \overline{L}_k(\bar
h_i)\,\partial_{\bar t_i'}$.  If we denote the complex gradient
$\nabla_t(r_j)$ by $([\partial r_j/\partial t_i](t,\bar
t))_{1\leq i\leq n}$, we may compute the gradient of both sides of the
vector identity $r'(h(t),\bar h(\bar t)) \equiv a(t,\bar t)\,
r(t,\bar t)$, which yields for $j=1,\dots,d$:
\def\theequation{3.2.9}\begin{equation}
\left\{
\aligned
\left(
\sum_{i'=1}^n\, 
{\partial h_{i'}\over \partial t_i}(t)\, 
{\partial r_j'\over \partial t_{i'}'}
(h(t),\bar h(\bar t))
\right)_{1\leq i\leq n}\equiv\\
\equiv \left(
\sum_{l=1}^d\, 
{\partial a_{j,l}\over \partial t_i}(t,\bar t)\, 
r_l(t,\bar t)+
\sum_{l=1}^d\, 
a_{j,l}(t,\bar t)\, 
{\partial r_l\over \partial t_i}(t,\bar t)
\right)_{1\leq i\leq n}.
\endaligned\right.
\end{equation}
Next, we apply the derivations $(\overline{L})^\beta$ to this
identity. Taking into account that the holomorphic terms are not
differentiated, using the definition of $\overline{L}_k'$ and using
the rule of Leibniz for the differentiation of a product, 
we get the following expression
\def\theequation{3.2.10}\begin{equation}
\left\{
\aligned
{} &
\left(
\sum_{i'=1}^n\, 
{\partial h_{i'}\over \partial t_i}(t)\, 
(\overline{L'}^\beta)\left({\partial r_j'\over \partial t_{i'}'}\right)
(h(t),\bar h(\bar t))
\right)_{1\leq i\leq n}\equiv\\
& \ \ \
\equiv \left(
\sum_{\gamma\leq \beta}\,
{\beta! \over \gamma! \
(\beta-\gamma)!}\,
\left(
\sum_{l=1}^d\, 
\left(
\overline{L}^\gamma\left(
{\partial a_{j,l}\over \partial t_i}
\right)(t,\bar t) \ \,
\overline{L}^{\beta-\gamma}(r_l)
(t,\bar t)+\right.\right.\right.\\
&
\left.\left.\left.
\ \ \ \ \ \ \ \ \ \ \ \ \ \ \ \ \ \ \ \ \ \ \ \
+
\overline{L}^\gamma(a_{j,l})(t,\bar t) \ \,
\overline{L}^{\beta-\gamma}\left(
{\partial r_l\over \partial t_i}
\right)(t,\bar t)
\right)
\right)
\right)_{1\leq i\leq n}.
\endaligned\right.
\end{equation}
Here, for a multiindex $\gamma\in\N^m$, we write $\gamma\leq \beta$ if
$\gamma_k\leq \beta_k$ for $k=1,\dots,m$. In this expression, we set
$t=0$.  Since the vector fields $\overline{L}_k$ are tangent to $M$,
all the expressions $\overline{L}^{\beta-\gamma}(r_l)(0,0)$ in the
right hand side vanish.  Using the invertibility of the Jacobian
matrix $(\partial h_{i'}/\partial t_i(0))_{1\leq i,i'\leq n}$, we then
see that the vectors $\left((\overline{L'})^\beta\left( \partial
r_j'/\partial t_{i'}'\right) (0,0)\right)_{1\leq i'\leq n}$ are linear
combinations with constant coefficients of the vectors
$\left((\overline{L}^{ \beta_1})\left( \partial
r_{j_1}/\partial t_i\right) (0,0)\right)_{1\leq i\leq n}$, where
$j_1=1,\dots,d$ and $\vert \beta_1\vert \leq \vert \beta \vert$.
This entails that we have the inclusion relation
\def\theequation{3.2.11}\begin{equation}
\left\{
\aligned
{\rm Span}_\C \{
(\overline{L'})^\beta [\nabla_{t'}(r_j')(t',\bar t')]\vert_{t'=0}:\, 
\vert \beta\vert \leq \ell_0, \, 
j=1,\dots,d\}\subset \\
\subset
{\rm Span}_\C \{
\overline{L}^\beta [\nabla_t(r_j)(t,\bar t)]\vert_{t=0}:\, 
\vert \beta\vert \leq \ell_0, \, 
j=1,\dots,d\}.
\endaligned\right.
\end{equation}
Let $t=h'(t')$ denote the inverse of $t'=h(t)$. Reasoning as above, we get the
opposite inclusion relation, hence an equality in~\thetag{3.2.11}, 
as desired. We notice that the preceding reasoning also shows that
for a fixed coordinate system, the condition~\thetag{3.2.7} does not
depend on the choice of $d$ defining equations for $M$.

Finally, we check that this definition of finite
nondegeneracy coincides with the first one given
in the beginning of \S3.2.3. As we have shown that the second
definition of $\ell_0$-finite nondegeneracy does not depend on the choice of
defining equations, we can assume that $r_j(t,\bar t) :=
\Theta_j(\bar z,t)-\bar w_j$. Then $\nabla_t(r_j)(t,\bar t)=
([\partial\Theta_j/\partial t_1](\bar z,t),\dots,
[\partial\Theta_j/\partial t_n](\bar z,t))$.  We can also assume that
the basis of $(0,1)$ vector fields tangent to $M$ is the usual one, 
as in Chapter~2, which is given by
$\overline{L}_k:=\partial_{\bar z_k}+ \sum_{j=1}^d\, \Theta_{j,\bar
z_k}(\bar z,t)\, \partial_{\bar w_j}$, for $k=1,\dots,m$. Then using
the development $\Theta_j(\bar z,t)=\sum_{\gamma\in\N^m}\, (\bar
z)^\gamma\, \Theta_{j,\gamma}(t)$, we may compute
\def\theequation{3.2.12}\begin{equation}
\overline{L}^\beta [\nabla_t(r_j)(t,\bar t)]\vert_{t=0}=
\beta! \, 
([\partial \Theta_{j,\beta}/\partial t_1](0),\dots,
[\partial \Theta_{j,\beta}/\partial t_n](0)).
\end{equation}
As the right hand side coincides up to a nonzero 
factor with the $(j,\beta)$-th column of the Jacobian 
matrix of the Segre mapping $\mathcal{Q}_{\ell_0}:
t\mapsto (\Theta_{j,\beta}(t))_{1\leq j\leq
d, \, \vert \beta\vert \leq \ell_0}$, we see immediately 
that $\mathcal{Q}_{\ell_0}$ is of rank $n$ at the origin 
if and only if~\thetag{3.2.7} holds.
The proof of Lemma~3.2.6 is complete.
\endproof

\def\theexample{3.2.13}\begin{example}
{\rm 
We provide some elementary examples of finitely nondegenerate
hypersurfaces at the origin:
\begin{itemize}
\item[{\bf (1)}]
$\bar w=w+i[z^5\bar z+\bar z^5z]$ in $\C^2$ is
$5$-finitely nondegenerate.
\item[{\bf (2)}]
$\bar w=w+i[z_1\bar z_1+z_1^2\bar z_2+\bar z_1^2z_2]$ in $\C^3$
is $2$-finitely nondegenerate.
\item[{\bf (3)}]
$\bar w=w+i[z_1\bar z_1+z_1^2\bar z_2+\bar z_1^2z_2+
z_1^3\bar z_3+\bar z_1^3z_3]$ in $\C^4$ is
$3$-finitely nondegenerate.
\end{itemize}
More generally,
let $\varphi_1(z),\dots,\varphi_m(z)$ be a collection of
holomorphic functions vanishing at the origin in $\C^m$ such that the mapping
$z\mapsto (\varphi_1(z),\dots,\varphi_m(z))$ is of rank $m$ at the
origin. Let $\psi_1(z),\dots,\psi_m(z)$ be an arbitrary collection
of nonconstant holomorphic functions with 
different order of vanishing at $0$ and let $\ell_0$ be the highest
order of vanishing of the $\psi_k$. Then the hypersurface
\def\theequation{3.2.14}\begin{equation}
\bar w= w+i\left[
\varphi_1(z)\overline{\psi_1(z)}+
\overline{\varphi_1(z)}\psi_1(z)+\cdots+
\varphi_m(z)\overline{\psi_m(z)}+
\overline{\varphi_m(z)}\psi_m(z)\right]
\end{equation} 
is $\ell_0$-finitely nondegenerate at the origin.
On the contrary, $\bar w = w+iz^2\bar z^2$ is not 
finitely nondegenerate at the origin (it is in
fact essentially finite, {\it see}~\S3.2.25 below).
}
\end{example}

Finite nondegeneracy of $M$ is not a gratuitous generalization of the
notion of Levi nondegeneracy, which would be simply something like a
folklore ``higher order Levi form''. On the contrary, it will appear
to be a very natural nondegeneracy condition. In particular, we shall
establish that an arbitrary real algebraic or analytic generic
submanifold $M$ is, locally in a neighborhood of a Zariski-generic
point $p\in M$, biholomorphic to a product $\underline{M}_p'\times
\Delta^{\kappa}$ of finitely nondegenerate generic submanifold
$\underline{M}_p'\subset \C^{n-\kappa}$ with a certain polydisc
$\Delta^\kappa$.  This property says that up to neglecting the ``flat
part'' $\Delta^\kappa$, every real algebraic or analytic generic
submanifold is finitely nondegenerate at a Zariski-generic point.
Furthermore, some classical CR manifolds (as the tube over
the two-dimensional light cone for
instance), are Levi degenerate at every point, but are finitely
nondegenerate.

\def\theexample{3.2.15}\begin{example}
{\rm 
The tube over the two-dimensional light cone $\Gamma_\C$ in $\C^3$ is
the singular hypersurface defined by $u^2=x_1^2+x_2^2$, where $u={\rm
Re}\, w$ and $x_k={\rm Re}\, z_k$, $k=1,2$. We consider the regular
part of $\Gamma_\C$, which coincides with $\Gamma_\C\cap
\{(x_1,x_2)\neq (0,0)\}$.  In a neighborhood of the smooth point
$(1,0,1)$, we check that $\Gamma_\C$ is Levi degenerate. The $(0,1)$
vector fields tangent to $\Gamma_\C$ are generated by
\def\theequation{3.2.16}\begin{equation}
\left\{
\aligned
\overline{L}_1:=
{\partial\over \partial \bar z_1}+
{x_1\over u}\, 
{\partial \over \partial \bar w},\\
\overline{L}_2:=
{\partial \over \partial \bar z_2}+
{x_2\over u}\, 
{\partial \over \partial \bar w}.
\endaligned\right.
\end{equation}
We observe that the dilatation vector field
\def\theequation{3.2.17}\begin{equation}
\overline{T}:=x_1\, 
{\partial \over \partial \bar z_1}+
x_2\,
{\partial \over \partial \bar z_2}+
u\, 
{\partial \over \partial \bar w},
\end{equation}
which coincides with $x_1\,\overline{L}_1+x_2\, \overline{L}_2$ on
$M$, lies in the kernel of the Levi form, since
$[\overline{L}_1,T]={1\over 2}\, L_1$ and $[\overline{L}_2,T]={1\over
2}\, L_2$. We also observe that, according to~[7], [11], $\Gamma_\C$
is necessarily foliated by complex curves. In fact, the regular locus
of $\Gamma_\C$ is globally foliated by complex lines as follows:
$z_1:=(r+is)\cos\theta+i\lambda$, $z_2:=(r+is)\sin\theta+i\mu$,
$z_3:=r+is$, where $r,s,\theta,\lambda,\mu$ are real
parameters. Finally, applying Lemma~3.2.6, we may check rapidly that
$\Gamma_\C$ is $2$-finitely nondegenerate at every point.
}
\end{example}

\def\theexample{3.2.18}\begin{example}
{\rm 
A generalization of this example is the regular part $M$ of M.
Freeman's cubic $x_1^3+x_2^3-u^3=0$, {\it cf.}~[11].
On $\{u\neq 0\}$, the $(0,1)$
vector fields tangent to $M$ are generated by
\def\theequation{3.2.19}\begin{equation}
\left\{
\aligned
\overline{L}_1:=
{\partial\over \partial \bar z_1}+
{x_1^2\over u^2}\, 
{\partial \over \partial \bar w},\\
\overline{L}_2:=
{\partial \over \partial \bar z_2}+
{x_2^2\over u^2}\, 
{\partial \over \partial \bar w}.
\endaligned\right.
\end{equation}
Again, the dilatation vector field lies in the kernel of the Levi
form, $M$ is foliated by complex lines, but $M$ is (2- or 3-) finitely
nondegenerate at every point.
}
\end{example}

\def\theexample{3.2.20}\begin{example}
{\rm 
Another example of everywhere Levi degenerate real
algebraic hypersurface in $\C^3$ is the hypersurface $M_0$ defined in
the domain
$\{(z_1,z_2,w)\in \C^3: \, \vert z_2 \vert < 1\}$ by the equation
\def\theequation{3.2.21}\begin{equation}
\bar w=w+i\left[
{2z_1\bar z_1+z_1^2\bar z_2
+\bar z_1^2z_2)\over (1-z_2\bar z_2)} 
\right].
\end{equation}
The $(0,1)$ vector fields tangent to $M_0$ are generated by
\def\theequation{3.2.22}\begin{equation}
\left\{
\aligned
\overline{L}_1:=
{\partial\over\partial \bar z_1}+i\left[
{2z_1+2\bar z_1 z_2\over 1
-z_2\bar z_2}\right]\,
{\partial\over\partial \bar w},\\
\overline{L}_2:=
{\partial\over\partial \bar z_2}+i\left[
{(z_1+\bar z_1 z_2)^2\over
(1-z_2\bar z_2)^2}\right]\,
{\partial\over\partial \bar w}.\\
\endaligned\right.
\end{equation}
The kernel of the Levi form is generated by the vector field
\def\theequation{3.2.23}\begin{equation}
\overline{T}:=
-\left[
{z_1+\bar z_1 z_2\over
1-z_2\bar z_2}\right]\,
{\partial\over\partial \bar z_1}+
{\partial\over\partial \bar z_2}-
i\left[{(z_1+\bar z_1 z_2)^2\over
(1-z_2\bar z_2)^2}\right]\,
{\partial\over\partial \bar w}.
\end{equation}
Indeed, we compute:
\def\theequation{3.2.24}\begin{equation}
[\overline{L}_1,T]=
-\left({1\over 1-z_2
\bar z_2}\right)\,
L_1, \ \ \ \ \ \ \ \ \
[\overline{L}_2,T]=-
\left({z_1+\bar z_1 z_2\over
(1-z_2\bar z_2)^2}\right)\,L_1.
\end{equation}
Finally, according to~[7], [11], $M_0$ is
necessarily foliated by complex curves. In fact, $M_0$ is foliated by
the complex lines $z_1:=z_0-\bar z_0\zeta$, $z_2:=\zeta$,
$w:=-iz_0\bar z_0+i\zeta\bar z_0^2+\lambda$, where
$z_0\in\C$, $\lambda\in\R$ and where the complex
variable $\zeta\in\C$ satisfies $\vert
\zeta\vert< 1$. Finally, by a direct application of Lemma~3.2.4, we
see that $M_0$ is $2$-finitely nondegenerate at every point.
}
\end{example}

\subsection*{3.2.25.~Essential finiteness}
More generally, we shall say that $M$ is {\sl essentially finite at
$p_0$} if there exists an integer $k$ such that the Segre mapping
$\mathcal{Q}_k$ is locally finite in a neighborhood of the origin. Of
course, this implies that $\mathcal{Q}_l$ is also locally finite for
all $l\geq k$. Moreover, finite nondegeneracy implies trivially
essential finiteness. For the moment, we shall admit that essential
finiteness is independent of the choice of coordinates and of defining
equations (this will be proved in Section~3.3.).  We say that $M$ is
{\sl $\ell_0$-essentially finite at $p_0$} if $\ell_0$ is the smallest
such integer.  By D.~Hilbert's Nullstellensatz, the following characterization is
immediate.

\def\thelemma{3.2.26}\begin{lemma}
The generic submanifold $M$ is $\ell_0$-essentially finite at the
origin if and only if the complex algebraic or analytic set defined by
\def\theequation{3.2.27}\begin{equation}
\Theta_{j,\beta}(t)-\Theta_{j,\beta}(0)=0, \ \ \ \ \ 
j=1,\dots,d, \ \ \vert \beta \vert \leq \ell_0,
\end{equation}
is zero-dimensional at $0$, but the same complex algebraic
or analytic subset defined with $\vert \beta \vert \leq \ell_0-1$
is positive-dimensional at $0$.
\end{lemma}

Classically, it is known that this property is equivalent to
the fact that the ideal generated by the functions in the left hand
side of~\thetag{3.2.27} is of finite codimension in $\mathcal{A}_\C\{t\}$
or in $\C\{t\}$. More precisely, we define the integers $\ell_1$ and
$\varepsilon_1$ as follows: $\ell_1\geq \ell_0$ is the smallest integer such that
the ideal generated by all $\Theta_{j, \beta}(t)-\Theta_{j, \beta}(0)$
coincides with the ideal generated by the $\Theta_{j, \beta}
(t)-\Theta_{ j,\beta}(0)$ with $\vert \beta \vert \leq k$.  This
integer exists, by noetherianity of $\mathcal{A}_\C\{t\}$ or of
$\C\{t\}$. Also, we define the integer $\varepsilon_1$ to be the
codimension of the ideal $\left<
\Theta_{j,\beta}(t)-\Theta_{j,\beta}(0) \right>_{1\leq j\leq d, \,
\vert \beta\vert \leq \ell_1}$ in $\mathcal{A}_\C\{t\}$ or in $\C\{t\}$.
By essential finiteness, $\varepsilon_1<\infty$. We shall also observe in
Section~3.3 below that $\varepsilon_1$ is a biholomorphic invariant of $M$ at
$p_0$. We call $\varepsilon_1$ the {\sl essential type of $M$ at 
$p_0$} and we denote it by ${\rm Ess\, Type} (M,p_0)$.

Essential finiteness of $M$ at the origin can be expressed
geometrically as follows. We introduce the locus of coincidence
of Segre varieties
\def\theequation{3.2.28}\begin{equation}
\A_{p_0}:=\{t\in \Delta_n(\rho_1): \, 
S_{\bar t}=S_{\bar p_0}\}.
\end{equation}

\def\thelemma{3.2.29}\begin{lemma} 
The set $\A_{p_0}$ is a complex
algebraic or analytic subset of a neighborhood of $p_0$ in $\C^n$
which is contained in $M$ and which is described in local coordinates
by $\A_{p_0}=\{t\in\Delta_n(\rho_1): \,
\Theta_{j,\beta}(t)=\Theta_{j,\beta}(0),\, j=1,\dots,d, \,
\beta\in\N^m\}$.
\end{lemma}

\proof
Let $t\in \A_{p_0}$. Since $p_0\in M$, we have $p_0\in S_{\bar p_0}$
by Lemma~2.2.9. So $p_0\in S_{\bar t}$, whence again by Lemma~2.2.9, $t\in
S_{\bar p_0}=S_{\bar t}$, whence $t\in M$. This shows that $\A_{p_0}$
is contained in $M$. Next, in coordinates, we represent $S_{\bar
p_0}=\{(z,w)\in \Delta_n(\rho_1): \, w_j=\sum_{\beta\in\N^m}\,
z^\beta\, \overline{\Theta}_{j,\beta}(0), \, j=1,\dots,d\}$ and $S_{\bar
t}=\{(z,w)\in \Delta_n(\rho_1): \, w_j=\sum_{\beta\in\N^m}\,
z^\beta\, \overline{\Theta}_{j,\beta}(\bar t), \, j=1,\dots,d\}$. The two
$m$-dimensional complex manifolds $S_{\bar p_0}$ and $S_{\bar t}$
coincide if and only if all the coefficients of their graphing
functions coincide, namely $\overline{\Theta}_{j,\beta}(\bar t)
=\overline{\Theta}_{j,\beta}(0)$
for all $j$ and all $\beta$, which completes the proof.
\endproof

\def\theexample{3.2.30}\begin{example}
{\rm 
We provide some elementary examples of essentially finite
hypersurfaces at $0$.
\begin{itemize}
\item[{\bf (1)}]
$\bar w=w+i[z^N\bar z^N]$ in $\C^2$ has essential type equal to $N$.
\item[{\bf (2)}]
$\bar w=w+i[z_1^3\bar z_1^3+z_2^4\bar z_2^4]$ in $\C^3$
has essential type equal to $12$.
\end{itemize}
More generally, let $\varphi_1(z), \dots, \varphi_m(z)$, be an
arbitrary collection of holomorphic functions with
\def\theequation{3.2.31}\begin{equation}
\dim\left(\C\{z\}/\left<
(\varphi_k(z))_{1\leq k\leq m}
\right>\right)=:
\varepsilon_1<\infty. 
\end{equation}
Then the essential type at the origin of $w=\bar
w+i[\varphi_1(z)\overline{\varphi_1(z)}+\cdots+
\varphi_m(z)\overline{\varphi_m(z)}]$ is equal to $\varepsilon_1$.  On
the contrary, the hypersurface $\bar w=w+i[z_1\bar z_1(1+z_2\bar
z_2)]$ is not essentially finite at the origin (it is in fact Segre
nondegenerate, {\it see}~\S3.2.32 just below).
}
\end{example}

\subsection*{3.2.32~Segre nondegeneracy}
More generally, we shall say that $M$ is {\sl Segre nondegenerate at}
$p_0$ if there exists an integer $k$ such that the restriction to
$S_{\bar p_0}$ of the $k$-th Segre mapping is of maximal possible generic rank
equal to $m$:
\def\theequation{3.2.33}\begin{equation}
{\rm genrk}_\C\left(S_{\bar p_0}\ni t\mapsto 
\mathcal{Q}_k(t)\right)=m,
\end{equation}
which means more precisely that
\def\theequation{3.2.34}\begin{equation}
{\rm genrk}_\C\left(
z\mapsto \left(\Theta_{j,\beta}(z,\overline{\Theta}(z,0))_{1\leq j\leq d, \, 
\vert \beta \vert \leq k}
\right)\right)=m.
\end{equation}
We say that $M$ is {\sl $\ell_0$-Segre nondegenerate at $p_0$}
if $\ell_0$ is the smallest such integer.
Then the following characterization is immediate

\def\thelemma{3.2.35}\begin{lemma}
The generic submanifold $M$ is $\ell_0$-Segre nondegenerate at the
origin if and only if there exist multiindices
$\beta_*^1,\dots,\beta_*^m\in \N^m$ with $\vert \beta_*^k \vert \leq
\ell_0$, $k=1,\dots,m$, and integers $j_*^1,\dots,j_*^m$ with $1\leq j_*^k\leq d$,
$k=1,\dots,m$, such that the determinant
\def\theequation{3.2.36}\begin{equation}
{\rm det}\, 
\left( \left(
[L_{k_1}\Theta_{j_*^{k_2},\beta_*^{k_2}}](
z,\overline{\Theta}(z,0))\right)_{1\leq k_1,k_2\leq m}
\right)
\end{equation}
does not vanish identically as a power series in $z$, but such 
a property is impossible for multiindices 
satisfying $\vert \beta_*^k \vert \leq \ell_0-1$.
\end{lemma}

\def\theexample{3.2.37}\begin{example}
{\rm 
Let $\varphi_1(z),\dots,\varphi_m(z)$ be a collection of holomorphic
functions defined in a neighborhood of the origin $\C^m$ such that the
generic rank of $z\mapsto ( \varphi_1(z), \dots, \varphi_m (z))$ is
equal to $m$.  Then the hypersurface $\bar w= w+i [\varphi_1 (z)
\overline{\varphi_1 (z)} +\cdots+ \varphi_m (z) \overline{ \varphi_m
(z)}]$ is Segre nondegenerate at the origin.  On the contrary, the
real algebraic hypersurface $M$ in $\C^3$ of equation
\def\theequation{3.2.38}\begin{equation}
{\rm Im}\, w={\vert z_1 \vert^2 \vert 1+ z_1 \bar{z}_2\vert^2 
\over 1+\hbox{Re}\, (z_1\bar{z}_2)} -{\rm Re}\, w \, 
{\hbox{Im}\,(z_1\bar{z}_2) \over
1 +\hbox{Re} (z_1\bar{z}_2)}
\end{equation}
is not Segre nondegenerate at the origin (it is in fact
holomorphically nondegenerate, {\it see}~\S3.2.40 just below).  Indeed, solving
with respect to $\bar w$, we may compute its complex defining
equation, which yields
\def\theequation{3.2.39}\begin{equation}
\bar w=-2i\, z_1\bar z_1(1+z_1\bar z_2)+
w\, (1+z_1\bar z_2)/(1+\bar z_1 z_2).
\end{equation}
Then $S_{\bar 0}=\{(z,0)\}$ and the mappings $\mathcal{Q}_k\vert_{
S_{\bar 0}}$, for $k\geq 2$ are (up to zero terms) equal to
$(z_1,z_2)\mapsto (-2iz_1,-2iz_1^2)$, hence of generic rank $1<2$.
}
\end{example}

\subsection*{3.2.40.~Holomorphic nondegeneracy}
More generally, we shall say that $M$ is {\sl holomorphically
nondegenerate at} $p_0$ if there exists an integer $k$ such that the
$k$-th Segre mapping is of maximal generic rank equal to $n$. This
means that
\def\theequation{3.2.41}\begin{equation}
{\rm genrk}_\C\left(t\mapsto (\Theta_{j,\beta}(t))_{1\leq j\leq d,\, 
\vert \beta \vert \leq k}\right)=n.
\end{equation}
We say that $M$ is $\ell_0$-holomorphically nondegenerate at $p_0$ if
$\ell_0$ is the smallest possible such integer.
Then the following characterization is immediate

\def\thelemma{3.2.42}\begin{lemma}
The generic submanifold $M$ is holomorphically nondegenerate at the
origin if and only if there exist multiindices $\beta_*^1, \dots,
\beta_*^n \in \N^n$ with $\vert \beta_*^i \vert \leq
\ell_0$, $i=1,\dots,n$, and integers $j_*^1, \dots,j_*^n$ with $1\leq
j_*^i \leq d$, $i=1,\dots,n$, such that the determinant
\def\theequation{3.2.43}\begin{equation}
{\rm det}\, 
\left( \left(
{\partial \Theta_{j_*^{i_1},\beta_*^{i_1}}\over
\partial t_{i_2}}(t)
\right)_{1\leq i_1,i_2\leq m}
\right)
\end{equation}
does not vanish identically as a power series in $t$, but such 
a property is impossible for multiindices satisfying $\vert \beta_*^i 
\vert \leq \ell_0-1$.
\end{lemma}

\subsection*{3.2.44.~Links between the five nondegeneracy conditions}
Finally, we shall establish the following hierarchy between
the five nondegeneracy conditions presented in this chapter.

\def\thetheorem{3.2.45}\begin{theorem}
Let $M$ be a real algebraic or analytic generic
submanifold in $\C^n$ and let $p_0\in M$. Then the following
four implications hold:
\def\theequation{3.2.46}\begin{equation}
\left\{
\aligned
{}
& \
M \ \text{\it is $\ell_0$-holomorphically nondegenerate at $p_0$} 
\Longleftarrow\\
\Longleftarrow
& \
M  \ \text{\it is $\ell_0$-Segre nondegenerate at $p_0$} 
\Longleftarrow\\
\Longleftarrow
& \
M  \ \text{\it is $\ell_0$-essentially finite at $p_0$}
\Longleftarrow\\
\Longleftarrow
& \
M  \ \text{\it is $\ell_0$-finitely nondegenerate at $p_0$}
\Longleftarrow\\
\Longleftarrow
& \
M  \ \text{\it is Levi nondegenerate at $p_0$}.
\endaligned\right.
\end{equation}
\end{theorem}

\proof
We prove the first implication, considering that the other three
follow from known results in local complex analytic geometry ({\it
see} especially Lemma~4.1.4 below). By specializing the functional
equations~\thetag{2.1.25}, we obtain
$\Theta(0,z,\overline{\Theta}(z,0,0))\equiv 0$ and $w\equiv
\overline{\Theta}(0,0,\Theta(0,0,w))$. Remind the notational
coincidence $\Theta_{j,0}(z,w)\equiv \Theta_j(0,z,w)$. It follows that
the mapping $(0,w)\mapsto (\Theta_{j,0}(0,w))_{1\leq j\leq d}$ is
already of rank $d$ at the origin.  Furthermore, the restriction to
$S_{\bar 0}$ of the zeroth Segre mapping is identically zero:
$\mathcal{Q}_0(z,\overline{\Theta}(z,0))= (\Theta_{j,0}(z,
\overline{\Theta}(z,0)))_{1\leq j\leq d}\equiv 0$.

Suppose now that $M$ is $\ell_0$-Segre nondegenerate, hence $\ell_0$
is the smallest integer such that the generic rank of the mapping
$z\mapsto (\Theta_{j,\beta}(z, \overline{\Theta}(z,0))_{1\leq j\leq
d,\ 1\leq \vert \beta\vert \leq \ell_0})$ is equal to $m$ (notice that
we have written $1\leq \vert \beta \vert \leq \ell_0$ and not $\vert
\beta \vert \leq \ell_0$). Let $k\in\N$ with $k\geq 1$.  It follows
immediately that the generic rank of the mapping
\def\theequation{3.2.47}\begin{equation}
(z,w)\mapsto \left((\Theta_{j,0}(0,w))_{1\leq j\leq d},  
(\Theta_{j,\beta}(z,
\overline{\Theta}(z,0)))_{
1\leq j\leq d,\ 1\leq \vert
\beta\vert \leq k}\right)
\end{equation}
is equal to $n$ if only if $k\geq \ell_0$, hence $M$ is 
$\ell_0$-holomorphically nondegenerate. This completes the proof of 
Theorem~3.2.45.
\endproof

\subsection*{3.2.48.~Expression of the five nondegeneracy conditions in 
normal coordinates} Assume now that the coordinates $(z,w)$ are
normal, namely we have $\Theta_j(0,z,w)\equiv w_j$ and $\Theta_j(\bar
z,0,w)\equiv w_j$, {\it cf.} Theorem~2.1.32. 
It follows that in the development $\bar
w_j=\sum_{\beta\in\N^m}\, (\bar z)^\beta\, \Theta_{j,\beta}(t)$ of the
defining equations of $M$, we have $\Theta_{j,0}(t)\equiv w_j$ and
$\Theta_{j,\beta}(0)=0$ for all $j$ and all $\beta$. Then we can
simplify a little bit the expression of the five nondegeneracy
conditions, which is sometimes useful in applications.

\def\thelemma{3.2.49}\begin{lemma}
In normal coordinates,  we have the following
characterizations:
\begin{itemize}
\item[{\bf (1)}]
$M$ is Levi nondegenerate at the origin 
if and only if there exist multiindices $\beta_*^1,\dots,\beta_*^m\in 
\N^m$ with $\vert \beta_*^k\vert =1$, $k=1,\dots,m$, and integers 
$j_*^1,\dots,j_*^m$ with $1\leq j_*^k\leq d$, $k=1,\dots,m$,
such that ${\rm det} \left(
[\partial \Theta_{j_*^{k_1},\beta_*^{k_1}}/
\partial z_{k_2}](0)\right)_{1\leq  k_1,k_2\leq m}\neq 0$.
\item[{\bf (2)}]
$M$ is $\ell_0$-finitely nondegenerate at the origin if and only if there exist
multiindices $\beta_*^1,\dots,\beta_*^m\in \N^m$ with $\vert
\beta_*^k\vert \leq\ell_0$, $k=1,\dots,m$, and integers
$j_*^1,\dots,j_*^m$ with $1\leq j_*^k\leq d$, $k=1,\dots,m$, such that
${\rm det} \left( [\partial \Theta_{j_*^{k_1},\beta_*^{k_1}}/ \partial
z_{k_2}](0)\right)_{1\leq k_1,k_2\leq m}\neq 0$, but such a property is impossible
for multiindices $\beta_*^k$ satisfying $\vert \beta_*^k \vert \leq \ell_0-1$.
\item[{\bf (3)}]
$M$ is essentially finite at the origin if and only if there exists an
integer $\ell_0$ such that the ideal generated by the
$\Theta_{j,\beta}(z,0)$ for $j=1,\dots,d$ and $\vert \beta \vert \leq
\ell_0$ is of finite codimension in
$\mathcal{A}_\C\{z\}$ or in $\C\{z\}$,
but the same ideal for $\vert \beta \vert \leq \ell_0$ is
of infinite codimension.
\item[{\bf (4)}]
$M$ is Segre nondegenerate at the origin if and only if there exist
multiindices $\beta_*^1,\dots,\beta_*^m\in \N^m$ with $\vert
\beta_*^k\vert \leq\ell_0$, $k=1,\dots,m$, and integers
$j_*^1,\dots,j_*^m$ with $1\leq j_*^k\leq d$, $k=1,\dots,m$, such that
${\rm det} \left( [\partial \Theta_{j_*^{k_1},\beta_*^{k_1}}/ \partial
z_{k_2}](z,0)\right)_{1\leq k_1,k_2\leq m}\not\equiv 0$ in
$\mathcal{A}_\C\{z\}$ or in $\C\{z\}$, but such a property is impossible
for multiindices $\beta_*^k$ satisfying $\vert \beta_*^k \vert \leq \ell_0-1$.
\item[{\bf (5)}]
$M$ is holomorphically nondegenerate at the origin if and only if there exist
multiindices $\beta_*^1,\dots,\beta_*^m\in \N^m$ with $\vert
\beta_*^k\vert \leq\ell_0$, $k=1,\dots,m$, and integers
$j_*^1,\dots,j_*^m$ with $1\leq j_*^k\leq d$, $k=1,\dots,m$, such that
${\rm det} \left( [\partial \Theta_{j_*^{k_1},\beta_*^{k_1}}/ \partial
z_{k_2}](z,w)\right)_{1\leq k_1,k_2\leq m}\not\equiv 0$ in
$\mathcal{A}_\C\{z,w\}$ or in $\C\{z,w\}$, but such a property is impossible
for multiindices $\beta_*^k$ satisfying $\vert \beta_*^k \vert \leq \ell_0-1$.
\end{itemize}
\end{lemma}

\proof 
Thanks to $\Theta_{j,0}(t)\equiv w_j$, we observe that the zero-th
order Segre mapping already provides a rank $d$ subset of power
series. Then each one of the five characterizations may be
checked easily.
\endproof

\section*{\S3.3.~Biholomorphic invariance} 

\subsection*{3.3.1.~Finite nondegeneracy}
Let $t'=h(t)$ be a change of coordinates centered at the origin and
let $M':=h(M)$ as in \S3.1.5. Our purpose is to verify that the above
five nondegeneracy conditions are biholomorphically invariant. Then
the following lemma is relatively crucial.

\def\thelemma{3.3.2}\begin{lemma}
For every $k\in \N$, the ranks at the origin of the
two Segre mappings $t\mapsto \mathcal{Q}_k(t)$ and
$t'\mapsto \mathcal{Q}_k'(t')$ are equal.
\end{lemma}

\proof 
Let ${\bf 1}_k=(0,\dots,0,1,0,\dots,0)\in\N^m$ denote the multiindex
with $1$ at the $k$-th place and zero elsewhere. By
differentiating~\thetag{3.1.10} with respect to $t_i$ at $t_i=0$, we
have the following relations for all $j=1,\dots,d$ and all
$\beta\in\N^m$:
\def\theequation{3.3.3}\begin{equation}
\left\{
\aligned
\sum_{i'=1}^n\, 
{\partial \Theta_{j,\beta}'\over
\partial t_{i'}'}(0)\, 
{\partial h_{i'}\over \partial t_i}(0)+
\sum_{k=1}^m\, 
\sum_{l=1}^d\,
(\beta_k+1)\,
{\partial \bar f_k\over \partial w_l}(0)\,
{\partial \Theta_{0,l}\over \partial t_i}(0)\, 
\Theta_{j,\beta+{\bf 1}_k}'(0)=\\
=\sum_{j_1=1}^d\, \sum_{\vert\beta_1\vert\leq
\vert \beta \vert}\, 
{\partial Q_{j,\beta}\over
\partial \Theta_{j_1,\beta_1}}(\{
\Theta_{j_1,\beta_1}(0)\}_{1\leq j_1\leq d, \, 
\vert \beta_1 \vert \leq \vert \beta \vert})\, 
{\partial \Theta_{j_1,\beta_1}\over \partial t_i}(0).
\endaligned\right.
\end{equation}
Since the Jacobian matrix $([\partial h_{i'}/\partial t_i](0))_{1\leq
i,i'\leq n}$ is
invertible, we deduce immediately from this relation that each partial
derivative $[\partial \Theta_{j,\beta}'/\partial t_{i'}'](0)$ is a
linear combination with constant coefficients of the partial
derivatives $[\partial \Theta_{j_1,\beta_1}/\partial t_i](0)$, with
$i=1,\dots,n$, $j_1=1,\dots,d$ and $\vert \beta_1\vert \leq \vert
\beta \vert$. Consequently, for all $k\in\N$, we have the following
inequality
\def\theequation{3.3.4}\begin{equation}
{\rm rk}_0\left(
t'\mapsto 
(\Theta_{j,\beta}'(t'))_{1\leq j\leq d,\, 
\vert \beta \vert \leq k}
\right)\leq
{\rm rk}_0\left(
t\mapsto 
(\Theta_{j,\beta}(t))_{1\leq j\leq d,\, 
\vert \beta \vert \leq k}
\right).
\end{equation}
Applying the same reasoning to the inverse transformation 
$t=h'(t')$, we also obtain the reverse inequality. In conclusion, 
we have the equality of ranks
\def\theequation{3.3.5}\begin{equation}
{\rm rk}_0\left(
t'\mapsto 
(\Theta_{j,\beta}'(t'))_{1\leq j\leq d,\, 
\vert \beta \vert \leq k}
\right)=
{\rm rk}_0\left(
t\mapsto 
(\Theta_{j,\beta}(t))_{1\leq j\leq d,\, 
\vert \beta \vert \leq k}
\right),
\end{equation}
which completes the proof of Lemma~3.3.2.
\endproof

In particular, it follows immediately that $M$ is
$\ell_0$-finitely nondegenerate at $p_0$ if and only if $M'$ is
$\ell_0$-finitely nondegenerate at $p_0'=h(p_0)$.  This proves that the
definition given in \S3.2.3 is invariant under complex algebraic or
analytic changes of coordinates.

\subsection*{3.3.6.~Levy multitype at the origin}
More generally, since the ranks of the Segre mappings $\mathcal{Q}_k$
are invariant, we may introduce the successive invariant integers
$\lambda_{k,0}$ such that the rank at the origin of $\mathcal{Q}_k$
equals $\lambda_{0,0}+\lambda_{1,0}+ \cdots+\lambda_{k,0}$. Since
the mapping $w\mapsto \Theta(0,0,w)$ is invertible, we have
$\lambda_{0,0}=d$. Since $M$ is $\ell_0$-finitely nondegenerate at the
origin, we have
\def\theequation{3.3.7}\begin{equation}
d+\lambda_{1,0}+\cdots+\lambda_{\ell_0,0}=n.
\end{equation} 
If $M$ is $\ell_0$-finitely nondegenerate at $p_0$, 
we call the multiplet $(d,\lambda_{1,0}, \dots,\lambda_{\ell_0,0})$ the
{\sl Levi multitype of $M$ at $p_0$}.

\subsection*{3.3.8.~Essential finiteness}
Now, we check that essential finiteness is a biholomorphically invariant
property. It suffices to prove the following lemma.

\def\thelemma{3.3.9}\begin{lemma}
If $p_0'=h(p_0)$, we have 
$h(\A_{p_0})=\A_{p_0'}'$.
\end{lemma}

\proof
Let $t\in\A_{p_0}$, namely $\Theta_{j,\beta}(t)= \Theta_{j,\beta}(0)$
for all $j=1,\dots,d$ and all $\beta\in \N^m$.  The recipe is again to
look at~\thetag{3.1.10}. As $\Theta_0(t)=\Theta_0(0)=0$ and as $\bar
f(0)=0$, we obtain from~\thetag{3.1.10}
\def\theequation{3.3.10}\begin{equation}
\left\{
\aligned
\Theta_{j,\beta}'(h(t))\equiv & \
R_{j,\beta}(\{
\Theta_{j_1,\beta_1}(t)\}_{
1\leq j_1\leq d, \, 
\vert \beta_1 \vert \leq 
\vert \beta \vert})\\
\equiv & \
R_{j,\beta}(\{
\Theta_{j_1,\beta_1}(0)\}_{
1\leq j_1\leq d, \, 
\vert \beta_1 \vert \leq 
\vert \beta \vert})\equiv 
\Theta_{j,\beta}'(h(0)),
\endaligned\right.
\end{equation}
so $h(t)\in \A_{p_0'}'$, where $p_0'=h(p_0)$ is the 
origin in the coordinates $t'$. In other words, we have
shown that $h(\A_{p_0})\subset \A_{p_0'}'$.
Since $h$ is invertible, we also get similarly 
$h'(\A_{p_0'}')\subset \A_{p_0}$. In conclusion, 
$h(\A_{p_0})=\A_{p_0'}'$, as desired. This completes the
proof.
\endproof
 
\subsection*{3.3.11.~Segre nondegeneracy}
Next, we check that Segre nondegeneracy is a biholomorphically
invariant property. At first, we observe that it follows from the first
functional equation in~\thetag{2.1.25} that $\Theta(0,z, \overline{
\Theta}(z,0)) \equiv 0$, or equivalently $\Theta_0(z, \overline{
\Theta}(z,0)) \equiv 0$. The recipe is again to
look at~\thetag{3.1.10}, replacing $t$ by 
$(z,\overline{\Theta}(z,0))$, which 
yields
\def\theequation{3.3.12}\begin{equation}
\Theta_{j,\beta}'(h(z,\overline{\Theta}(z,0)))\equiv
R_{j,\beta}(\{
\Theta_{j_1,\beta_1}(z,\overline{\Theta}(z,0))\}_{
1\leq j_1\leq d, \, 
\vert \beta_1 \vert \leq 
\vert \beta \vert}).
\end{equation}
Since $h$ is invertible and maps the Segre
variety $S_{\bar p_0}$ onto the Segre variety
$S_{\bar p_0'}'$, we get the inequality
\def\theequation{3.3.13}\begin{equation}
{\rm genrk}_\C\left(z'\mapsto 
\mathcal{Q}_k'(z',\overline{\Theta}'(z',0))\right)\leq
{\rm genrk}_\C\left(z\mapsto 
\mathcal{Q}_k(z,\overline{\Theta}(z,0))\right).
\end{equation}
Reversing the r\^oles of $t$ and of $t'$, we also
get the opposite inequality, hence an equality.
In particular, $M$ is Segre nondegenerate at $p_0$
if and only if $M'$ is Segre nondegenerate at $p_0'$.

\subsection*{3.3.14.~Holomorphic nondegeneracy} 
Using again~\thetag{3.1.10} and
the analogous relation in which the r\^oles of $t$ and of $t'$ are
reversed, we may establish that for $k\in \N$, we have
\def\theequation{3.3.15}\begin{equation}
{\rm genrk}_\C\left(t'\mapsto \mathcal{Q}_k'(t')\right)= {\rm
genrk}_\C\left(z\mapsto \mathcal{Q}_k(t)\right).
\end{equation}
In particular, $M$ is holomorphically nondegenerate at $p_0$ 
if and only if $M'$ is holomorphically nondegenerate at $p_0'$.

\section*{\S3.4.~Manifolds without nondegeneracy conditions}

\subsection*{3.4.1.~Reflection mapping} It is also interesting to study generic
submanifolds without requiring any nondegeneracy condition on them, as we
shall see in our study of the generalized reflection principle 
in Part~II of this memoir. Let $t'=h(t)$ be a local complex algebraic or analytic
equivalence defined in a neighborhood of $p_0\in M$ satisfying the
same conditions as in \S3.1.5. We do not require that $h$ is invertible.  Let
$\bar w_j'=\Theta_j'(\bar z',t')$, $j=1,\dots,d$ be the equations of
$M'$.  Let $\bar\nu':=(\bar\lambda',\bar\mu')\in\C^m\times\C^d$. We
define the {\sl reflection mapping} associated to such defining
equations to be the vectorial power series
\def\theequation{3.4.2}\begin{equation}
\mathcal{R}_h'(t,\bar \nu'):=
\left(\bar \mu_j'-\Theta_j'
(\bar \lambda',h(t))\right)_{1\leq j\leq d},
\end{equation}
which belongs to $\mathcal{A}_\C\{t,\bar\nu'\}^d$ or to
$\C\{t,\bar\nu'\}^d$. By developing the right hand side in powers of
$\bar\lambda'$, we can write more explicitely
\def\theequation{3.4.3}\begin{equation}
\mathcal{R}_h'(t,\bar \nu')=
\left(\bar\mu_j'-\sum_{\beta\in\N^m}\, 
(\bar\lambda')^\beta\, 
\Theta_{j,\beta}'(h(t))
\right)_{1\leq j\leq d}.
\end{equation}
The datum of $\mathcal{R}_h'$ is essentially equivalent to the datum
of the infinite collection of complex algebraic or analytic functions
$\Theta_{j,\beta}'(h(t))$, or equivalently to the composition of the
infinite Segre mapping of $M'$ with $h$, namely the mapping $t\mapsto
\mathcal{Q}_\infty'(h(t))$. 

We observe that $\mathcal{R}_h'$ is biholomorphic invariant in the
following sense. Let $t''=h'(t')$ be a second local mapping, which we
assume to be invertible (do not confuse here $h'$ with the notation
used in \S3.1.5 for the inverse of $h$).  Let $M'':=h'(M')$ and assume
that its equations are $\bar w_j''= \Theta_j''(\bar z'',z'',w'')$,
$j=1,\dots,d$. We consider the composition $t''=h'(h(t))$. Applying
Theorem~3.1.9 and~\thetag{3.1.11}, we know that
$\Theta_{j,\beta}''(h'(t'))\equiv R_{j,\beta}'(
\{\Theta_{j,\beta}'(t')\}_{1\leq j\leq d,\, \beta\in\N^m})$, where the
$R_{j,\beta}'$ are certain algebraic or analytic expressions which
depend only on $h'$. It follows that we can write
\def\theequation{3.4.4}\begin{equation}
\left\{
\aligned
\mathcal{R}_{h'\circ h}''(t,\bar \nu'')= 
& \
\bar \mu''-\sum_{\beta\in\N^m}\, 
(\bar \lambda'')^\beta \, 
\Theta_{j,\beta}''(h'(h(t)))\\
= 
& \
\bar \mu''-\sum_{\beta\in\N^m}\, 
(\bar \lambda'')^\beta \, 
R_{j,\beta}'(
\{\Theta_{j,\beta}'(h(t))\}_{1\leq j\leq d,\, 
\beta\in\N^m}).
\endaligned\right.
\end{equation}
So we can essentially express the reflection mapping
$\mathcal{R}_{h'\circ h}''(t,\bar \nu'')$ by means of the reflection
mapping $\mathcal{R}_h'(t,\bar \nu')$, modulo algebraic or analytic
expressions $R_{j,\beta}'$ which only depend on the change of
coordinates $t''=h'(t')$. In particular, the $R_{j,\beta'}'$ are
algebraic if $M'$, $h'$ and $M''$ are algebraic and the
relation~\thetag{3.4.4} shows that if $\mathcal{R}_h'(t,\bar \nu')$
was also algebraic at the beginning, then its algebraicity is
preserved after the algebraic change of coordinates $t''=h'(t')$. It
is in this sense that we say that $\mathcal{R}_h'$ is a biholomorphic
invariant object. We shall study the reflection mapping thoroughly in
Part~II of this memoir.

\subsection*{3.4.5.~Tangent holomorphic vector fields}
It is desirable to obtain a rank property analogous to~\thetag{3.3.5}
for points $t_p$ close to the origin. Whereas for a fixed integer $k$,
it is in general untrue that~\thetag{3.3.5} holds with the first ${\rm
rk}_0$ replaced by ${\rm rk}_{h(t_p)}$ and the second ${\rm rk}_0$
replaced by ${\rm rk}_{t_p}$ ({\it see} Example~3.5.16
below), the corresponding property for
$k=\infty$ is true.

\def\thecorollary{3.4.6}\begin{corollary}
For all $t_p$ in a neighborhood of the origin,
the rank at $t_p$ of the infinite Segre mapping
$t\mapsto \mathcal{Q}_\infty(t)$ of $M$ coincides with the 
rank at $t_{p'}'=h(t_p)$ of the infinite Segre mapping
$t'\mapsto \mathcal{Q}_\infty'(t')$ of $M'$.
\end{corollary}

\proof
There exists an integer $\ell_p$ such that the rank $n_p$ at $t_p$ of
the infinite Segre mapping $t \mapsto \mathcal{Q}_\infty(t)$ coincides
with the rank at $t_p$ of the $\ell_p$-th Segre mapping $t\mapsto
\mathcal{ Q}_{ \ell_p}(t)$.  Hence for every $l=1,\dots,d$ and every
$\vert \beta \vert \geq \ell_p+1$, the gradient of $\Theta_{
j,\beta}(t)$ at $t_p$ is a linear combination of the columns of the
Jacobian matrix of ${\rm Jac}\,\mathcal{ Q}_{ \ell_p}(t_p)$. Using
this fact, using the invertibilty of $h$ and differentiating the two
sides of~\thetag{3.1.10}, we deduce that the rank $n_{p'}'$ at
$t_{p'}'$ of the infinite Segre mapping $t' \mapsto
\mathcal{Q}_\infty'(t')$ is less than or equal to $n_p$, namely
$n_{p'}' \leq n_p$. Since $h$ is invertible, by considering the
inverse $t=h'(t')$, we can reverse the r\^oles of $t_p$ and of
$t_{p'}'$ and we get the opposite inequality $n_p\leq n_{p'}'$, which
completes the proof.
\endproof

In particular, the generic rank $n_M=\max_{p\in M}\, n_p$ of the
infinite Segre mapping of $M$ is a biholomorphic invariant of $M$. We
call this integer the {\sl essential holomorphic dimension of} $M$.
Concretely, $n_M$ is the smallest integer such that there exists a
$n_M\times n_M$ minor of the infinite Jacobian matrix
\def\theequation{3.4.7}\begin{equation}
{\rm Jac}\, 
\mathcal{Q}_\infty(t)=
\left([\partial 
\Theta_{j,\beta}/
\partial t_i](t)\right)_{
1\leq j\leq d, \, 
\beta\in\N^m}^{1\leq i\leq n}
\end{equation}
that does not vanish identically, but all its $(n_M+1)\times
(n_M+1)$ minors do vanish identically.

\def\thetheorem{3.4.8}\begin{theorem}
With this integer $n_M$,
there exist $(n-n_M)$ holomorphic vector fields 
\def\theequation{3.4.9}\begin{equation}
\mathcal{T}_k=\sum_{i=1}^n\, 
a_i(t)\, {\partial\over
\partial t_i}, \ \ \ \ \ \ \
k=1,\dots, n-n_M, 
\end{equation}
which have complex algebraic or analytic
coefficients $a_i(t)$, which are defined in a neighborhood $V_0$ of
the origin, which are tangent to $M\cap V_0$ and which are linearly
independent at a Zariski-generic point $p\in V_0$. Conversely, the integer
$(n-n_M)$ is the maximal number of holomorphic vector fields with
complex algebraic or analytic coefficients defined in a neighborhood
$V_0$ of the origin which are tangent to $M\cap V_0$ and linearly
independent at a Zariski-generic point.
\end{theorem}

\proof
We choose integers $j_*^1, \dots, j_*^{n_M}$ with $1\leq j_*^i\leq d$, $i=1,
\dots, n_M$, and multiindices $\beta_*^1,\dots,\beta_*^{n_M} \in\N^m$
such that the generic rank of the mapping $t\mapsto (\Theta_{j_*^l,
\beta_*^l}(t))_{ 1\leq l\leq n_M}$ is equal to $n_M$. In other words,
the Jacobian matrix $([\partial\Theta_{j_*^l, \beta_*^l}(t)/\partial
t_i](t))_{1\leq l\leq n_M, \, 1\leq i\leq n}$ possesses a $n_M\times
n_M$ minor which does not vanish identically.  If $n-n_M>0$, applying
classical linear algebra (Cramer's rule for solving systems of linear
equations), we see that there exist $(n-n_M)$ independent
power series vectorial solutions
$(a_{k,1}(t),\dots,a_{k,n}(t))$, $k=1,\dots,n-n_M$, of the system of
$n_M$ equations
\def\theequation{3.4.10}\begin{equation}
\sum_{i=1}^n\, 
a_{k,i}(t)\, 
{\partial \Theta_{j_*^l, 
\beta_*^l}\over \partial t_i}(t)\equiv 0, \ \ \ \ \
l=1,\dots,n_M.
\end{equation}
Equivalently, the $(n-n_M)$ vector fields
\def\theequation{3.4.11}\begin{equation}
\mathcal{T}_k:=\sum_{i=1}^n\, 
a_{k,i}(t)\, {\partial\over \partial t_i},
\end{equation}
$k=1,\dots,n-n_M$, are linearly independent at a Zariski-generic point
of $V_0$ and they satisfy $\mathcal{T}_k\Theta_{j_*^l,
\beta_*^l}(t)\equiv 0$ for $k=1,\dots,n-n_M$ and $l=1,\dots,
n_M$. Since $M$ is generic, the restriction to $M$ of the vector fields
$\mathcal{T}_k$ are also linearly independent at a Zariski-generic
point.

Let now $(j,\beta)\neq (j_*^l,\beta_*^l)$, for
$l=1,\dots,n_M$. By assumption, we also have
\def\theequation{3.4.12}\begin{equation}
{\rm genrk}_\C \left(t\longmapsto
\left(
(\Theta_{j_*^l,\beta_*^l}(t))_{1\leq l\leq n_M},\,
\Theta_{j,\beta}(t)
\right)\right)=n_M.
\end{equation}
In a neighborhood $V_p$ of a point $t_p\in V_0$ at which this rank is
equal to its maximum $n_M$, there exists a complex algebraic or
analytic mapping $R_{j,\beta}$ such that we can write
\def\theequation{3.4.13}\begin{equation}
\Theta_{j,\beta}(t)\equiv 
R_{j,\beta}( 
\Theta_{j_*^1,\beta_*^1}(t),\dots,
\Theta_{j_*^{n_M},\beta_*^{n_M}}(t))
\end{equation}
for all $t\in V_p$. Since $\mathcal{T}_k \Theta_{j_*^l, \beta_*^l}(t)
\equiv 0$, it follows that $\mathcal{T}_k \Theta_{ j,\beta} (t) \equiv
0$ for $t\in V_p$, hence for all $t\in V_0$ thanks to the principle of
analytic continuation. In summary, we have shown that $\mathcal{T}_k
\Theta_{j, \beta}(t)\equiv 0$ for all $k=1,\dots,n-n_M$, all $j=1,
\dots,d$ and all $\beta\in\N^m$.  We conclude immediately that the
$\mathcal{T}_k$ are tangent to $M$, since
\def\theequation{3.4.14}\begin{equation}
\mathcal{T}_k(\bar w_j-
\Theta_j(\bar z,t))\equiv
\sum_{\beta\in\N^m}\, 
(\bar z)^\beta\, 
\mathcal{T}_k\Theta_{j,
\beta}(t)\equiv 0.
\end{equation}

Conversely, suppose that there exist $\chi$ holomorphic vector fields
$\mathcal{T}_k$, $k=1,\dots,\chi$, like~\thetag{3.4.9} with complex
algebraic or analytic coefficients which are linearly independent at a
Zariski-generic point of $V_0$ and such that $\mathcal{T}_k$ is
tangent to $M\cap V_0$. By the tangency condition~\thetag{3.4.14}, we get
$\mathcal{T}_k\Theta_{j,\beta}(t)\equiv 0$ for all $k=1,\dots,\chi$, 
all $j=1,\dots,d$ and
all $\beta \in\N^m$.  By considering these equations at a point at
which the vector fields $\mathcal{T}_k$ are linearly independent, we
deduce the inequality $\chi\leq n-n_M$.  The proof of Theorem~3.4.8 is
complete.
\endproof

\def\thecorollary{3.4.15}\begin{corollary}
The real algebraic or analytic generic submanifold $M$ is
holomorphically nondegenerate at $0$, {\it i.e.}
$n_M=n$, if and only if there does not exist 
any nonzero holomorphic vector field defined in a neighborhood 
$V_0$ of $0$ which is tangent to $M\cap V_0$.
\end{corollary}

This property may be considered as an equivalent definition of
holomorphic nondegeneracy, as was done by N. Stanton in~[28].
However, we believe that the previous definition in terms of the
generic rank of the Segre mapping is more adequate.

\subsection*{3.4.16.~Exceptional locus of $M$} We define the
{\sl extrinsic exceptional locus of $M$} to be the proper complex
analytic set $\mathcal{E}^{\rm exc}$ which is the zero locus of all
$n_M\times n_M$ minors of the Jacobian matrix ${\rm Jac}\,
\mathcal{Q}_\infty(t)$. By Corollary~3.4.6, $\mathcal{E}^{\rm exc}$ is an
invariant complex algebraic or analytic set which is independent of
coordinates. We define the {\sl intrinsic exceptional locus} of $M$ to
the proper real algebraic or analytic subset $M\cap \mathcal{E}^{\rm
exc}$.  By definition, the rank of $t\mapsto\mathcal{Q}_\infty(t)$ at
$t_p\in V_0$ equals $n_M$ if and only if $t_p\in M\backslash
\mathcal{E}^{\rm exc}$.  We shall come back to $\mathcal{E}^{\rm exc}$
in Corollary~3.5.53 below.

\section*{\S3.5.~Jets of Segre varieties and global nondegeneracy conditions}

\subsection*{3.5.1.~Fundamental definitions}
Let $\mathcal{M}$ be the extrinsic complexification of $M$ given by
the equations $\xi_j=\Theta_j(\zeta,t)$, $j=1,\dots,d$. We shall
assume that $\mathcal{M}$ is complex algebraic or analytic and
defined in the polydisc $\Delta_{2n}(\rho_1)$, as in Definition~2.1.44. 
Recall that the
conjugate complexified Segre variety $\underline{\mathcal{S}}_t$ is the
$m$-dimensional complex algebraic or analytic submanifold defined by
the equations $\xi_j=\Theta_j(\zeta,t)$, $j=1,\dots,d$, where $t$ is fixed.  We
consider the mapping of $k$-th order jets of 
$\underline{\mathcal{S}}_t$ at one of its point 
$(\zeta,\Theta(\zeta,t))$ which is defined by 
\def\theequation{3.5.2}\begin{equation}
J_\tau^k\underline{\mathcal{S}}_t:=\left(
\zeta,\left({1\over\beta!}\, 
\partial_\zeta^\beta \Theta_j(\zeta,t)\right)_{1\leq j\leq d,\, 
\vert \beta \vert \leq k}\right).
\end{equation} 
In this section, we shall study the mapping $J_\tau^k\underline{\mathcal{S}}_t$
thoroughly. It is complex algebraic or analytic with values in
$\C^{m+N_{d,m,k}}$. If $k_2\geq k_1$ and if $\pi_{k_2,k_1}$ denotes
the canonical projection $\C^{m+ N_{d,m,k_2}} \to \C^{
m+N_{d,m,k_1}}$, we obviously have $\pi_{ k_2, k_1}( J_\tau^{k_2}
\underline{ \mathcal{S}}_t)= J_\tau^{ k_1} \underline{ \mathcal{
S}}_t$.

Compared to the Segre mapping introduced in \S3.1.1, here we notice that
the term $\zeta$ is present and we may obviously identify $J_\tau^k
\underline{ \mathcal{S}}_t \vert_{\zeta=0}$ with $\mathcal{Q}_k(t)$.

We shall sometimes denote this mapping by $(t,\tau)\mapsto J_\tau^k
\underline{\mathcal{S}}_t$, where we implicitely mean that
$(t,\tau)\in \mathcal{M}$.

Analogously, we may also consider the mapping of $k$-th order jets of
the complexified Segre variety $\mathcal{S}_\tau$ defined by the
equations $w_j=\overline{\Theta}_j(z,\tau)$, $j=1,\dots,d$, where $\tau$ is fixed.
Its explicit expression is similar:
\def\theequation{3.5.3}\begin{equation}
J_t^k\mathcal{S}_\tau:=\left(
z,\left({1\over\beta!}\, 
\partial_z^\beta 
\overline{\Theta}_j(z,\tau)\right)_{1\leq j\leq d,\, 
\vert \beta \vert \leq k}\right).
\end{equation}
The link between these two mappings is very simple:
\def\theequation{3.5.4}\begin{equation}
J_{\bar t}^k\underline{\mathcal{S}}_{\bar \tau} \equiv
\overline{
J_t^k\mathcal{S}_\tau
}
\end{equation}
In other words, the following diagram is commutative.
$$
\diagram \mathcal{M} \rto^{\sigma} \dto_{J_\bullet^k\mathcal{S}_\bullet} 
& \mathcal{M} 
\dto^{J_\bullet^k\underline{\mathcal{S}}_\bullet}\\
\C^{m+N_{d,m,k}} \rto^{(\overline{\bullet})} & \C^{m+N_{d,m,k}}
\enddiagram,
$$
where $(\overline{\bullet})$ denotes the complex conjugation operator.
Since the two jet mappings are therefore essentially equivalent, we
shall only study the nondegeneracy conditions for the mapping
$J_\tau^k\underline{\mathcal{S}}_t$.

\subsection*{3.5.5.~Invariance under changes of coordinates}
As in \S3.1.5, let $t'=h(t)$ be a change of complex algebraic or
analytic coordinates. We shall prove the following theorem in Section~3.6 below.  
Notice that for $\zeta=0$, we recover Theorem~3.1.9.

\def\thetheorem{3.5.6}\begin{theorem}
For every $j=1,\dots,d$ and every $\beta\in\N^m$, there exists a
complex algebraic or analytic mapping in its variables $Q_{j,\beta}$
such that
\def\theequation{3.5.7}\begin{equation}
{1\over \beta!}\, 
{\partial^{\vert \beta\vert} \Theta_j'\over
\partial (\zeta')^\beta}\left(
\bar f(\zeta,\Theta(\zeta,t)),
h(t)\right)\equiv
Q_{j,\beta}\left(\zeta,\ 
\left({1\over \beta_1!}\,\partial_\zeta^{\beta_1}\Theta_{j_1}(\zeta,t)
\right)_{1\leq j_1\leq d, \,
\vert \beta_1\vert \leq \vert \beta\vert}\right).
\end{equation}
\end{theorem}

Here, the points $(\zeta,t)$ belong to a neighborhood of the origin, say to
the polydisc
$\Delta_{2m+d}(\rho_1)$.  Fix $p=(t_p,\tau_p)\in\mathcal{M}$ which we
identify with $(\zeta_p,t_p)\in\C^{2m+d}$ and denote $t_{p'}':=h(t_p)$
and $\zeta_{p'}':=\bar f(\tau_p)$. Locally in a neighborhood of
$(\zeta_p,t_p)$, the mapping 
\def\theequation{3.5.8}\begin{equation}
(\zeta,t)\mapsto (\bar f(\zeta,\Theta(\zeta,t)),h(t))=:(\zeta',t')
\end{equation}
is invertible by assumption (remind that $\zeta\mapsto 
\bar f(\zeta,\Theta(\zeta,0))$ is invertible at $\zeta=0$). 
Applying Theorem~3.5.6, we
may deduce:

\def\thecorollary{3.5.9}\begin{corollary}
For every $k\in\N$, the following equality of ranks holds
\def\theequation{3.5.10}\begin{equation}
\left\{
\aligned
{}
&
{\rm rk}_{(\zeta_p,t_p)}\left(
(\zeta,t)\mapsto  (\zeta,\, (1/\beta!)\, 
(\partial_\zeta^\beta\Theta_j(\zeta,t))_{1\leq j\leq d,\, 
\vert\beta\vert \leq k})
\right)=\\
&
\ \ \ \ \ \ \ \ \ \ \ 
={\rm rk}_{(\zeta_{p'}',t_{p'}')}\left(
(\zeta',t')\mapsto (\zeta', \, (1/\beta!) \,
(\partial_{\zeta'}^\beta
\Theta_j'(\zeta',t'))_{1\leq j\leq d,\, 
\vert \beta\vert \leq k})
\right).
\endaligned\right.
\end{equation}
\end{corollary}

\proof 
Indeed, using the change of coordinates~\thetag{3.5.8}, it suffices to
show that the rank at $(\zeta_p,t_p)$ of the mapping $(\zeta,t)\mapsto
\left(\bar f(\zeta,\Theta(\zeta,t)), \, 
{1\over \beta!}\, \partial_{\zeta'}^\beta\Theta_j'(\bar
f(\zeta,\Theta(\zeta,t)),h(t)) \right)_{1\leq j\leq d,\, \vert \beta
\vert \leq k}$ is less than or equal to the rank at $(\zeta_p,t_p)$ of
the mapping $(\zeta,t)\mapsto \left(\zeta, \,  {1\over \beta!}\,
\partial_\zeta^\beta\Theta_j(\zeta,t)\right)_{1\leq j\leq d, \, \vert
\beta\vert \leq k}$, because after reversing the r\^oles of $t$ and
of $t'$, we also get the opposite inequality. But this inequality
follows directly by differentiating the two sides of~\thetag{3.5.7}
with respect to $(\zeta,t)$ at $(\zeta_p,t_p)$.
\endproof

Let $p\in M$. The ranks at $(p,\bar p)$ of the mappings
$(t,\tau)\mapsto J_\tau^k\underline{\mathcal{S}}_t$ are invariant
under changes of coordinates. Thus, we may introduce several pointwise
invariants of $M$ at $p$ as follows. We denote by $m+n_p\leq m+n$ the
maximal rank of the mapping $(t,\tau)\mapsto
J_\tau^k\underline{\mathcal{S}}_t$ at $(t_p,\bar t_p)$ for
$k=0,1,\dots$ and by $\ell_p$ the smallest integer $k$ such that the
rank at $(t_p, \bar t_p)$ of the mapping $(t,\tau)\mapsto
J_\tau^k\underline{\mathcal{S}}_t$ is equal to $m+n_p$. More
generally, for $k=0,\dots,\ell_M$, we denote by $\lambda_{k,p}$ the
nonnegative integers satisfying
\def\theequation{3.5.11}\begin{equation}
{\rm rk}_{(t_p,\bar t_p)}\left(
(t,\tau)\mapsto J_\tau^k\underline{\mathcal{S}}_t
\right)=m+\lambda_{0,p}+\cdots+\lambda_{k,p}.
\end{equation}
Obviously, the functions $p\mapsto n_p$, $p\mapsto \ell_p$, $p\mapsto
\lambda_{k,p}$ are lower semicontinuous in the Zariski topology.

\subsection*{3.5.12.~Generic ranks}
To begin with, we observe that it follows from Corollary~3.5.9 that the
generic ranks of the mappings $(t,\tau)\mapsto
J_\tau^k\underline{\mathcal{S}}_t$, which increase with $k$, are invariant
under changes of coordinates. We need the following stabilization
result.

\def\thelemma{3.5.13}\begin{lemma}
If ${\rm genrk}_\C\left((t,\tau)\mapsto
J_\tau^{k+1}\underline{\mathcal{S}}_t\right)={\rm genrk}_\C
\left((t,\tau)\mapsto J_\tau^k\underline{\mathcal{S}}_t\right)$, then
for all $l\geq 1$, we also have ${\rm genrk}_\C\left((t,\tau)\mapsto
J_\tau^{k+l}\underline{\mathcal{S}}_t\right)={\rm genrk}_\C
\left((t,\tau)\mapsto J_\tau^k\underline{\mathcal{S}}_t\right)$.
\end{lemma}

\proof
By assumption, in a neighborhood $\mathcal{V}_p$ of a point
$(t_p,\tau_p)\in\mathcal{M}$ at which the ranks of the first two
mappings are equal to their generic rank, hence maximal and locally
constant, it follows from the
constant rank theorem that for every $j=1,\dots,d$ and for every multiindex
$\beta$ with $\vert \beta \vert=k+1$, there exists a complex algebraic
or analytic function $R_{j,\beta}$ such that we can write
\def\theequation{3.5.14}\begin{equation}
{1\over \beta!}\,
\partial_\zeta^\beta
\Theta_j(\zeta,t)=R_{j,\beta}\left(\zeta, \, 
\left({1\over \beta_1!}\,
\partial_\zeta^{\beta_1}
\Theta_{j_1}(\zeta,t)\right)_{1\leq j_1\leq d, \, 
\vert \beta_1\vert \leq \vert k \vert}\right),
\end{equation}
for all $(t,\tau)\in \mathcal{V}_p$. Differentiating these relations
with respect to $\zeta$ and making substitutions, we obtain that for
every $j=1,\dots,d$ and every multiindex $\beta$ with $\vert
\beta\vert =k+l$, $l\geq 1$, there exists a complex algebraic or
analytic function $R_{j,\beta}$ satisfying a relation 
like~\thetag{3.5.14}. This implies that the generic rank of the
mapping $\mathcal{V}_p \ni(t,\tau) \mapsto J_\tau^{k+l} \underline{
\mathcal{S}}_t$ is the same as the generic rank of the mapping $\mathcal{V}_p
\ni(t, \tau) \mapsto J_\tau^k \underline{ \mathcal{S}}_t$. As the
generic rank propagates by the principle of analytic continuation, the
lemma follows.
\endproof

\subsection*{3.5.15.~Reformulation of the five nondegeneracy conditions}
We shall now observe that the five nondegeneracy conditions introduced in
Section~3.2 may also be expressed by means of the morphism of jets of
Segre varieties.  This formulation will be much better than the
formulation in terms of the Segre mapping given in Section~3.2,
because it will be valuable not only for the central point 
$p_0$, but also for an arbitrary point $p$ varying in a 
neighborhood of $p_0$. Before stating the theorem, 
we observe that, on the contrary, the $k$-th Segre mappings $t\mapsto 
\mathcal{Q}_k(t)$ are not appropriate to express the nondegeneracy 
conditions for other points than the origin.

\def\theexample{3.5.16}\begin{example}
{\rm 
The first Segre mapping $\mathcal{Q}_1$ of
the real algebraic hypersurface $M$ of
$\C^3$ given by the equation (which is the cubic
tangent to the Example~3.2.20)
\def\theequation{3.5.17}\begin{equation}
\bar w=w+i[2z_1\bar z_1+z_1^2\bar z_2+\bar z_1^2z_2]
\end{equation}
is the map
\def\theequation{3.5.18}\begin{equation}
(z_1,z_2,w)\longmapsto 
(w, \, 2iz_1, \, iz_1^2),
\end{equation}
which is only of rank $2$ at every point. This shows that the rank of
the Levi form of $M$ at the origin is equal to $1$, which is
true. Does this imply that the rank of the Levi form of $M$ is equal
to $1$ at every point~?  Of course not, because the Levi matrix
\def\theequation{3.5.19}\begin{equation}
\mathcal{H}(\varphi)(z,\bar z)=
\left(
\begin{array}{ll}
{\partial^2 \varphi \over \partial z_1 \partial \bar z_1} & 
{\partial^2 \varphi \over \partial z_2 \partial \bar z_1} 
 \\
{\partial^2 \varphi \over \partial z_1 \partial \bar z_2} &
{\partial^2 \varphi \over \partial z_2 \partial \bar z_2} 
\end{array}
\right)
=
\left(
\begin{array}{ll}
2 & 2z_1 \\
2\bar z_1 & 0
\end{array}
\right)
\end{equation}
is of rank $2$ at every point with $z_1\neq 0$, hence 
$M$ is Levi nondegenerate outside $\{z_1=0\}$.
This example shows that the $k$-th Segre mappings 
$\mathcal{Q}_k$ are appropriate to define nondegeneracy 
conditions only at the origin.
}
\end{example}

We shall therefore make some translations.
Let $p=t_p=(z_p,w_p)$ be a point varying in 
a neighborhood of the central point $p_0$. To express
the five nondegeneracy conditions using the original definitions
given in Section~3.2, we must choose coordinates vanishing
at $p$. Since we have already argued that the 
nondegeneracy conditions are then independent of
the choice of coordinates vanishing at $p$, we can 
simply make a translation of coordinates by setting
\def\theequation{3.5.20}\begin{equation}
t_1:= t-t_p, \ \ \ \ \ \ \
\text{\rm or equivalently} \ \ \ \ \ \ \ 
t=t_p+t_1.
\end{equation}
We shall assume that $t_p$ belongs to $M$, hence
$p^c=(t_p,\bar t_p)$ belongs to 
$\mathcal{M}$. The precise changes of coordinates will
therefore be
\def\theequation{3.5.21}\begin{equation}
\left\{
\aligned
z_1
& \ = & \
z-z_p, 
\ \ \ \ \ \ \ \ \ \ \
z
& \ = & \
z_1+z_p,\\
w_1
& \ = & \
w-w_p, 
\ \ \ \ \ \ \ \ \ \ \
w
& \ = & \
w_1+w_p,\\
\zeta_1
& \ = & \
\zeta-\bar z_p, 
\ \ \ \ \ \ \ \ \ \ \
\zeta
& \ = & \
\zeta_1+\bar z_p,\\
\xi_1
& \ = & \
\xi-\bar w_p, 
\ \ \ \ \ \ \ \ \ \ \
\xi
& \ = & \
\xi_1+\bar w_p.
\endaligned
\right.
\end{equation}
In the coordinates $t_1=(z_1,w_1)$ vanishing at $t_p$, we can 
represent the translation $\mathcal{M}_1$ of $\mathcal{M}$ by 
complex defining equations
\def\theequation{3.5.22}\begin{equation}
\xi_{1,j}=\Theta_{1,j}(\zeta_1,t_1), \ \ \ \ \ \ \ 
j=1,\dots,d.
\end{equation}
Of course, we may compute the $\Theta_{1,j}(\zeta_1,t_1)$ in terms 
of the $\Theta_j(\zeta,t)$ as follows. Since
$\bar w_{j,p}=\Theta_j(\bar z_p,t_p)$, we have
\def\theequation{3.5.23}\begin{equation}
\Theta_{1,j}(\zeta_1,t_1)=\xi_{1,j}=\xi_j
-\bar w_{j,p}=\Theta_j(\zeta,t)-\Theta_j(\bar z_p,t_p),
\end{equation}
which yields for $j=1,\dots,d$:
\def\theequation{3.5.24}\begin{equation}
\Theta_{1,j}(\zeta_1,t_1)=\Theta_j(\zeta_1+\bar z_p, t_1+t_p)-
\Theta_j(\bar z_p, t_p).
\end{equation}
If we now develope $\Theta_{1,j}(\zeta_1,t_1)$ in powers of 
$\zeta_1$ (as we did for $\Theta(\zeta,t)$)
\def\theequation{3.5.25}\begin{equation}
\Theta_{1,j}(\zeta_1,t_1)=\sum_{\beta\in\N^m}\, 
(\zeta_1)^\beta \, 
\Theta_{1,j,\beta}(t_1),
\end{equation}
then putting $\zeta_1=0$ in~\thetag{3.5.24}, we obtain for
$\beta=0\in\N^m$ the formula
\def\theequation{3.5.26}\begin{equation}
\Theta_{1,j,0}(t_1)=\sum_{\beta\in\N^m}\, 
(\bar z_p)^\beta \, 
[\Theta_{j,\beta}(t_1+t_p)-\Theta_{j,\beta}(t_p)],
\end{equation}
and differentiating~\thetag{3.5.24} at $\zeta_1=0$, we also 
obtain for all nonzero multiindices
$\beta\in\N^m\backslash \{0\}$ the important general explicit formulas
\def\theequation{3.5.27}\begin{equation}
\left\{
\aligned
\Theta_{1,j,\beta}(t_1)
& \
=\sum_{\gamma\in\N^m}\, 
{(\beta +\gamma)! \over \beta! \ \gamma !}\, 
(\bar z_p)^\gamma \, 
\Theta_{j,\beta+\gamma}(t_1+t_p)\\
& \
={1\over \beta!} \, 
\left[
\partial_{\zeta}^\beta \Theta_j(\zeta,t)\right]_{\zeta=\bar z_p, \, 
t=t_1+t_p}.
\endaligned\right.
\end{equation}
By means of these formulas, we have thus expressed $\Theta_{1,j}
(\zeta_1,t_1)$ in terms of $\Theta(\zeta,t)$.

Consequently, we can introduce the $k$-th 
Segre mapping in the coordinates $t_1$
\def\theequation{3.5.28}\begin{equation}
\mathcal{Q}_{1,k}: \ \ \C^n\ni t_1 \longmapsto
(\Theta_{1,j,\beta}(t_1))_{1\leq j\leq d,\, 
\vert\beta\vert\leq k}\in \C^{N_{d,n,k}},
\end{equation}
and speak of the five nondegeneracy conditions in terms of
$\mathcal{Q}_{1,k}$, since the coordinates $t_1$ are centered at
$t_p$.

However, a better way of considering the nondegeneracy conditions
at points $p$ in a neighborhood of $p_0$ would be to express
them in a single system of coordinates.

The following theorem provides the desired characterization of the five
nondegeneracy conditions by means of the morphism of jets of Segre varieties, 
expressed in a single system of coordinates.

\def\thetheorem{3.5.29}\begin{theorem}
Let $M$ be a real algebraic or analytic local 
generic submanifold of $\C^n$ given as usual by the equations~\thetag{3.1.2}
and let $J_\tau^k \underline{\mathcal{S}}_t$ be the morphism
of $k$-th jets of conjugate Segre varieties given explicitely by 
\def\theequation{3.5.30}\begin{equation}
\left\{
\aligned
(\zeta,t)\longmapsto 
J_\tau^k\underline{\mathcal{S}}_t:=
& \
\left(
\zeta, \ \left({1\over\beta!}\, 
\partial_\zeta^\beta \Theta_j(\zeta,t)\right)_{1\leq j\leq d,\, 
\vert \beta \vert \leq k}\right)\\
= & \
\left(
\zeta, \ \left(
\sum_{\gamma\in\N^m}\, 
{(\beta+\gamma)! \over \beta! \ \gamma !}\, 
(\zeta)^\gamma \, \Theta_{j,\beta+\gamma}(t)
\right)_{1\leq j\leq d,\, 
\vert \beta \vert \leq k}\right),
\endaligned
\right.
\end{equation}
which is a complex algebraic or analytic map defined in
$\Delta_{m+n}(\rho_1)$ with values in $\C^{m+N_{d,m,k}}$.  Let $t_p\in
M$ with $\vert t_p \vert < \rho_1$, let $(t_p,\bar t_p)\in\mathcal{M}$
and let $\ell_0\in\N$ with $\ell_0\geq 1$. Then
\begin{itemize}
\item[{\bf (1)}]
$M$ is Levi nondegenerate at $t_p$ if and only if $J_\tau^1
\underline{\mathcal{S}}_t$ is of rank equal to $m+n$ at $(\bar z_p, t_p)$.
\item[{\bf (2)}]
$M$ is $\ell_0$-finitely nondegenerate at $t_p$ if and only if 
$\ell_0$ is the smallest integer $k$ such that
$J_\tau^k\underline{\mathcal{S}}_t$ is
of rank equal to $m+n$ at $(\bar z_p,t_p)$.
\item[{\bf (3)}]
$M$ is $\ell_0$-essentially finite at $t_p$ if and only if 
$\ell_0$ is the smallest integer $k$ such that
$J_\tau^k\underline{\mathcal{S}}_t$ is a locally finite
complex algebraic or analytic map in a neighborhood of
$(\bar z_p, t_p)$.
\item[{\bf (4)}]
$M$ is $\ell_0$-Segre nondegenerate at $t_p$ if and only if 
$\ell_0$ is the smallest integer $k$ such that the restriction of
$J_\tau^k\underline{\mathcal{S}}_t$ to the second
Segre chain, namely
\def\theequation{3.5.31}\begin{equation}
\left\{
\aligned
\mathcal{S}_{\bar t_p}^2=
& \
\left\{\left(
z_1+z_p, \,
\overline{\Theta}(z_1+z_p, 
\bar t_p), \, 
\zeta_1+\bar z_p, \, \right. \right. \\
& \ \ \ \
\left.
\Theta(\zeta_1+\bar z_p,z_1+z_p, \overline{\Theta}(z_1+z_p,\bar t_p))
\right)\in \Delta_{2n}(\rho_1): \, z_1\in \C^m,\, \zeta_1\in\C^m\}
\endaligned\right.
\end{equation} 
is of generic rank equal to $2m$ at $(z_1,\zeta_1)=(0,0)$, hence
all over $\mathcal{S}_{\bar t_p}^2$.
\item[{\bf (5)}]
$M$ is $\ell_0$-holomorphically nondegenerate at $t_p$ if and only if 
$\ell_0$ is the smallest integer $k$ such that
$J_\tau^k\underline{\mathcal{S}}_t$ is of generic 
rank $m+n$ in a neighborhood of $(\bar z_p,t_p)$, hence all
over $\Delta_{m+n}(\rho_1)$.
\end{itemize}
\end{theorem}

\proof
To establish this theorem, it suffices to inspect the five
definitions given in Section~3.2 in the convenient system of coordinate
$t_1=t-t_p$ vanishing at $t_p$ which was introduced before
stating the theorem.

Indeed, thanks to the expressions~\thetag{3.5.26} and~\thetag{3.5.27}  of
$\Theta_{1,j,\beta}(t_1)$, we observe that except for $\beta=0$,
the components of $J_\tau^k\underline{\mathcal{S}}_t$ in~\thetag{3.5.30}
after the first $m$ components $(\zeta_1,\dots,\zeta_m)$ coincide with
the components $\Theta_{1,j,\beta}(t_1)$
of $\mathcal{Q}_{1,k}(t_1)$, after $(\zeta)^\gamma$ has
been replaced by $(\bar z_p)^\gamma$. For $\beta=0$, according
to~\thetag{3.5.26}, the difference between the
two mappings is only the constant $-\sum_{\beta\in\N^m}
\, (\bar z_p)^\beta \, \Theta_{j,\beta}(t_p)=-\bar w_p$, which 
disappears by differentiation.

Consequently, we verify the following relation between the Jacobian
matrix of $J_\tau^k \underline{\mathcal{S}}_t$ at $(\zeta,t)=(\bar
z_p, t_p)$ and the Jacobian matrix of $\mathcal{Q}_{1,k}$ at $t_1=0$,
{\it i.e.} at $t=t_p$:
\def\theequation{3.5.32}\begin{equation}
{\rm  Jac} (J_\tau^k \underline{\mathcal{S}}_t)(\bar z_p,t_p)=
\left(
\begin{array}{ll}
{\rm I}_{m\times m} & 
0 
 \\
\ast \ast \ast &
{\rm Jac}\, \mathcal{Q}_{1,k} (0)
\end{array}
\right)
\end{equation}
where ${\rm I}_{m\times m}$ denotes the identity $m\times m$ matrix
and $\ast \ast \ast$ some terms which we need not to compute.
We deduce at once that
\def\theequation{3.5.33}\begin{equation}
{\rm rk}_{\bar z_p,t_p}\left((\zeta,t)\mapsto
J_\tau^k\underline{\mathcal{S}}_t\right) = m+
{\rm rk}_0\left( t_1\mapsto 
\mathcal{Q}_{1,k}(t_1)\right).
\end{equation}

By expressing Levi nondegeneracy and $\ell_0$-finite nondegeneracy in
terms of $\mathcal{Q}_{1,k}(t_1)$ in the coordinates $t_1$ vanishing
at $t_p$, we immediately get characterizations {\bf (1)} and {\bf (2)}
of Theorem~3.5.29.

Since more generally, at a point $(z_1+\bar z_p, t_1+t_p)$ varying in
a neighborhood of $(\bar z_p, t_p)$, we have
\def\theequation{3.5.34}\begin{equation}
{\rm  Jac} (J_\tau^k \underline{
\mathcal{S}}_t)(z_1+\bar z_p,t_1+t_p)=
\left(
\begin{array}{ll}
{\rm I}_{m\times m} & 
0 
 \\
\ast \ast \ast &
{\rm Jac}\, \mathcal{Q}_{1,k} (t_1)
\end{array}
\right),
\end{equation}
we deduce at once that
\def\theequation{3.5.35}\begin{equation}
{\rm genrk}_\C\left(
(\zeta,t)\mapsto J_\tau^k
\underline{\mathcal{S}}_t
\right)=m+
{\rm genrk}_\C \left(t_1\mapsto 
\mathcal{Q}_{1,k}(t_1)\right).
\end{equation}
By expressing holomorphic nondegeneracy in terms of
$\mathcal{Q}_{1,k}(t_1)$ in the coordinates $t_1$ vanishing at $t_p$,
we immediately get the characterization {\bf (5)} of Theorem~3.5.29.

Next, we also observe that
\def\theequation{3.5.36}\begin{equation}
\left\{
\aligned
\dim_\C \, \C\{t_1\}/
& \
\left<
\Theta_{1,j,\beta}(t_1)
\right>_{1\leq j\leq d, \, 
\vert \beta \vert \leq k}=\\
& \
=
\dim_\C \, \C\{\zeta-\bar z_p, \, 
t-t_p\}\left/ 
\left< 
j_\tau^k\underline{\mathcal{S}}_t-
j_{\bar t_p}^k\underline{\mathcal{S}}_{t_p}
\right>.\right.
\endaligned
\right.
\end{equation}
We deduce the characterization {\bf (3)} of
Theorem~3.5.29.

For the last case {\bf (4)} to be considered, we notice that the restriction of 
$\mathcal{Q}_{1,k}(t_1)$ to the Segre variety of $M_1$
passing through the origin in coordinates $t_1$, which is 
by definition the map
\def\theequation{3.5.37}\begin{equation}
z_1\longmapsto \left(
\Theta_{1,j,\beta}(z_1,\overline{\Theta}_1(z_1,0))
\right)_{1\leq j\leq d, \, 
\vert \beta \vert \leq k},
\end{equation}
coincides thanks to~\thetag{3.5.27} (neglecting the constant
$-\bar w_p$ which appears for $\beta=0$) with 
\def\theequation{3.5.38}\begin{equation}
z_1\longmapsto 
\left(
\sum_{\gamma\in\N^m}\, 
{(\beta+\gamma)! \over
\beta! \ \gamma!}\, 
(\bar z_p)^\gamma \, 
\Theta_{j,\beta+\gamma}(
z_1+z_p, 
\overline{\Theta}_1(z_1,0)+w_p)
\right)_{1\leq j\leq d,\, 
\vert \beta \vert \leq k}.
\end{equation}
But since we have by~\thetag{3.5.23}
\def\theequation{3.5.39}\begin{equation}
\overline{\Theta}_1(z_1,0)+w_p=\overline{\Theta}(z_1+z_p, \bar t_p),
\end{equation}
we can rewrite~\thetag{3.5.38} as
\def\theequation{3.5.40}\begin{equation}
z_1\longmapsto 
\left(
\sum_{\gamma\in\N^m}\, 
{(\beta+\gamma)! \over
\beta! \ \gamma!}\, 
(\bar z_p)^\gamma \, 
\Theta_{j,\beta+\gamma}(
z_1+z_p,
\overline{\Theta}(z_1+z_p,\bar t_p))
\right)_{1\leq j\leq d,\, 
\vert \beta \vert \leq k}.
\end{equation}
We claim that
this mapping coincides with the last components of the restriction of
the mapping $J_\tau^k\underline{\mathcal{S}}_t$ to the second Segre
chain~\thetag{3.5.31}. Indeed, computing explicitely this restriction, 
we get exactly
\def\theequation{3.5.41}\begin{equation}
\left\{
\aligned
(z_1,\zeta_1)\longmapsto 
\left(
\zeta_1+\bar z_p, 
\left(
\sum_{\gamma\in\N^m}\, 
{(\beta+\gamma)! \over
\beta! \ \gamma!}\, 
(\zeta_1+\bar z_p)^\gamma \, \right.\right. \\
\ \ \ \ \ \ \ \ 
\left.
\left.
\Theta_{j,\beta+\gamma}(
z_1+z_p,
\overline{\Theta}(z_1+z_p,\bar t_p))
\right)_{1\leq j\leq d, \, 
\vert \beta \vert < k}
\right).
\endaligned\right.
\end{equation}
Consequently, we deduce
\def\theequation{3.5.42}\begin{equation}
{\rm genrk}_\C\left(
(z_1,\zeta_1)\longmapsto 
J_\tau^k\underline{\mathcal{S}}_t\vert_{\mathcal{S}_{\bar t_p}^2}
\right)=m+
{\rm genrk}_\C \left(
z_1 \mapsto \mathcal{Q}_{1,k}(z_1,\overline{\Theta}_1(z_1,0))
\right).
\end{equation}
By expressing Segre nondegeneracy in terms of
$\mathcal{Q}_{1,k}(t_1)$ in the coordinates $t_1$ vanishing at $t_p$,
we immediately get the characterization {\bf (4)}
of Theorem~3.5.29.

The proof of Theorem~3.5.29. is complete.
\endproof

\subsection*{3.5.43.~Essential holomorphic
dimension and Levi multitype} If $M$ is a local piece of generic
submanifold as above, we denote by $\ell_M$ the smallest integer $k$
such that Lemma~3.5.13 holds and we call it the {\sl Levi type of $M$}. We
denote by $m+n_M\leq m+n$ the generic rank of the mapping
$(t,\tau)\mapsto J_\tau^{\ell_M}\underline{\mathcal{S}}_t$ and we call
$n_M$ the {\sl essential holomorphic dimension of $M$}.  This
terminology is justified by the fact that locally in a neighborhood of
a Zariski-generic point $p\in M$, then $M$ is biholomorphically
equivalent to a product $\underline{M}_p'\times \Delta^{n-n_M}$, where
$\underline{M}_p'$ is a generic submanifold of codimension $d$ of
$\C^{n_M}$ ({\it see}~Theorem~3.5.48 below). 

By a specialization of
the second functional equation~\thetag{2.1.25}, which yields 
$w\equiv \overline{\Theta}(0,0,
\Theta(0,0,w))$, we see that the mapping
$w\mapsto \Theta(0,0,w)$ is invertible. Consequently, the rank and the
generic rank of the zeroth jet mapping ${\rm genrk}_\C \left((t,\tau)
\mapsto J_\tau^0 \underline{ \mathcal{S}}_t \right)=
(\zeta,\Theta(\zeta,z,w))$ is equal to $m+d$. Thus, the 
integer $n_M$ always satisfies the inequalities $d\leq n_M\leq n$.

Generally speaking, we may define $\lambda_{0,M}:={\rm genrk}_\C
\left((t,\tau) \mapsto J_\tau^0 \underline{ \mathcal{S}}_t \right)=d$
and for every $k=1,\dots,\ell_M$,
\def\theequation{3.5.44}\begin{equation}
\lambda_{k,M}:={\rm genrk}_\C\left(
(t,\tau)\mapsto 
J_\tau^k\underline{\mathcal{S}}_t
\right)-{\rm genrk}_\C\left(
(t,\tau)\mapsto 
J_\tau^{k-1}\underline{\mathcal{S}}_t
\right).
\end{equation} 
By Lemma~3.5.13, we have $\lambda_{1,M}\geq 1,\dots,
\lambda_{\ell_M,M}\geq 1$.
With these definitions, we have the relations
\def\theequation{3.5.45}\begin{equation}
{\rm genrk}_\C\left(
(t,\tau)\mapsto 
J_\tau^k\underline{\mathcal{S}}_t
\right)=
d+\lambda_{1,M}+\cdots+\lambda_{k,M},
\end{equation}
for $k=0,1,\dots,\ell_M$ and
\def\theequation{3.5.46}\begin{equation}
{\rm genrk}_\C\left(
(t,\tau)\mapsto 
J_\tau^k\underline{\mathcal{S}}_t
\right)=n_M=
d+\lambda_{1,M}+\cdots+\lambda_{\ell_M,M},
\end{equation}
for all $k\geq \ell_M$. It follows that 
we have the inequality
\def\theequation{3.5.47}\begin{equation}
\ell_M\leq \lambda_{1,M}+\cdots+\lambda_{\ell_M,M}=
n_M-d.
\end{equation}
We are now in position to state and to prove the main theorem of this
chapter in which we remind all the essential assumptions. Up to now,
we have worked locally in a neighborhood of a point $p_0\in M$. In the
following theorem, we observe that we may easily globalize our 
constructions, provided that $M$ is connected.

\def\thetheorem{3.5.48}\begin{theorem} Let $M$ be a connected real
algebraic or analytic generic submanifold in $\C^n$ of codimension
$d\geq 1$ and of CR dimension $m=n-d\geq 1$. Then there exist well
defined integers $n_M$, $\ell_M$ and $\lambda_{0,M},
\lambda_{1,M},\dots,\lambda_{\ell_M,M}$ and a proper real algebraic or analytic
subvariety $E$ of $M$ such that for every point $p\in M\backslash E$ and for
every system of coordinates $(z,w)$ vanishing at $p$ in which $M$ is
represented by defining equations $\bar w_j=\Theta_j(\bar z,t)$,
$j=1,\dots,d$, then the following four properties hold{\rm :}
\begin{itemize}
\item[{\bf (1)}]
$\lambda_{0,M}=d$, 
$d\leq n_M\leq n$ and $\ell_M\leq n_M-d$.
\item[{\bf (2)}]
For every $k=0,1,\dots,\ell_M$, the mapping of $k$-th order jets of
the conjugate complexified Segre varieties $(t,\tau)\mapsto
J_\tau^k\underline{\mathcal{S}}_t$ is of rank equal to
$m+\lambda_{0,M}+\cdots+\lambda_{k,M}$ at $(t_p,\bar t_p)=(0,0)$.
\item[{\bf (3)}]
$n_M=d+ \lambda_{1,M} +\cdots+ \lambda_{ \ell_M,M}$ and for every
$k\geq \ell_M$, the mapping of $k$-th order jets of the conjugate
complexified Segre varieties $(t,\tau) \mapsto J_\tau^k \underline{
\mathcal{S}}_t$ is of rank equal to $n_M$ at $(t_p, \bar t_p)=(0,0)$.
\item[{\bf (4)}]
There exists a complex algebraic or analytic change of coordinates
$t'=h(t)$ vanishing at $p$ and defined in a neighborhood of $p$ such
that the image $M_p':=h(M)$ is the product $\underline{M}_p'\times
\Delta^{n-n_M}$ of a real algebraic or analytic generic submanifold of
codimension $d$ in $\C^{n_M}$ by a complex polydisc $\Delta^{n-n_M}$.
Furthermore, at the central point $\underline{p}' \in \underline{M}_p'
\subset \C^{n_M}$, the generic submanifold $\underline{M}_p'$ is 
$\ell_M$-finitely nondegenerate, hence in particular its essential
holomorphic dimension $n_{\underline{M}_p'}$ coincides with 
$n_M$.
\end{itemize}
\end{theorem}

\proof
We fix $p_0\in M$ and coordinates $(z,w)$ as above vanishing at
$p_0$. Let $V_{p_0}$ be small neighborhood of $p_0$ in $M$. We first
define $E\cap V_{p_0}$: it consists of the set of points $p\in
V_{p_0}$ at which $\ell_p$ is not minimal, $n_p$ is not maximal and
the $\lambda_{k,p}$ are not maximal. Clearly, this set may be
described by the vanishing of a collection of minors of the Jacobian
matrix of the jet mapping $J_\tau^k \underline{ \mathcal{S}}_t$, so it
is a proper real algebraic or analytic subvariety $E_{p_0}$ of
$V_{p_0}$.  Next, we verify that the various $E_{p_0}$ glue together.
Indeed, let assume that $V_{p_0}$ and $V_{q_0}$ overlap.  Let $p\in
V_{p_0}\cap V_{q_0}$. In the intersection $V_{p_0}\cap V_{q_0}$, we
have to compare three real algebraic or analytic subvarieties
$E_{p_0}$, $E_p$ and $E_{q_0}$ defined in terms of three morphisms of
$k$-th order jets of conjugate Segre varieties. Using the important
relation given by Theorem~3.5.6 ({\it cf.} also Corollary~3.5.9), and using
an explicit description of the above mentioned collection of minors we
may establish easily that $E_{p_0}$ and $E_{q_0}$ coincide with $E_p$
in $V_{p_0}\cap V_{q_0}$.  Consequently, the various $E_{p_0}$ glue
together. Taking account of the considerations which precede the
statement of Theorem~3.5.48, this proves properties {\bf (1)}, {\bf (2)}
and {\bf (3)}.

Let us now prove {\bf (4)}. Let $p\in M\backslash E$.  We choose
coordinates $t=(z,w)$ vanishing at $p$. By assumption, for every
$k\geq \ell_M$, the $k$-th order Segre mapping $(\zeta,t)\mapsto
(\zeta, {1\over \beta!}\, (\partial_\zeta^\beta
\Theta_j(\zeta,t))_{1\leq j\leq d, \, \vert \beta\vert \leq k})$ is of
constant rank $m+n_M$ in a neighborhood of the origin (the point $p$) 
in $\C^{m+n}$.
In particular, this entails that at every point of the form $(0,t_p)$
in a neighborhood of the origin, the mapping $t\mapsto
(\Theta_{j,\beta}(t))_{1\leq j\leq d, \, \vert \beta \vert \leq k}$ is
of constant rank $n_M$. It follows from the constant rank theorem that
there exists an open neighborhood $V_0$ of the origin in $\C^n$ such
that the union of level sets $\mathcal{F}_q:=\{t\in V_0: \,
\Theta_{j,\beta}(t)=\Theta_{j,\beta}(q), \, j=1,\dots,d,\,
\beta\in\N^m\}$, for $q$ running in $V_0$, constitutes a local complex
algebraic or analytic foliation of $V_0$ by $(n-n_M)$-dimensional
complex manifolds. We can straighten this foliation to a product
$\Delta^{n_M}\times \Delta^{n-n_M}$, where the second factor
corresponds to the leaves of this foliation. Let $t'=h(t)$ denote such 
a straightening change of coordinates. Let $M_0':=h(M)$ be of equation
$\bar w'=\Theta'(\bar z',t')$. Thanks to Theorem~3.1.9, we observe that
this foliation is again defined by the level sets of the
functions $\Theta_{j,\beta}'(t')$, namely $\mathcal{F}_{p'}=
\{t'\in V_0': \, \Theta_{j,\beta}'(t')=\Theta_{j,\beta}'(p'), \, 
j=1,\dots,d, \, \beta\in \N^m\}$. To conclude that in these
coordinates, $M_0'$ is a product $\underline{M}_0'\times \Delta^{n-n_M}$, 
it suffices to establish the following lemma.

\def\thelemma{3.5.49}\begin{lemma}
If a point $p'\in V_0'$ belongs to $M_0'$,
then its leaf $\mathcal{F}_{p'}$ is entirely contained in $M_0'$.
\end{lemma}

\proof
Indeed, let $q'\in\mathcal{F}_{p'}$, {\it i.e.}
$\Theta_{j,\beta}'(t_{q'}')=\Theta_{j,\beta}'(t_{p'}')$ for all
$j$ and all $\beta$. It follows first that
\def\theequation{3.5.50}\begin{equation}
0=\bar w_{p'}'-\Theta'(\bar z_{p'}',t_{p'}')=
\bar w_{p'}'-\Theta'(\bar z_{p'}',t_{q'}').
\end{equation}
Next, thanks to the reality of $M_0'$, by Lemma~2.1.27, there exists an
invertible $d \times d$ matrix of power series $a' (t',\tau')$ such
that $w'- \overline{ \Theta}'(z',\tau')\equiv a'(t',\tau')\, [\xi'-
\Theta'( \zeta',t')]$, so we deduce $0=w_{ q'}'- \overline{ \Theta}'(
z_{q'}',\bar t_{p'}')$ and by conjugating
\def\theequation{3.5.51}\begin{equation}
0=\bar w_{q'}'-\Theta'(\bar z_{q'}',t_{p'}').
\end{equation}
Finally, using again the property
$\Theta_{j,\beta}'(t_{q'}')=\Theta_{j,\beta}'(t_{p'}')$ for all $j$ and all
$\beta$, we deduce
\def\theequation{3.5.52}\begin{equation}
0=\bar w_{q'}'-\Theta'(\bar z_{q'}',t_{q'}'),
\end{equation}
which shows that $q'\in M_0'$. This completes the proof of Lemma~3.5.49.
\endproof

The proof of Theorem~3.5.48 is complete.
\endproof

At a Zariski-generic point, the generic submanifold $M$ incorporates a
factor $\Delta^{n-n_M}$ which is in a certain sense ``flat'' with
respect to the point of view of CR geometry. After dropping this
innocuous factor, we come down to the study of a finitely
nondegenerate generic submanifold. Thus, in a certain sense, finitely
nondegenerate submanifolds $M$ are the ``generic'' models. This is why
it is interesting to state a direct corollary of Theorem~3.5.48 about
holomorphically nondegenerate submanifolds.

\def\thecorollary{3.5.53}\begin{corollary} 
Let $M$ be a connected real algebraic or analytic generic submanifold
in $\C^n$ of codimension $d\geq 1$ and of CR dimension $m=n-d\geq 1$.
Assume that $M$ is holomorphically nondegenerate. Then
\begin{itemize}
\item[{\bf (1)}]
There exists an integer $\ell_M$ with $1\leq \ell_M\leq m$ and a
proper real algebraic or analytic subset $E$ of $M$ such that $M$ is
$\ell_M$-finitely nondegenerate at every point of $M\backslash E$.
\item[{\bf (2)}]
There exists a proper complex algebraic or analytic subset
$\mathcal{E}^{\rm exc}$ defined in a neighborhood of $M$ in $\C^n$ which
depends only on $M$ such that $M$ is finitely nondegenerate at a point
$p$ if and only if $p\in M\backslash \mathcal{E}^{\rm exc}$.
\item[{\bf (3)}]
In general, the inclusion $(\mathcal{E}^{\rm exc}\cap M)\subset E$ is
strict.
\end{itemize}
\end{corollary} 

\proof
In a polydisc neighborhood $V_0\subset
\C^n$ of an arbitrary point $p_0\in M$, we have defined
in \S3.4.16 a local exceptional locus $\mathcal{E}_{p_0}^{\rm exc}\subset V_0$
such that $M\cap \mathcal{E}_{p_0}^{\rm exc}$ consists exactly of finitely
degenerate points. Thanks to their biholomorphic invariance, these local
complex algebraic or analytic subsets $\mathcal{E}_{p_0}^{\rm exc}$ glue
together in a well defined global exceptional locus $\mathcal{E}^{\rm
exc}$ defined in a neighborhood of $M$ in $\C^n$.  Finally, $E\backslash
(\mathcal{E}^{\rm exc}\cap M)$ consists of points which are
$k$-finitely nondegenerate for some $k\geq \ell_M+1$, and so is clearly
nonempty in general. This completes the proof of
Corollary~3.5.53.
\endproof

\section*{\S3.6.~Transformation rules for jets of Segre varieties}

We now establish the biholomorphic invariance of the mapping of jets
of Segre varieties. As in \S3.1.5, let $t'=h(t)$ with 
inverse $t=h'(t')$ be a complex algebraic or
analytic local biholomorphism fixing the origin 
and let $M':=h(M)$. Theorems~3.1.9 follows from 
the two relations~\thetag{8.5.2} and~\thetag{8.5.3}
after putting $\zeta=0$, taking account of the fact that
$\Theta_j(0,t)$ coincides with $\Theta_{j,0}(t)$ (in the notation 
$\Theta_{j,\beta}(t)$).
Theorem~3.5.6 follows immediately from
the two relations~\thetag{8.5.2} and~\thetag{8.5.3} just below, if we
decide to consider the last argument of $Q_{j,\beta}$ simply as
functions of $(\zeta,(\Theta_j(\zeta,t))_{1\leq j\leq d})$.

\def\thetheorem{3.6.1}\begin{theorem}
For every $j=1,\dots,d$ and every $\beta\in\N^m$, there exists a
universal rational mapping in its variables $Q_{j,\beta}$ such that
\def\theequation{3.6.2}\begin{equation}
\left\{
\aligned
{}
&
{1\over \beta!}\, 
{\partial^{\vert \beta\vert} \Theta_j'\over
\partial (\zeta')^\beta}\left(
\bar f(\zeta,\Theta(\zeta,t)),
h(t)\right)\equiv\\
&
\ \ \ \ \
\equiv
Q_{j,\beta}\left(
\left(\partial_\zeta^{\beta_1}\Theta_{j_1}(\zeta,t)
\right)_{1\leq j_1\leq d, \,
\vert \beta_1\vert \leq \vert \beta\vert}, \,
\left(\partial_t^{\alpha_1}
\bar h_{i_1}(\zeta,\Theta(\zeta,t))
\right)_{1\leq i_1\leq n,\, 
\vert \alpha_1 \vert\leq 
\vert \beta \vert}\right).
\endaligned\right.
\end{equation}
Here, the $Q_{j,\beta}$ are algebraic or analytic in a neighborhood of
the constant jet $((\partial_\zeta^{\beta_1}\Theta_{j_1}(0,0)
)_{1\leq j_1\leq d, \, \vert \beta_1\vert \leq \vert \beta\vert},
\,(\partial_t^{\alpha_1} \bar h_{i_1}(0,0) )_{1\leq
i_1\leq n,\, \vert \alpha_1 \vert\leq \vert \beta \vert})$.
Equivalently, we have the relations for all $j=1,\dots, d$ and all
$\beta\in\N^m$:
\def\theequation{3.6.3}\begin{equation}
\left\{
\aligned
{}
&
\Theta_{j,\beta}'(h(t))\equiv
\sum_{\gamma\in\N^m}\, 
(-1)^\gamma \, 
{(\beta+\gamma)!\over \beta! \ \gamma !}\, 
\bar f(\zeta,\Theta(\zeta,t))^\gamma\, 
\\
&
Q_{j,\beta+\gamma}\left(
\left(\partial_\zeta^{\beta_1}\Theta_{j_1}(\zeta,t)
\right)_{1\leq j_1\leq d, \,
\vert \beta_1\vert \leq \vert \beta\vert+\vert \gamma \vert}, \,
\left(\partial_t^{\alpha_1}
\bar h_{i_1}(\zeta,\Theta(\zeta,t))
\right)_{1\leq i_1\leq n,\, 
\vert \alpha_1 \vert\leq 
\vert \beta \vert+\vert \gamma \vert}
\right).
\endaligned\right.
\end{equation}
\end{theorem}

\proof
As in \S3.1.5, we may assume that the complex
defining equations of $M'$ are of the form 
$\bar w_j'=\Theta_j'(\bar z', t')$, $j=1,\dots,d$ in 
coordinates $t'=(z',w')\in\C^m\times \C^d$. Geometrically 
speaking, this means that the linear mapping $\pi_{z'}'\circ
dh: T_0^cM \to \C_{z'}^m$ is submersive, where $\pi_{z'}': \C^n\to 
\C_{z'}^m$ is the natural projection onto the $z'$-space.
We decompose the mapping $h(t)=(f(t),g(t))\in\C^m\times \C^d$. 
By
complexifying the fundamental relations $\bar g_j(\tau)=\Theta_j'(\bar
f(\tau), h(t))$, $j=1,\dots,d$, which express that $h$ maps $M$ into
$M'$ and by replacing $\xi$ by $\Theta(\zeta,\tau)$, we obtain the
following power series identities
\def\theequation{3.6.4}\begin{equation}
\bar g_j(\zeta,\Theta(\zeta,t))\equiv \Theta_j'(\bar f(\zeta,\Theta(\zeta,t)),
h(t)), \ \ \ \ \ \ j=1,\dots,d.
\end{equation}
We differentiate this relation with respect to $\zeta_k$, for $k=1,\dots,m$.
Remembering that the explicit expression of the natural 
basis of complexified
$(0,1)$-vector fields is given by
\def\theequation{3.6.5}\begin{equation}
\underline{\mathcal{L}}_k=
{\partial \over \partial \zeta_k}+\sum_{j=1}^1 \, 
{\partial \Theta_j\over \partial \zeta_k}(\zeta,t)\, 
{\partial \over \partial \xi_j},
\end{equation}
for $k=1,\dots,m$, we immediately see that differentiation with respect to 
$\zeta_k$ of a power series $\psi(\zeta,\Theta(\zeta,t))$ is equivalent
to applying the vector field $\underline{\mathcal{L}}_k$ to $\psi$, viewed
as a derivation. So we get by the chain rule
\def\theequation{3.6.6}\begin{equation}
\underline{\mathcal{L}}_k \bar g_j(\zeta,\Theta(\zeta,t))\equiv
\sum_{l=1}^m\, \underline{\mathcal{L}}_k
\bar f_l(\zeta,\Theta(\zeta,t))\, 
{\partial \Theta_j'\over \partial \zeta_l'}(
\bar f(\zeta,\Theta(\zeta,t)), h(t)).
\end{equation}
Since $\pi_{z'}'\circ dh: T_0^cM\to \C_{z'}^m$ is submersive, we have
the nonvanishing of the following determinant
\def\theequation{3.6.7}\begin{equation}
{\rm det}\, 
\left(\underline{\mathcal{L}}_{k_1} \bar f_{k_2}(0)\right)_{
1\leq k_1,k_2\leq m}\neq 0.
\end{equation}
Hence we can divide locally for $(\zeta,t)$ in a neighborhood of the
origin by the determinant
\def\theequation{3.6.8}\begin{equation}
\mathcal{D}(\zeta,t):=
{\rm det}\, 
\left(\underline{\mathcal{L}}_{k_1} \bar f_{k_2}(\zeta,\Theta(\zeta,t))\right)_{
1\leq k_1,k_2\leq m}.
\end{equation}
Viewing~\thetag{3.6.6} as an inhomogeneous linear system and
using the classical rule of Cramer, for every $j=1,\dots,d$, we can 
solve the first partial derivatives $\partial \Theta_j'/\partial \zeta_k'$
with respect to the other terms, which yields
expressions of the form
\def\theequation{3.6.9}\begin{equation}
{\partial \Theta_j'\over \partial \zeta_k'}(
\bar f(\zeta,\Theta(\zeta,t)), h(t))\equiv
{R_{j,k}\left((\underline{\mathcal{L}}_{k_1'} \bar h_{i_1}(\zeta,
\Theta(\zeta,t)))_{1\leq i_1\leq n, \, 
1\leq k_1'\leq m}\right)\over
\mathcal{D}(\zeta,t)}.
\end{equation}
Here, by the very application of Cramer's rule, it follows that the
terms $R_{j,k}$ are certain universal polyomials of
determinant type (some minors). 
By differentiating again~\thetag{3.6.9} with respect to
the variables $\zeta_k$, using again Cramer's rule, we get that for every
pair of integers $(k_1, k_2)$ with $1\leq k_1,\,k_2\leq m$ and for
every $j=1,\dots,d$, there exist a universal polynomial
$R_{j,k_1,k_2}$ such that we can write
\def\theequation{3.6.10}\begin{equation}
{\partial^2\Theta_j'\over \partial \zeta_{k_1}'\partial \zeta_{k_2}'}
\left(\bar f(\zeta,\Theta(\zeta,t)),h(t)\right)\equiv
{R_{j,k_1,k_2}\left((\underline{\mathcal{L}}_{k_1',k_2'} \bar h_{i_1}(\zeta,
\Theta(\zeta,t)))_{
1\leq i_1\leq n, \, 1\leq k_1',k_2'\leq m}\right)\over
\mathcal{D}(\zeta,t)^3}.
\end{equation}
The reader should notice the exponent $3$ in the denominator, with the
decomposition ``$3$''$=$``$2$''$+$``$1$'' 
where ``$2$'' comes from the derivatives of
the quotient $R_{j,k}/\mathcal{D}$ in~\thetag{3.6.9} and where ``$1$'' comes
from the second application of Cramer's rule.

Remind that for an arbitrary multi-index $\beta=(\beta_1,\dots,\beta_m)
\in\N^m$, we denote by $\underline{\mathcal{L}}^\beta$ the antiholomorphic
derivation of order $\vert \beta\vert$ defined by
$(\underline{\mathcal{L}}_1)^{\beta_1}\cdots (\underline{\mathcal{L}}_m)^{\beta_m}$.

Differentiating more generally the relations~\thetag{3.6.4} with
respect to $\zeta^\beta=\zeta_1^{\beta_1}\cdots \zeta_m^{\beta_m}$,
we see by an easy induction that for every multindex $\beta\in\N^m$
and for every $j=1,\dots,d$, there exists a complicated but universal
polynomial $R_{j,\beta}$ such that the following identity holds:
\def\theequation{3.6.11}\begin{equation}
{1\over \beta!}\, 
{\partial^{\vert \beta\vert} \Theta_j'\over
\partial (\zeta')^\beta}
\left(\bar 
f(\zeta,\Theta(\zeta,t)),h(t)\right)
\equiv
{R_{j,\beta}\left((\underline{\mathcal{L}}^{\beta_1}
\bar h_{i_1}(\zeta,\Theta(\zeta,t)))_{1\leq i_1\leq n,\,
1\leq \vert \beta_1 \vert \leq \vert \beta\vert}\right)\over
[\mathcal{D}(\zeta,t)]^{2\vert \beta\vert -1}}.
\end{equation} 
An important observation
is in order. The composed derivations $\underline{\mathcal{L}}^{\beta_1}$ are
certain differential operators with nonconstant coefficients. Using
the explicit expressions of the $\underline{\mathcal{L}}_k$, we see that all
these coefficients are certain universal polynomials of the collection
of partial derivatives $(\partial^{\vert \beta_2\vert}
\Theta_{j_2}(\zeta,t)/\partial z^{\beta_2})_{1\leq j_2\leq d, \,
1\leq \vert \beta_2\vert \leq \vert \beta_1\vert}$. Thus the numerator
of~\thetag{3.6.11} becomes a certain universal (computable
by means of combinatorial formulas) algebraic
or analytic function of the collection
\def\theequation{3.6.12}\begin{equation}
\left(
\left(\partial_\zeta^{\beta_1}\Theta_{j_1}(\zeta,t)
\right)_{1\leq j_1\leq d, \,
\vert \beta_1\vert \leq \vert \beta\vert}, \,
\left(\partial_t^{\alpha_1}
\bar h_{i_1}(\zeta,\Theta(\zeta,t))
\right)_{1\leq i_1\leq n,\, 
\vert \alpha_1 \vert\leq 
\vert \beta \vert}\right)
\end{equation} 
A similar property holds for the denominator. In conclusion, 
we have constructed the rational mapping $Q_{j,\beta}$ satisfying~\thetag{3.6.2}. 

For the second part of Theorem~3.6.1, let us rewrite the
relations~\thetag{3.6.2} in the following more explicit form, simply
obtained by developing the left hand side with respect to the powers
$(\bar f)^\gamma$:
\def\theequation{3.6.13}\begin{equation}
\left\{
\aligned
\sum_{\gamma\in\N^m}\, 
{(\beta+\gamma)! \over \beta! \ \gamma!}\, 
& \
\bar f(\zeta,\Theta(\zeta,t))^\gamma \, 
\Theta_{j,\beta+\gamma}'(h(t))\equiv \\
\equiv 
& \
Q_{j,\beta}\left(
\left(\partial_\zeta^{\beta_1}\Theta_{j_1}(\zeta,t)
\right)_{1\leq j_1\leq d, \,
\vert \beta_1\vert \leq \vert \beta\vert}, \,
\left(\partial_t^{\alpha_1}
\bar h_{i_1}(\zeta,\Theta(\zeta,t))
\right)_{1\leq i_1\leq n,\, 
\vert \alpha_1 \vert\leq 
\vert \beta \vert}\right).
\endaligned\right.
\end{equation}
We may interpret this infinite collection of
identities as an infinite upper triangular inhomogeneous linear system
with unknowns being the $\Theta_{j,\beta}'(h(t))$. The inversion of this
infinite triangular matrix is in fact very elementary. 
Indeed, by interpreting Taylor's formula at a purely formal level, 
we see that if we are given an infinite collection
of equalities with complex coefficients and with $\zeta\in\C^m$ which 
is of the form
\def\theequation{3.6.14}\begin{equation}
\sum_{\gamma\in\N^m}\, 
{(\beta+\gamma)! \over \beta! \ \gamma!}\,
\zeta^\gamma \, 
\Theta_{j,\beta+\gamma}'=Q_{j,\beta},
\end{equation}
for all $j=1,\dots,d$ and all $\beta\in\N^m$, then we can solve the unknowns
$\Theta_{j,\beta}'$ in terms of the right hand side terms $Q_{j,\beta}$
by means of a totally similar formula, except for signs:
\def\theequation{3.6.15}\begin{equation}
\sum_{\gamma\in\N^m}\,
(-1)^\gamma\, 
{(\beta+\gamma)! \over \beta! \ \gamma!}\, 
\zeta^\gamma \,
Q_{j,\beta+\gamma}=\Theta_{j,\beta}', 
\end{equation}
for all $j=1,\dots,d$ and all $\beta\in\N^m$.  Applying this
observation to~\thetag{3.6.13}, we deduce the
representation~\thetag{3.6.3} in Theorem~3.6.1, which completes the
proof.
\endproof

\section*{\S3.7.~Local geometry of CR submanifolds at a Zariski-generic 
point} 

Let $M\subset \C^n$ be a {\it not necessarily generic}\, connected real
algebraic or analytic CR submanifold of codimension $d$, of CR
dimension $m$ and of holomorphic codimension $c=d-n+m$. Combining
Theorem~2.1.9, Theorem~2.1.32, Corollary~2.8.6 and Theorem~3.5.48, we
obtain the following local explicit coordinate representation of $M$
locally in a neighborhood of a Zariski-generic point. This theorem
will be useful in Part~II of this memoir.

\def\thetheorem{3.7.1}\begin{theorem}
There exists a proper real algebraic or analytic subset $E$ of $M$ and
integers $m_1$, $m_2$, $d_1$, $d_2$, $c$ which depend only on
$M$ and which satisfy
\def\theequation{3.7.2}\begin{equation}
\left\{
\aligned
d=
& \
d_1+d_2+2c,\\
m=
& \
m_1+m_2, \\
\endaligned\right.
\end{equation}
and in the case where $m_1\geq 1$, there exist moreover two integers $\ell_M$
and $\nu_M$ which satisfy
\def\theequation{3.7.3}\begin{equation}
\left\{
\aligned
\ell_M 
& \
\leq m_1,\\
\nu_M 
& \
\leq d_1+1
\endaligned\right.
\end{equation}
such that for every point $p_0\in 
M\backslash E$, there exist local
complex analytic or algebraic normal coordinates
\def\theequation{3.7.4}\begin{equation}
(z_1,z_2,w_1,w_2,w_3)\in 
\C^{m_1}\times 
\C^{m_2}\times
\C^{d_1}\times
\C^{d_2}\times
\C^c
\end{equation}
vanishing at $p_0$ and complex algebraic or analytic defining functions
$\Theta_{1,j_1}(\bar z_1,z_1, w_1,w_2)$, $j_1=1,\dots,d_1$ which
converge normally in $\Delta_{2m_1+d_1+d_2}(2\rho_1)$ for some
$\rho_1>0$, which satisfy 
$\Theta_{1,j_1}(0,z_1,w_1,w_2)\equiv 0$ for
$j_1=1,\dots,d_1$, and which are independent of $(\bar z_2,z_2)$ such that $M$
is represented locally in a neighborhood of $p_0$ by the complex
defining equations
\def\theequation{3.7.5}\begin{equation}
\left\{
\aligned
0=
& \
w_3, \\
0=
& \
\bar w_2-w_2,\\
0=
& \
\bar w_1-\Theta_1(\bar z_1,z_1,w_1,w_2),
\endaligned\right.
\end{equation}
in the polydic $\Delta_n(\rho_1)$ and such that, moreover,
for every constant $u_{2,q}\in \R^{d_2}$, the
generic submanifold $M_{1,u_{2_q}}$ of $\C^{m_1}\times \C^{d_1}$ defined 
by the complex equations 
\def\theequation{3.7.6}\begin{equation}
0=\bar w_1-\Theta_1(\bar z_1,z_1,w_1,u_{2,q}),
\end{equation}
which identifies with the intersection of $M$ with the complex
subspace $\{w_2= u_{2,q}=ct., \, w_3=0\}$, is minimal of Segre type
$\nu_M$ and $\ell_M$-finitely nondegenerate at $(z_1,w)=(0,0)$.  In the
case where $m_1=0$, the third equation in~\thetag{3.7.5} should be
replaced by the simpler vectorial equation $\bar w_1=w_1$, hence in
this case $M$ identifies in a neighborhood of $p_0$ with the
intersection of $\Delta_n (\rho_1)$ with the Levi-flat product
$\C^{m_2} \times \R^{d_1} \times \R^{d_2} \times \{0\}$ in $\C^{m_2}
\times \C^{d_1} \times \C^{d_2} \times \C^c$.
\end{theorem}

\bigskip

\noindent
{\Large \bf Chapter~4: Nondegeneracy conditions for power series CR mappings}

\section*{\S4.1.~CR-horizontal nondegeneracy conditions for power series 
CR mappings}

\subsection*{\S4.1.1.~Nondegeneracy conditions for power series mappings}
The datum of a formal, or complex algebraic or complex analytic
mapping $h: \C^n\to \C^{n'}$ with $h(0)=0$ is equivalent to the datum
of a collection of $n'$ power series $(h_1(t),\dots,h_{n'}(t))$, where
$t\in \C^n$, with the $h_{i'}(t)$ being scalar power series vanishing
at the origin and belonging to $\C\dl t\dr$, to $\C\{t\}$
or to $\mathcal{A}_\C\{t\}$.  We introduce five classical nondegeneracy
conditions, which we formulate in the case where $h(t)\in\C\dl
t\dr^{n'}$.

\def\thedefinition{4.1.2}\begin{definition}
{\rm 
A formal power series mapping $h(t)=(h_1(t),\dots,
h_{n'}(t))$, with components 
$h_{i'}(t)\in\C\dl t\dr$, $i'=1,\dots,n'$, is
called
\begin{itemize}
\item[{\bf (1)}]
{\sl Invertible} if $n'=n$ and ${\rm det}\, 
([\partial h_{i_1}/\partial t_{i_2}](0))_{1\leq i_1,i_2\leq n}\neq 0$.
\item[{\bf (2)}]
{\sl Submersive} if $n\geq n'$ and there exist integers
$1\leq i(1)< \cdots < i(n')\leq  n$  such that ${\rm det}\, 
([\partial h_{i_1'}/\partial t_{i(i_2')}](0))_{1\leq i_1',
i_2'\leq n'}\neq 0$.
\item[{\bf (3)}]
{\sl Finite} if the ideal generated by the components
$h_1(t),\dots,h_{n'}(t)$ is of finite codimension in $\C\dl t\dr$.
This implies $n'\geq n$.
\item[{\bf (4)}]
{\sl Dominating} if $n\geq n'$ and there exist integers 
$1\leq i(1)< \cdots < i({n'})\leq n$ such that
${\rm det}([\partial h_{i'}/\partial t_{i(i_2')}](t))_{1\leq
i_1',i_2'\leq n'}\not \equiv 0$ in $\C\dl t\dr$.
\item[{\bf (5)}]
{\sl Transversal} if there does not exist a nonzero power series
$G(t_1',\dots,t_{n'}')\in \C\dl t_1',\dots,t_{n'}'\dr$ such that
$G(h_1(t),\dots,h_{n'}(t))\equiv 0$ in $\C\dl t\dr$.
\end{itemize}
}
\end{definition}
It is elementary to see that invertibility implies submersiveness
which implies domination. Furthermore, we have:

\def\thelemma{4.1.3}\begin{lemma}
If a formal power series is either invertible, submersive
or dominating, then it is transversal.
\end{lemma}

\proof
It suffices to prove the statement in the case where $h$ is
dominating.  Suppose on the contrary that there exists a nonzero power
series $G(t_1',\dots,t_{n'}')\in \C\dl t_1',\dots,t_{n'}'\dr$ such
that $G(h_1(t),\dots, h_{n'}(t))\equiv 0$ in $\C\dl t \dr$. By
differentiating this identity with respect to the $t_i$ and looking at
the linear homogeneous system so obtained, we deduce that all the
$n'\times n'$ minors (of maximal dimension) of the Jacobian matrix
$([\partial h_{i'}/\partial t_i](t))_{1\leq i\leq n, \, 1\leq i' \leq
n'})$ vanish identically in $\C\dl t \dr$, contradicting domination.
\endproof

\def\thelemma{4.1.4}\begin{lemma}
In the equidimensional case $n'=n$, invertibility implies {\rm (}and
is in fact equivalent to{\rm )} submersiveness, which implies
finiteness, which implies domination, which finally implies
transversality.
\end{lemma}

\proof
Classical result from local complex analytic geometry, which may
easily reproved by the reader or found for instance~[1], [6] or in the
references therein.
\endproof

\subsection*{4.1.5. Power series CR mappings and CR-horizontal 
nondegeneracy conditions} Let $M$ and $M'$ be two real algebraic or
analytic generic submanifolds of $\C^n$ and of $\C^{n'}$. Let
$r_j(t,\bar t):=\bar w_j-\Theta_j(\bar z,t)$, $j=1,\dots,d$ and
$r_{j'}'(t',\bar t'):=\bar w_{j'}'- \Theta_{j'}'(\bar z',t')$,
$j'=1,\dots,d'$, be defining equations for $M$ and for $M'$.
Following the definition given in \S2.1.5, we say that {\sl $h$ is a
{\rm (}local{\rm )} power series CR mapping from $M$ to $M'$} if there
exists a $d'\times d$ matrix $a(t,\bar t)$ of formal, analytic or
algebraic power series such that, in vectorial notation
\def\theequation{4.1.6}\begin{equation}
r'(h(t),\bar h(\bar t))\equiv a(t,\bar t)\,
r(t,\bar t).  
\end{equation}
Setting $\bar t=0$ in this matrix identity, we get
\def\theequation{4.1.7}\begin{equation}
r'(h(t),0)\equiv a(t,0)\, r(t,0).
\end{equation}
The set defined by the equations $r_j(t,0)=0$ is of
course the Segre variety $S_{\bar p_0}$ passing through the
origin. Similarly, the set defined by $r_{j'}'(t',0)=0$,
$j'=1,\dots,d'$, is the Segre variety $S_{\bar p_0'}$.
Then~\thetag{4.1.7} shows that $h$ induces a power series mapping
from the Segre variety $S_{\bar p_0}$ to $S_{\bar p_0'}'$. We can thus
apply Definition~4.1.2 at the level of Segre varieties.

\def\thedefinition{4.1.8}\begin{definition}
{\rm
A power series CR mapping $h: M\to M'$ is called
\begin{itemize}
\item[{\bf (1)}]
{\sl CR-invertible at $p_0$} if
$m=m'$ and if the induced formal mapping $h\vert_{S_{\bar p_0}}: 
S_{\bar p_0}\to S_{\bar p_0'}'$ between Segre varieties passing
through the origin is a formal invertible mapping at $0$.
\item[{\bf (2)}]
{\sl CR-submersive at $p_0$} 
if $m\geq m'$ and if the induced formal mapping $h\vert_{S_{\bar p_0}}: 
S_{\bar p_0}\to S_{\bar p_0'}'$ between Segre varieties passing
through the origin is a formal submersion at $0$.
\item[{\bf (3)}]
{\sl CR-finite at $p_0$} if the induced formal mapping
$h\vert_{S_{\bar p_0}}: S_{\bar p_0}\to S_{\bar p_0'}'$ between Segre
varieties passing through the origin is finite at $0$.
\item[{\bf (4)}]
{\sl CR-dominating at $p_0$} 
if the induced formal mapping $h\vert_{S_{\bar p_0}}: 
S_{\bar p_0}\to S_{\bar p_0'}'$ between Segre varieties passing
through the origin is dominating at $0$.
\item[{\bf (5)}]
{\sl CR-transversal at $p_0$} 
if the induced formal mapping $h\vert_{S_{\bar p_0}}: 
S_{\bar p_0}\to S_{\bar p_0'}'$ between Segre varieties passing
through the origin is transversal at $0$.
\end{itemize}
}
\end{definition}

\subsection*{4.1.9.~Complexified mapping}
We may reformulate these definitions in a more concrete way as follows.
First of all, we may of course
complexify the defining identities~\thetag{4.1.6}, which
yields $r'(h(t),\bar h(\tau))\equiv
a(t,\tau) \, r(t,\tau)$. These complexified
identities show that the {\sl 
complexified mapping} 
\def\theequation{4.1.10}\begin{equation}
h^c(t,\tau):=(h(t),\bar h(\tau))\in\C\dl t\dr^n\times \C\dl \tau \dr^n.
\end{equation}
induces a power series mapping from $\mathcal{M}$ to 
$\mathcal{M}'$. As explained in \S1.2.2, the complexified
mapping $h^c$ stabilizes the two pairs
of foliations $(\mathcal{F}, \underline{\mathcal{F}})$ and
$(\mathcal{F}', \underline{\mathcal{F}}')$, namely 
$h^c$ sends (conjugate) complexified Segre varieties of $\mathcal{M}$ to 
(conjugate) complexified Segre varieties of 
$\mathcal{M}'$. 

After replacing $\xi$ by $\Theta(\zeta,t)$ in~\thetag{4.1.10}, we obtain
the following power series identities in $\C\dl \zeta,t\dr$:
\def\theequation{4.1.11}\begin{equation}
\bar g_{j'}(\zeta,\Theta(\zeta,t))\equiv
\Theta_{j'}'(\bar f(\zeta,\Theta(\zeta,t)),h(t)), \ \ \ \ \ \
j'=1,\dots,d'.
\end{equation}
In fact, as may be easily established, these identities are equivalent
to the existence of a $d'\times d$ matrix of formal power series
$a(t,\tau)$ which satisfies~\thetag{4.1.6}.  However, {\it throughout
this memoir, we shall work only with the more convenient fundamental
power series identities}~\thetag{4.1.11}.

Of course, the Segre variety passing through $p_0$ is represented by
\def\theequation{4.1.12}\begin{equation}
S_{\bar p_0}: \ \ \ \ w_j=\overline{\Theta}_j(z,0), 
\ \ \ \ \ j=1,\dots,d.
\end{equation}
Similarly, the Segre variety passing through 
$p_0'$ is represented by
\def\theequation{4.1.13}\begin{equation}
S_{\bar p_0'}': \ \ \ \ w_{j'}'=\overline{\Theta}_{j'}'(z',0), 
\ \ \ \ \ j'=1,\dots,d'.
\end{equation}
If we split $h=(f,g)\in\C^{m'}\times \C^{d'}$, then the restriction of
$h$ to $S_{\bar p_0}$ coincide with the power series 
CR mapping
\def\theequation{4.1.14}\begin{equation}
\C^n\ni z\longmapsto 
\left(f(z,\overline{\Theta}(z,0)), \, 
\overline{\Theta}'(f(z,\overline{\Theta}'(z,0)), \, 0)\right)\in \C^{m'}\times
\C^{d'}.
\end{equation}
By projection onto the $\C^{m'}\times \{0\}$, we may of course
identify this mapping with the {\sl CR-horizontal part} of the mapping
defined by
\def\theequation{4.1.15}\begin{equation}
\C^n\ni z\longmapsto 
f(z,\overline{\Theta}(z,0))\in \C^{m'}.
\end{equation}
In the special case where $M'$ is given in normal coordinates, we have
$\Theta'(z',0)\equiv 0$, hence the last $d'$ terms in~\thetag{4.1.7} all
vanish and the identification of $h\vert_{S_{\bar p_0}}$ with its
CR-horizontal part is trivial in this case (but we shall avoid for the
moment to introduce normal coordinates). With such notations, we can
reformulate the above five nondegeneracy conditions more concretely.

\def\thedefinition{4.1.16}\begin{definition}
{\rm
Such a power series CR mapping $h: M\to M'$ is 
\begin{itemize}
\item[{\bf (1)}]
{\sl CR-invertible at $p_0$} 
if its CR-horizontal part is a formal invertible map at $0$.
\item[{\bf (2)}]
{\sl CR-submersive at $p_0$} 
if its CR-horizontal part is a formal submersion at $0$.
\item[{\bf (3)}]
{\sl CR-finite at $p_0$} if its CR-horizontal part
is finite at $0$.
\item[{\bf (4)}]
{\sl CR-dominating at $p_0$} 
if its CR-horizontal part is dominating at $0$.
\item[{\bf (5)}]
{\sl CR-transversal at $p_0$} if
its CR-horizontal part is transversal at $0$.
\end{itemize}
}
\end{definition}

To conclude, as direct corollaries of
Lemmas~4.1.3 and~4.1.4, we have

\def\thelemma{4.1.17}\begin{lemma}
If a formal power series CR mapping $h: M\to M'$ is either
CR-invertible, CR-submersive or CR-dominating, then it is
CR-transversal. Furthermore, in the CR-equidimensional case $m'=m$,
CR-invertibility implies {\rm (}and is in fact equivalent to{\rm )} CR
submersiveness, which implies CR-finiteness, which implies
CR-domination, which finally implies CR-transversality.
\end{lemma}

\section*{\S4.2.~Segre nondegeneracy conditions for power series
CR mappings} 

\subsection*{4.2.1.~Preliminaries}
After having introduced the five CR-horizontal nondegeneracy
conditions on the power series CR mapping $h$, we may introduce {\sl
Segre nondegeneracy conditions} on $h$ which are related to the
reflection principle, which will be studied in the next chapters of
Part~II of this memoir. Consequently, {\it this Section~4.2 is of
utmost importance to understand the whole memoir}.

As in \S4.1.9, let $h^c: \mathcal{M}\to
\mathcal{M}'$ be a complexified power series CR mapping, namely start
with the fundamental power series identities in $\C\dl \zeta,t\dr$:
\def\theequation{4.2.2}\begin{equation}
\bar g_{j'}(\zeta,\Theta(\zeta,t))-
\Theta_{j'}'(\bar f(\zeta,\Theta(\zeta,t)),h(t))\equiv 0,
\end{equation} 
for $j'=1,\dots,d'$.  As in~\thetag{2.3.8}, let
$\underline{\mathcal{L}}_1,\dots, \underline{\mathcal{L}}_m$ be the
basis of complexified $(0,1)$ vector fields tangent to $\mathcal{M}$
given by
\def\theequation{4.2.3}\begin{equation}
\underline{\mathcal{L}}_k=
{\partial\over \partial \zeta_k}+ \sum_{j=1}^d\, 
{\partial \Theta_j \over \partial \zeta_k} 
(\zeta,t)\, 
{\partial \over \partial \xi_j}. 
\end{equation}
Let $\beta=(\beta_1, \dots, \beta_m) \in \N^m$.  Applying the
derivations $\underline{ \mathcal{ L}}^\beta= \underline{ \mathcal{
L}}_1^{ \beta_1} \cdots \underline{ \mathcal{ L}}_m^{ \beta_m}$
to~\thetag{4.2.2} and without writing the
arguments, we get
\def\theequation{4.2.4}\begin{equation}
\underline{\mathcal{L}}^\beta\bar g_{j'}-
\sum_{\gamma'\in\N^{m'}}\, 
\underline{\mathcal{L}}^\beta(\bar f^\gamma)\, 
\Theta_{j',\gamma'}'(h)\equiv 0, 
\end{equation}
for all $\beta \in \N^m$,
for $j'=1,\dots,d'$ and for $(t,\tau)\in\mathcal{M}$.  By developping
the derivations $\underline{\mathcal{L}}^\beta$, we observe that

\def\thelemma{4.2.5}\begin{lemma}
For every $i'=1,\dots,n'$ and every $\beta\in\N^m$, there exists a
polynomial $P_{i',\beta}$ in the jet $J_\tau^{\vert \beta
\vert} \bar h(\tau)$ with coefficients being power series in
$(t,\tau)$ which depend only on the defining functions
$\xi_j-\Theta_j(\zeta,t)$ of $\mathcal{M}$ and which can be computed
by means of some combinatorial formula, such that
\def\theequation{4.2.6}\begin{equation}
\underline{\mathcal{L}}^\beta \bar h_{i'}(\tau)\equiv
P_{i',\beta}(t,\tau,J_\tau^{\vert \beta \vert} \bar h(\tau)).
\end{equation}
\end{lemma}

\proof
For $\vert \beta \vert =1$, this follows by inspecting the
coefficients of the vector fields $\underline{ \mathcal{ L}}_k$
in~\thetag{4.2.4}. By induction, assuming that such a formula holds for
all $\beta\in\N^m$ with $\vert \beta \vert =k\geq 1$, applying the
vector fields $\underline{ \mathcal{ L}}_1,\dots,
\underline{\mathcal{L}}_m$ to~\thetag{4.2.6} and using
the chain rule, we get a similar type of formula for all $\beta\in\N^m$ with
$\vert \beta \vert = k+1$. Clearly, the coefficients of $P_{i',\beta}$
depend on the partial derivatives of the $\Theta_j(\zeta,t)$ and one
can refine this statement by providing an explicit combinatorial
formula (which we shall not need).
\endproof

In formula~\thetag{4.2.6}, we have written the first two arguments of
$P_{i',\beta}$ to be $(t,\tau)$. In fact, by inspecting the
coefficients of the vector fields $\underline{\mathcal{L}}_k$, these
first two arguments are $(\zeta,t)$. However, since we are always
considering the variables $(t,\tau)$ to belong to $\mathcal{M}$, we
have to replace everywhere $\tau$ by $\Theta(\zeta,t)$ or $w$ by
$\overline{\Theta}(z,t)$. Consequently, we
can identify a function of $(t,\tau)$ with
a function of $(\zeta,t)$ or with a function 
of $(z,\tau)$. Before going further, we shall make the
following convention which will allow us to make some slight abuse
of notation on occasion.

\def\theconvention{4.2.7}\begin{convention}
{\rm 
Let $k,l\in\N$. On the complexification $\mathcal{M}$, we identify
(notationally) a power series written under the complete form
\def\theequation{4.2.8}\begin{equation}
R(t,\tau,J^k
h(t),J^l\bar h(\tau)),
\end{equation}
with a power series written under one
of the following four forms
\def\theequation{4.2.9}\begin{equation}
\left\{
\aligned
{\bf (1)} & \ \ 
R\left(t,\zeta,\Theta(\zeta,t),J^k h(t), J^l\bar
h(\zeta,\Theta(\zeta,t))\right),\\
{\bf (2)} & \ \ 
R\left(t,\zeta,J^k h(t),
J^l\bar h(\zeta,\Theta(\zeta,t)\right),\\
{\bf (3)} & \ \
R\left(z,\overline{\Theta}(z,\tau), J^kh(z,\overline{\Theta}(z,\tau)),
J^l\bar h(\tau)\right),\\
{\bf (4)} & \ \
R\left(z,\tau, J^kh(z,\overline{\Theta}(z,\tau)),
J^l\bar h(\tau)\right).
\endaligned\right.
\end{equation}
}
\end{convention}

Admitting this convention, applying Lemma~4.2.5 and using the chain rule,
we deduce that for every $j'=1,\dots,d'$ and every $\beta\in\N^m$,
there exist formal power series $R_{j',\beta}'=
R_{j',\beta}'(t,\tau,J_\tau^{\vert \beta \vert} : \, t')$ which depend only
on the defining equations of $\mathcal{M}$ and on the defining
equations of $\mathcal{M}'$ such that we can rewrite the
equations~\thetag{4.2.4} under the general form
{\large 
\def\theequation{4.2.10}\begin{equation}
\underline{\mathcal{L}}^\beta [\bar g_{j'}(\tau)-
\Theta_{j'}(\bar f(\tau), h(t))]=:
R_{j',\beta}'(t,\tau,J_\tau^{\vert \beta \vert} \bar h(\tau): \, 
h(t))\equiv 0,
\end{equation}}
for $j'=1,\dots,d'$, where the identity ``$\equiv 0$'' is understood
``on $\mathcal{M}$'', namely as a formal power series identity in
$\C\dl \zeta, t\dr$ after replacing $\xi$ by $\Theta(\zeta,t)$ or
equivalently, as a formal power series identity in $\C\dl z , \tau\dr$
after replacing $w$ by $\overline{\Theta}(z,\tau)$. We shall
constantly refer to these identities~\thetag{4.2.10} in the sequel.

{\it Importantly, the smoothness of the power series $R_{j',\beta}'$
is the minimum of the two smoothesses of $M$ and of $M'$. This crucial
remark will be the basis of all the various formal reflection
principles developed in the next chapters of Part~II}.

For instance, the power series $R_{j',\beta}'$ are all complex analytic
if $M$ is real analytic and if $M'$ is real algebraic, even if the power
series CR mapping $h(t)$ was assumed to be purely formal and non convergent.

By a careful inspection of the application of the chain rule in the
development of ~\thetag{4.2.5}, we even see that each $R_{j',\beta}'$ is
relatively polynomial with respect to the derivatives of positive
order $(\partial_\tau^\alpha \bar h(\tau))_{1\leq \vert \alpha \vert
\leq \vert \beta \vert}$.

\subsection*{4.2.11.~Segre nondegeneracy conditions for CR mappings} 
We are now ready to define nondegeneracy conditions for power series
CR mappings which generalize the nondegeneracy conditions for generic
submanifolds introduced in Chapter~3.  In the equations~\thetag{4.2.10},
we replace $h(t)$ by a new independent variable $t'\in\C^{n'}$, we set
$(t,\tau)=(0,0)$, and we define the following collection of power
series
\def\theequation{4.2.12}\begin{equation}
\Psi_{j',\beta}'(t'):=
\left[\underline{\mathcal{L}}^\beta\bar g_{j'}-
\sum_{\gamma\in\N^m}\, 
\underline{\mathcal{L}}^\beta(\bar f^{\gamma'})\, 
\Theta_{j',\beta}'(t')\right]_{t=\tau=0}, 
\end{equation}
for $j'=1,\dots,d'$ and $\beta\in\N^m$. Here, for $\beta=0$, we
mean that $\Psi_{j',0}'(t')=-\Theta_{j'}'(0,t')$.
Clearly, an equivalent way of defining $\Psi_{j',\beta}'(t')$ is 
as follows
\def\theequation{4.2.13}\begin{equation}
\Psi_{j',\beta}'(t'):=
R_{j',\beta}'(0,0,J_\tau^{\vert \beta \vert} \bar h(0) : \, t').
\end{equation}
Now, just before introducing the desired five definitions, we make the
following crucial heuristic remark.  When $n=n'$, $M=M'$ and $h={\rm
Id}$, we drop the dashes and we denote by $T$ (instead of $t'$) the new
independent variable, hence we may compute
\def\theequation{4.2.14}\begin{equation}
\left\{
\aligned
\Psi_{j,\beta}(T)
& \
=\left[\underline{
\mathcal{L}}^\beta\xi_j-
\sum_{\gamma\in\N^m}\, 
\underline{\mathcal{
L}}^\beta(\zeta)^\gamma\, 
\Theta_{j,\beta}(T)
\right]_{t=\tau=0}\\
& \ 
=\left[\underline{
\mathcal{L}}^\beta \Theta_j(\zeta,t)-
\beta! \, \Theta_{j,
\beta}(T)\right]_{t=\tau=0}\\
& \
=
\beta! \, [\Theta_{j,\beta}(0)-\Theta_{j,\beta}(T)].
\endaligned
\right.
\end{equation}
Consequently, we see that up to an affine combination, we recover with
$\Psi_{j,\beta}(T)$ the components of the infinite Segre mapping~\thetag{3.1.3}.
Furthermore, using the computation ~\thetag{4.2.14}, we may check easily
that the following five nondegeneracy conditions just below are a
generalization to CR mappings of the five nondegeneracy conditions
introduced in Section~3.2 for generic submanifolds of $\C^n$.

\def\thedefinition{4.2.15}\begin{definition}
{\rm 
The formal, algebraic or analytic power series CR mapping $h$ is called
\begin{itemize}
\item[{\bf (1)}]
{\sl Levi-nondegenerate} at the origin if the mapping
\def\theequation{4.2.16}\begin{equation}
t'\mapsto \left(
R_{j',\beta}'(0,0,J_\tau^{\vert \beta \vert} \bar h(0) :\,
t')
\right)_{1\leq j'\leq d', \, \vert
\beta \vert \leq 1}
\end{equation}
is of rank $n'$ at $t'=0$. This condition requires the
dimensional inequality $d'(m+1)\leq n'$.
\item[{\bf (2)}]
{\sl $\ell_1$-finitely nondegenerate} at the origin
if there exists an integer $\ell$ such that
the mapping
\def\theequation{4.2.17}\begin{equation}
t'\mapsto \left(
R_{j',\beta}'(0,0,J_\tau^{\vert 
\beta \vert} \bar h(0) :\, t')
\right)_{1\leq j'\leq d', \, \vert
\beta \vert \leq \ell}
\end{equation}
is of rank $n'$ at $t'=0$, and if $\ell_1$ is the smallest such
integer $\ell$.
\item[{\bf (3)}]
{\sl $\ell_1$-Segre finite} at the origin if there exists an integer
$\ell$ such that the mapping
\def\theequation{4.2.18}\begin{equation}
t'\mapsto \left(
R_{j',\beta}'(0,0,J_\tau^{\vert \beta \vert} \bar h(0) :\,
t')
\right)_{1\leq j'\leq d', \, \vert
\beta \vert \leq \ell}
\end{equation}
is locally finite at $t'=0$,
and if $\ell_1$ is the smallest such integer.
\item[{\bf (4)}]
{\sl $\ell_1$-Segre nondegenerate} if there exist an integer $\ell$,
integers ${j_*'}^1,\dots, {j_*'}^{n'}$ with $1\leq {j_*'}^{i'}\leq d'$
for $i'=1,\dots,n'$ and multiindices $\beta_*^1,\dots,\beta_*^{n'}$
with $\vert \beta_*^{i'} \vert \leq \ell$ for $i'=1,\dots,n'$, such
that the determinant
\def\theequation{4.2.19}\begin{equation}
{\rm det}\, \left(
{\partial R_{{j_*'}^{i_1'},\beta_*^{i_1'}}'\over
\partial t_{i_2'}'}\left(
z,\overline{\Theta}(z,0),0,0,
J^{\vert \beta \vert} \bar h(0) :\,
h(z,\overline{\Theta}(z,0))
\right)
\right)_{
1\leq i_1',i_2'\leq n'}
\end{equation}
does not vanish identically in $\C\dl z \dr$,
and if $\ell_1$ is the smallest such integer $\ell$.
\item[{\bf (5)}]
{\sl $\ell_1$-holomorphically nondegenerate} at the origin if there
exists an integer $\ell$, integers ${j_*'}^1,\dots, {j_*'}^{n'}$ with
$1\leq {j_*'}^{i'}\leq d'$ for $i'=1,\dots,n'$ and multiindices
$\beta_*^1,\dots,\beta_*^{n'}$ with $\vert \beta_*^{i'} \vert \leq
\ell$ for $i'=1,\dots,n'$, such that the determinant
\def\theequation{4.2.20}\begin{equation}
{\rm det}\, \left(
{\partial R_{{j_*'}^{i_1'},\beta_*^{i_1'}}'\over
\partial t_{i_2'}'}\left(
0,0,0,0,
J^{\vert \beta \vert} \bar h(0) :\, 
h(t)
\right)
\right)_{
1\leq i_1',i_2'\leq n'}
\end{equation}
does not vanish identically in $\C\dl t \dr$, and if $\ell_1$ is the
smallest such integer $\ell$.
\end{itemize}
}
\end{definition}

These five nondegeneracy conditions for power series CR mappings are
of utmost importance for the reflection principle and they will be
studied in Part~II of this memoir. Some of them are suggested by
K.~Diederich and S.M.~Webster in~[9].  The notion of $\ell_1$-finite
nondegeneracy unifies conditions which appear in the works~[1], [8],
[14], [17], [22], [23], [27], [30], [34], [35], [36].

We notice that Levi nondegeneracy implies finite nondegeneracy which
implies Segre finiteness. However, Segre finiteness and Segre
nondegeneracy are totally independent conditions,
as shown by the following two examples.

\def\theexample{4.2.21}\begin{example}
{\rm 
We already know that essential finiteness of $M$ does not imply Segre
nondegeneracy, hence simply by considering the identity map of the
hypersurface $v=z_1\bar z_1(1+z_1\bar z_2)$ of $\C^3$, we see that
{\bf (3)} does not imply {\bf (4)} in Definition~4.2.15 just above.

On the other hand, it is easy to check that the mapping
$(z_1,w)\mapsto (z_1,0,w)$ from $M: \bar w=w+i\, z_1^2\bar z_1^2$ to
$M': \bar w'=w'+i[ {z_1'}^2{\bar{z_1'}}^2+ z_1'{\bar{z_2'}}^2+ \bar
z_1' {z_2'}^2]$ is Segre finite at the origin but it is not
Segre nondegenerate.
}
\end{example}

We also mention that the above five nondegeneracy conditions are
meaningful for sufficiently smooth local CR mappings, by considering
the Taylor series of $M$ at $p_0$, of $h$ at $p_0$ and
of $M'$ at $p_0'$.

\subsection*{4.2.22.~Necessary conditions}
Coming back to~\thetag{4.2.12}, noticing that $\bar f(0)=0$, we see
that the constant $\underline{\mathcal{L}}^\beta(\bar f^\gamma)
\vert_{t=\tau=0}$ vanishes for $\vert \gamma \vert > \vert \beta
\vert$, whence every $\Psi_{j',\beta}'(t')$ is an affine combination
with constant coefficients of some (but not all) of the power series
$\Theta_{j_1',\beta_1'}'(t')$, for $1\leq j_1'\leq d'$ and $\vert
\beta_1' \vert\leq \vert \beta \vert$. Based on this observation, we
deduce immediately properties {\bf (1)}, {\bf (2)}, {\bf (3)} and {\bf
(5)} of Lemma~4.2.23 just below, which states necessary conditions for
$h$ to be nondegenerate in each one of the above five senses. The
proof of {\bf (4)} is also elementary.

\def\thelemma{4.2.23}\begin{lemma}
Let $h: M\to M'$ be a power series CR mapping as above.
\begin{itemize}
\item[{\bf (1)}]
If $h$ is Levi-nondegenerate at the origin, then
$M'$ is Levi-nondegenerate at the origin. 
\item[{\bf (2)}]
If $h$ is $\ell_1$-finitely nondegenerate at the origin, then 
$M'$ is $\ell_0'$-finitely nondegenerate at the origin for some
$\ell_0'\leq \ell_1$.
\item[{\bf (3)}]
If $h$ is $\ell_1$-essentially finite at the origin, then 
$M'$ is $\ell_0'$-essentially finite at the origin for some
$\ell_0'\leq \ell_1$.
\item[{\bf (4)}]
If $h$ is $\ell_1$-Segre nondegenerate at the origin, then 
$M'$ is $\ell_0'$-Segre nondegenerate at the origin for some
$\ell_0'\leq \ell_1$.
\item[{\bf (5)}]
If $h$ is $\ell_1$-holomorphically nondegenerate at the origin, then 
$M'$ is $\ell_0'$-holomorphically nondegenerate at the origin for some
$\ell_0'\leq \ell_1$.
\end{itemize}
\end{lemma}

\section*{\S4.3.~Study of CR-transversal power series CR mappings}

We now show that CR-transversality of the mapping $h$ insures that it
enjoys exactly the same nondegeneracy condition as the target
$M'$. The condition of CR transversality is much more general than the
condition of domination, because for instance it does not impose any
dimensional inequality between $m$ and $m'$ or between $n$ and
$n'$. We shall end this section with the proof of the following
theorem, which is quite long and technical, but of central
importance. Here, we assume that $M$ and $M'$ are either algebraic or
analytic, and that $h$ is algebraic, analytic or even
formal. According to Lemma~4.1.17, the following lemma applies to many
situtations.

\def\thetheorem{4.3.1}\begin{theorem}
Assume that $h$ is CR-transversal at $p_0$. Then the following
five properties hold:
\begin{itemize}
\item[{\bf (1)}]
If $M'$ is Levi nondegenerate at $p_0'$, then
$h$ is $\ell_1$-finitely nondegenerate at $p_0$ for some $\ell_1\geq 1$.
\item[{\bf (2)}]
If $M'$ is $\ell_0'$-finitely nondegenerate at $p_0'$, then 
$h$ is $\ell_1$-finitely nondegenerate at $p_0$ for some
$\ell_1\geq \ell_0'$.
\item[{\bf (3)}]
If $M'$ is $\ell_0'$-essentially finite at $p_0'$, then 
$h$ is $\ell_1$-Segre finite at $p_0$ for some
$\ell_1\geq \ell_0'$.
\item[{\bf (4)}]
If $M'$ is $\ell_0'$-Segre nondegenerate at $p_0'$, then 
$h$ is $\ell_1$-Segre nondegenerate at $p_0$ for
some $\ell_1\geq \ell_0'$.
\item[{\bf (5)}]
If $M'$ is $\ell_0'$-holomorphically nondegenerate at 
$p_0'$, and if moreover $h$ is transversal at 
$p_0$, then $h$ is $\ell_1$-holomorphically nondegenerate at
$p_0$ for some $\ell_1\geq \ell_0'$.
\end{itemize}
\end{theorem}

\proof
In order to simplify a little bit the notations and the computations, we
shall assume that the coordinates $(z,w)$ for $M$ near $p_0$ and
$(z',w')$ for $M'$ near $p_0'$ are normal. Thus, the Segre variety
$S_{\bar p_0}$ is given by $\{(z,0)\}$ instead of
$\{(z,\overline{\Theta}(z,0))\}$ and similarly for $S_{\bar p_0'}'$,
which will slightly simplify the presentation of the formal calculations below.

We remind that by the fundamental definition \thetag{4.2.10}, the
functions $R_{j',\beta}'(t,\tau, J^{\vert \beta \vert} \bar h(\tau)
:\, t')$ are the power series development of
$\underline{\mathcal{L}}^\beta r_{j'}'(t',\bar h(\tau))$. Here, the
integer $j'$ satisfies $1\leq j'\leq d'$ and the multiindex $\beta$
belongs to $\N^m$.

We consider the gradient of $R_{j',\beta}'$ with respect to the
distinguished variable $t'$:
\def\theequation{4.3.2}\begin{equation}
\nabla_{t'} \, R_{j',\beta}'(t,\tau, 
J^{\vert \beta \vert}\bar h(\tau) :\,
t'):=
\left(
{\partial R_{j',\beta}'\over \partial t_i'}
(t,\tau,J^{\vert \beta \vert}\bar h(\tau) :\, t')
\right)_{1\leq i'\leq n'},
\end{equation}
considered as a vertical vector, {\it i.e.} as a 
$n' \times 1$ matrix. Also, we shall work in the sequel with 
a fixed $j'$ and we shall consider the expression 
\def\theequation{4.3.3}\begin{equation}
\left(
\nabla_{t'} \, R_{j',\beta}'
(t,\tau,J^{\vert \beta \vert}\bar h(\tau) :\, t')
\right)_{\beta \in \N^m}
\end{equation}
as an $n'\times \infty$ matrix. We introduce a new notation.  Let
$\nu\in\N$ with $\nu\geq 1$, let $x=(x_1,\dots,x_\nu)\in\K^\nu$, let
$\mu\in\N$ with $\mu\geq 1$ and let $A(x)$ by a $ \mu\times \infty$
matrix of power series.  By ${\rm genMrk} \, A(x)$, we denote the
generic rank of the matrix $A(x)$, which is defined to be the largest
integer $\kappa\leq \mu$ such that there exists a $\kappa\times
\kappa$ minor of $A(x)$ which does not vanish identically as a power
series in $x$. The letter ``M'' in {\rm genMrk} stands for the word
``Matrix''. The notation ${\rm genMrk}$ should not be confused with the
notation ${\rm genrk}_\C$ introduced in \S2.1.5.

In the sequel, we shall in fact put $(t,\tau):=(0,0)$ in
$R_{j',\beta}'$. Then, the generic rank of the matrix
\def\theequation{4.3.4}\begin{equation}
{\rm genMrk}\, \left(
\nabla_{t'} \, R_{j',\beta}'
(0,0,J^{\vert \beta \vert}\bar h(0) :\, t')
\right)_{\beta \in \N^m}
\end{equation}
is of course the largest integer $\kappa'$ such that there exists a
$\kappa'\times \kappa'$ minor of~\thetag{4.3.4} which does not vanish
identically as a power series in $t'$.  In the case where we put
$t'=0$, we even do not have to speak of generic rank, hence we simply
denote by
\def\theequation{4.3.5}\begin{equation}
{\rm Mrk}\, \left(
\nabla_{t'} \, R_{j',\beta}'
(0,0,J^{\vert \beta \vert}\bar h(0) :\, 0)
\right)_{\beta \in \N^m}
\end{equation}
the rank of a $n'\times \infty$ constant matrix. Of course, we use the same
notations ${\rm genMrk}$ and ${\rm Mrk}$ for the truncated matrices
where we allow only $\vert \beta \vert\leq k$, for some integer $k\in
\N$.

First of all we prove part {\bf (2)} of Theorem~4.3.1, which contains of course part
{\bf (1)}. We state a technical lemma which holds essentially in the
case of codimension $d=1$.

\def\thelemma{4.3.6}\begin{lemma}
Fix $j'$ with $1\leq j'\leq d'$.  Let $n_{j'}'$ be the integer defined
by
\def\theequation{4.3.7}\begin{equation}
n_{j'}':= {\rm Mrk}\, 
\left(
\nabla_{t'} \, \Theta_{j',\gamma'}'(0)
\right)_{\gamma'\in\N^{m'}}.
\end{equation}
Assume that $h$ is CR-transversal at $p_0$. Then we also have
\def\theequation{4.3.8}\begin{equation}
{\rm Mrk}\, 
\left(
\nabla_{t'} \, R_{j',\beta}'(0,0,J^{\vert \beta \vert} \bar h(0) :\,
0)
\right)_{\beta\in\N^m}=n_{j'}'.
\end{equation}
\end{lemma}

We first show that this technical lemma implies part {\bf (2)} of the
theorem. Indeed, by the definition of $R_{j',\beta}'$ ({\it see}\, the
expression~\thetag{4.2.5}) $\nabla_{t'} \, R_{j',\beta}'(0,0,J^{\vert
\beta \vert} \bar h(0) :\, 0)$ is a linear combination with complex
coefficients of the vectors $\nabla_{t'} \, \Theta_{j',\gamma'}'(0)$,
for $\vert \gamma'\vert \leq \vert \beta \vert$, namely more precisely
\def\theequation{4.3.9}\begin{equation}
\nabla_{t'}\, R_{j',\beta}'(0,0,
J^{\vert \beta \vert} \bar h(0) :\, 0)
=
-\sum_{\gamma'\in\N^{m'}}\, 
\underline{\mathcal{L}}^\beta (\bar f^{\gamma'})(0) \, 
[\nabla_{t'}\, \Theta_{j',\gamma'}'(0)].
\end{equation}
We notice that this sum is in fact truncated, with 
$\vert \gamma' \vert \leq \vert \beta \vert$.
Consequently, we deduce that
\def\theequation{4.3.10}\begin{equation}
{\rm Span}_\C\, 
\left\{\nabla_{t'} \, 
R_{j',\beta}'(0,0,J^{\vert \beta \vert} \bar h(0) :\,
0) : \, \beta \in \N^m\right\}
\end{equation}
is automatically contained in 
\def\theequation{4.3.11}\begin{equation}
{\rm Span}_\C\, 
\left\{
\nabla_{t'} \, \Theta_{j',\gamma'}'(0)
 : \, \gamma' \in \N^{m'}\right\}.
\end{equation}
But thanks to the rank condition stated in Lemma~4.3.6 (to be proved
below), we deduce
that these two subspaces must coincide. In other words, for
each $j'=1,\dots,d'$, there
exists integers $\ell_{1,j'}$ and $\ell_{0,j'}'$ with 
$\ell_{1,j'}\geq \ell_{0,j'}'$ such that 
\def\theequation{4.3.12}\begin{equation}
\left\{
\aligned
{\rm Span}_\C\, 
\left\{\nabla_{t'} \, 
R_{j',\beta}'(0,0,J^{\vert \beta \vert} \bar h(0) :\,
0) : \, \vert \beta \vert \leq \ell_{1,j'} \right\}=\\
\ \ \ \ \ \ \ \ \ 
{\rm Span}_\C\, 
\left\{
\nabla_{t'} \, \Theta_{j',\gamma'}'(0)
 : \, \vert \gamma' \vert 
\leq \ell_{0,j'}'\right\}.
\endaligned\right.
\end{equation}
Finally, by assumption of $\ell_0'$-finite nondegeneracy of $M'$ at
$p_0'$, the vector space spanned by the collection of terms in the
right hand sides of~\thetag{4.3.12}, when $j'=1,\dots,d'$, is equal to the
whole of $\C^{n'}$, with of course $\ell_0'=\max
(\ell_{0,1}',\dots,\ell_{0,n'}')$.  It follows that the vector space
spanned by the collection of terms in the left hand sides
of~\thetag{4.3.12}, when $j'=1,\dots,d'$, is also equal to the whole of
$\C^{n'}$. In conclusion, setting of course $\ell_1=\max
(\ell_{1,1},\dots,\ell_{1,n'})\geq \ell_0'$, we have shown that
$h$ is $\ell_1$-finitely nondegenerate at $p_0$.
It remains to establish the technical lemma.

\proof
So we fix $j'$. We proceed by contradiction. Assume that 
the rank $\kappa'$ of the $n'\times \infty$ matrix 
\def\theequation{4.3.13}\begin{equation}
\left(
\nabla_{t'}\, R_{j',\beta}'(0,0,J^{\vert \beta \vert} \bar h(0) :\,
0)\right)_{\beta\in\N^m}
\end{equation}
is stricly smaller than $n_{j'}'$, namely $\kappa'\leq n_{j'}'-1$.
Equivalently, there exist multiindices $\beta_1,\dots,\beta_{\kappa'}\in
\N^m$ such that for every multiindex $\beta\in\N^m$ different
from $\beta_1,\dots,\beta_{\kappa'}$, there exist constanst
$\Lambda_\beta^1,\dots,\Lambda_\beta^{\kappa'}$ such that
we can write
\def\theequation{4.3.14}\begin{equation}
\left\{
\aligned
\nabla_{t'} \, R_{j',\beta}'(0,0,J^{\vert
\beta \vert} \bar h(0) :\, 0)=
& \
\Lambda_\beta^1 \left[\nabla_{t'}\, 
R_{j',\beta_1}'(0,0, J^{\vert \beta \vert} \bar h(0) :\, 
0)\right]+\cdots+
\\
& \
+\Lambda_\beta^{\kappa'} \left[ \nabla_{t'}\, 
R_{j',\beta_{\kappa'}}'(0,0, 
J^{\vert \beta \vert} \bar h(0) :\, 
0)\right].
\endaligned\right.
\end{equation}
Now, replacing in these linear combinations the terms $\nabla_{t'} \,
R_{j',\beta}'$ by their explicit expression~\thetag{4.3.9}, we obtain
\def\theequation{4.3.15}\begin{equation}
\left\{
\aligned
\sum_{\gamma'\in\N^{m'}}\, 
\underline{\mathcal{L}}^\beta (\bar f^{\gamma'})(0) \, 
& \
[\nabla_{t'} \, \Theta_{j',\gamma'}'(0)]=\\
& \
\Lambda_\beta^1 \left( 
\sum_{\gamma'\in\N^{m'}}\, 
\underline{\mathcal{L}}^{\beta_1} (\bar f^{\gamma'})(0)\, 
[\nabla_{t'} \, \Theta_{j',\gamma'}'(0)]\right)+\cdots+\\
& \
\Lambda_\beta^{\kappa'}\left(\sum_{\gamma'\in\N^{m'}}\, 
\underline{\mathcal{L}}^{\beta_{\kappa'}} (\bar f^{\gamma'})(0)\, 
[\nabla_{t'} \, \Theta_{j',\gamma'}'(0)]\right).
\endaligned\right.
\end{equation}
As we are in normal coordinates, the conjugate complexified
Segre variety passing through the origin is given by 
$\underline{\mathcal{S}}_0=\{(0,0,\zeta,0)\}$ in $(z,w,\zeta,\xi)$
coordinates. Let $\gamma'\in\N^{m'}$. The restriction 
of $\bar f^{\gamma'}$ to $\underline{\mathcal{S}}_0$ is given 
by $\bar f^{\gamma'}(\zeta,0)$. We develope its Taylor series
with respect to the powers of $\zeta$ as follows
\def\theequation{4.3.16}\begin{equation}
\bar f^{\gamma'}(\zeta,0)=
\sum_{\beta\in\N^m}\, 
\bar f_{\gamma',\beta}\, 
{\zeta^\beta \over \beta!},
\end{equation}
where the $\bar f_{\gamma',\beta}$ are constants in $\C$. Since the
coordinates are normal, with the usual vector fields
$\underline{\mathcal{L}}_k$ given by
\def\theequation{4.3.17}\begin{equation}
\underline{\mathcal{L}}_k=
{\partial \over \partial \zeta_k}+
\sum_{j=1}^d\, 
{\partial \Theta_j \over \partial \zeta_k}(\zeta,t)\, 
{\partial \over \partial \xi_j},
\end{equation}
taking into account that $\Theta_j(\zeta,0,w)\equiv \Theta_j(0,z,w)\equiv 0$
({\it cf.} Theorem~2.1.32),
we immediately see that the constants $\bar f_{\gamma',\beta}$ are
simply given by
\def\theequation{4.3.18}\begin{equation}
\bar f_{\gamma',\beta} =\underline{\mathcal{L}}^\beta
(\bar f^{\gamma'})(0,0,0,0).
\end{equation}
We may now rewrite the expression~\thetag{4.3.15} as follows
\def\theequation{4.3.19}\begin{equation}
\sum_{\gamma'\in\N^{m'}}\, 
\left(
\bar f_{\gamma',\beta}-\Lambda_\beta^1 \, 
\bar f_{\gamma',\beta_1}-\cdots-
\Lambda_\beta^{\kappa'}\, 
\bar f_{\gamma', \beta_{\kappa'}}
\right)\, 
\left[
\nabla_{t'}\, 
\Theta_{j',\gamma'}'(0)
\right]
=0.
\end{equation}
Let us rewrite this expression temporarily as a linear system of the
form
\def\theequation{4.3.20}\begin{equation}
\sum_{\gamma'\in\N^{m'}}\, 
\overline{F}_{\beta,\gamma'}\, 
\left[
\nabla_{t'}\, 
\Theta_{j',\gamma'}'(0)
\right]
=0,
\end{equation}
where the $\overline{F}_{\beta,\gamma'}$ are complex constant.  Now, we use
the assumption~\thetag{4.3.7}. Since the rank of the $n'\times \infty$
matrix $(\nabla_{t'}\, \Theta_{j', \gamma'}'(0))_{\gamma'\in \N^{m'}}$
is equal to $n_{j'}'$ it follows that after making some linear
combinations between the lines of the system~\thetag{4.3.20} that there
exist $n_{j'}'$ distinct multiindices $\gamma_1', \dots,
\gamma_{n_{j'}'}'\in\N^{m'}$ such that for $i'=1,\dots,n_{j'}'$, we
can solve
\def\theequation{4.3.21}\begin{equation}
\overline{F}_{\beta,\gamma_{i'}'}=
\sum_{\gamma'\neq \gamma_1',\dots,
\gamma_{n_{j'}'}'}\,
A_{i',\gamma'}' \, 
\overline{F}_{\beta,\gamma'},
\end{equation}
where the $A_{i',\gamma'}'$ are complex constants.

We now replace the $\overline{F}_{\beta,\gamma'}$ by their values, which yields
for $\beta\neq \beta_1,\dots,\beta_{\kappa'}$ and for
$i'=1,\dots,n_{j'}'$ the following equalities
\def\theequation{4.3.22}\begin{equation}
\left\{
\aligned
{}
&
\bar f_{\gamma_{i'}',\beta}-\Lambda_\beta^1\, 
\bar f_{\gamma_{i'}',\beta_1}-\cdots-\Lambda_\beta^{\kappa'}\, 
\bar f_{\gamma_{i'}',\beta_{\kappa'}}=
\\
& \
=\sum_{\gamma'\neq \gamma_1',\dots,\gamma_{n_{j'}'}'}\, 
A_{i',\gamma'}'\left[
\bar f_{\gamma',\beta}-
\Lambda_\beta^1\, \bar f_{\gamma',\beta_1} -\cdots -
\Lambda_\beta^{\kappa'} \bar f_{\gamma',\beta_{\kappa'}}
\right].
\endaligned\right.
\end{equation}
We can now mix linearly the left and the right hand sides to obtain
for $i'=1,\dots,n_{j'}'$ and for $\beta\neq \beta_1,\dots,\beta_{\kappa'}$
\def\theequation{4.3.23}\begin{equation}
\left\{
\aligned
\bar f_{\gamma_{i'}',\beta}-
\sum_{\gamma'\neq \gamma_1',\dots,\gamma_{n_{j'}'}'}\, 
A_{i',\gamma'}'\, 
\bar f_{\gamma',\beta}
=
& \
\Lambda_\beta^1\left[
\bar f_{\gamma_{i'}',\beta_1}-\sum_{\gamma'\neq \gamma_1',\dots,
\gamma_{n_{j'}'}'}\, 
A_{i',\gamma'}' \, \bar f_{\gamma',\beta_1}
\right]\\
& \
+\cdots+\\
& \
\Lambda_\beta^{\kappa'}\left[
\bar f_{\gamma_{i'}',\beta_{\kappa'}}-\sum_{\gamma'\neq \gamma_1',\dots,
\gamma_{n_{j'}'}'}\, 
A_{i',\gamma'}' \, \bar f_{\gamma',\beta_{\kappa'}}
\right].
\endaligned\right.
\end{equation}
On the other hand, by the definition~\thetag{4.3.16} of the constants
$\bar f_{\gamma',\beta}$, we can develope the following expression 
in power series of $\zeta$
\def\theequation{4.3.24}\begin{equation}
\left\{
\aligned
\bar f^{\gamma_i'}(\zeta,0)
-
& \
\sum_{\gamma'\neq \gamma_1',\dots,\gamma_{n_{j'}'}'}\, 
A_{i',\gamma'}'\, \bar f^{\gamma'}(\zeta,0)=\\
& \
\sum_{\beta\neq \beta_1,\dots,\beta_{\kappa'}}\left(
\bar f_{\gamma_{i'}',\beta}-\sum_{\gamma'\neq \gamma_1',\dots,
\gamma_{n_{j'}'}'}\, A_{i',\gamma'}'\, \bar f_{\gamma',\beta}
\right)\, {\zeta^\beta\over \beta!}+
\\
& \ 
+\left( \bar f_{\gamma_{i'}',\beta_1}-
\sum_{\gamma'\neq \gamma_1',\dots,
\gamma_{n_{j'}'}'}\, 
A_{i',\gamma'}' \, \bar f_{\gamma',
\beta_1}\right)\, 
{\zeta^{\beta_1}\over \beta_1!}+\cdots+\\
& \
+\left( \bar f_{\gamma_{i'}',\beta_{
\kappa'}}-\sum_{\gamma'\neq \gamma_1',\dots,
\gamma_{n_{j'}'}'}\, 
A_{i',\gamma'}' \, \bar f_{\gamma',\beta_{\kappa'}}\right)\, 
{\zeta^{\beta_{\kappa'}}\over \beta_{\kappa'}!}.
\\
\endaligned\right.
\end{equation}
Now, we introduce the following new complex constants for $i'=1,\dots,n_{j'}'$
\def\theequation{4.3.25}\begin{equation}
\left\{
\aligned
\Pi_{i',1}:=
& \
\bar f_{\gamma_{i'}',\beta_1}-\sum_{\gamma'\neq
\gamma_1',\dots,
\gamma_{n_{j'}'}'}\, 
A_{i',\gamma'}'\, \bar f_{\gamma',\beta_1},\\
& \
\cdots \cdots \cdots\\
\Pi_{i',\kappa'}:=
& \
\bar f_{\gamma_{i'}',\beta_{\kappa'}}-\sum_{\gamma'\neq
\gamma_1',\dots,
\gamma_{n_{j'}'}'}\, 
A_{i',\gamma'}'\, \bar f_{\gamma',\beta_{\kappa'}},\\
\endaligned\right.
\end{equation}
and we use the preceding relations~\thetag{4.3.23} to rewrite~\thetag{4.3.24}
more simply as follows, for $i'=1,\dots,n_{j'}'$:
\def\theequation{4.3.26}\begin{equation}
\left\{
\aligned
\bar f^{\gamma_{i'}'}
& \
(\zeta,
0)-
\sum_{\gamma'\neq
\gamma_1',\dots,\gamma_{n_{j'}'}'}\, 
A_{i',\gamma'}' \, 
\bar f^{\gamma'}(\zeta,0)=\\
& \
=\sum_{\beta\neq \beta_1,\dots,\beta_{\kappa'}}\left[
\Lambda_\beta^1\, \Pi_{i',1}+\cdots+
\Lambda_\beta^{\kappa'} \, \Pi_{i',\kappa'}
\right]\, 
{\zeta^\beta \over \beta!}+
\Pi_{i',1}\, {\zeta^{\beta_1}\over \beta_1!}+
\cdots+
\Pi_{i',\kappa'}\, {\zeta^{\beta_{\kappa'}}\over
\beta_{\kappa'}!}\\
& \ 
=\Pi_{i',1}\left(
{\zeta^{\beta_1}\over \beta_1!}+
\sum_{\beta\neq \beta_1,\dots,\beta_{\kappa'}}\, 
\Lambda_\beta^1\, {\zeta^\beta \over \beta!}\right)+\cdots+
\Pi_{i',\kappa'}\left(
{\zeta^{\beta_{\kappa'}}\over \beta_{\kappa'}!}+
\sum_{\beta\neq \beta_1,\dots,\beta_{\kappa'}}\, 
\Lambda_\beta^{\kappa'}\, 
{\zeta^\beta \over \beta!}\right)\\
& \
=:
\Pi_{i',1}\, G_1(\zeta)+\cdots+
\Pi_{i',\kappa'}\, G_{\kappa'}(\zeta),
\endaligned\right.
\end{equation}
where the power series $G_1(\zeta),\dots,G_{\kappa'}(\zeta)$ are
defined by the last equality. Now, we use the assumption 
$\kappa'\leq n_{j'}'-1$. Since there are striclty less than
$n_{j'}'$ functions $G$ in~\thetag{4.3.26}, it follows that there
exist complex constants $\mu_1,\dots,\mu_{n_{j'}'}$ not
all zero such that
\def\theequation{4.3.27}\begin{equation}
\left\{
\aligned
0\equiv 
\mu_1
& \
\left(
\bar f^{\gamma_1'}(\zeta,0)
-\sum_{\gamma'\neq \gamma_1',\dots,
\gamma_{n_{j'}'}'} \, A_{1,\gamma'}\, 
\bar f^{\gamma'}(\zeta,0)\right)+\cdots+\\
& \
+\mu_{n_{j'}'} \, \left(
\bar f^{\gamma_{n_{j'}'}'}(\zeta,0)
-\sum_{\gamma'\neq \gamma_1',\dots,
\gamma_{n_{j'}'}'} \, A_{n_{j'}',\gamma'}\, 
\bar f^{\gamma'}(\zeta,0)\right),
\endaligned\right.
\end{equation}
and this last equality clearly contradicts the assumption that
$h$ is CR-transversal at $p_0$. 

The proofs of Lemma~4.3.6 and of parts {\bf (1)} and {\bf (2)} of
Theorem~4.3.1 are complete.
\endproof

We now prove part {\bf (3)} of Theorem~4.3.1.  We also proceed by
contradiction. Suppose that the ideal
\def\theequation{4.3.28}\begin{equation}
\left<
R_{j',\beta}'(0,0,
J^{\vert \beta \vert} \bar h(0) :\, t')
\right>_{1\leq j'\leq d', \, 
\beta \in\N^m}
\end{equation}
is of infinite codimension in $\mathcal{A}_\C\{t'\}$
or in $\C\{ t'\}$.  By the Nullstellensatz,
it follows that there exists a nonzero algebraic or analytic piece of
curve passing through the origin
\def\theequation{4.3.29}\begin{equation}
\C\ni s'\mapsto t'(s')\in \C^{n'},
\end{equation}
namely $t'(s')\in\mathcal{A}_\C\{s'\}^{n'}$ or $t'(s')\in\C\{s'\}^{n'}$, such that
\def\theequation{4.3.30}\begin{equation}
R_{j',\beta}'(0,0,J^{\vert \beta \vert} \bar h(0) :\, t'(s'))\equiv 0,
\end{equation}
for all $j'=1,\dots, d'$ and all $\beta\in\N^m$. Replacing the
$R_{j',\beta}'$ by their definition~\thetag{4.2.10}, this means that
\def\theequation{4.3.31}\begin{equation}
\underline{\mathcal{L}}^\beta \, \bar g_{j'}(0)-\sum_{\gamma'\in\N^{m'}}\, 
\underline{\mathcal{L}}^\beta (\bar f^{\gamma'})(0) \, 
\Theta_{j',\gamma'}'(t'(s'))\equiv 0,
\end{equation}
for $j'=1,\dots,d'$.
But since the coordinates $(z,w)$ and $(z',w')$ are
normal, we claim that we have
\def\theequation{4.3.32}\begin{equation}
\underline{\mathcal{L}}^\beta \, \bar g_{j'}(0)=0
\end{equation}
for all $\beta\in\N^m$ and all $j'=1,\dots,d'$.
Indeed, setting $t=0$ in the fundamental power series identities
\def\theequation{4.3.33}\begin{equation}
\bar g_{j'}(\zeta,\Theta(\zeta,t))\equiv
\Theta_{j'}'(\bar f(\zeta,\Theta(\zeta,t)),
h(t)), 
\end{equation}
we get thanks to the normality of coordinates
\def\theequation{4.3.34}\begin{equation}
\bar g_{j'}(\zeta,0)\equiv \Theta_{j'}'(\bar f(\zeta,0),0)\equiv 0,
\end{equation}
hence 
\def\theequation{4.3.35}\begin{equation}
\underline{\mathcal{L}}^\beta \, \bar g_{j'}(0)=
\partial_\zeta^\beta \bar g_{j'}(\zeta,0)\vert_{\zeta=0}=0,
\end{equation}
as claimed. Developing 
\def\theequation{4.3.36}\begin{equation}
\bar f^{\gamma'}(\zeta,0)=\sum_{\beta\in\N^m}\, 
\bar f_{\gamma',\beta}\, 
{\zeta^\beta\over \beta!},
\end{equation}
where (again thanks to normal coordinates)
\def\theequation{4.3.37}\begin{equation}
\bar f_{\gamma',\beta}=\underline{\mathcal{L}}^\beta 
\bar f^{\gamma'}(0), 
\end{equation}
we can now rewrite the expression~\thetag{4.3.31} in under the simpler form
\def\theequation{4.3.38}\begin{equation}
\sum_{\gamma'\in\N^{m'}}\, 
\bar f_{\gamma',\beta}\, 
\Theta_{j',\gamma'}'(t'(s'))\equiv 0.
\end{equation}
Because $M'$ is essentially finite at $p_0'$, there exists at smallest
one integer $j_0'$ with $1\leq j_0'\leq d'$ such that not all
$\Theta_{j_0',\gamma'}'(t'(s'))$ vanish identically, for
$\gamma'\in\N^{m'}$. Hence there exists $s_0'\in\C$ arbitrarily close
to the origin such that the infinite family of complex constants
$(\Theta_{j_0',\gamma'}'(t'(s_0')))_{\gamma'\in\N^{m'}}$ are not
all zero. We put 
\def\theequation{4.3.39}\begin{equation}
\theta_{\gamma'}':=\Theta_{j_0',\gamma'}'(t'(s_0'))\in \C.
\end{equation}
Setting $s':=s_0'$ in~\thetag{4.3.38}, we therefore get
for all  $\beta\in\N^m$ a relation of the form
\def\theequation{4.3.40}\begin{equation}
\sum_{\gamma'\in\N^{m'}}\, 
\bar f_{\gamma',\beta}\, 
\theta_{\gamma'}'=0,
\end{equation}
which yields after integrating with respect to $\zeta$ the
formal identity
\def\theequation{4.3.41}\begin{equation}
\sum_{\gamma'\in\N^{m'}}\,  \theta_{\gamma'}' \,
\bar f^{\gamma'}(\zeta,0) \equiv 0.
\end{equation}
But this clearly contradicts the CR-transversality of $h$ at $p_0$.

The proof of part {\bf (3)} of Theorem~4.3.1 is complete.

\smallskip

We now prove part {\bf (4)} of Theorem~4.3.1. The proof
has some similarities with the proof of 
part {\bf (2)}, but is a little bit more technical.
We use the notations introduced in the beginning of the proof
of Theorem~4.3.1. 

It follows directly 
from Section~3.2 ({\it cf.} especially~\thetag{3.2.47}
and Lemma~3.2.49) that in normal coordinates, and using 
the gradient notation $\nabla_{t'}$, then $M'$ is Segre
nondegenerate if and only if 
\def\theequation{4.3.42}\begin{equation}
{\rm genMrk}\, \left(
\nabla_{t'}\, \Theta_{j',\gamma'}'(z',0)
\right)_{1\leq j'\leq d',\, 
\gamma'\in\N^{m'}}=n'.
\end{equation}
Notice that here we let the integer $j'$ vary from $1$ to 
$d'$, but in the following lemma, we shall fix $j'$.

\def\thelemma{4.3.43}\begin{lemma}
Fix $j'$ with $1\leq j'\leq d'$. Let $n_{j'}'$ be the
integer defined by
\def\theequation{4.3.44}\begin{equation}
n_{j'}':=
{\rm genMrk}\, 
\left(
\nabla_{t'}\, \Theta_{j',\gamma'}'(z',0)
\right)_{\gamma'\in\N^{m'}}.
\end{equation}
Assume that $h$ is CR-transversal at $p_0$. Then we also 
have
\def\theequation{4.3.45}\begin{equation}
{\rm genMrk}\, \left(
\nabla_{t'}\, R_{j',\beta}'(z,0,0,0,J^{\vert \beta \vert}\bar h(0) 
:\, h(z,0))
\right)_{\beta\in\N^m}=n_{j'}'.
\end{equation}
\end{lemma}

We first show that this technical lemma implies part {\bf (4)} of
Theorem~4.3.1. Indeed, by the definitions~\thetag{4.2.5}
and~\thetag{4.2.10} of $R_{j',\beta}'$, after specifiying on the Segre
chain $\mathcal{S}_0=\{(z,0,0,0)\}$, we have
\def\theequation{4.3.46}\begin{equation}
\left\{
\aligned
\nabla_{t'} \, R_{j',\beta}'
& \
(z,0,0,0,J^{\vert \beta \vert}\bar
h(0) :\, h(z,0))=\\
& \
=-\sum_{\gamma'\in\N^{m'}}\, 
\underline{\mathcal{L}}^\beta (\bar f^{\gamma'})(z,0,0,0)\, 
[\nabla_{t'} \, \Theta_{j',\gamma'}'(h(z,0))].
\endaligned\right.
\end{equation}
Notice that we write here the first arguments of $R_{j', \beta}'$ and
the arguments of the differentiated expression $\underline{ \mathcal{
L}}^\beta ( \bar f^{ \gamma'})$ as $(z,w, \zeta, \xi)$, and not as
$(t,\tau)$. Here, in the arguments $(z,0,0,0)$ of this expression
$\underline{ \mathcal{ L}}^\beta ( \bar f^{ \gamma'})$, the term $z$
comes from the coefficients of the vector fields
$\underline{\mathcal{L}}_k$, but $(\zeta,\xi)=(0,0)$.  It follows that
the sum in~\thetag{4.3.46} is in fact truncated with $\vert \gamma'
\vert \leq \vert \beta \vert$, since $\bar f^{\gamma'}$ vanishes to
order $\vert \gamma' \vert$ at the origin.  Thus, the columns of the
matrix $(\nabla_{t'} R_{j',\beta}')_{\beta\in\N^m}$ on the left hand
side of~\thetag{4.3.46} are obtained as a linear combination (with
coefficients being certain formal power series in $z$) of the columns
of the matrix $(\nabla_{t'}\Theta_{j',\gamma'}')_{\gamma'\in \N^{m'}}$
on the right hand side.  Thanks to Lemma~4.3.43 (to be proved below),
we deduce that the formal linear space spanned by the columns of the
matrix $(\nabla_{t'} R_{j',\beta}')_{\beta\in\N^m}$ on the left
coincides with the formal linear space spanned by the columns of
matrix $(\nabla_{t'}\Theta_{j',\gamma'}')_{\gamma'\in \N^{m'}}$ on the
right. Finally, thanks to the Segre nondegeneracy
assumption~\thetag{4.3.42}, we deduce
\def\theequation{4.3.47}\begin{equation}
{\rm genMrk}\, \left(
\nabla_{t'} R_{j',\beta}'
(z,0,0,0,J^{\vert \beta \vert}\bar
h(0) :\, h(z,0))
\right)_{1\leq j'\leq d', \, 
\beta\in\N^m}=n',
\end{equation}
which shows that $h$ is $\ell_1$-Segre nondegenerate, for some
$\ell_1\geq \ell_0'$. It remains to establish the technical Lemma~4.3.43.

\proof
So we fix $j'$. Again, we proceed by contradiction.  Assume that the
generic rank $\kappa'$ of the $n'\times \infty$ matrix
\def\theequation{4.3.48}\begin{equation}
\left(
\nabla_{t'}\, R_{j',\beta}'(z,0,0,0,J^{\vert \beta \vert}\bar h(0) :\,
h(z,0))
\right)_{\beta\in\N^m}
\end{equation}
is strictly smaller than $n_{j'}'$, namely $\kappa'\leq n_{j'}'-1$. We
choose $\kappa'$ distinct multiindices
$\beta_1,\dots,\beta_{\kappa'}\in\N^m$ such that the generic rank of
the $n'\times \kappa'$ extracted matrix
\def\theequation{4.3.49}\begin{equation}
\left(
\nabla_{t'}\, R_{j',\beta_{i'}}'(z,0,0,0,
J^{\vert \beta \vert}\bar h(0) :\,
h(z,0))
\right)_{1\leq i'\leq \kappa'}
\end{equation}
is equal to $\kappa'$. We let $\Lambda(z) \in \C\dl z \dr$ denote a
not identically zero $\kappa'\times \kappa'$ minor of this matrix. It
then follows from Cramer's rule and from the rank
assumption that for every multiindex $\beta\in\N^m$ different from
$\beta_1, \dots, \beta_{ \kappa'}$, there exist formal power series
$\Lambda_\beta^1 (z),\dots, \Lambda_\beta^{ \kappa'}(z) \in\C\dl z\dr$
such that we can write
\def\theequation{4.3.50}\begin{equation}
\left\{
\aligned
\Lambda(z)\, 
\nabla_{t'}\, R_{j',\beta}'
& \
(z,0,0,0,J^{\vert \beta 
\vert} \bar h(0) :\,
h(z,0))\equiv \\
& \
\equiv\Lambda_\beta^1(z) \, 
\nabla_{t'} \, R_{j',\beta_1}'(z,0,0,0,
J^{\vert \beta \vert} \bar h(0) :\,
h(z,0))+\cdots+\\
& \ \ \ \
+\Lambda_\beta^{\kappa'}(z) \, 
\nabla_{t'} \, R_{j',\beta_{\kappa'}}'(z,0,0,0,
J^{\vert \beta \vert} \bar h(0) :\,
h(z,0)).
\endaligned\right.
\end{equation}
Replacing the $R_{j',\beta}'$ by their values given by~\thetag{4.3.46}, 
we get
\def\theequation{4.3.51}\begin{equation}
\left\{
\aligned
\Lambda(z)\, 
\sum_{\gamma'\in\N^{m'}}
& \ 
\underline{\mathcal{L}}^\beta (\bar f^{\gamma'})
(z,0,0,0)\,
[\nabla_{t'}\, \Theta_{j',\gamma'}'(h(z,0))]
\equiv \\
\equiv
& \
\sum_{\gamma'\in\N^{m'}}\left(
\Lambda_\beta^1(z)\, 
\underline{\mathcal{L}}^{\beta_1}(\bar f^{\gamma'})(z,0,0,0)+
\cdots+ \right. \\
& \
\left.
+\Lambda_\beta^{\kappa'}(z)\, 
\underline{\mathcal{L}}^{\beta_{
\kappa'}} (\bar f^{\gamma'})(z,0,0,0)
\right)
\left[
\nabla_{t'} \, \Theta_{j',\gamma'}'(h(z,0))
\right].
\endaligned\right.
\end{equation}
We remind that in normal coordinates, we have
$g(z,0)\equiv 0$, hence
\def\theequation{4.3.52}\begin{equation}
h(z,0)\equiv (f(z,0),0).
\end{equation}
Before going further, we need the following elementary observation.

\def\thelemma{4.3.53}\begin{lemma}
Assume that $h$ is CR-transversal at $p_0$ and fix $j'$ as in 
Lemma~4.3.43. Then 
\def\theequation{4.3.54}\begin{equation}
{\rm genMrk} \left(
\nabla_{t'} \,\Theta_{j',\gamma'}'(f(z,0),0)
\right)_{\gamma'\in\N^{m'}}=n_{j'}'.
\end{equation}
\end{lemma}

\proof
Suppose on the contrary that this generic matrix rank is 
equal to an integer $\kappa'\leq n_{j'}'-1$. By the definition~\thetag{4.3.45}
of $n_{j'}'$, there exists a $n_{j'}'\times n_{j'}'$ minor
$d'(z')$ of the $n'\times \infty$ matrix 
$(\nabla_{t'}\, \Theta_{j',\gamma'}'(z',0))_{\gamma'\in\N^{m'}}$
which does not vanish identically as a power series of $z'$.
We deduce that $d'(f(z,0))\equiv 0$, contradicting the 
CR-transversality assumption on $h$, which 
completes the proof.
\endproof

After making a linear combination,
we can rewrite temporarily the relations~\thetag{4.3.51} under the
form
\def\theequation{4.3.55}\begin{equation}
\sum_{\gamma'\in\N^{m'}}\, 
\overline{F}_{\beta,\gamma'}(z)\, 
\left[
\nabla_{t'}\,\Theta_{j',\gamma'}'(h(z,0))
\right]\equiv 0.
\end{equation}
By Lemma~4.3.53, there exist 
multiindices $\gamma_1',\dots,\gamma_{n_{j'}'}'$ such that
a $n_{j'}'\times n_{j'}'$ minor $A(z)$
of the $n'\times n_{j'}'$ matrix 
\def\theequation{4.3.56}\begin{equation}
\left(
\nabla_{t'} \, \Theta_{j',\gamma_{i'}'}'(h(z,0))
\right)_{1\leq i'\leq n_{j'}'}
\end{equation}
does not vanish identically as a power series of $z$.  We deduce that
after making some linear combinations between the lines of the
system~\thetag{4.3.55} that for every $i'=1,\dots,n_{j'}'$, there
exist formal power series $A_{i',\gamma'}'(z)$ such that we can solve
\def\theequation{4.3.57}\begin{equation}
A(z)\, 
\overline{F}_{\beta,\gamma_{i'}'}(z)\equiv
\sum_{\gamma'\neq \gamma_1',\dots,\gamma_{n_{j'}'}'}\, 
A_{i',\gamma'}'(z)\, 
\overline{F}_{\beta,\gamma'}(z).
\end{equation}
We now replace the $\overline{F}_{\beta,\gamma'}$ by their
values given by~\thetag{4.3.51} and~\thetag{4.3.55}. We obtain
the following formal
equalities, valuable for $i'=1,\dots,n_{j'}'$ and for
$\beta\neq \beta_1,\dots,\beta_{\kappa'}$:
\def\theequation{4.3.58}\begin{equation}
\left\{
\aligned
A(z)\, \Lambda(z) \,
& \
\underline{\mathcal{L}}^\beta 
(\bar f^{\gamma_{i'}'})(z,0,0,0)-
A(z) \,
\Lambda_\beta^1(z)\, 
\underline{\mathcal{L}}^{\beta_1}(
\bar f^{\gamma_{i'}'})(z,0,0,0)-\cdots-\\
& \
-A(z) \,
\Lambda_\beta^{\kappa'}(z)\, 
\underline{\mathcal{L}}^{\beta_{
\kappa'}}(\bar f^{\gamma_{i'}'})(z,0,0,0)
\equiv\\
\equiv 
& \
\sum_{\gamma'\neq \gamma_1',\dots,\gamma_{n_{j'}'}'}\, 
\left(
A_{i',\gamma'}'(z) \, \Lambda(z)\, 
\underline{\mathcal{L}}^{
\beta}(\bar f^{\gamma'})(z,0,0,0)  - \right. \\
& \
-A_{i',\gamma'}'(z) \, \Lambda_\beta^1(z)\, 
\underline{\mathcal{L}}^{
\beta_1}(\bar f^{\gamma'})(z,0,0,0) - \cdots - \\
& \
\left.
-A_{i',\gamma'}'(z) \, \Lambda_\beta^1(z)\, 
\underline{\mathcal{L}}^{
\beta_{\kappa'}}(\bar f^{\gamma'})(z,0,0,0)
\right).
\endaligned\right.
\end{equation}
After making some linear combinations, we can rewrite this
identity as
\def\theequation{4.3.59}\begin{equation}
\left\{
\aligned
A(z)\, \Lambda(z) \,
& \
\underline{\mathcal{L}}^\beta (\bar f^{\gamma_{i'}'})(z,0,0,0)-
\sum_{\gamma'\neq \gamma_1',\dots,\gamma_{n_{j'}'}'}\,
A_{i',\gamma'}'(z)\, \Lambda(z)\, 
\underline{\mathcal{L}}^\beta (\bar f^{\gamma'})(z,0,0,0)\equiv\\
\equiv 
& \
\Lambda_\beta^1(z)\left(
A(z)\, \underline{\mathcal{L}}^{\beta_1}
(\bar f^{\gamma_{i'}'})(z,0,0,0)- \right.
\\
& \ \ \ \ \ \ \ \ \ 
-\left. \sum_{\gamma'\neq
\gamma_1',\dots,\gamma_{n_{j'}'}'}\, 
A_{i',\gamma'}'(z)\, 
\underline{\mathcal{L}}^{
\beta_1} (\bar f^{\gamma'})(z,0,0,0)
\right)+\cdots+\\
& \
+ \Lambda_\beta^{\kappa'}(z)\left(
A(z)\, \underline{\mathcal{
L}}^{\beta_{\kappa'}}
(\bar f^{\gamma_{i'}'})
(z,0,0,0)- \right.
\\
& \ \ \ \ \ \ \ \ \ 
-\left. \sum_{\gamma'\neq
\gamma_1',\dots,\gamma_{n_{j'}'}'}\, 
A_{i',\gamma'}'(z)\, 
\underline{\mathcal{L}}^{\beta_{\kappa'}} 
(\bar f^{\gamma'})(z,0,0,0)
\right).
\endaligned\right.
\end{equation}
Next, we introduce the following new notation
\def\theequation{4.3.60}\begin{equation}
\bar f_{\gamma',\beta}(z):=
\underline{\mathcal{L}}^\beta (\bar f^{\gamma'})(z,0,0,0).
\end{equation}
Of course, we have
\def\theequation{4.3.61}\begin{equation}
\bar f^{\gamma'}(\zeta,0)\equiv \sum_{\beta\in\N^m}\, \bar
f_{\gamma',\beta}(0)\, {\zeta^\beta \over \beta!}
\end{equation}
and 
\def\theequation{4.3.62}\begin{equation}
\left\{
\aligned
\bar f^{\gamma'}(\zeta,\Theta(\zeta,z,0))
& \
\equiv \sum_{\beta\in\N^m}\, {\zeta^\beta \over \beta!}\, \left[
{\partial^{\vert \beta \vert} \over \partial \zeta^\beta} \left( \bar
f^{\gamma'}(\zeta,\Theta(\zeta,z,0)) \right)\right]_{\zeta=0}\\
& \ 
\equiv 
\underline{\mathcal{L}}^\beta (\bar f^{\gamma'})(z,0,0,0)\\
& \
\equiv
\sum_{\beta\in\N^m}\, 
{\zeta^\beta\over \beta!}\, 
\bar f_{\gamma',\beta}(z).
\endaligned\right.
\end{equation}
With this notation, we can rewrite~\thetag{4.3.59} as follows, 
where $i'=1,\dots,n_{j'}'$ and
$\beta\neq \beta_1,\dots,\beta_{\kappa'}$:
\def\theequation{4.3.63}\begin{equation}
\left\{
\aligned
A(z)\, \Lambda(z)\, 
& \
\bar f_{\gamma_{i'}'}(z)-
\sum_{\gamma'\neq \gamma_1',\dots,\gamma_{n_{j'}'}'}\, 
\Lambda(z) \, A_{i',\gamma'}'(z)\, 
\bar f_{\gamma',\beta}(z)\equiv\\
& \
\equiv
\Lambda_\beta^1(z)\left(
A(z)\, \bar f_{\gamma_{i'}',\beta}(z)-
\sum_{\gamma'\neq \gamma_1',
\dots,\gamma_{n_{j'}'}'}\, 
A_{i',\gamma'}'(z)\, 
\bar f_{\gamma',\beta}(z)
\right)+\cdots+\\
& \ \ \
+\Lambda_\beta^1(z)\left(
A(z)\, \bar f_{\gamma_{i'}',\beta}(z)-
\sum_{\gamma'\neq \gamma_1',
\dots,\gamma_{n_{j'}'}'}\, 
A_{i',\gamma'}'(z)\, 
\bar f_{\gamma',\beta}(z)
\right).
\endaligned\right.
\end{equation}
On the other hand, taking the definition~\thetag{4.3.60} of 
$\bar f_{\gamma',\beta}(z)$ and~\thetag{4.3.62} into account, we have the
power series development
\def\theequation{4.3.64}\begin{equation}
\left\{
\aligned
\sum_{\beta\in\N^m}\, 
& \
{\zeta^\beta \over \beta!}\,
\left(
A(z)\, \Lambda(z)\, 
\bar f_{\gamma_{i'}',\beta}(z)-\sum_{\gamma'\neq
\gamma_1',\dots,\gamma_{n_{j'}'}'}\, 
\Lambda(z)\, A_{i',\gamma'}'(z)\, 
\bar f^{\gamma',\beta}(z)\right)\equiv\\
& \
\equiv
A(z)\, \Lambda(z)\, 
\bar f^{\gamma_{i'}'}(\zeta,\Theta(\zeta,z,0))-\sum_{\gamma'\neq
\gamma_1',\dots,\gamma_{n_{j'}'}'}\, 
\Lambda(z)\, 
A_{i',\gamma'}'(z)\, 
\bar f^{\gamma'}(\zeta,\Theta(\zeta,z,0)).
\endaligned\right.
\end{equation}
In this identity, we decompose the sum $\sum_{\beta\in\N^m}$ as the
sum $\sum_{\beta\neq \beta_1,\dots,\beta_{\kappa'}}$ plus the
$\kappa'$ remaining terms corresponding to
$\beta=\beta_1,\dots,\beta_{\kappa'}$ and we substitute the
expressions obtained just previously in~\thetag{4.3.63}, which yields
\def\theequation{4.3.65}\begin{equation}
\left\{
\aligned
\sum_{\beta\in\N^m}\, 
{\zeta^\beta \over \beta!}\,
& \ 
\left(
A(z)\, \Lambda(z)\, 
\bar f_{\gamma_{i'}',\beta}(z)-\sum_{\gamma'\neq
\gamma_1',\dots,\gamma_{n_{j'}'}'}\, 
\Lambda(z)\, A_{i',\gamma'}'(z)\, 
\bar f_{\gamma',\beta}(z)\right)\equiv\\
\equiv
\sum_{\beta\neq\beta_1,\dots,\beta_{\kappa'}}
& \
{\zeta^\beta \over \beta!}\left(
A(z)\, \Lambda(z)\, 
\bar f_{\gamma_{i'}',\beta}(z)
-\sum_{\gamma'\neq \gamma_1',\dots,\gamma_{n_{j'}'}'}\, 
\Lambda(z)\, A_{i',\gamma'}'(z)\, 
\bar f_{\gamma',\beta}(z)
\right)+\\
& \
+ {\zeta^{\beta_1}\over \beta_1!}\left(A(z)\, \Lambda(z)\, 
\bar f_{\gamma_{i'}',\beta_1}(z)
-\sum_{\gamma'\neq \gamma_1',\dots,\gamma_{n_{j'}'}'}\, 
\Lambda(z)\, A_{i',\gamma'}'(z)\, 
\bar f_{\gamma',\beta_1}(z)
\right)+\\
& \
+\cdots+\\
& \
+ {\zeta^{\beta_{\kappa'}}\over \beta_{
\kappa'}!}\left(A(z)\, \Lambda(z)\, 
\bar f_{\gamma_{i'}',\beta_{\kappa'}}(z)
-\sum_{\gamma'\neq \gamma_1',\dots,\gamma_{n_{j'}'}'}\, 
\Lambda(z)\, A_{i',\gamma'}'(z)\, 
\bar f_{\gamma',\beta_{\kappa'}}(z)
\right).
\endaligned\right.
\end{equation}
Now, we make some linear combinations, which yields
the following representation of the right hand side of~\thetag{4.3.65}
\def\theequation{4.3.66}\begin{equation}
\left\{
\aligned
A(z)\, \Lambda(z)\,
& \ 
\bar f^{\gamma_{i'}'}(\zeta,\Theta(\zeta,z,0))-\sum_{\gamma'\neq
\gamma_1',\dots,\gamma_{n_{j'}'}'}\, 
\Lambda(z)\, 
A_{i',\gamma'}'(z)\, 
\bar f^{\gamma'}(\zeta,\Theta(\zeta,z,0))
\equiv\\
\equiv
\sum_{\beta\neq \beta_1,\dots,\beta_{\kappa'}}\, 
{\zeta^\beta\over \beta!}
& \
\left[
\Lambda_\beta^1(z)\left(
A(z)\, \bar f_{\gamma_{i'}',\beta_1}(z)-
\sum_{\gamma'\neq 
\gamma_1',\dots,\gamma_{n_{j'}'}'}\, 
A_{i',\gamma'}'(z)\,
\bar f_{\gamma',\beta_1}(z)
\right)+ \right.\\
& \
+\cdots+\\
& \
\left.
+
\Lambda_\beta^{\kappa'}(z)\left(
A(z)\, \bar f_{\gamma_{i'}',\beta_{\kappa'}}(z)-
\sum_{\gamma'\neq \gamma_1',\dots,\gamma_{n_{j'}'}'}\, 
A_{i',\gamma'}'(z)\,
\bar f_{\gamma',\beta_{\kappa'}}(z)
\right)
\right]+\\
& \
+{\zeta^{\beta_1}\over \beta_1!}\left[
A(z)\, \bar f_{\gamma_{i'}',\beta_1}(z)-
\sum_{\gamma'\neq 
\gamma_1',\dots,\gamma_{n_{j'}'}'}\, 
A_{i',\gamma'}'(z)\,
\bar f_{\gamma',\beta_1}(z)
\right]+\\
& \
+\cdots+\\
& \
+{\zeta^{\beta_{\kappa'}}\over \beta_{\kappa'}!}\left[
A(z)\, \bar f_{\gamma_{i'}',\beta_{\kappa'}}(z)-
\sum_{\gamma'\neq 
\gamma_1',\dots,\gamma_{n_{j'}'}'}\, 
A_{i',\gamma'}'(z)\,
\bar f_{\gamma',\beta_{\kappa'}}(z)
\right].
\endaligned\right.
\end{equation}
If we now set for $i'=1,\dots,n_{j'}'$:
\def\theequation{4.3.67}\begin{equation}
\left\{
\aligned
\Pi_{i',1}(z):=
& \
A(z)\, \bar f_{\gamma_{i'}',\beta_1}(z)-
\sum_{\gamma'\neq 
\gamma_1',\dots,\gamma_{n_{j'}'}'}\, 
A_{i',\gamma'}'(z)\,
\bar f_{\gamma',\beta_1}(z),\\
& \
\cdots\cdots\cdots \\
\Pi_{i',\kappa'}(z):=
& \
A(z)\, \bar f_{\gamma_{i'}',\beta_{\kappa'}}(z)-
\sum_{\gamma'\neq 
\gamma_1',\dots,\gamma_{n_{j'}'}'}\, 
A_{i',\gamma'}'(z)\,
\bar f_{\gamma',\beta_{\kappa'}}(z),
\endaligned\right.
\end{equation}
then we can rewrite~\thetag{4.3.66} as follows
\def\theequation{4.3.68}\begin{equation}
\left\{
\aligned
A(z)\, \Lambda(z)\,
& \ 
\bar f^{\gamma_{i'}'}(\zeta,\Theta(\zeta,z,0))-\sum_{\gamma'\neq
\gamma_1',\dots,\gamma_{n_{j'}'}'}\, 
\Lambda(z)\, 
A_{i',\gamma'}'(z)\, 
\bar f^{\gamma'}(\zeta,\Theta(\zeta,z,0))\equiv\\
\equiv
& \
\Pi_{i',1}(z)
\left(
{\zeta^{\beta_1}\over \beta_1!}\, 
\Lambda(z)+\sum_{\beta\neq \beta_1,\dots,\beta_{\kappa'}}\, 
\Lambda_\beta^1(z)\, 
{\zeta^\beta\over \beta!}
\right)+\\
& \
+\cdots+\\
& \
+\Pi_{i',\kappa'}(z)
\left(
{\zeta^{\beta_{\kappa'}}\over \beta_{\kappa'}!}\, 
\Lambda(z)+\sum_{\beta\neq \beta_1,\dots,\beta_{\kappa'}}\, 
\Lambda_\beta^{\kappa'}(z)\, 
{\zeta^\beta\over \beta!}\right)=: \\
=:
& \
\Pi_{i',1}(z)\, G_1(z,\zeta)+\cdots+
\Pi_{i',\kappa'}(z)\, 
G_{\kappa'}(z,\zeta).
\endaligned\right.
\end{equation}
We now set $C(z):= A(z)\, \Lambda(z)$ and
for $i'=1,\dots,n_{j'}'$ and
$\gamma'\neq \gamma_1',\dots,\gamma_{n_{j'}'}'$, 
\def\theequation{4.3.69}\begin{equation}
B_{i',\gamma'}:= \Lambda(z)\, A_{i',\gamma'}'(z).
\end{equation}
In summary, we have obtained that for $i'=1,\dots,
n_{j'}'$, we can write
\def\theequation{4.3.70}\begin{equation}
\left\{
\aligned
C(z)\, 
\bar f^{\gamma_{i'}'}(\zeta,\Theta(\zeta,z,0))-
& \
\sum_{\gamma'\neq
\gamma_1',\dots,\gamma_{n_{j'}'}'}\, 
B_{i',\gamma'}(z)\, 
\bar f^{\gamma'}(\zeta,\Theta(\zeta,z,0))\equiv\\
& \
\equiv
\Pi_{i',1}(z)\, G_1(z,\zeta)+\cdots+
\Pi_{i',\kappa'}(z)\, 
G_{\kappa'}(z,\zeta).
\endaligned\right.
\end{equation}
Because $\kappa'\leq n_{j'}'-1$, hence there are less
than $n_{j'}'$ functions $G_{i'}(z,\zeta)$ in the
right hand side of~\thetag{4.3.70}, it follows
that there exist power series $\mu_1(z),\dots,
\mu_{n_{j'}'}(z)\in\C\dl z\dr$, not all zero, such that
\def\theequation{4.3.71}\begin{equation}
\left\{
\aligned
0\equiv 
& \
\mu_1(z)\, \left(
C(z)\, \bar f^{\gamma_1'}(\zeta,\Theta(\zeta,z,0))-\sum_{\gamma'\neq
\gamma_1',\dots,\gamma_{n_{j'}'}'}\, 
B_{1,\gamma'}(z)\, 
\bar f^{\gamma'}(\zeta,\Theta(\zeta,z,0))
\right)+\\
& \
+\cdots+\\
& \ 
+\mu_{n_{j'}'}(z)\, \left(
C(z)\, \bar f^{\gamma_{n_{j'}'}'}(\zeta,\Theta(\zeta,z,0))-\sum_{\gamma'\neq
\gamma_1',\dots,\gamma_{n_{j'}'}'}\, 
B_{n_{j'}',\gamma'}(z)\, 
\bar f^{\gamma'}(\zeta,\Theta(\zeta,z,0))
\right). 
\endaligned\right.
\end{equation}
Finally, we simplify a little bit this expression by 
writing it under the form
\def\theequation{4.3.72}\begin{equation}
\left\{
\aligned
0\equiv 
& \
\mu_1(z)\, C(z)\, 
\bar f^{\gamma_1'}(\zeta,\Theta(\zeta,z,0))+\cdots+
\mu_{n_{j'}'}(z)\, 
C(z)\, 
\bar f^{\gamma_{n_{j'}'}'} (\zeta,\Theta(\zeta,z,0))+\\
& \
+\sum_{\gamma'\neq \gamma_1',\dots,\gamma_{n_{j'}'}'}\, 
E_{\gamma'}(z)\, 
\bar f^{\gamma'} (\zeta,\Theta(\zeta,z,0)).
\endaligned\right.
\end{equation}
where the $E_{\gamma'}(z)$ are formal power series with 
respect to $z$. We are now in position to conclude the
proof of Lemma~4.3.53, namely to come to an absurd as announced
in the beginning of the proof.

Indeed, since $C(z)\not\equiv 0$ by construction, and since there
exists at smallest one power series $\mu_{i'}(z)$ which does not vanish
identically, we can apply the following elementary lemma to conclude
that $\bar f(\zeta,0)$ satisfies a nontrivial power series identity, which 
clearly contradicts the CR-transversality assumption.

\def\thelemma{4.3.73}\begin{lemma}
Assume that there exists a collection of power series $E_{\gamma'}(z)$
indexed by $\gamma'\in\N^{m'}$ with the property that there exists at smallest one
multiindex $\gamma_0'\in\N^{m'}$ with $E_{\gamma_0'}(z)\not \equiv 0$
in $\C\dl z\dr$ and assume that $\bar f(\zeta,\Theta(\zeta,z,0))$
satisfies the formal power series identity
\def\theequation{4.3.74}\begin{equation}
\sum_{\gamma'\in\N^{m'}}\, 
E_{\gamma'}(z)\, \bar f^{\gamma'}(\zeta,\Theta(\zeta,z,0))\equiv 0
\end{equation}
in $\C\dl z,\zeta \dr$. Then there exists a collection of
constants $F_{\gamma'}\in\C$ indexed by $\gamma'\in\N^{m'}$ {\rm which 
do not all vanish} such that 
\def\theequation{4.3.75}\begin{equation}
\sum_{\gamma'\in\N^{m'}}\, 
F_{\gamma'}\, \bar f^{\gamma'}(\zeta,0)\equiv 0
\end{equation}
in $\C\dl \zeta \dr$. In other words, the formal mapping $\zeta\mapsto
\bar f(\zeta,0)$ is {\rm not} transversal, in the sense of Definition~4.1.2
{\bf (5)}.
\end{lemma}

\proof
We put $z=0$ in~\thetag{4.3.74}. If there exists a multiindex
$\gamma'\in\N^{m'}$ such that $E_{\gamma'}(0)\neq 0$, we are done.
Otherwise, we differentiate~\thetag{4.3.74} with respect to 
$z$ and we put $z=0$. If there exist an integer 
$k$ with $1\leq k\leq m$ and a multiindex $\gamma'\in\N^{m'}$
such that $[\partial E_{\gamma'}(0)/\partial z_k](0)\neq 0$,
we are done. Otherwise, we again differentiate with respect to 
$z$ and put $z=0$. Since $E_{\gamma_0'}(z)\not\equiv 0$
in $\C\dl z \dr$, this process converges towards the conclusion
after a finite number of steps, which completes the proof. 
\endproof

The proofs of Lemma~4.3.53 and of part {\bf (4)} of Theorem~4.3.1  are now
complete.
\endproof

We now prove part {\bf (5)} of Theorem~4.3.1. Proceeding as for parts
{\bf (2)} and {\bf (4)}, we shall essentially reason by contradiction,
but we shall summarize the main part of the proof (Lemma~4.3.79 below),
because it is very similar to the proof of Lemma~4.3.53.

By assumption of $\ell_0'$-holomorphic nondegeneracy 
of $M'$ at $p_0'$, we have
\def\theequation{4.3.76}\begin{equation}
{\rm genMrk} \left(
\nabla_{t'} \, \Theta_{j',\gamma'}'(t')
\right)_{1\leq j' \leq d', \, 
\gamma'\in \N^{m'}}=n'.
\end{equation}
Notice that in {\bf (5)}, we make one supplementary assumption, by 
requiring in addition that $h$ is transversal at $p_0$.
This to insure that the following holds.

\def\thelemma{4.3.77}\begin{lemma}
Assume that $h$ is transversal at $p_0$ and that 
$M'$ is holomorphically nondegenerate at $p_0$. Then 
\def\theequation{4.3.78}\begin{equation}
{\rm genMrk} \left(
\nabla_{t'} \, \Theta_{j',\gamma'}'(h(t))
\right)_{1\leq j' \leq d', \, 
\gamma'\in \N^{m'}}=n'.
\end{equation}
\end{lemma}

\proof
Suppose on the contrary that this generic matrix
rank if equal to an integer $\kappa'\leq n'-1$.
By hypothesis, there exists a $n'\times n'$ minor
$d'(t')$ of the $n'\times \infty$ matrix~\thetag{4.3.76}
which does not vanish identically as a power series
of $t'$. We deduce $d'(h(t))\equiv 0$, contradicting the
transversality of $h$ at $p_0$, which 
completes the proof.
\endproof

To establish part {\bf (5)} of Theorem~4.3.1, we need
the following technical lemma.

\def\thelemma{4.3.79}\begin{lemma}
Fix $j'$ with $1\leq j'\leq d'$. Let
$n_{j'}'$ be the integer defined by
\def\theequation{4.3.80}\begin{equation}
n_{j'}':=
{\rm genMrk} \left(
\nabla_{t'} \, \Theta_{j',\gamma'}'(t')
\right)_{\gamma'\in\N^{m'}}.
\end{equation}
Assume that $h$ is CR-transversal at $p_0$. Then 
we also have
\def\theequation{4.3.81}\begin{equation}
{\rm genMrk} \left(
\nabla_{t'}\, 
R_{j',\beta}(0,0,0,0,J^{\vert \beta \vert} \bar h(0) :\,
h(t))
\right)_{\beta \in \N^m}=n_{j'}'.
\end{equation}
\end{lemma}

As in the paragraph after the statement of Lemma~4.3.53, 
we can easily verify that this technical lemma implies
part {\bf (5)} of Theorem~4.3.1. It remains to 
prove Lemma~4.3.79.

\proof
So we fix $j'$. Again, we proceed by contradiction.  We summarize the
proof.  Assuming that the generic rank of the $n'\times \infty$ matrix
\def\theequation{4.3.82}\begin{equation}
\left(
\nabla_{t'}\, 
R_{j',\beta}'(0,0,0,0,J^{\vert \beta \vert} \bar h(0) :\,
h(t))_{\beta\in\N^m}
\right)
\end{equation}
is equal to an integer $\kappa'\leq n_{j'}'-1$ and taking into 
account that 
\def\theequation{4.3.83}\begin{equation}
{\rm genMrk}\, 
\left(
\Theta_{j',\gamma'}'(h(t))
\right)_{\gamma'\in\N^{m'}}=n_{j'}'
\end{equation}
(by a slight generalization of Lemma~4.3.77) since $h$ is transversal at
$p_0$, we deduce that there exist $\kappa'$ distinct multiindices
$\beta_1,\dots,\beta_{\kappa'}\in\N^m$ and a not identically zero
power series $\Lambda(t)\not \equiv 0$ and for every multiindex
$\beta\neq \beta_1,\dots,\beta_{\kappa'}$, power series
$\Lambda_\beta^1(t),\dots,\Lambda_\beta^{\kappa'}(t)$ such that we can
write
\def\theequation{4.3.84}\begin{equation}
\left\{
\aligned
0\equiv 
\sum_{\gamma'\in\N^{m'}}\, 
& \
\left(
\Lambda(t)\, 
\underline{\mathcal{L}}^\beta(\bar f^{\gamma'})(0) \, 
-\Lambda_\beta^1(t)\, 
\underline{\mathcal{L}}^{\beta_1}(\bar f^{\gamma'})(0) \, - \cdots -\right.\\
& \
\left.
- \Lambda_\beta^{\kappa'}(t) \, 
\underline{\mathcal{L}}^{\beta_{\kappa'}}(
\bar f^{\gamma'})(0)
\right)\, \left[
\nabla_{t'} \, \Theta_{j',\gamma'}'(h(t))
\right].
\endaligned\right.
\end{equation}
Reasoning as in the proof of Lemma~4.3.53, we deduce that there exist
$n_{j'}'$ distinct multiindices $\gamma_1',\dots,
\gamma_{j'}'\in\N^{m'}$, that there exists a not identically zero
power series $A(t)\not\equiv 0$, that for $i'=1,\dots,n_{j'}'$ and for
$\gamma'\neq \gamma_1',\dots, \gamma_{n_{j'}'}'$, there exist power
series $A_{i',\gamma'}'(t)\in\C\dl t\dr$, that there exist power
series $\Pi_{i',1}(t),\dots,\Pi_{i',\kappa'}(t)\in\C\dl t\dr$ and that
there exist power series $G_1(t,\zeta),\dots,G_{\kappa'}(t,\zeta)\in
\C\dl t,\zeta\dr$ such that we can write for $i'=1,\dots,
n_{j'}'$:
\def\theequation{4.3.85}\begin{equation}
\left\{
\aligned
A(t)\, 
\Lambda(t)\,
& \
\bar f^{\gamma_{i'}'}(\zeta,0)-\sum_{\gamma'\neq
\gamma_1',\dots,\gamma_{n_{j'}'}'}\, 
\Lambda(t)\, 
A_{i',\gamma'}'(t)\, 
\bar f^{\gamma'}(\zeta,0)\equiv\\
& \
\equiv
\Pi_{i',1}(t)\, G_1(t,\zeta) +\cdots+
\Pi_{i',\kappa'}(t)\, G_{\kappa'}(t,\zeta).
\endaligned\right.
\end{equation}
Since $\kappa'\leq n_{j'}'-1$, we deduce that there  exist
power series $\mu_1(t),\dots,\mu_{n_{j'}'}(t)\in\C\dl t\dr$ not
all zero such that we can write, after setting
$C(t):= A(t)\, \Lambda(t)$ and 
$B_{i',\gamma'}(t):=\Lambda(t)\, 
A_{i',\gamma'}'(t)$:
\def\theequation{4.3.86}\begin{equation}
\left\{
\aligned
0\equiv 
& \
\mu_1(t)\left(
C(t)\, \bar f^{\gamma_1'}(\zeta,0)-
\sum_{\gamma'\neq \gamma_1',\dots,\gamma_{n_{j'}'}'}\, 
B_{1,\gamma'}(t)\, 
\bar f^{\gamma'}(\zeta,0)
\right)+\\
& \
+\cdots+\\
& \
+\mu_{n_{j'}'}(t)\left(
C(t)\, \bar f^{\gamma_{n_{j'}'}'}(\zeta,0)-
\sum_{\gamma'\neq \gamma_1',\dots,\gamma_{n_{j'}'}'}\, 
B_{n_{j'}',\gamma'}(t)\, 
\bar f^{\gamma'}(\zeta,0)
\right).
\endaligned\right.
\end{equation}
Simplifying a little bit this expression by writing it
under the form
\def\theequation{4.3.87}\begin{equation}
\left\{
\aligned
0\equiv 
& \
\mu_1(t) \, C(t)\, 
\bar f^{\gamma_1'}(\zeta,0)+\cdots+
\mu_{n_{j'}'}(t) \, 
C(t) \, 
\bar f^{\gamma_{n_{j'}'}'} (\zeta,0)+ \\
& \
+\sum_{\gamma'\neq \gamma_1',\dots,\gamma_{n_{j'}'}'}\, 
E_{\gamma'}(t)\, 
\bar f^{\gamma'}(\zeta,0),
\endaligned\right.
\end{equation}
we deduce by applying the same reasoning as in Lemma~4.3.73 that there
exists a not identically zero power series relation
\def\theequation{4.3.88}\begin{equation}
\sum_{\gamma'\in\N^{m'}}\, 
F_{\gamma'} \, 
\bar f^{\gamma'}(\zeta,0)\equiv 0, 
\end{equation}
contradicting the assumption of CR-transversality, which 
completes the proof of Lemma~4.3.79.
\endproof

This completes the proof of part {\bf (5)}. In conclusion, the
proof of Theorem~4.3.1 is complete.
\endproof

\bigskip

\noindent
{\Large \bf Bibliography for Part~I}

\vfill
\end{document}